\documentclass[12pt,twoside]{report}
\usepackage{subfiles}
\usepackage[a4paper,asymmetric,hmargin={3cm,3cm},vmargin=3cm]{geometry}
\usepackage{setspace}
\usepackage{fancyhdr,graphicx,lastpage}
\usepackage{listings}
\usepackage{amsmath}
\usepackage{amssymb}
\usepackage{amsthm}
\usepackage{enumerate}
\usepackage{parskip}
\usepackage{multirow}
\usepackage{url}
\usepackage{makeidx}
\usepackage{mathrsfs}
\usepackage{mathtools}
\usepackage{pdfpages}
\usepackage{bookmark}
\usepackage[backend=bibtex,style=alphabetic,citestyle=alphabetic]{biblatex}
\makeindex
\newtheoremstyle{krisjan}
{}
{}
{}
{}
{\bf}
{.}
{3pt}
{}

\theoremstyle{krisjan}
\newtheorem{thm}{Theorem}[chapter]
\newtheorem{prop}[thm]{Proposition}
\newtheorem{lem}[thm]{Lemma}
\newtheorem{cor}[thm]{Corollary}
\newtheorem{defn}[thm]{Definition}

\begin{document}
\begin{titlepage}
	\centering
	{\large \bfseries Sums and products of square-zero matrices \par}
	\vspace{1.5 cm}
	{by}\\
	\vspace{1.5 cm}
	{CHRISTIAAN JOHANNES HATTINGH\par}
	\vspace{1.5cm}
	{submitted in accordance with the requirements\\
for the degree of} \\
	\vspace{1.5cm}
	{MASTER OF SCIENCE\par}
	\vspace{1cm}
	{in the subject}\\
	\vspace{1cm}
	{MATHEMATICS\par}
	\vspace{1.5cm}
	{at the} \\
	\vspace{1cm}
	{UNIVERSITY OF SOUTH AFRICA\par}
	\vspace{1.5cm}
	{SUPERVISOR: PROF J D BOTHA \par} 
	\vspace{2cm}
	{MARCH 2018\par}

	
\end{titlepage}
%



\thispagestyle{empty}
\onehalfspacing
\begin{abstract}
Which matrices can be written as sums or products of square-zero matrices? This question is the central premise of this dissertation. Over the past 25 years a significant body of research on products and linear combinations of square-zero matrices has developed, and it is the aim of this study to present this body of research in a consolidated, holistic format, that could serve as a theoretical introduction to the subject. 

The content of the research is presented in three parts: first results within the broader context of sums and products of nilpotent matrices are discussed, then products of square-zero matrices, and finally sums of square-zero matrices.

\end{abstract}
\addtocounter{page}{-1}
\thispagestyle{empty}

\hspace{0pt}
\vfill
This work is based on research supported in part by the National Research Foundation of South Africa (Grant Number 103932).
\vfill
\hspace{0pt}

\tableofcontents

\newpage

\pagestyle{fancy}
\fancyhf{}
\fancyhead[LO]{\small{\leftmark}}
\fancyhead[RE]{\small{\rightmark}}
\fancyfoot{}
\cfoot{\arabic{page}/\pageref{LastPage}}
\raggedbottom
\chapter{Abstract and Statement of the Research Problem} 
Which matrices can be written as sums or products of square-zero matrices? This question is the central premise of this dissertation. Over the past 25 years a significant body of research on products and linear combinations of square-zero matrices has developed, and it is the aim of this study to present this body of research in a consolidated, holistic format, that could serve as a theoretical introduction to the subject. 

The current body of research arose mainly from two related fields: functional analysis, and linear algebra and matrix theory. Consequently some results are confined to bounded operators on a complex Hilbert space, so referring to square matrices acting as linear transformations on vector spaces that are endowed with an inner product structure \cite{ww} \cite{novak} \cite{novak2} \cite{takahashi} \cite{bukovsek}. These results usually assume the underlying field to be algebraically closed, or fixed to the field of complex numbers. Some of these results extend to the case where the vector space is infinite dimensional. The second source of results are within the context of ``pure" linear algebra and matrix theory: results applicable to vector spaces over an arbitrary field \cite{botha2} \cite{botha3} \cite{pazzis2}. These results are usually confined to the finite dimensional case, and non-square matrices are also considered in some results. 

In light of the above, I would like to state the bounds of the research presented here. Results in this text is framed within a matrix theory and linear algebra perspective. The focus is on finite dimensional vector spaces over an arbitrary field, and as such infinite dimensional results on Hilbert spaces will not be considered as a rule. Furthermore, within this frame of reference, extension of results previously confined to vector spaces over algebraically closed fields will always be considered. Similarly the generalization of results relying on an inner product will be considered.   

The question: ``Which matrices can be written as products of square-zero matrices?" was addressed in a unifying theory, the core result being necessary and sufficient conditions for two matrices over an arbitrary field to have a square-zero quotient, and whenever this is possible, what is the bounds on the rank of such a quotient. This result is then used as a basis for results pertaining to decomposition of a square matrix into square-zero factors: necessary and sufficient conditions for factorization into two, three or arbitrarily many square-zero factors. These results are due to Botha \cite{botha2} and also include formulae to construct a square-zero quotient (whenever two given matrices have such a quotient). 

It is worth noting that before Botha's results, Novak \cite{novak} presented necessary and sufficient conditions for an operator to be a product of two and three square-zero operators, respectively. These results are restricted to square matrices over an algebraically closed field. Some results for the infinite dimensional case were also presented in Novak's article, but will not be considered in this dissertation.

Sums of square-zero matrices were first considered by Wang and Wu \cite{ww}. Results were restricted to bounded linear operators on a complex Hilbert space. The main results are necessary and sufficient conditions for a finite matrix (operator) to be a sum of two square-zero operators, necessary conditions as well as some sufficient conditions for a finite operator to be a sum of three square-zero operators, and necessary and sufficient conditions for a (finite) matrix to be written as a sum of four square-zero matrices. Some results are also provided for the infinite dimensional case. The results by Wang and Wu served as a basis for all further work on the topic. 

First, Takahashi \cite{takahashi} extended the results for sums of three square-zero matrices. The results were constrained to the complex number field, and included necessary and sufficient conditions for some special cases (i.e. diagonalizable matrices with minimum polynomial of degree 2), as well as an extension of sufficient conditions, based on the size of the largest eigenspace and the number of invariant factors of size two. 

Novak \cite{novak} extended results to the special case of sums of two commuting square-zero matrices, providing necessary and sufficient conditions for such a decomposition. Botha \cite{botha3} extended the result on sums of two square-zero matrices: necessary and sufficient conditions for a matrix over an arbitrary field to be written as a sum of two square-zero matrices. 

Finally, de Seguins Pazzis \cite{pazzis2} extended the results on sums of three square-zero matrices to include: necessary and sufficient conditions for a matrix over a field of characteristic two to be a sum of three square-zero matrices. In the same article the result by Wang and Wu \cite{ww} on sums of four square-zero matrices was generalized to the case where the underlying field is arbitrary. Some results for sums of three square-zero matrices over a field not of characteristic two, where the matrix may be augmented by rows or columns of zeros, were also presented. The last-mentioned results are of interest specifically for infinite matrices, and will not be considered as part of this research. 

As part of the research I will also present a comprehensive treatment of results within the more general context: sums and products of nilpotent matrices. Results on products of nilpotent matrices due to Wu \cite{wu2}, Laffey \cite{laffey} and Sourour \cite{sour} will be discussed. The original results contain some errata: these are discussed and corrected in comprehensive proofs. Results on sums of nilpotent matrices due to Wang and Wu \cite{ww} were extended by Breaz and C{\u a}lug{\u a}ƒreanu \cite{breaz}. The relevant result in the text by Breaz and C{\u a}lug{\u a}ƒreanu \cite{breaz} is Proposition 3, and makes use of a special case of Theorem 2 in the same text, which pertains to matrices over a commutative ring. This result is presented here in the special case (that is, for matrices over a field) as a self-contained result.

This concludes a statement of the research problem. 
\chapter{Introduction}
\section{Definitions and Notation}
Table \ref{matnot} provides a summary of the mathematical notation used in this text. Further key notational conventions are described below.

\begin{table}[!ht]
\small
\begin{center}
{\renewcommand{\arraystretch}{1.1}
\renewcommand{\tabcolsep}{0.4cm}
\begin{tabular}{rl}
$\mathscr{F}$ & Arbitrary field. \\
$\mathscr{F}^n$ & Vector space consisting of all $n$-tuples, or column vectors with \\
& $n$ components, over the field $\mathscr{F}$. \\ 
$\mathscr{F}[x]$ & The (infinite dimensional) vector space of all polynomials of\\ 
& finite degree in $x$ over the field $\mathscr{F}$. \\
$\boldsymbol \Phi_{\alpha \beta}(A)$ & The matrix representation of the linear transformation \\
& $A:\mathscr{F}^n \rightarrow \mathscr{F}^m$ with respect to the basis $\alpha$ for $\mathscr{F}^n$, and the \\
& basis $\beta$ for $\mathscr{F}^m$.\\
$A^T$ & The transpose of $A$. \\
$A^*$ & The conjugate transpose of $A=(a_{ij})$, defined as $(a_{ij})^* = (\overline{a_{ji}})$. \\
$A^R$ & Right inverse of $A$. \\
$\text{row}(A)$ & Row space of the matrix $A$. \\
$\text{col}(A)$ & Column space of the matrix $A$. \\
$\text{N}(A)$ & Null space, or kernel, of $A$. \\
$\text{R}(A)$ & The range, or image, of $A$.\\
$\dim(\mathcal{V})$ & Dimension of the vector space $\mathcal{V}$. \\
$\text{n}(A)$ & The dimension of the null space, also known as the nullity, of $A$. \\
$\text{r}(A)$ & The dimension of the column space or range, also known\\
& as the rank of $A$. \\
$\text{tr}(A)$ & The trace of $A$, defined as the sum of entries on the diagonal \\
& of a square matrix $A$. \\
$\det(A)$ & The determinant of the square matrix $A$. \\
$M_{m \times n}(\mathscr{F})$ & The set of all matrices of order $m \times n$ with entries from $\mathscr{F}$. \\
$M_{m}(\mathscr{F})$ & The set of all square matrices of order $m \times m$ with entries \\
& from $\mathscr{F}$. \\
$J_{m}(\lambda)$ & Square matrix of order $m$ with entries: $\lambda$ on the diagonal,\\ 
& one on the subdiagonal, and zero elsewhere.\\
$C(p(x))$ & The companion matrix of $p(x)$ (a detailed definition \\
&follows below). \\
$H(p(x)^e)$ & The hypercompanion matrix of $p(x)^e$ (a detailed definition \\
& follows below). \\
$A \restriction \mathcal{V}$ & The restriction of $A$ to the domain $\mathcal{V}$. \\
$\sim$ & Matrix equivalence. \\
$\stackrel{R}{\sim}$ & Row equivalence. \\
$\approx$ & Matrix similarity. \\
$\oplus$ & Direct sum of vector subspaces.\\
$\text{Dg}[A,B]$ & Diagonal block matrix with $A$ in upper left block,\\
& and $B$ in the lower right block.
\end{tabular}}
\caption{Mathematical notation used in this text. \label{matnot}}
\end{center}
\end{table}

Capital letters such as $A$ indicate a matrix, or a linear transformation defined as multiplication on the left by $A$: so for $A \in M_{m \times n}(\mathscr{F})$ we have the corresponding linear transformation \[A:\mathscr{F}^n \rightarrow \mathscr{F}^m \text{, defined by } A(v)=Av,\] for every $v \in \mathscr{F}^n$. The difference will always be clear within the context of use. Occasionally a matrix may be indicated as $A=(a_{ij})$, where $a_{ij}$ is the entry in row $i$ and column $j$, to place emphasis on a generalized entry in the matrix. A bold font $\mathbf{0}$ is used to indicate a matrix or vector which has only zero entries, emphasizing that the relevant matrix is not of order $1 \times 1$.

Lower case letters such as $v$ generally indicate either column vectors, or row vectors, including scalars. Again, the meaning will be clear from the context. A standard basis vector, usually assumed to be a column vector, will be denoted by $e_i$: explicitly, $e_i$ is the vector with a one in coordinate $i$ and zeros elsewhere. 

Lower case greek letters such as $\alpha$ generally indicate a basis, and as such is a set of linearly independent vectors that span a particular vector space. An exception is made for $\lambda$, which customarily indicates an eigenvalue, or diagonal entry in a matrix. Calligraphic capitals such as $\mathcal{V}$ indicate a vector space or subspace, depending on the context. 

The order of a matrix $A$ is the amount of rows and columns of $A$, for example if $A$ has order $m \times n$, $A$ has $m$ rows and $n$ columns. In general I have followed the convention of indicating a diagonal block matrix consisting of blocks $A$ and $B$ as $\text{Dg}[A,B]$, whereas $\mathcal{A} \oplus \mathcal{B}$ indicates the direct sum of vector subspaces, but the meaning should also be clear from the context.

A matrix $A \in M_n(\mathscr{F})$, where $\mathscr{F}$ is an arbitrary field, is defined as \emph{scalar} if $A = \lambda I_n$ for some $\lambda \in \mathscr{F}$. The \emph{determinant function}, indicated as $\det$, is defined as the unique function $\det: M_n(\mathscr{F}) \rightarrow \mathscr{F}$ with the properties
\begin{enumerate}
\item $\det(AB) = \det(A) \det(B)$ for all $A,B \in M_n(\mathscr{F})$.
\item Let $E_k$ be the elementary matrix obtained from $I_n$ by multiplying any row with $k \neq 0$. Then $\det(E_k) = k$.
\end{enumerate}
The theory related to the determinant function can be referenced in the texts by Cullen \cite[Chapter 3]{cullen} and Golan \cite[Chapter 11]{golan}.

The following definitions are adapted from Cullen \cite[Chapters 5-7]{cullen}.
A simple Jordan block, indicated as $J_{m}(\lambda)$ is defined as the $m \times m$ matrix with $\lambda$ on the diagonal, ones on the subdiagonal, and zeros elsewhere, for example
\[J_3(2) = \begin{bmatrix} 2 & 0 & 0 \\
1& 2 & 0 \\
0 & 1 & 2\end{bmatrix}\]

Define $\mathscr{F}[x]$ as the infinite-dimensional vector space of all polynomials of finite degree in $x$ over the field $\mathscr{F}$. Let $p(x) \in \mathscr{F}[x]$ be a monic polynomial of degree $m$: 
\[p(x)=x^m+a_{m-1}x^{m-1}+ \cdots + a_{1}x + a_0 = x^m-(-a_{m-1}x^{m-1}- \cdots - a_{1}x - a_0). \]
The \emph{companion matrix} $C(p(x))$ of $p(x)$ is the $m \times m$ matrix
\[ C(p(x)) = 
\begin{bmatrix} 
0 & 0 &  \cdots & 0 & -a_0 \\
1 & 0 &  \cdots & 0 & -a_1 \\
0 & 1 & \ddots & \vdots & \vdots \\
\vdots & \ddots & \ddots &0& -a_{m-2}\\
0& \cdots & 0 & 1 & -a_{m-1} 
\end{bmatrix}.\]

The \emph{minimum polynomial} of a matrix $A \in M_n(\mathscr{F})$ is the monic polynomial $p(x)$ of least degree such that $p(A)=\mathbf{0}$. In some texts this polynomial is referred to as the \emph{minimal polynomial} of $A$. If the minimum polynomial of $A$ is the same as the characteristic polynomial of $A$, we say that $A$ is \emph{nonderogatory}. In some texts a matrix with this property is referred to as a \emph{cyclic} matrix. It is easy to show that a companion matrix is always nonderogatory \cite[Theorem 7.1]{cullen}.  

Let $N$ be an $m \times m$ matrix whose only nonzero element is a 1 in entry $(1,m)$. The \emph{hypercompanion matrix} $H(p(x)^e)$ of $p(x)^e$ (where $p(x)$ is some monic polynomial of degree $m$) is the $(e \cdot m) \times (e \cdot m)$ matrix 
\[ H(p(x)^e) = 
\begin{bmatrix} 
C(p(x)) & \mathbf{0} & \mathbf{0} & \cdots & \mathbf{0} \\
N & C(p(x)) & \mathbf{0} & \cdots & \mathbf{0} \\
\mathbf{0} & N & \ddots &  & \vdots \\
\vdots &  & \ddots &&\mathbf{0} \\
\mathbf{0}& \cdots & \mathbf{0} & N & C(p(x)) 
\end{bmatrix}.\]

Suppose  $\mathscr{F}$ is an arbitrary field and $A, P_1 \in M_n(\mathscr{F}), B, P_2 \in M_k(\mathscr{F})$ with $P_1, P_2$ non-singular. Then $A$ is similar to $P_1^{-1}AP_1$, $B$ is similar to $P_2^{-1}BP_2$ and it is easy to verify that $\text{Dg}[A,B]$ is similar to \[\begin{bmatrix}P_1^{-1} & \mathbf{0} \\ \mathbf{0} & P_2^{-1} \end{bmatrix} \begin{bmatrix}A & \mathbf{0} \\ \mathbf{0} & B \end{bmatrix}\begin{bmatrix} P_1 & \mathbf{0} \\ \mathbf{0} & P_2\end{bmatrix}.\] This fact will often be used to simplify a proof, as it is then sufficient to prove that the individual diagonal blocks are similar to some matrix of interest, or if the blocks share some property, it might be sufficient to prove the result for only one of the diagonal blocks.

\section{Significance of the Research}
A square-zero matrix $A$, is a nilpotent matrix with order of nilpotence 2, that is, $A^2=\mathbf{0}$. Square-zero matrices have a well-known structure, and therefore a factorization or linear decomposition producing such matrices could be beneficial in understanding the structure of the matrix being factored, or could result in computational simplification.

Matrix factorization could be described in the simplest sense as writing a matrix as a product of (other) matrices. Sometimes, to distinguish among the resultant factors it is useful to define matrix factorization in terms of matrix division. Since matrix multiplication is not commutative we can define a factorization both in terms of left division and right division, depending on the resultant factors. A formal definition is provided in Section \ref{prelim}. 

A sum decomposition of matrices, is the expression of one matrix as a linear combination of other matrices. The set of all matrices in $M_{m \times n}(\mathscr{F})$ is a vector space over the field $\mathscr{F}$. It is often useful to express a vector as a linear combination of simpler vectors, or vectors that possess particular properties. 

In this section I will primarily focus on the significance of research into the aspects as mentioned above, with reference to the historical context where applicable. I will not attempt to provide a detailed chronological and historical account of developments. I divide the discussion into three subsections: matrix factorization, matrix sums and square-zero matrices.
\subsection{Matrix factorization}
An apt description employed by Halmos, that captures the essential purpose of matrix factorization is ``bad products of good matrices" \cite{halmos}. Given a ``bad" matrix, is it possible to factor it into a product of ``good" matrices? 

Within the field of matrix theory and linear algebra a primary concern is to understand the properties, or structure, of a matrix or the linear transformation it represents. A ``bad" matrix within this context could be described as a matrix with structure or properties that are not readily perceivable. The primary mechanism whereby the properties of such a matrix may be understood is through its equivalence (related via an equivalence relation) to a canonical form (the ``good matrix"). These equivalence relations are in most cases defined in terms of a factorization. Here there are numerous examples: row equivalence, column equivalence, (matrix) equivalence, similarity (which is just a special case of equivalence), etcetera. 

In this sense matrix factorization is fundamental to matrix theory and linear algebra, and it is worth noting that linear algebra's growing importance in undergraduate mathematics \cite{tucker} renders its importance to a variety of fields outside of theoretical mathematics. Case in point, much of undergraduate tuition in linear algebra, from students in engineering to economic sciences, at least includes row reduction (as a method to solve a system of linear equations) and the diagonalization problem. Row reduction could be seen as implicit factorization into a row equivalent right factor with right quotient an invertible matrix, and diagonalization is a similarity factorization. These are typical examples of the many ``great matrix factorizations traditionally written in the form $A=()()$ or $A=()()()$" \cite[p.684]{lee} - indicating that factorization is most useful if there are three or less factors. 

Within applied fields it has been suggested that there are two historical reasons for matrix factorization \cite{hubert}, and the split is best described in terms of two major application fields of matrix theory. The first is numerical linear algebra, which is closely associated with the development of digital computers. In most instances within this field the primary reason for matrix factorization is computational simplification - a typical example is solving a system of linear equations in the form of the matrix equation $Ax=B$. Rather than solving this equation directly, it is sometimes computationally less expensive to first factor: $A=LU$ where $L$ is lower triangular with ones on the diagonal and $U$ is upper triangular \cite[p.68]{hubert}. The $LU$-factorization of a matrix was introduced by Alan Turing in 1948 \cite[p.7]{tucker}.

Whereas the purpose for factorization mentioned above does not imply any significant attribution of meaning to the factors themselves, there are other fields where the factors have significant meaning in terms of the linear transformations represented. Essentially here, the reason for factorization is similar to understanding linear operators through their relation to canonical forms: a ``hidden" structure is revealed that allow insight into the context, or simplifies its complexity. 

In data mining and statistical applications, datasets are represented by matrices, where the rows represent observations and the columns variables or attributes. Matrix factorization within this context can reveal significant underlying structure: generally such factorizations are in the form $A=CWF$ (again, more than three factors are seldom seen within applied fields), where the columns of $C$ represent latent attributes, the rows of $F$ latent observations, and $W$ is usually diagonal, each entry of which indicates the relative importance of the corresponding underlying latent factors \cite[p.24]{skillicorn}. Within this class of factorizations is singular value decomposition, mostly used for filtering out ``noise", and closely related are principle component and factor analysis. 

The last-mentioned two methods amount to orthogonal diagonalization of a symmetric matrix (a correlation or covariance matrix) and are mainly used for removing multi-collinearity (in linear algebra terms: columns that are linearly dependent or nearly linearly dependent). There are many other examples within this category of matrix factorization in data mining applications in particular: semidiscrete decomposition, graph analysis with adjacency matrices, non-negative matrix factorizations, etcetera. For a comprehensive treatment the reader is referred to a text such as \emph{Understanding Complex Datasets, Data Mining with Matrix Decompositions} by Skillicorn \cite{skillicorn}.

\subsection{Matrix sums}
As with matrix factorization, decomposition into a sum of simpler matrices can provide a deeper understanding of the structure of a matrix. The spectral decomposition theorem \cite[p.248]{cullen} \cite[p.429]{golan} shows that a diagonalizable matrix can be expressed as a linear combination of projections onto independent subspaces, namely the eigenspaces of the relevant matrix. A direct application of spectral decomposition within applied science (optics in this case) can be seen in \emph{Sum decomposition of Mueller-matrix images and spectra of beetle cuticles} \cite{arwin}, as an example. 

Spectral decomposition makes it possible to select only a few components of interest, that in effect summarize the essential characteristics of the original linear transformation. Components here refer to the eigenspaces of the original transformation. This method forms the core of statistical techniques such as principal component analysis: a correlation matrix (which is symmetric, and hence always orthogonally diagonalizable) is decomposed by spectral decomposition, and then only the components associated with the largest eigenvalues are selected. These components represent the components in the data that explain most of the observed variation (in a statistical sense) in the data. Typically all components with eigenvalues less than 1 are discarded. In many implementations of such techniques the sum decomposition is implicit, and not necessarily shown as part of the procedure.

A closely related example is singular value decomposition, which is a matrix factorization also mentioned in the previous section. This decomposition can also be expressed as a sum of rank one matrices \cite[p.57]{hopcroft}. For a matrix $A$, each term in this sum represents the projection of the rows of $A$ onto the subspace spanned by some singular vector of $A$ (hence these spaces are independent). A singular vector in this case refers to a vector associated with an eigenvalue of $A^*A$. The singular value decomposition can be written as 
\[ A = \sum_{i=1}^n \sqrt{\lambda_i} u_i v_i^T, \] where $\lambda_1, \lambda_2, \ldots, \lambda_n$ are the non-zero eigenvalues of $A^*A$, $v_1, v_2, \ldots ,v_n$ the orthonormal set of associated eigenvectors of $A^*A$, and $u_i = (1/\sqrt{\lambda_i})Av_i$ for each $i=1,2,\ldots,n$. The sum is arranged in such a way that $\sqrt{\lambda_1} \geq \sqrt{\lambda_2} \geq \cdots \geq \sqrt{\lambda_n}$.
When presented in this way the first $k$ terms of the sum is the best rank $k$ approximation of the transformation $A$, and herein lies vast applications within statistics and computer science: data compression, dimensionality reduction, noise reduction, etcetera. As a reference for the interested reader, the theory from a linear algebra perspective is presented in Proposition 18.15 in the textbook by Golan \cite{golan}.

In terms of a typical graduate course within the sciences, the spectral decomposition theorem is encountered in the intermediate to advanced undergraduate curriculum, and singular value decomposition is often only presented at advanced level or in postgraduate qualifications of certain applied fields such as statistics. One of the first encounters a student in the sciences may have with sum decomposition, would be to write a matrix as a sum of a symmetric and skew-symmetric matrix (when the underlying field of the entries is the real numbers), or the sum of a Hermitian and skew-Hermitian matrix (in the case where the underlying field is the complex numbers). 

Most students may know this decomposition simply as an interesting exercise in matrix algebra, but it does indeed also have important practical implications. One example is within computational fluid dynamics, where a specialisation of the mentioned decomposition can be applied to the matrix representing the pressure and continuity terms of the Navier-Stokes equations for incompressible laminar flow. Due to the particular structure of this matrix, it may be written as the sum of a positive semi-definite symmetric matrix and a skew-symmetric matrix under certain assumptions, and this fact is used to formulate a stable representation of the Navier-Stokes equations. Stable in this case refers to solvability via a systematic iterative method such as Gauss-Seidel or Jacobi iterative algorithms. I refer here to research by Roscoe \cite{roscoe}.

I conclude this section with a brief discussion on the significance of square-zero matrices.
\subsection{Square-zero matrices}
In the preceding discussion we have seen that factorization or sum decomposition reveals the structure of a matrix by exposing its linear components (in the case of sums), or factors (in the case of factorization). These components (or factors) are simpler matrices with a known structure and effect (when considered as linear transformations). Let us consider the structure and transformative effect of a square-zero matrix.
 
A square-zero matrix which is not the zero matrix itself has minimum polynomial $x^2$, and is therefore similar to a matrix with diagonal blocks consisting of $H(x^2)=J_2(0)$,\footnote{$H(x^2)$ indicates the hypercompanion matrix of $x^2$ \cite[Definition 7.3, p.244]{cullen}} and $H(x)=J_1(0)$, so that it is a matrix with a known structure and properties which are particularly accessible, even when compared to nilpotent matrices with a higher order of nilpotence. Applying a square-zero transformation twice to any vector, results in the zero-vector, that is, it ``annihilates" every vector in its image.

One example where the order 2 nilpotency of a square-zero matrix is particularly useful is in calculating a matrix exponential. Let $A$ be square-zero, then 
\[e^A = \sum_{k=0}^\infty \frac{1}{k!}A^k = I + A.\] Therefore if we can decompose a matrix into a sum of square-zero matrices, there is potential for simplifying calculation of its exponential.

Nilpotent operators and matrices in general, are also of importance in other ways in linear algebra and matrix theory. One example is the development of the theory pertaining to the Jordan canonical form, where nilpotent operators play a central role \cite[Chapter 5]{cullen} \cite[Chapter 13]{golan}. In the historical memoir by Peirce titled \emph{Linear Associative Algebra} he regards nilpotent expressions as ``an essential element of the calculus of linear algebras. Unwillingness to accept them has retarded the progress of discovery and the investigation of quantitative algebras." \cite[p.188]{peirce}. In this text Peirce devotes a large section to the development and investigation of properties of nilpotent and idempotent elements within linear associative algebra.

This concludes the discussion on significance of the research.
\chapter{Sums and Products of Nilpotent Matrices}
This chapter is concerned with the broader context: the set of all square-zero matrices is a subset of the set of all nilpotent matrices, and results proved here will therefore be applicable to square-zero matrices.
\section{Preliminaries}
I start by presenting supporting results that are used extensively in what is to follow. First I briefly mention results that are well known.

Let $\mathscr{F}$ be an arbitrary field. A matrix $A \in M_n(\mathscr{F})$ is nilpotent if and only if its characteristic polynomial is $x^n$. Necessity follows directly from the Cayley-Hamilton theorem \cite[Proposition 12.16]{golan} \cite[Theorem 5.1]{cullen}, and sufficiency from the fact that a nilpotent matrix is similar to a matrix in Jordan canonical form, where each simple Jordan block on the diagonal is of the form $J_i(0)$ where $i$ is some positive integer \cite[Proposition 13.4]{golan} \cite[Theorem 5.9]{cullen}. By modifying the necessary part it is easy to show that a similar result holds in terms of the minimum polynomial. Since similar matrices have the same characteristic equation, it also follows that nilpotency is preserved under similarity.

From the above it follows that the characteristic polynomial of a nilpotent matrix splits over an arbitrary field, and its only characteristic value is zero. Now suppose some power of a matrix $A$, say $A^k$, is similar to a matrix with characteristic polynomial $x^n$. By invariance of the characteristic polynomial under similarity it follows that $A^k$ is nilpotent. So there exists a positive integer $m$ so that $(A^k)^m$ is the zero matrix. But then $A^{k\cdot m}$ is the zero matrix and it follows that $A$ must also be nilpotent. These facts will be used implicitly in some of the proofs that follow.

I now present results supporting the main theorems in this section. The first two results form the core constructive argument employed by both Fillmore \cite{fillmore} and Laffey \cite{laffey}, in proofs involving induction on the order of square matrices.
\begin{quote}
\begin{lem}
Let $\mathscr{F}$ be an arbitrary field. If $A$ is a non-scalar matrix in $M_n(\mathscr{F})$, then there exists $x \in M_{n \times 1}(\mathscr{F})$ such that $Ax$ and $x$ are linearly independent.
\label{lem:nonscalar1}
\end{lem}
\end{quote}
{\bf Proof.} I consider two mutually exclusive cases: first suppose that $A$ is a diagonal matrix with entries $\lambda_1, \lambda_2, \ldots, \lambda_n$. Since $A$ is not scalar there exist positive integers $i,j$ such that $\lambda_i \neq \lambda_j$. Let $x \in M_{n \times 1}(\mathscr{F})$ be the vector with $1$ in each entry. Then $Ax = \begin{bmatrix} \lambda_1 & \lambda_2 & \cdots & \lambda_n \end{bmatrix}^T$, and since $\lambda_i \neq \lambda_j$ it follows that $Ax$ is not a scalar multiple of $x$, that is, $Ax$ and $x$ are linearly independent.

Now suppose that $A$ is not diagonal, so that there exist positive integers $i,j$ such that $i \neq j$ and $a_{ij} \neq 0$. Let $x \in M_{n \times 1}(\mathscr{F})$ have a $1$ in entry $j$ and zeros elsewhere. Then entry $i$ of $Ax$ is $a_{ij}$, but entry $i$ of $x$ is zero, and it follows that $Ax$ and $x$ are linearly independent. \hfill $\square$

\begin{quote}
\begin{prop}
Let $\mathscr{F}$ be an arbitrary field and let $\lambda \in \mathscr{F}$. The non-scalar matrix $A \in M_n(\mathscr{F})$, $n \geq 2$ is similar to
\[ A' = \begin{bmatrix} \lambda & c \\ 
e_1 & B \end{bmatrix}\]
 where $c \in M_{1 \times (n-1)}(\mathscr{F})$, $e_1 \in M_{(n-1) \times 1}(\mathscr{F})$ and $e_1$ has a $1$ in entry $(1,1)$ and zero elsewhere. Furthermore the matrix $B \in M_{(n-1) \times (n-1)}(\mathscr{F})$ has trace $\text{tr}(A)-\lambda$, and is non-scalar if $n \geq 3$.
\label{prop:specifyDiagonal}
\end{prop}
\end{quote}
{\bf Proof.} Since $A$ is non-scalar, by Lemma \ref{lem:nonscalar1}, we can find a vector $x$ such that $Ax$ is not a scalar multiple of $x$, then it follows that $\{x, Ax - \lambda x\}$ is a set of linearly independent vectors. Complete this set to a basis $\{x,Ax-\lambda x, x_2, \ldots, x_{(n-1)}\}$. Representing $A$ with respect to this basis we have:
\[ A' = \begin{bmatrix} \lambda & c \\ 
e_1 & B \end{bmatrix}\]
  $c \in M_{1 \times (n-1)}(\mathscr{F})$, $e_1 \in M_{(n-1) \times 1}(\mathscr{F})$ and $e_1$ has a $1$ in entry $(1,1)$ and zero elsewhere, and $B \in M_{(n-1) \times (n-1)}(\mathscr{F})$. Note that $A'$ and $A$ have the same trace since they are similar. This fact follows from a result which is easily proved by routine computation: for all matrices $A,B \in M_n(\mathscr{F})$ we have $\text{tr}(AB) = \text{tr}(BA)$ \cite[p.318]{golan}.  Since $A'$ and $A$ have the same trace it follows that $\text{tr}(B) = \text{tr}(A)-\lambda$. Now if $n=2$ or $B$ is non-scalar the proof is complete.
  
Suppose $n>2$ and $B$ is scalar, so that we have
\[A' = \begin{bmatrix} \lambda & c_1 & \cdots && c_{n-1} \\ 
1 & \lambda_1 & 0 &\cdots &0 \\
0 &0&\ddots&&\vdots \\
\vdots & \vdots & & \ddots & 0 \\
0 & 0 & \cdots & 0 & \lambda_1
 \end{bmatrix},\] where $\lambda_1 \in \mathscr{F}$ is some scalar such that $(n-1)\lambda_1 = \text{tr}(A)-\lambda$. I will now show that in such a case $A'$ is similar to a matrix 
 \[ \begin{bmatrix} \lambda & c' \\ e_1 & B' \end{bmatrix}\] where $B'$ is non-scalar. This similarity transform is achieved by constructing a change-of-basis matrix 
\[ P_1 = \begin{bmatrix} 1 & e_2^T \\ \mathbf{0} & I_{n-1} \end{bmatrix} \] where $e_2^T$ is the $1 \times (n-1)$ matrix with a 1 in entry $(1,2)$ and zero elsewhere. It then follows that
\[P_1^{-1} A' P_1 = \begin{bmatrix} \lambda & c_1 & d_1 & c_3 & \cdots & c_{n-1} \\ 
1 & \lambda_1 & 1 &0&\cdots &0 \\
0 &0&\lambda_1&0&\cdots&0 \\
\vdots & \vdots & & \ddots& & \vdots\\
& & & & & 0 \\
0 & 0 & \cdots & & 0 & \lambda_1
 \end{bmatrix},\] where $d_1 = \lambda + c_2 - \lambda_1$. So we still have $\lambda$ in entry $(1,1)$ and the $(n-1) \times (n-1)$ submatrix in the lower right is non-scalar, completing the proof. \hfill $\square$
 
The following result is due to Fillmore, and forms the core of the argument by Wang and Wu \cite{ww}, and de Seguins Pazzis \cite{pazzis2}, in proving necessary and sufficient conditions for a non-scalar matrix to be written firstly as a sum of two nilpotent matrices and secondly as a sum of four square-zero matrices. Since the result is conveyed very economically in the original text, and since it is the principal device underpinning important proofs in the body of research on square-zero sums, it warrants a detailed description here. 

\begin{quote}
\begin{prop}[Theorem 2 \cite{fillmore}]
Let $\mathscr{F}$ be an arbitrary field. The non-scalar matrix $A \in M_n(\mathscr{F})$ is similar to a matrix with main diagonal $\lambda_1, \lambda_2, \ldots, \lambda_n$, if and only if $\lambda_1 + \lambda_2 + \cdots + \lambda_n = \text{tr}(A)$.
\label{prop:trace}
\end{prop}\end{quote}
{\bf Proof.} First, necessity follows from the fact that similar matrices have the same trace, which as mentioned before easily follows from the fact that for all matrices $A,B \in M_n(\mathscr{F})$ we have $\text{tr}(AB) = \text{tr}(BA)$ \cite[p.318]{golan}.

The proof of the converse is by induction. Since, by hypothesis we are restricted to non-scalar square matrices, a proof by induction on $n \times n$, the order of a square matrix will start with initial step $n=2$: let $A \in M_2(\mathscr{F})$ be a non-scalar matrix. Then the desired result follows directly from the initial construction in the proof of Proposition \ref{prop:specifyDiagonal}, proving the initial step. Explicitly, having selected a pair $\lambda_1, \lambda_2$ such that $\text{tr}(A) = \lambda_1 + \lambda_2$, by Proposition \ref{prop:specifyDiagonal} $A$ is then similar to \[\begin{bmatrix}\lambda_1 & y \\ 1 & \text{tr}(A)-\lambda_1 \end{bmatrix}= \begin{bmatrix}\lambda_1 & y \\ 1 & \lambda_2 \end{bmatrix}\] where $y \in \mathscr{F}$. 

The induction hypothesis is formulated as follows: assume that the result holds for a matrix $A$ of order $k \times k$. We need to prove that the result holds for a matrix of order $(k+1) \times (k +1)$. First, select $\lambda_1, \lambda_2, \ldots, \lambda_{k+1}$ such that $\text{tr}(A) = \lambda_1 + \lambda_2 + \cdots + \lambda_{k+1}$. By Proposition \ref{prop:specifyDiagonal} $A$ is similar to 
\[ \begin{bmatrix} \lambda_1 & c \\ 
e_1 & B \end{bmatrix},\]
where $c \in M_{1 \times k}(\mathscr{F})$, $e_1 \in M_{k \times 1}(\mathscr{F})$ and $e_1$ has a $1$ in entry $(1,1)$ and zero elsewhere, and $B \in M_{k \times k}(\mathscr{F})$ is non-scalar, and such that $\text{tr}(B) = \text{tr}(A)-\lambda_1=\lambda_2 + \cdots + \lambda_{k+1}$. By the induction hypothesis it follows that $B$ is similar to a matrix with diagonal $\lambda_2, \lambda_3, \ldots, \lambda_{k+1}$. Let $P^{-1}BP$ be a matrix with this property. It follows that 
\[A \approx \begin{bmatrix} 1 & \mathbf{0} \\ 
\mathbf{0} & P^{-1} \end{bmatrix}\begin{bmatrix} \lambda_1 & c \\ 
e_1 & B \end{bmatrix}\begin{bmatrix} 1 & \mathbf{0} \\ 
\mathbf{0} & P \end{bmatrix},\] and the product on the right-hand side is a matrix with diagonal $\lambda_1,\lambda_2, \ldots, \lambda_{k+1}$, completing the proof. \hfill $\square$

The following result was proved independently by Laffey \cite{laffey} and Sourour \cite{sour2} and was then used by Laffey to prove an important result concerning products of nilpotent matrices.

\begin{quote}
\begin{prop}[Theorem 1.1 \cite{laffey}]
Let $\mathscr{F}$ be an arbitrary field. Let $A \in M_n(\mathscr{F})$ be non-scalar and invertible. Let \[x_1,x_2,\ldots,x_n,y_1,y_2, \ldots,y_n\] be arbitrary elements of $\mathscr{F}$ subject to the constraint \[\det (A) = x_1x_2\cdots  x_ny_1y_2 \cdots y_n.\] Then $A$ is similar to $LU$, where $L$ is lower triangular with diagonal entries $x_1,x_2,\ldots,x_n$, and $U$ is upper triangular with diagonal entries $y_1,y_2, \ldots,y_n$.
\label{prop:laffeyTriangle}
\end{prop}\end{quote}
{\bf Proof.}
Let $z_i = x_iy_i$ where $i \in \{1,2,\ldots,n\}$. Notice that it is sufficient to prove that a matrix $A$ as specified in the statement of the theorem is similar to $L_1U_1$ where $L_1$ is lower triangular with only ones on the diagonal and $U_1$ is upper triangular with diagonal $z_1,z_2, \ldots, z_n$. Let $D_1 = \text{Dg}[x_1,x_2, \ldots, x_n]$, $L=L_1D_1$ and $U=D_1^{-1}U_1$. Then $LU = L_1U_1$ and $L$ and $U$ have the desired diagonal elements. I proceed with a proof by induction.
 
Since we are dealing with non-scalar matrices, the initial step of induction is the case $n=2$. By Proposition \ref{prop:specifyDiagonal} $A$ is similar to
\[ \begin{bmatrix} z_1 & y \\ 1 & x \end{bmatrix} = \begin{bmatrix} 1 & 0 \\ z_1^{-1} & 1 \end{bmatrix} \begin{bmatrix} z_1 & y \\ 0 & x - z_1^{-1}y \end{bmatrix}, \] where $x, y \in \mathscr{F}$, and the factorization is possible since $A$ is non-singular, and $z_1$ is therefore not zero. Now since $\det(A) = z_1z_2$, it follows that $z_2 = x-z_1^{-1}y$.

Now suppose the result holds for $n=k$, that is any invertible non-scalar matrix of order $k \times k$ is similar to $L_1U_1$. I proceed to prove that the result then holds for any invertible non-scalar $A \in M_{k+1}(\mathscr{F})$. By Proposition \ref{prop:specifyDiagonal} $A$ is similar to 
\begin{equation} A'=\begin{bmatrix} z_1 & c \\ e_1 & B \end{bmatrix} = \begin{bmatrix} 1 & \mathbf{0} \\ z_1^{-1}e_1 & I_k \end{bmatrix} \begin{bmatrix} z_1 & c \\ \mathbf{0} & B - z_1^{-1}e_1c \end{bmatrix},\label{eq:luFactor} \end{equation} where $B \in M_{k}(\mathscr{F})$, $c \in M_{1 \times k}(\mathscr{F})$ and $e_1 \in M_{k \times 1}(\mathscr{F})$ has a 1 in entry $(1,1)$ and zeros elsewhere. Now to be able to apply the induction hypothesis to $B-z_1^{-1}e_1c$ (let this matrix be $B_1$) two requirements must be fullfilled. First, $B_1$ must be non-scalar. Suppose this requirement is not fulfilled, so that $B_1=bI_k$ for some scalar $b \in \mathscr{F}$. Notice that 
\[B_1 = B-z_1^{-1}e_1c = B-z_1^{-1}\begin{bmatrix} c \\ \mathbf{0} \end{bmatrix},\] that is, $B_1$ is $B$ with the first row modified by subtracting $z_1^{-1}c$.

Let \[P = \begin{bmatrix} 1 & e_2^T \\ \mathbf{0} & I_k\end{bmatrix},\] where $e_2^T \in M_{1 \times k}(\mathscr{F})$ has a one in entry $(1,2)$ and zeros elsewhere. Then 
\[ PA'P^{-1}= \begin{bmatrix} 1 & e_2^T \\ \mathbf{0} & I_k\end{bmatrix}
\begin{bmatrix} z_1 & c \\ e_1 & B \end{bmatrix}
 \begin{bmatrix} 1 & -e_2^T \\ \mathbf{0} & I_k\end{bmatrix} = 
  \begin{bmatrix} z_1 & c' \\ e_1 & B-e_1e_2^T \end{bmatrix},
\]
where $c' = -z_1e_2^T + c + e_2^TB$. Now we can write
\[
 \begin{bmatrix} z_1 & c' \\ 
 e_1 & B-e_1e_2^T \end{bmatrix} = 
 \begin{bmatrix} 1 & \mathbf{0} \\ 
 z_1^{-1}e_1 & I_k \end{bmatrix} 
 \begin{bmatrix} z_1 & c' \\ 
 \mathbf{0} & B_2 \end{bmatrix},
\]
where $B_2 = B - e_1e_2^T - z_1^{-1}e_1c'$. Simplifying $B_2$ we have \begin{eqnarray}
\nonumber B - e_1e_2^T - z_1^{-1}e_1c' &=& B - e_1e_2^T - z_1^{-1}e_1(-z_1e_2^T + c + e_2^TB) \\
\nonumber &=& B- e_1e_2^T + e_1e_2^T - z_1^{-1}e_1c - z_1^{-1}e_1e_2^TB \\
\nonumber &=& B_1 - z_1^{-1}\begin{bmatrix} e_2^T \\ \mathbf{0} \end{bmatrix}B \\
\nonumber &=& bI_{k} - z_1^{-1}\begin{bmatrix} \text{row}_2(B) \\ \mathbf{0} \end{bmatrix} \\
\nonumber &=& \begin{bmatrix} b & -e_1' \\ \mathbf{0} & bI_{k-1} \end{bmatrix},
\end{eqnarray}
where $e_1'  = \begin{bmatrix} z_1^{-1}b & 0 & \cdots & 0 \end{bmatrix}$, which follows from the fact that row 2 of $B$ and $B_1$ are the same. It follows that $B_2$ is non-scalar, proving that we can always satisfy the first requirement, by applying an additional similarity transformation step as above when necessary, replacing $A'$ with $PA'P^{-1}$ and, consequentially $B_1$ with $B_2$. 

At this point the induction hypothesis can be applied to the matrix $B_1$, provided that we can find a similarity transformation such that $z_1$ will remain in entry $(1,1)$, and therefore factorization over the entire matrix remains valid (this is the second requirement). Let $P_1$ be a matrix such that $P_1^{-1}B_1P_1=L_1U_1$ as per the induction hypothesis. Then \[\begin{bmatrix} 1 & \mathbf{0} \\ \mathbf{0} & P_1^{-1} \end{bmatrix}
\begin{bmatrix} z_1 & c \\ e_1 & B \end{bmatrix}
\begin{bmatrix} 1 & \mathbf{0} \\ \mathbf{0} & P_1 \end{bmatrix} = 
\begin{bmatrix} z_1 & cP_1 \\ P_1^{-1}e_1 & P_1^{-1}BP_1 \end{bmatrix}, \] and then applying the factorization \eqref{eq:luFactor}, the matrix on the right-hand side can be written as
\begin{eqnarray} \nonumber \begin{bmatrix} z_1 & cP_1 \\ P_1^{-1}e_1 & P_1^{-1}BP_1 \end{bmatrix} &=& \begin{bmatrix} 1 & \mathbf{0} \\ z_1^{-1}P_1^{-1}e_1 & I_k \end{bmatrix} \begin{bmatrix} z_1 & cP_1 \\ \mathbf{0} & P_1^{-1}BP_1 - z_1^{-1}P_1^{-1}e_1cP_1 \end{bmatrix} \\
\nonumber &=& \begin{bmatrix} 1 & \mathbf{0} \\ z_1^{-1}P_1^{-1}e_1 & I_k \end{bmatrix} \begin{bmatrix} z_1 & cP_1 \\ \mathbf{0} & P_1^{-1}B_1P_1 \end{bmatrix} \\
\nonumber &=& \begin{bmatrix} 1 & \mathbf{0} \\ z_1^{-1}P_1^{-1}e_1 & I_k \end{bmatrix} \begin{bmatrix} z_1 & cP_1 \\ \mathbf{0} & L_1U_1 \end{bmatrix} \\
\nonumber &=&  \begin{bmatrix} 1 & \mathbf{0} \\ z_1^{-1}P_1^{-1}e_1 & L_1 \end{bmatrix} \begin{bmatrix} z_1 & cP_1 \\ \mathbf{0} & U_1 \end{bmatrix},
\end{eqnarray} which is a factorization of the desired form, completing the proof. \hfill $\square$

The following result is the main device used by Sourour \cite{sour} in an inductive proof of the main result on products of nilpotent matrices. Note that there is an implicit assumption made in the proof of this proposition by Sourour \cite[p.304]{sour}. On p.305 \cite{sour} the statement is made that $\text{N}(PAP) = \mathcal{M}_1 \oplus \text{N}(A)$ (I make use of the exact notation used by Sourour here). By the construction given in the proof, it is assumed that the vector $e_0$ is not in $\mathcal{M}_2 = \text{R}(P)$. Suppose $e_0$ is in $\mathcal{M}_2 = \text{R}(P)$: then we will have $N(PAP) = \mathcal{M}_1 \oplus \text{N}(A) \oplus \text{span}\{e_0\}$. 

Now assuming $e_0 \in \mathcal{M}_2$, I present a simple counter-example to this part of the proof: Consider $A=J_4(0)$, which is not square-zero. Select as $e_0$ the standard basis vector with a one in the first entry and zeros elsewhere (the vector denoted as $e_1$ in this dissertation). Then by the given construction $\mathcal{M}_1 = \text{span}\{(0,1,0,0)^T\}$, and $\text{N}(J_4(0)) = \text{span}\{(0,0,0,1)^T\} \subseteq \mathcal{M}_2$. Furthermore, if we assume we may select $e_0$ as one of the basis vectors in $\mathcal{M}_2$, then we will have
\[\text{N}(PAP) = \text{span}\{(0,1,0,0)^T\} \oplus \text{span}\{(0,0,0,1)^T\} \oplus \text{span}\{e_0\}.\] 
Therefore the rank of $PAP$ will be less than the rank of $AP$ and $PA$, and the result of the proposition cannot be achieved.

I now present a proof of the sufficient part of Sourour's proposition where the situation as shown above is explicitly avoided. Sufficiency is all that is required for the proof of the main theorem on nilpotent factorization.
\begin{quote}
\begin{prop}[Proposition \cite{sour}]
Let $\mathscr{F}$ be an arbitrary field. If $A \in M_n(\mathscr{F})$ is not square-zero, then it is similar to a matrix of the form 
\[\begin{bmatrix} \lambda & c \\ b & D \end{bmatrix}\] with $\lambda \in \mathscr{F}$, $D \in M_{(n-1)}(\mathscr{F})$, $\text{r}(D) = \text{r}(A) - 1$, $b \in \text{R}(D)$ and $c \in \text{row}(D)$.
\label{prop:sourInductSingular}
\end{prop}\end{quote}
{\bf Proof.} First suppose that $A$ is a scalar matrix. In this case the result is immediate, as $A=\text{Dg}[\lambda, \lambda I_{(n-1)}]$ is already in the desired form.

Suppose that $A$ is not square-zero, and $A$ is not scalar. Then we can find a vector $x_0$ such that $Ax_0$ and $x_0$ are linearly independent, and $A^2x_0 \neq \mathbf{0}$. Explicitly, if $A$ is full rank, this result is easy to see, as the fact that $A$ is not scalar ensures that there exists a vector $x_0$ such that $Ax_0$ and $x_0$ are linearly independent, and then $A^2x_0 \neq \mathbf{0}$ since $A$ is full rank. Suppose $A$ is not full rank. Since $A$ is not square-zero there exists $x_1 \in \text{R}(A)$ such that $x_1 \notin \text{N}(A)$. Let $Ax_0=x_1$. Now if $x_0$ and $Ax_0 = x_1$ are linearly independent, the claim holds. If this is not the case then $\lambda x_0 = Ax_0$. Replace $x_0$ with $x_0 + n$ where $n$ is a nonzero vector in the null space of $A$: then $A(x_0 + n) = Ax_0 = \lambda x_0$ and since $x_0$ and $x_0 + n$ are linearly independent the claim holds.

Having established the existence of a suitable vector $x_0$, construct a basis for $\mathscr{F}^n$ in the following way 
\[ \alpha_1 = \{x_1, x_1-x_0, n_1, n_2, \ldots, n_k, u_1, u_2, \ldots, u_{(n-k-2)}\},\] where $n_1, n_2, \ldots n_k$ are basis vectors for the null space of $A$, and $u_1, u_2, \ldots, u_{(n-k-2)}$ are arbitrary subject to $\alpha_1$ being a linearly independent set. Let \[\mathcal{V}_1 = \text{span}\{x_1\} \text{ and}\] \[\mathcal{V}_2 = \text{span}\{x_1-x_0,n_1, n_2, \ldots, n_k,u_1,u_2,\ldots,u_{(n-k-2)}\},\] so that $\mathscr{F}^n = \mathcal{V}_1 \oplus \mathcal{V}_2$. Let $P$ be the projection along $\mathcal{V}_1$ onto $\mathcal{V}_2$. Now the following two results hold:
\[ \text{N}(PA) = \text{span}\{x_0\} \oplus \text{N}(A), \]
\[ \text{N}(AP) = \mathcal{V}_1 \oplus \text{N}(A).\]
The results hold, since by virtue of the construction of $\alpha_1$ and $P$ we have that the null space of $P$ is contained in the range of $A$, and the null space of $A$ is contained in the range of $P$.

Now the key to obtaining the required result is that we must have
\[\text{N}(PAP) = \text{N}(AP) =  \mathcal{V}_1 \oplus \text{N}(A).\]
For this result to be true, we must have $x_1 \notin \text{R}(AP)$. I will now prove this fact. By virtue of the construction of $\alpha_1$ we have $x_0 \notin \text{R}(P)$, for if it was then $x_0 + (x_1-x_0) = x_1 \in \text{R}(P) = \mathcal{V}_2$ which is a contradiction. Since $x_0 \notin \text{R}(P)$ it follows that $Ax_0 = x_1 \notin \text{R}(AP)$, as required.

Applying the rank-nullity theorem to each of the products above we have 
\begin{equation} \text{r}(PA) = \text{r}(AP) = \text{r}(PAP) = \text{r}(A) - 1. \label{eq:sourRank} \end{equation}

Now relative to the basis $\alpha_1$ the matrix representations of the products above are
\[PA \approx \begin{bmatrix} 0 & \mathbf{0} \\ b & D \end{bmatrix},\]
\[AP \approx \begin{bmatrix} 0 & c \\ \mathbf{0} & D \end{bmatrix},\]
\[PAP \approx \begin{bmatrix} 0 & \mathbf{0} \\ \mathbf{0} & D \end{bmatrix}.\]

Combining this result with \eqref{eq:sourRank}, it is immediately apparent that $b$ is in the column space of $D$ and $c$ is in the row space of $D$, as required. 
%
%
\hfill $\square$

\section{Products of Nilpotent Matrices}
The first result on products of nilpotent matrices is due to Wu \cite{wu2}, but was restricted to $M_n(\mathbb{C})$, the set of square matrices over the complex field. This result was extended to $M_n(\mathscr{F})$ where $\mathscr{F}$ is an arbitrary field, independently by Sourour \cite{sour} and Laffey \cite{laffey}. The statement of the theorem remained essentially the same as formulated by Wu. As the approaches of Laffey and Sourour are quite different I will present these two proofs. 
\subsection{Laffey's proof}
Laffey's proof requires some of Wu's results which are valid over an arbitrary field (since it deals only with nilpotent matrices). 

Now I provide some remarks regarding the proofs of the following two results below: specifically the method of construction employed in finding most of the change-of-basis matrices in these results (please note that this method was not used to find every change-of-basis matrix listed in the proofs). By Theorem 5.7, 5.8 and 5.9 in the textbook by Cullen \cite{cullen}, for a nilpotent matrix $A \in M_n(\mathscr{F})$, there exist vectors $\xi_1, \ldots, \xi_k$ and invariant subspaces $\mathcal{V}_1, \ldots, \mathcal{V}_k$ of $\mathscr{F}^n$ such that $\mathcal{V}_1 \oplus \cdots \oplus \mathcal{V}_k = \mathscr{F}^n$ and for each $i \in \{1,2,\ldots,k\}$ the set $\{\xi_i, A\xi_i, \ldots, A^{p_i}\xi_i\}$ is a basis for $\mathcal{V}_i$. Note that $p_i$ is less than or equal to the index of nilpotency of $A$. 

Now since most of the factors listed in Lemma \ref{lem:j0nilpotent} and Lemma \ref{lem:j0nilpotent2} below consist entirely of columns that are either zero vectors or standard basis vectors, it is trivial to select the vectors $\xi_i$: suppose row $j$ of the nilpotent factor $B$ is a zero vector, then let $\xi_i = e_j$, and then by repeatedly applying $B, B^2, \ldots, B^{p_i}$ to $e_j$ (until it is annihilated), the basis for the relevant subspace may be constructed. The process is repeated for all remaining zero rows of $B$. 

Here follows a short example: take $N_1$ as listed in \eqref{eq:jk0_k_odd} which consists only of standard basis vectors and zero columns. Since row $1$ is the zero vector we know that $e_{1}$ is not in the range of $N_1$, therefore if it is in a basis $\{\xi_1, N_1\xi_1, \ldots, N_1^{p}\xi_1\}$ it must be $\xi_1$ itself. And indeed notice that 
\[N_1e_{1} = e_{3}, N_1e_{3} = e_5,\ldots,N_1e_{k} = e_2, N_1e_{2} = e_4,\ldots,N_1e_{k-3} = e_{k-1}.\] The sequence $e_1,e_3, \ldots, e_k, e_2, e_4, \ldots, e_{k-1}$ forms a basis, relative to which the transformation $N_1$ shifts each basis vector to the next basis vector, which is exactly the effect of $J_k(0)$ on the standard basis ($J_k(0)e_1 = e_2, J_k(0)e_2=e_3,\ldots$). 

I now proceed with a proof of Wu's results.
\begin{quote}
\begin{lem}[Lemma 2 \cite{wu2}]
For any positive integer $k \neq 2$ the simple Jordan block $J_k(0)$ is the product of two nilpotent matrices with the same rank as $J_k(0)$.
\label{lem:j0nilpotent}
\end{lem}
\end{quote}
{\bf Proof.} In what is to follow $E_{(i,j)}$ denotes the square matrix with a 1 in entry $(i,j)$ and zero elsewhere. The case where $k=1$ is trivial, so assume $k \geq 3$. Now either $k$ is odd, or $k$ is even. Suppose first that $k$ is odd, then 

\begin{equation} J_k(0) = \left [ {\renewcommand{\arraystretch}{1.2}\begin{array}{c|c|c} 
\mathbf{0} &  \begin{matrix} 0 \\ 0 \end{matrix} & \begin{matrix} 0 \\ 1 \end{matrix}\\
\hline 
I_{k-2} &  \mathbf{0} &  \mathbf{0} \end{array}} \right ]
\left [ {\renewcommand{\arraystretch}{1.2}\begin{array}{c|c|c} 
\mathbf{0} &  I_{k-2} & \mathbf{0} \\
\hline 
0 & \mathbf{0} & 0 \\
\hline
1 &  \mathbf{0} &  0 \end{array}} \right ]=N_1N_2.
\label{eq:jk0_k_odd}
\end{equation}
Consider the right-hand side. Both factors have rank $k-1$ as desired. Denote the factor on the left as $N_1$. Let 
\begin{equation} Q_1=[e_1,e_3,\ldots,e_{k},e_2,e_4,\ldots,e_{k-1}] \label{eq:jk0_k_odd_q1} \end{equation} where $e_i$ denotes the standard basis  (column) vector with a one in coordinate $i$ and zeros elsewhere. Then $Q_1$ is a change-of-basis matrix and
\[Q_1^{-1}N_1Q_1 = J_k(0),\] that is, $N_1$ is similar to $J_k(0)$ and is therefore nilpotent.

Now consider the right factor: denote this factor by $N_2$. Let 
\[P_1 = \begin{bmatrix} I'_{k-1} & \mathbf{0} \\ \mathbf{0} & 1 \end{bmatrix},\] where $I'$ denotes the square matrix with all ones on the anti-diagonal, and zeros elsewhere. Then $P_1^{-1} N_2 P_1 = J_k(0)$, confirming that $N_2$ is nilpotent.

Now suppose $k$ is even, $k \geq 4$. Let $N_2$ be defined as above. Then
\begin{equation}
J_k(0) = \left [ {\renewcommand{\arraystretch}{1.2}\begin{array}{c|c|c} 
\mathbf{0} & \mathbf{0} & \begin{matrix} 0 & 0 \\ 1 & 1 \end{matrix}\\
\hline 
\begin{matrix} 1 & -1\\ 0 & 1 \end{matrix} &  \mathbf{0} &  \mathbf{0} \\
\hline
\mathbf{0} & I_{k-4} & \mathbf{0} \end{array}} \right ]
(N_2 + E_{(1,3)}) = N_3N_4.
\label{eq:jk0_k_even}
\end{equation} 
Again, we can directly observe that the factors on the right-hand side are each of the same rank as $J_k(0)$, and it remains to verify that each factor is nilpotent. Consider first the left factor (denoted as $N_3$): notice that the bottom row partition vanishes if $k=4$. Let 
\begin{equation}Q_2 = [e_1,e_3,\ldots,e_{k-1},e_2,e_4-e_3,e_6-e_5,\ldots,e_k-e_{k-1}]. 
\label{eq:jk0_k_even_q2} 
\end{equation} Then $Q_2$ is a change-of-basis matrix and 
\[Q_2^{-1}N_3Q_2 = J_k(0), \] and therefore $N_3$ is nilpotent.

Now observe that $P_1^{-1}(N_2 + E_{(1,3)})P_1 = J_k(0) + E_{(k-1,k-3)}$, where $P_1$ is the change of basis matrix as defined above. As this matrix is lower triangular with all zeros on the diagonal it is nilpotent, as desired, concluding the proof. \hfill $\square$

\begin{quote}
\begin{lem}[Lemma 3 \cite{wu2}]
Let $\mathscr{F}$ be an arbitrary field. Any nilpotent matrix $N \in M_n(\mathscr{F})$, with $n \neq 2$, is the product of two nilpotent matrices, each with the same rank as $N$.
\label{lem:j0nilpotent2}
\end{lem}
\end{quote}
{\bf Proof.} Note that, as in the proof of Lemma \ref{lem:j0nilpotent}, $E_{(i,j)}$ denotes the square matrix with a 1 in entry $(i,j)$ and zero elsewhere, and $e_i$ the standard basis vector which has a one in coordinate $i$ and zero elsewhere.

As mentioned in the preliminary material of this chapter, any nilpotent matrix is similar to a matrix in Jordan canonical form, where each simple Jordan block on the diagonal is of the form $J_i(0)$ ($i$ is some positive integer). Furthermore nilpotency is preserved under similarity, and therefore it is sufficient to prove this result for a matrix in Jordan canonical form with all zeros on the diagonal. Let us assume $N$ is in this form. 

Now by Lemma \ref{lem:j0nilpotent} if $N$ contains no simple Jordan blocks $J_2(0)$, the proof is complete, as each simple Jordan block on the diagonal is a product of two nilpotent matrices. So suppose $N$ contains one or more blocks $J_2(0)$. Notice that it is sufficient to prove the result for each of the following matrices: $\text{Dg}[J_k(0),J_2(0)]$ where $k \geq 1$ and $\text{Dg}[J_2(0),J_2(0),J_2(0)]$. All other $N \in M_n(\mathscr{F})$ with $n \neq 2$ containing $J_2(0)$ can be constructed from a combination of these two and a submatrix containing only diagonal blocks of the form $J_k(0)$ where $k \neq 2$, and therefore combining Lemma \ref{lem:j0nilpotent} with factorizations for the aforementioned forms will yield the desired result. 

Let us first consider $\text{Dg}[J_k(0),J_2(0)]$. First consider the case $k=1$, then 
\begin{equation} \text{Dg}[J_1(0),J_2(0)] = E_{(3,1)}E_{(1,2)}, \label{eq:nilpotent_jkj2_1} \end{equation} and it is immediately apparent that both factors are nilpotent and rank 1 as required.
Now suppose $k \geq 2$. Then either $k$ is even, or $k$ is odd. Suppose first that $k$ is even. Then 
\begin{equation} \text{Dg}[J_k(0),J_2(0)] = \left [ {\renewcommand{\arraystretch}{1.2}\begin{array}{c|c} 
\mathbf{0} & J_k(0) \\ 
\hline
\begin{matrix} 0 & 0 \\ 0 & 1 \end{matrix} &  \mathbf{0} \end{array}} \right ] 
\left [ {\renewcommand{\arraystretch}{1.2}\begin{array}{cc|c} 
\mathbf{0} & \mathbf{0} & J_2(0) \\ 
\hline
 I_{k-1} & \mathbf{0}  & \mathbf{0} \\ 
\mathbf{0} & 0 &  0 \end{array}} \right ]=N_1N_2. \label{eq:nilpotent_jkj2_even} \end{equation}
It is easy to verify that each factor has the same rank as $\text{Dg}[J_k(0),J_2(0)]$. Consider the factor on the left (denoted by $N_1$). Applying a similarity transform using the change-of-basis matrix $P_1 = \text{Dg}[I'_{k+1},1]$, where $I'$ denotes the matrix with ones on the anti-diagonal, the resultant matrix is lower triangular with only zeros on the diagonal, and therefore the left factor is nilpotent. Explicitly
\[P_1^{-1}N_1P_1 = \begin{bmatrix} J_k(0) & \mathbf{0} \\ 
E_{(2,k)} & \mathbf{0} \end{bmatrix},\] where $E_{(2,k)}$ indicates a matrix with a one in entry $(2,k)$ and zeros elsewhere.

To verify that the right factor is nilpotent (denoted by $N_2$), we can make use of the change-of-basis matrix $Q_1=[e_1,e_3,\ldots,e_{k+1},e_2,e_4,\ldots,e_{k+2}]$ to show that $N_2$ is similar to $\text{Dg}[J_{k+1}(0),J_1(0)]$.

Suppose $k$ is odd, then  
\begin{equation} \text{Dg}[J_k(0),J_2(0)] = \left [ {\renewcommand{\arraystretch}{1.2}\begin{array}{c|c} 
\mathbf{0} &  A_1 \\
\hline
J_2(0) & \mathbf{0} \end{array}} \right ] 
\left [ {\renewcommand{\arraystretch}{1.2}\begin{array}{c|c} 
\mathbf{0} &  \begin{matrix} 1 & 0 \\ 0 & 0 \end{matrix} \\ 
\hline
 A_2 & \mathbf{0} \end{array}} \right ]=N_3N_4, \label{eq:nilpotent_jkj2_odd} \end{equation}
 where 
 \[A_1 = [e_2,\mathbf{0},e_4,e_3,e_6,e_5,\ldots,e_{k-1},e_{k-2},e_k]\] and
 \[A_2 = [e_1,e_4,e_3,e_6,e_5,\ldots,e_{k-1},e_{k-2},e_k,\mathbf{0}].\]
 Notice in particular that for the case $k=3$ we have $A_1=[e_2,\mathbf{0},e_3]$ and $A_2 = [e_1,e_3,\mathbf{0}]$.
 Now the factor $N_3$ has rank $k$ which is the same as the rank of $\text{Dg}[J_k(0),J_2(0)]$. Let 
\[ Q_2 = [e_1,e_{k+2},e_{k},e_{k-1},e_{k-4},e_{k-5},e_{k-8},\ldots \qquad \qquad \qquad \qquad \qquad\]
 \begin{equation} \qquad \ldots,e_{3-(-1)^{(k-3)/2}},e_{k+1},e_{k-2},e_{k-3},\ldots,e_{3+(-1)^{(k-3)/2}}]. \label{eq:nilpotent_jkj2_odd_q} \end{equation} Then $Q_2$ is a change-of-basis matrix and 
 \[Q_2^{-1}N_3Q_2 = \text{Dg}[J_{k-2 \cdot (\lfloor (k-3)/4 \rfloor)+(-1)^{(k-3)/2}}(0),J_{k-2 \cdot (1+\lfloor (k-3)/4 \rfloor)}(0)],\] so that $N_3$ is nilpotent.
 
 Finally, it is easy to verify that $N_4$ also has rank $k$. Let 
 \[Q_3 = [e_4,e_8,\ldots,e_{k+1},e_1,e_3, \ldots,e_{k},e_2,e_6, \ldots,e_{k-1},e_{k+2} ] \] when $(k-3)/2$ is even, and 
 \[Q_3 = [e_2,e_6, \ldots, e_{k+1},e_1,e_3, \ldots,e_{k},e_4,e_8,\ldots,e_{k-1},e_{k+2}] \] when $(k-3)/2$ is odd. Then 
 \[Q_3^{-1} N_4 Q_3 = \text{Dg} [J_{k-\lfloor (k-3)/4 \rfloor}(0),J_{2+\lfloor (k-3)/4 \rfloor}(0)],\] confirming that $N_4$ is nilpotent.
 
 Finally,
 \begin{eqnarray} \nonumber \text{Dg}[J_2(0),J_2(0),J_2(0)] &=& \left [ {\renewcommand{\arraystretch}{1.2}\begin{array}{c|c|c} 
\mathbf{0} &  \mathbf{0} & \begin{matrix} 0 & 0 \\ 0 & 1 \end{matrix} \\ 
\hline
J_2(0) & \mathbf{0} & \mathbf{0} \\
\hline
 \mathbf{0} & J_2(0) & \mathbf{0} \end{array}} \right ]
 \left [ {\renewcommand{\arraystretch}{1.2}\begin{array}{c|c|c} 
\mathbf{0} & \begin{matrix} 1 & 0 \\ 0 & 0 \end{matrix} &  \mathbf{0} \\ 
\hline
\mathbf{0} & \mathbf{0} &  \begin{matrix} 1 & 0 \\ 0 & 0 \end{matrix}  \\
\hline
  J_2(0) & \mathbf{0} & \mathbf{0} \end{array}} \right ] \\
  &=& N_5N_6. \label{eq:j2j2j2} \end{eqnarray}
Again, it is easy to show that the rank of each factor is the same as the rank of $\text{Dg}[J_2(0),J_2(0),J_2(0)]$. Let \begin{equation} Q_4 = [e_1,e_4,e_3,e_6,e_2,e_5], \label{eq:j2j2j2_q} \end{equation} then
\[Q_4^{-1} N_5 Q_4  = \text{Dg}[J_2(0),J_3(0),J_1(0)],\] which shows that $N_5$ is nilpotent. Let $Q_5 = [e_2,e_4,e_5,e_3,e_1,e_6]$, then 
\[Q_5^{-1} N_6 Q_5  = \text{Dg}[J_1(0),J_1(0),J_4(0)],\] which shows that $N_6$ is nilpotent, concluding the proof. \hfill $\square$ 

Note that the original proof for the case $\text{Dg}[J_k(0),J_2(0)]$ where $k \geq 1$ is odd \cite[Lemma 3]{wu2} is invalid, since the matrix 
\[ \begin{bmatrix} \mathbf{0} & J_2(0) \\ J_k(0) & \mathbf{0} \end{bmatrix} \] is not nilpotent when $k=7$ (the first odd integer greater than 1 for which this matrix is not nilpotent)\footnote{For $k=7$ it can be checked that $A^3 = A^6 = \cdots = A^{3j}$ where $j > 2$ is an integer.}.

\begin{quote}
\begin{cor}
Let $\mathscr{F}$ be an arbitrary field. Any nilpotent matrix $N \in M_n(\mathscr{F})$, with $n \neq 2$, is similar to the product of two nilpotent matrices, where the first row and last column of the first factor (left factor) is the zero vector, and the last row of the second factor (right factor) is either the zero vector or $e_1^T$ (the vector with a one in the first entry and zeros elsewhere).
\label{cor:j0nilpotent2}
\end{cor}
\end{quote}
{\bf Proof.} The previous results present exhaustive configurations of $N$, and factorizations into nilpotent factors for each of these. Now we need to prove that each of these factorizations is of the desired form.

First, suppose $J_k(0)$ with $k \neq 2$ is the last simple Jordan block of $N$ on the diagonal (i.e. this block is not paired with a block $J_2(0)$). If $k=1$ the result is immediate. Suppose $k \geq 2$: now we need to consider the factorizations \eqref{eq:jk0_k_odd} and \eqref{eq:jk0_k_even}. For \eqref{eq:jk0_k_odd} let $Q_1$ be as defined in \eqref{eq:jk0_k_odd_q1}, then it is easy to verify that both factors in 
\[Q_1^{-1} J_k(0) Q_1 = (Q_1^{-1} N_1 Q_1) (Q_1^{-1} N_2 Q_1) \] are of the desired form. In particular $Q_1^{-1} N_2 Q_1$ has the zero vector as its last row. Now for \eqref{eq:jk0_k_even}, let $Q_2$ be as defined in \eqref{eq:jk0_k_even_q2}, then both factors in
\[Q_2^{-1} J_k(0) Q_2 = (Q_2^{-1} N_3 Q_2) (Q_2^{-1} N_4 Q_2) \] are of the desired form. In particular $Q_2^{-1} N_4 Q_2$ has as its last row the vector $e_1^T$. 

Now suppose the last diagonal block of $N$ is $J_2(0)$. Suppose first that we have $N = \text{Dg}[J,J_k(0),J_2(0)]$ where $J$ is some nilpotent matrix in Jordan canonical form. If $k=1$ the result is immediate by the factorization \eqref{eq:nilpotent_jkj2_1}. If $k \geq 2$ and $k$ is even the result is also immediate by the factorization \eqref{eq:nilpotent_jkj2_even}. Suppose $k \geq 2$ and $k$ is odd: we can make use of the factorization as defined in \eqref{eq:nilpotent_jkj2_odd}. By applying the change-of-basis matrix $Q_2$ as defined in \eqref{eq:nilpotent_jkj2_odd_q}, both factors in 
\[Q_2^{-1} \text{Dg}[J_k(0),J_2(0)] Q_2 = (Q_2^{-1} N_3 Q_2) (Q_2^{-1} N_4 Q_2) \] are of the desired form. In particular note that the last row of $Q_2^{-1} N_4 Q_2$ is the zero vector.

Finally suppose $N = \text{Dg}[J,J_2(0),J_2(0),J_2(0)]$. We can make use of the factorization~\eqref{eq:j2j2j2} and the change-of-basis matrix $Q_4$ as defined in \eqref{eq:j2j2j2_q}, then both factors in 
\[ Q_4^{-1} \text{Dg}[J_2(0),J_2(0),J_2(0)] Q_4 = (Q_4^{-1} N_5 Q_4) (Q_4^{-1} N_6 Q_4) \] are of the desired form, in particular note that the last row of $Q_4^{-1} N_6 Q_4$ is the zero vector. \hfill $\square$

On to the main result of this section, as proved by Laffey. I modify the proof to include the converse, for completeness.
\begin{quote}
\begin{thm}[Theorem 1.3 \cite{laffey}]
Let $\mathscr{F}$ be an arbitrary field. The matrix $A \in M_n(\mathscr{F})$ is a product of two nilpotent matrices if and only if it is singular, except when $A \in M_2(\mathscr{F})$ and $A$ is nilpotent and nonzero.
\label{thm:nilpotentProductLaffey}
\end{thm}
\end{quote}
{\bf Proof.} First, suppose $A$ is the product of two nilpotent matrices $N_1,N_2$. Now a nilpotent matrix cannot be full rank, since it has zero as a characteristic value. Furthermore $\text{r}(A) \leq \min(\text{r}(N_1),\text{r}(N_2))$ \cite[Proposition 6.11]{golan}, and it follows that $A$ cannot be full rank, and is therefore singular. 

Now suppose $A$ is singular. If $A$ is nilpotent then the result follows directly by Lemma \ref{lem:j0nilpotent2}. Suppose therefore that $A$ is not nilpotent. By Fitting's lemma \cite[Theorem 5.10]{cullen} $A \approx \text{Dg}[A_0,A_1]$ where $A_0$ is nilpotent and $A_1$ is invertible. Let us consider three mutually exclusive cases in terms of the matrix $A_0$: first suppose $A_0 = J_1(0)$. Let $k=n-1$. Since $A_1$ is similar to $LU$ where $L=(l_{ij})$ is lower triangular and $U=(u_{ij})$ is upper triangular (Proposition \ref{prop:laffeyTriangle}) we have 
\[A = \begin{bmatrix} 0 & 0 & \cdots & & 0 \\
 l_{11} & 0 & & & \vdots \\
 l_{21} & l_{22} & \ddots & & \vdots \\
 \vdots & & \ddots & 0 & 0 \\
 l_{k1} & \cdots & & l_{kk} & 0
 \end{bmatrix} 
 \begin{bmatrix} 0 &  u_{11} & u_{12} & \cdots & u_{1k} \\
 0 & 0 & u_{21} & & \vdots \\
 \vdots &  & \ddots & \ddots & \vdots \\
 \vdots & &  & 0 & u_{kk} \\
 0 & \cdots & &  & 0
 \end{bmatrix}.\]
 Note that each factor on the right-hand side has characteristic polynomial $x^n$ so that it is nilpotent. Notice also that each factor has the same rank as $A$, since each of the triangular blocks in the factors are full rank.
 
Now suppose that $A_0$ is similar to $J_2(0)$. Note that $J_2(0)$ is similar to 
\[\begin{bmatrix} 0 & 1 \\ 0 & 0 \end{bmatrix},\] which is the Jordan form preferred by some texts. Let $k=n-2$, and let
\begin{eqnarray} \nonumber A_2 &=& \begin{bmatrix} 
1 	& 0 		& \cdots 	& 		& 0 		& 1\\
0 	& 0 		& \cdots 	& 		& 0 		& 0\\
0 	& l_{11} 	& 0 		& 		& \vdots	& \vdots \\
0	& l_{21} 	& l_{22} 	& \ddots 	& 		&  \\
\vdots & \vdots 	& 		& \ddots 	& 0 		& 0 \\
-1 	& l_{k1} 	& \cdots 	& 		& l_{kk} 	& -1
 \end{bmatrix}
 \begin{bmatrix} 
0 	& 1 		& 0	 	& \cdots		&  		& 0\\
0 	& 0 		&  u_{11} 	& u_{12} 	& \cdots 	& u_{1k}\\
0 	& 0 		& 0 		& u_{21} 	& 		& \vdots \\
\vdots	& \vdots 	&  		& \ddots 	& \ddots 	& \vdots\\
 &  	& 		&  		& 0 		& u_{kk} \\
0 	& 0 		& \cdots 	& 		& 	0 	& 0
 \end{bmatrix}\\
 \nonumber &=&
 \left [ {\renewcommand{\arraystretch}{1.2}\begin{array}{c|c} 
\begin{matrix} 0 & 1 \\ 0 & 0 \end{matrix} &  \mathbf{0} \\ 
\hline
-E_{(k,2)} & A_1 \end{array}} \right ].
 \end{eqnarray}
By Roth's theorem $A_2$ is similar to $A$ \cite{roth}. Explicitly, 
\[ \begin{bmatrix} I & \mathbf{0} \\ -X & I\end{bmatrix}  \left [ {\renewcommand{\arraystretch}{1.2}\begin{array}{c|c} 
\begin{matrix} 0 & 1 \\ 0 & 0 \end{matrix} &  \mathbf{0} \\ 
\hline
-E_{(k,2)} & A_1 \end{array}} \right ] \begin{bmatrix} I & \mathbf{0} \\ X & I\end{bmatrix} =  \left [ {\renewcommand{\arraystretch}{1.2}\begin{array}{c|c} 
\begin{matrix} 0 & 1 \\ 0 & 0 \end{matrix} &  \mathbf{0} \\ 
\hline
\mathbf{0} & A_1 \end{array}} \right ]\] only if there exists a solution to 
\[-E_{(k,2)} = X \begin{bmatrix} 0 & 1 \\ 0 & 0 \end{bmatrix} - A_1X.\] 
Now since $A_1$ is full rank there exists a vector $x \in \mathscr{F}^k$ such that $A_1x = e_k$. It follows that $X = [\mathbf{0},x]$ is a matrix that will satisfy the requirements, proving that the given factorization is valid.

 Let the first factor in the factorization above be $N_1$ and the second $N_2$. Since $N_2$ is upper triangular with only 0 on the diagonal, it is immediately apparent that it is nilpotent. Now consider $N_1$: Let $P=[e_1,e_1-e_n,e_2,e_3,\ldots,e_{n-1}]$, then
\[ P^{-1} N_1 P = \left [ {\renewcommand{\arraystretch}{1.2}\begin{array}{c|c} 
\begin{matrix} 0 & 0 \\ 1 & 0 \end{matrix} &  \begin{matrix} l_{k1} & l_{k2} & \cdots & \qquad & l_{kk} \\ 
-l_{k1} & -l_{k2} & \cdots & & -l_{kk} \end{matrix} \\ 
\hline
\mathbf{0} & \begin{matrix} 
0 & 0 & \cdots & & 0 \\
l_{11} & 0 &  & & 0 \\ 
l_{21} & l_{22} & \ddots & & \vdots \\
\vdots & & \ddots & 0 & 0 \\
l_{1(k-1)} & \cdots & & l_{(k-1)(k-1)} & 0
\end{matrix} \end{array}} \right ]. \]
Since the determinant of a block triangular matrix is the product of the determinants of its diagonal blocks \cite[Proposition 11.12]{golan}, the characteristic polynomial of $N_1$ is the product of the characteristic polynomials of its diagonal blocks. It is easy to verify that the characteristic polynomial of $N_1$ is therefore $x^2 \cdot x^k = x^n$, confirming that it is nilpotent, and concluding the proof of the result in this case. Again, we can directly observe that each factor has the same rank as $A$.

It remains to prove the result for the case where $A_1$ is of order $k \times k$ where $k \leq n-3$, so that $A_0$ is of order $3 \times 3$ or larger. Note that in the original proof by Laffey \cite{laffey} there is an error in the factorization given at the top of page 99. If the last column of $N_1$ is not zero, then multiplication of the factors result in a block matrix \[ \begin{bmatrix} A_0 & C \\ B & A_1 \end{bmatrix} \] where the matrix $C$ is not necessarily the zero matrix, as indicated in the original proof. Therefore we cannot apply Roth's theorem \cite{roth} to prove similarity of the factorization with the original matrix $A$, as originally claimed. I will now show that with a slight modification of the factorization given in the original proof the result remains valid. For this purpose I make use of Corollary \ref{cor:j0nilpotent2} whereby we can assume $A_0 = N_1N_2$, where $N_1$ is nilpotent and its last column is the zero vector, and $N_2$ is nilpotent and its last row is either the zero vector or $e_1^T$.

Now we have
\[ \begin{bmatrix} A_0 & \mathbf{0} \\ B & A_1\end{bmatrix} =  \left [ {\renewcommand{\arraystretch}{1.2}\begin{array}{c|c} 
N_1 & \mathbf{0} \\
\hline
\begin{matrix} 
0 & \cdots & 0 & l_{11} \\
0 &  & \vdots &  l_{21} \\ 
\vdots &  & & \vdots \\
0 & \cdots & 0 & l_{k1} 
\end{matrix} & \begin{matrix} 
0 & 0 & \cdots & 0 \\
l_{22} & 0 &  &  0 \\ 
\vdots & \ddots & \ddots & \vdots \\
l_{k2} & \cdots & l_{kk} & 0 
\end{matrix}
\end{array}} \right ]
\left [ {\renewcommand{\arraystretch}{1.2}\begin{array}{c|c} 
N_2 &\begin{matrix} 
0 & \cdots &   & 0  \\
\vdots &  &    & \vdots \\ 
0 & \cdots  &  & 0 \\
u_{11} & u_{12} &  \cdots & u_{1k} 
\end{matrix}  \\
\hline
\mathbf{0} &
 \begin{matrix}
0  & u_{22} &  \cdots & u_{2k} \\
\vdots & \ddots & \ddots    &  \vdots \\ 
 \quad \, &  & 0 & u_{kk} \\
0 & \cdots & 0 & 0 
\end{matrix}
\end{array}} \right ].\]
If the last row of $N_2$ is zero then $B=\mathbf{0}$ and it follows that $A$ is similar to the given factorization, which is the desired result. Suppose the last row of $N_2$ is $e_1^T$, then   
\[ B= \begin{bmatrix} l_{11} &0 & \cdots & 0 \\
l_{21} &0 & \cdots & 0 \\
\vdots  & \vdots & & \vdots \\
l_{k1}&0 & \cdots & 0 \end{bmatrix}. \] Now we may use Roth's theorem to prove the result \cite{roth}. Explicitly, 
\[ \begin{bmatrix} I & \mathbf{0} \\ -X & I\end{bmatrix} \begin{bmatrix} A_0 & \mathbf{0} \\ B & A_1\end{bmatrix} \begin{bmatrix} I & \mathbf{0} \\ X & I\end{bmatrix} = \begin{bmatrix} A_0 & \mathbf{0} \\ 0 & A_1\end{bmatrix}\] only if there exists a solution to $B = XA_0 - A_1X$. Let $X$ be the matrix with $-1/u_{11}$ in entry $(1,1)$ and zero elsewhere. Then $XA_0 = \mathbf{0}$ since the first row of $A_0$ is zero (by corollary \ref{cor:j0nilpotent2} the first row of $N_1$ is the zero vector), and $-A_1X = B$, which yields the desired result. This proves that $A$ is similar to the given factorization.

It remains to prove that the factors are nilpotent. But now, as mentioned before, since the determinant of a block triangular matrix is the product of the determinants of its diagonal blocks \cite[Proposition 11.12]{golan}, the characteristic polynomial of each factor is the product of the characteristic polynomials of its diagonal blocks. Notice that for both these factors the diagonal blocks are all nilpotent, and therefore both factors are nilpotent, which concludes the proof: $A$ is the product of two nilpotent matrices as desired. \hfill $\square$

Notice again in the last case above, since the triangular blocks in each factorization are full rank, and by the results of Wu, each of the nilpotent matrices $N_1$ and $N_2$ has the same rank as $A_0$, it follows that each factor in this factorization has the same rank as $A$. 

The results in this section does not give any indication whether it is possible to vary the ranks of the nilpotent factors, but there are results available on \emph{Factorization of a singular matrix into nilpotent matrices with prescribed ranks} \cite{botha4}.

\subsection{Sourour's proof}
I now present Sourour's proof of the main result. Sourour's proof is divided into two parts, the first part proves the result for square-zero matrices and the second all matrices that are not square-zero. The second part consists of an inductive proof which requires Proposition \ref{prop:sourInductSingular}. The inductive step is however not valid for one exceptional case: every nonzero nilpotent matrix in $M_2(\mathscr{F})$ is square-zero since such a matrix cannot have order of nilpotency higher than two. Observe that every nonzero $2 \times 2$ nilpotent matrix is similar to 
\[\begin{bmatrix} 0 & 0 \\ 1 & 0 \end{bmatrix}.\] This matrix is square-zero, and furthermore by Theorem \ref{thm:botha2FactorSZ2} (page \pageref{thm:botha2FactorSZ2}) it is not a product of two square-zero matrices. It follows that a nonzero nilpotent matrix in $M_2(\mathscr{F})$ is not a product of two nilpotent matrices. The following lemma is required to prove the result of Theorem \ref{thm:nilpotentSourour} in this exceptional case:
\begin{quote}
\begin{lem}[Lemma 2 \cite{sour}]
Let $\mathscr{F}$ be an arbitrary field. If the characteristic equation of $A \in M_3(\mathscr{F})$ is $x^2(x-\lambda)$ then $A$ is a product of two nilpotent matrices, each of the same rank as $A$.
\label{lem:nilpotentSourour}
\end{lem}
\end{quote}
{\bf Proof.} If the minimum polynomial of $A$ is $x$, then $A$ is the zero matrix and the result is trivially true. If the minimum polynomial of $A$ is $x^2$ then $A$ is similar to
\[ \text{Dg}[J_1(0),J_2(0)] = [e_3,\mathbf{0},\mathbf{0}] [\mathbf{0},e_1,\mathbf{0}]. \]
If the minimum polynomial of $A$ is $x^3$ then $A$ is similar to
\[ J_3(0) =  [e_3,\mathbf{0},e_2][e_3,e_1,\mathbf{0}].\]
If the minimum polynomial of $A$ is $x(x-\lambda)$ ($\lambda \neq 0$) then $A$ is similar to
\[ \text{Dg}[\lambda,J_1(0),J_1(0)] = \lambda [\mathbf{0},\mathbf{0},e_1] [e_3,\mathbf{0},\mathbf{0}] .\]
Finally, if the minimum polynomial of $A$ is $x^2(x-\lambda)$ ($\lambda \neq 0$) then $A$ is similar to
\[ \text{Dg}[\lambda,J_2(0)] = \begin{bmatrix} \lambda & -1 & -\lambda \\
0 & 0 & 0 \\
\lambda & -2 & -\lambda \end{bmatrix}
 \begin{bmatrix} 1 & -1/\lambda & 0 \\
\lambda & -1 & 0 \\
-1 & 0 & 0 \end{bmatrix}. \]
It is easy to verify that the factorizations of each case above is valid, and for every case each factor is nilpotent and of the same rank as $A$. \hfill $\square$  

The statement of the main result is essentially the same as the statement of Theorem~\ref{thm:nilpotentSourour}, but in this case I will only present the sufficient part, and the ranks of the factors are explicitly stated in the proof (an argument based on rank forms the core of the proof by contradiction in the last part of the induction part).
\begin{quote}
\begin{thm}[Theorem \cite{sour}]
Let $\mathscr{F}$ be an arbitrary field. If the matrix $A \in M_n(\mathscr{F})$ is singular, but not a nonzero nilpotent matrix of order $2 \times 2$, then it is a product of two nilpotent matrices, each of which has the same rank as $A$.
\label{thm:nilpotentSourour}
\end{thm}
\end{quote}
{\bf Proof.} Let $A \in M_n(\mathscr{F})$ be square-zero, and $n \neq 2$. If $A$ is the zero matrix the result is trivially true, as $A=\mathbf{0}^2$. If $A$ is nonzero then $A$ is nilpotent and the result follows easily by Lemma \ref{lem:j0nilpotent2}. For completeness I also present here the proof for the case of square-zero $A$ due to Sourour, which is an elegant proof independent from Lemma \ref{lem:j0nilpotent2}. Suppose $A$ is square-zero with rank $1$: then $A$ is similar to a direct sum of a $3 \times 3$ nilpotent matrix with a zero matrix (notice that we must have $n \geq 3$, and if $n=3$ the zero matrix vanishes). So in this case the desired factorization of the nilpotent block follows directly from Lemma \ref{lem:nilpotentSourour}, and if a zero block is present its factorization is trivial, so that the desired result is achieved.

Now suppose the rank of $A$ is $k \geq 2$. Then it is easy to verify that $A$ is similar to a direct sum of a zero matrix of order $n-2k$ and the $2k \times 2k$ matrix 
\[\begin{bmatrix} \mathbf{0} & I_k \\ \mathbf{0} & \mathbf{0} \end{bmatrix}.\] To see this, notice that for given $n$ there is a single similarity class for a square-zero matrix of order $n$ and rank $k$, because every matrix in this class must be similar to a matrix which consists of a direct sum of $k$ blocks $J_2(0)$ and a zero block of order $(n-2k) \times (n-2k)$. Therefore since 
\[\text{Dg}\left [\mathbf{0}_{(n-2k)\times (n-2k)},\begin{bmatrix} \mathbf{0} & I_k \\ \mathbf{0} & \mathbf{0} \end{bmatrix}\right ]\] is square-zero and rank $k$ it is similar to any other square-zero matrix of order $n \times n$ and rank $k$, and it will be sufficient to prove the result for the matrix 
\[\begin{bmatrix} \mathbf{0} & I_k \\ \mathbf{0} & \mathbf{0} \end{bmatrix}.\] Notice that 
\[\begin{bmatrix} \mathbf{0} & I \\ \mathbf{0} & \mathbf{0} \end{bmatrix} = \begin{bmatrix} J_k(0)^T & J_k(0) \\ \mathbf{0} & \mathbf{0} \end{bmatrix} \begin{bmatrix} \mathbf{0} & E_{(2,1)} \\ \mathbf{0} & J_k(0)^T \end{bmatrix}. \] Each of these factors has rank $k$ and is nilpotent, as desired, which concludes the proof for the case where $A$ is square-zero.

Suppose $A$ is not square-zero (so $A$ is also not the zero matrix). I proceed by induction on the order $n$ of $A$. Since $A$ is a nonzero singular matrix, the base step of the inductive proof is $n=2$. By Proposition \ref{prop:sourInductSingular} a singular matrix $A \in M_2(\mathscr{F})$ that is not square-zero is similar to 
\[\begin{bmatrix} \lambda & 0 \\ 0 & 0 \end{bmatrix} = \begin{bmatrix} 0 & \lambda \\ 0 & 0 \end{bmatrix}\begin{bmatrix} 0 & 0 \\ 1 & 0 \end{bmatrix},\] which proves that $A$ is a product of two nilpotent matrices. 

Now suppose that every singular matrix $A \in M_{n-1}(\mathscr{F})$ that is not square-zero is a product of two nilpotent matrices of the same rank as $A$. We have to prove that every singular matrix $A \in M_n(\mathscr{F})$ that is not square-zero is the product of two nilpotent matrices with the same rank as $A$. Now by Proposition \ref{prop:sourInductSingular} $A$ is similar to
\[\begin{bmatrix} \lambda & c \\ b & D \end{bmatrix},\] with $\lambda \in \mathscr{F}$, $D \in M_{(n-1)}(\mathscr{F})$, $\text{r}(D) = \text{r}(A) - 1$, $b \in \text{R}(D)$ and $c \in \text{row}(D)$. Since the rank of $D$ is one less than the rank of $A$, the matrix $D$ is also singular. Now if $D$ is not square-zero it follows by the induction hypothesis that $D=N_1N_2$ where $N_1$ and $N_2$ are nilpotent and have the same rank as $D$. If $D$ is square-zero the same result holds by the argument at the beginning of the proof, except if $D$ is nonzero, square-zero and in $M_2(\mathscr{F})$: we will deal with this case separately. Now since $c \in \text{row}(D)$ and $b \in \text{R}(D)$ there exist $b_0, c_0^T \in \mathscr{F}^{n-1}$ such that $c_0D = c$ and $Db_0 = b$. Now I propose that we can find $n_1, n_2^T \in \mathscr{F}^{n-1}$ such that 
\begin{equation} \begin{bmatrix} \lambda & c \\ b & D \end{bmatrix} = \begin{bmatrix} 0 & c_0N_1 + n_2 \\ \mathbf{0} & N_1 \end{bmatrix}\begin{bmatrix} 0 & \mathbf{0} \\ N_2b_0 + n_1 & N_2 \end{bmatrix}. \label{sourMainInduct}\end{equation}
The equation above imposes the following constraints in terms of the anti-diagonal blocks of the matrix on the left-hand side:
\[c = (c_0N_1+n_2)N_2 = c_0D + n_2N_2 = c + n_2N_2 \text{, and}\]
\[b = N_1(N_2b_0+n_1) = Db_0 + N_1n_1 = b + N_1n_1,\] and it follows that we must have 
\begin{equation} N_1n_1 = \mathbf{0} \text{ and } n_2N_2 = \mathbf{0}.\label{sourarg1} \end{equation}
Furthermore the following constraint is imposed in terms of the block in entry $(1,1)$:
\[\lambda = (c_0N_1+n_2)(N_2b_0 + n_1) = c_0Db_0 + c_0N_1n_1 + n_2N_2b_0 + n_2n_1.\]
By the requirement \eqref{sourarg1} this equation reduces to
\begin{equation} \lambda - c_0Db_0 = n_2n_1.\label{sourarg2} \end{equation} 
Now $n_1$ and $n_2$ must both be nonzero, for else each of the factors in \eqref{sourMainInduct} have rank less than $A$. Note that $\lambda - c_0Db_0 \neq 0$ since otherwise $\lambda = c_0Db_0$ and so
\[A \approx \begin{bmatrix} \lambda & c \\ b & D \end{bmatrix} = \begin{bmatrix} 1 & c_0 \\ \mathbf{0} & I\end{bmatrix}\begin{bmatrix} 0 & \mathbf{0} \\ Db_0 & D \end{bmatrix},\] from which it follows that the rank of $A$ is the same as the rank of $D$, which is a contradiction.

Essentially the requirement \eqref{sourarg2} then reduces to the requirement that $n_2n_1 \neq 0$, since supposing $n_2n_1 = z \neq 0$, we can replace $n_2$ with \[\frac{(\lambda - c_0Db_0)}{z}n_2,\] ensuring that \eqref{sourarg2} is satisfied. This vector will also satisfy \eqref{sourarg1} if $n_2$ satisfied \eqref{sourarg1} since it is a scalar multiple of $n_2$. 

Thus we require that $n_2n_1 \neq 0$, and $n_1 \in \text{N}(N_1)$ and $n_2 \in \text{R}(N_2)^{\perp_f}$ where $f: \mathscr{F}^n \times \mathscr{F}^n \rightarrow \mathscr{F}$ is defined by $f(v,w)=v^Tw$. \footnote{It is easy to see that the function $f$ is a non-degenerate symmetric bilinear form, and $n_2$ is in the $f$-orthogonal complement of $\text{R}(N_2)$ \cite[p.453 - 458]{golan}. See Proposition \ref{prop:matrixDivision}, the proof of $b \Rightarrow c$ for a complete discussion.}

I proceed with a proof by contradiction. Suppose that we cannot find $n_1, n_2$ with the required properties. Certainly, since $N_1, N_2$ are nilpotent and hence singular, there exist nonzero vectors $x_1, x_2^T \in \mathscr{F}^n$ such that $N_1x_1 = \mathbf{0}$ and $x_2N_2 = \mathbf{0}$. If $x_2x_1=0$ for all $x_1, x_2^T \in \mathscr{F}^n$ we conclude that $\text{R}(N_2)^{\perp_f} \subseteq \text{N}(N_1)^{\perp_f}$, so that $\text{N}(N_1) \subseteq \text{R}(N_2)$ \cite[Proposition 20.3]{golan}. So the null space of $D=N_1N_2$ will consist of each vector in the null space of $N_2$ as well as each vector $N_2 u$ which is in the null space of $N_1$. Since $N_1$ is singular there is at least one such a vector which is nonzero. It follows that $\text{n}(D) = \text{n}(N_1) + \text{n}(N_2)$ so that the rank of $D$ is less than the rank of $N_1$ and $N_2$ which is a contradiction. Therefore we can find $n_1, n_2$ with the required properties, and the factorization \eqref{sourMainInduct} is valid. Since each of the factors in \eqref{sourMainInduct} are block triangular matrices with nilpotent diagonal blocks, it follows that $A$ is a product of nilpotent matrices as required \cite[Proposition 11.12]{golan}.

It remains to prove the exceptional case, where the matrix $D$ is a $2 \times 2$ nonzero nilpotent matrix. In this case we have 
\[A \approx \begin{bmatrix} \lambda & c_1 & 0 \\
0 & 0 & 0 \\
b_1 & 1 & 0 \end{bmatrix},\] so that the characteristic polynomial of $A$ is $x^2(x-\lambda)$ and therefore we can apply Lemma \ref{lem:nilpotentSourour}, and it follows that $A$ is a product of two nilpotent matrices, each with the same rank as $A$.
\hfill $\square$
  
\section{Sums of Nilpotent Matrices}
The first result on sums of nilpotent matrices appears in Theorem 3.6 of the article by Wang and Wu \cite{ww}, even though it is not the main result of the theorem, or mentioned explicitly in the statement of the theorem. Later it was explicitly mentioned by Wu \cite[Proposition 5.17]{wu}. Even though the context of the result is limited to the complex numbers, the result is valid over an arbitrary field since the proof relies on Proposition \ref{prop:trace}. This proposition is however limited to non-scalar matrices, so that a further result is required to complete the question regarding sums of nilpotent matrices. Such an extension was provided by Breaz and C{\u a}lug{\u a}ƒreanu \cite[Proposition 3]{breaz}. I present this result here, with the proof modified to show only a self-contained version for the special case of matrices over fields (the original result requires a special case of a result for matrices over a commutative ring \cite[Thoerem 2]{breaz}).

\begin{quote}
\begin{thm}[Proposition 3 \cite{breaz}]
Let $\mathscr{F}$ be an arbitrary field. The matrix $A \in M_n(\mathscr{F})$ is a sum of nilpotent matrices if and only if $\text{tr}(A)=0$. Furthermore, 
\begin{enumerate}
\item if $\text{tr}(A) =0$ and $A = \lambda I_n$, $\lambda \neq 0$, then $A$ is a sum of three nilpotent matrices, which is the best possible. 
\item if $\text{tr}(A) = 0$ and $A$ is the zero matrix, or non-scalar, then $A$ is a sum of two nilpotent matrices.
\end{enumerate}
\label{thm:nilpotentSums}
\end{thm}
\end{quote}
{\bf Proof.} First suppose that $A$ is a sum of nilpotent matrices. Since similar matrices have the same trace, and a nilpotent matrix is similar to a matrix in Jordan canonical form with only zero on the diagonal, a nilpotent matrix has trace zero. Therefore the sum of nilpotent matrices also has trace zero, and hence it follows that $\text{tr}(A)=0$.

For the converse, suppose first that $\text{tr}(A) =0$ and $A = \lambda I_n$, $\lambda \neq 0$. Combining these requirements we must have $n \lambda = 0$, i.e. $n$ is a multiple of the characteristic of the underlying field. Notice that  
\[ \lambda I_n = -\lambda \begin{bmatrix*}
0 & 0 & \cdots   & 0 \\
1 & 0 &  & 0 \\
\vdots & \ddots & \ddots  & \vdots \\
1 & \cdots & 1  & 0
\end{bmatrix*} + 
\lambda \begin{bmatrix*}
1 & 1 & \cdots &   1 \\
1 & 1 &  & \vdots \\
\vdots &  & \ddots  & \\
1 & \cdots & & 1
\end{bmatrix*}
-\lambda \begin{bmatrix*}
0 & 1 & \cdots   & 1 \\
0 & 0 & \ddots & \vdots \\
\vdots &  & \ddots  & 1 \\
0 & 0 & \cdots & 0
\end{bmatrix*}.
\] To be clear, $\lambda I_n = -\lambda T_1 + \lambda B - \lambda T_2$, where $T_1$ is a matrix that is lower triangular with only zero on the diagonal, and only ones in each entry below the diagonal, $T_2$ is upper triangular with only zero on the diagonal and ones in each entry above the diagonal. The matrix $B$ is the matrix with ones in each of its entries. Now it is easy to see that each of these summands are nilpotent. Explicitly, since the only characteristic value of $T_1$ and $T_2$ is zero, both $T_1$ and $T_2$ are nilpotent. Furthermore, since by hypothesis $n \lambda = 0$, and therefore $n \lambda^2=0$, it follows that $\text{row}_i(B) \text{col}_j(B) = 0$ for any $i,j \in \{1,2, \ldots, n\}$, so that $B$ is square-zero. It follows that $A$ is a sum of three nilpotent matrices.

For this case it remains to prove that $A$ is not a sum of two nilpotent matrices. Suppose to the contrary that $A = N_1 + N_2$ where $N_1$ and $N_2$ are both nilpotent. Now consider $N_2 = A - N_1$: let $P$ be a change of basis matrix such that $P^{-1}N_1P$ is in Jordan canonical form. Then
\[ P^{-1}N_2P = P^{-1}(A-N_1)P = P^{-1}(\lambda I_n)P-P^{-1}N_1P = \lambda I_n - J\] where $J$ is the Jordan form of $N_1$ and therefore has only zero on the diagonal. It follows that $N_2$ is similar to a matrix in Jordan canonical form which has only $\lambda$ on the diagonal, and therefore it cannot be nilpotent, which is a contradiction. So $A$ is not a sum of two nilpotent matrices.

Finally, let $A$ be a matrix such that $\text{tr}(A) = 0$ and $A= \mathbf{0}$ or $A$ is non-scalar. If $A=\mathbf{0}$ the result is obvious, so suppose $A$ is non-scalar. By Proposition \ref{prop:trace}  $A$ is similar to a matrix which has only zero on the diagonal, let this matrix be $A_1=(a_{ij})$. Then
\[A_1 = \begin{bmatrix*}[c] 0 & 0 & \cdots & 0 \\
a_{21} & 0 & & 0 \\
\vdots & \ddots & \ddots & \vdots \\
a_{n1} & \cdots & a_{n(n-1)}&0\end{bmatrix*} + 
\begin{bmatrix}
0 & a_{12} & \cdots & a_{1n}\\
0 & 0 & \ddots & \vdots \\
\vdots &  & \ddots & a_{(n-1)n} \\
0 & 0 & \cdots & 0
\end{bmatrix}, \] and since each of the summands above are triangular matrices with only zero as characteristic value, $A$ is the sum of two nilpotent matrices, which completes the proof. \hfill $\square$
\section{Concluding Remarks and Open Problems}
Which matrices are sums or products of nilpotent matrices? The answers presented here satisfy these questions with simple and intuitive solutions. For the special case of square-zero matrices the answers to similar questions are  in many cases more intricate and less perceptible. Additionally, the proofs presented here are of an adequately constructive nature to facilitate construction of the factors or summands involved, whenever it is possible to write a matrix as a sum or product of nilpotent matrices.

Can anything be added to these results? Certainly, in terms of factorization an interesting question remains:
\begin{quote} Let $\mathscr{F}$ be a field. Which matrices in $M_{m \times n}(\mathscr{F})$ have a nilpotent divisor or quotient? \end{quote}
In other words, when can $G,F \in M_{m \times n}(\mathscr{F})$ be written as $G=NF$ where $N \in M_{m \times m}$ is nilpotent, or $G=FN$ where $N \in M_{n \times n}(\mathscr{F})$ is nilpotent? For the special case of square-zero matrices this question was answered by Botha \cite{botha2}, but the more general question regarding nilpotent matrices remain open.

Various other questions could be posed in terms of specific subsets of matrices, for example in terms of commuting nilpotent matrices \cite{bukovsek}, however the most general questions have been answered by the content presented here. This concludes my discussion on the broader context of sums and products of nilpotent matrices.
 
\chapter{Products of Square-Zero Matrices}

Results on products of square-zero matrices can be presented as a unified theory when approached from the framework of matrix division. In this chapter I will first present results as they became available chronologically, and then I will present a unified development of results based on matrix division with a square-zero divisor. I start with preliminary results, which constitute much of the underlying theory to matrix division with a square-zero quotient.
\section{Preliminary Results: Matrix Division \label{prelim}}
Matrix division can be defined as follows (\cite[p.71]{botha1}): Let $G \in M_{m \times n}(\mathscr{F})$, $F \in M_{k \times n}(\mathscr{F})$ where $\mathscr{F}$ is a field. Then $F$ divides $G$ on the right if there exists $H \in M_{m \times k}(\mathscr{F})$ such that $G=HF$. In such a case $F$ is called a right divisor (or right factor) of $G$ and $H$ a right quotient of $G$ and $F$. As with integer division, divisibility implies a zero remainder. A similar definition can be given for left division.

\subsection{Necessary and sufficient conditions for matrix division}
I will now discuss, in detail, necessary and sufficient conditions for a matrix $F \in M_{k \times n}(\mathscr{F})$ to divide a matrix $G \in M_{m \times n}(\mathscr{F})$ on the right, following a proposition and proof by Botha \cite[p.72]{botha1}. 
\begin{quote}
\begin{prop}[Proposition 1 \cite{botha1}] Let $G \in M_{m \times n}(\mathscr{F})$ and $F \in M_{k \times n}(\mathscr{F})$ where $\mathscr{F}$ is a field. The following statements are equivalent:
\begin{enumerate}[(a)]
\item There exists a matrix $H \in M_{m \times k}(\mathscr{F})$ such that $G=HF$.
\item $\text{row}(G) \subseteq \text{row}(F).$
\item $\text{N}(F) \subseteq \text{N}(G).$
\item $\text{r}(F)=\text{r}\left ( \begin{bmatrix} G \\ F \end{bmatrix} \right ).$
\end{enumerate}
\begin{equation} \label{prop1pt1} \end{equation}
If any of the statements above are true, then there exists a quotient $H$ of each rank such that
\begin{equation}\text{r}(G)\leq \text{r}(H ) \leq \text{r}(G) + \min(\text{n}(G^T),\text{n}(F^T)). \label{prop1pt2}\end{equation}
\label{prop:matrixDivision}\end{prop}
\end{quote}
{\bf Proof.} For the proof of \eqref{prop1pt1} I will provide the details of the equivalence of (a) and (b), and the equivalence of (b) and (d) which are not explicitly provided in the original proof. Then I will provide an alternate proof for the equivalence of (b) and (c) using the theory of non-degenerate symmetric bilinear forms.

$a \Rightarrow b$. Suppose there exists a matrix $H \in M_{m \times k}(\mathscr{F})$ such that $G=HF$. Now following a proof in Cullen \cite[Theorem 2.5, p.77]{cullen}, we have $v^T \in$ row$(G)$ if and only if \[v^T = \sum_{i=1}^m c_i \text{row}_i (G)= \begin{bmatrix} c_1 & c_2 & \cdots & c_m \end{bmatrix} G \in \{c^T G: c \in \mathscr{F}^m\}, \]
so that row$(G) = \{c^T G: c \in \mathscr{F}^m\}$. It follows that 
\begin{eqnarray} \nonumber \text{row}(G)&=&\text{row}(HF)\\
\nonumber &=& \{c^T (HF): c \in \mathscr{F}^m\} \\
\nonumber &=&\{(c^T H)F: c \in \mathscr{F}^m\} \\
\nonumber  &\subseteq& \{d^T F: d \in \mathscr{F}^k\} \\
\nonumber &=& \text{row}(F).\end{eqnarray}

$b \Rightarrow a$. Suppose row$(G) \subseteq$ row$(F)$. Then we can write each row of $G$ as \[\text{row}_i(G)=\sum_{j=1}^k h_{ij}\text{row}_j(F), \] where $h_{ij} \in \mathscr{F}$, for $i \in \{1,2,\ldots,m\}$. Let $h_{ij}$ be the entry in row $i$ and column $j$ of a matrix $H$. Then $H=(h_{ij}) \in M_{m \times k}(\mathscr{F})$, and \[ G= \begin{bmatrix} \text{row}_1(G) \\ \text{row}_2(G) \\ \vdots \\ \text{row}_m(G) \end{bmatrix} = (h_{ij}) \begin{bmatrix} \text{row}_1(F) \\ \text{row}_2(F) \\ \vdots \\ \text{row}_k(F) \end{bmatrix} = HF. \]

$b \Rightarrow c$. Consider the function $f: \mathscr{F}^n \times \mathscr{F}^n \rightarrow \mathscr{F}$, defined by $f(v,w)=v^Tw$. It is easy to see that:
\begin{itemize}
\item $f_{w_0}:\mathscr{F}^n \rightarrow \mathscr{F}$ defined by $f_{w_0}(v)=f(v,w_0)$ for any given $w_0 \in \mathscr{F}^n$ is a linear functional and therefore belongs to the dual space of $\mathscr{F}^n$.\footnote{A linear functional is a linear transformation from a vector space $\mathcal{V}$ (over $\mathscr{F}$) to $\mathscr{F}$ (considered as a vector space of dimension 1). The vector space of all linear functionals from $\mathcal{V}$ to $\mathscr{F}$ is called the dual space of $\mathcal{V}$. \cite[p.317]{golan}} 
\item Similarly $f_{v_0}:\mathscr{F}^n \rightarrow \mathscr{F}$ defined by $f_{v_0}(w)=f(v_0,w)$ for any given $v_0 \in \mathscr{F}^n$ is a linear functional and therefore belongs to the dual space of $\mathscr{F}^n$. 
\end{itemize} 
By the properties above, $f$ is a bilinear form, and furthermore it is a symmetric bilinear form since $f(v,w)=f(w,v)$ for any $v,w \in \mathscr{F}^n$. A symmetric bilinear form $g:\mathcal{V} \times \mathcal{V} \rightarrow \mathscr{F}$ is non-degenerate if and only if, for every nonzero $v \in \mathcal{V}$ there exists $w \in \mathcal{V}$ so that $g(v,w) \neq 0$ \cite[p.455]{golan}. Let $v=(v_1,v_2,\ldots,v_n)^T \in \mathscr{F}^n$ be nonzero. Suppose $i \in \{1,2,\ldots,n\}$ is the least integer so that $v_i$ is nonzero. Let $w=(0,\ldots,0,v_i,0,\ldots,0)^T \in \mathscr{F}^n$ (with $v_i$ in entry $i$ of $w$). Then \[f(v,w)=v^Tw=v_i^2 \neq 0,\] since a field contains no divisors of zero. This proves that $f$ is non-degenerate for any field $\mathscr{F}$.

Now suppose row$(G) \subseteq$ row$(F)$. The right $f$-orthogonal complement \cite[p.458]{golan} of row$(G)$ is \begin{eqnarray}
\nonumber \text{row}(G)^{\perp_f}&=&\{w \in \mathscr{F}^n:f(v,w)=0 \text{ for all } v \in \text{row}(G) \} \\
\nonumber &=& \{w \in \mathscr{F}^n:Gw=\mathbf{0} \} \\
 &=& \text{N}(G). \label{eqp12}
\end{eqnarray}      
Similarly $\text{row}(F)^{\perp_f}=$ N$(F)$. Now if (b) holds it follows directly from Proposition 20.3, Golan \cite[p.458]{golan}, that \[\text{N}(F) \subseteq \text{N}(G).\]

$c \Rightarrow b$. Now suppose $\text{N}(F) \subseteq \text{N}(G)$. For any $v \in$ row$(G)$, by \eqref{eqp12} and the symmetry of $f$, $v \in \{v \in \mathscr{F}^n: f(w,v)=0$ for all $w \in$ N$(G)\}$, implying that row$(G) \subseteq$ N$(G)^{\perp_f}$. Furthermore, since $f$ is non-degenerate, it follows from Proposition 20.4, Golan \cite[p.459]{golan} that\footnote{The subspace $(\mathscr{F}^n)^{\perp_f}$ is trivial when $f$ is non-degenerate, and therefore satisfies the conditions of Proposition 20.4, Golan. See p.458 Golan for details \cite{golan}.} \[\text{n}(G)+\dim(\text{N}(G)^{\perp_f})=n.\] Now since we also have n$(G)+$ r$(G)=n$ it follows that row$(G)=$ N$(G)^{\perp_f}$. By a similar argument row$(F)=$ N$(F)^{\perp_f}$. Now again, given (c) holds, it follows directly from Proposition 20.3, Golan \cite[p.458]{golan}:
\[\text{row}(G) \subseteq \text{row}(F).\]

$b \Leftrightarrow d$. \begin{eqnarray} \nonumber \text{row}(G) \subseteq \text{row}(F) &\Leftrightarrow& \text{row}(F) = \text{row} \left ( \begin{bmatrix} G \\ F \end{bmatrix} \right ) \\
\nonumber &\Leftrightarrow& \text{r}(F)= \text{r} \left ( \begin{bmatrix} G \\ F \end{bmatrix} \right ).
\end{eqnarray}

The proof of \eqref{prop1pt2} by Botha \cite{botha1} consists of existential quantification, followed by verification of universality. Existential quantification in this case amounts to establishing the existence of a class of matrices satisfying \eqref{prop1pt1} and with ranks varying from the lower to upper bound as stated in \eqref{prop1pt2}, and verifying its universality of showing that there can be no matrix $H$ satisfying \eqref{prop1pt1} with rank outside of \eqref{prop1pt2}. 

The existence part is constructive, in that a matrix $H$ satisfying $G=HF$ is constructed from two components $H_1,H_2$ with disjoint domains and codomains, and known ranks. Essentially $H$ is then equivalent to a generalized diagonal block matrix where the diagonal blocks are matrices equivalent to $H_1$ and $H_2$ and therefore the rank of $H$ is the sum of the rank of $H_1$ and $H_2$. For $H_1, H_2$ to have known ranks they are defined in terms of the actions of $H$ on vectors in the range of $F$ (defined as the transformation $H_1$), and vectors outside of the range of $F$ (the transformation~$H_2$). So decomposing the domain and codomain of $H$:
\begin{itemize}
\item For the domain, let $\mathcal{V} \subseteq \mathscr{F}^k$ be such that $\mathscr{F}^k=$ R$(F)\oplus \mathcal{V}$.
\item For the codomain, let $\mathcal{W} \subseteq \mathscr{F}^m$ be such that $\mathscr{F}^m=$ R$(G)\oplus \mathcal{W}$.
\end{itemize}
Then let $H \restriction$ R$(F)=H_1:$ R$(F) \rightarrow$ R$(G)$ be defined as the linear transformation sending $Fv$ to $Gv$ for each $v \in \mathscr{F}^n$. The transformation $H$ must fulfill this requirement in order for $G=HF$ to be true, and furthermore, in satisfying the conditions of \eqref{prop1pt1} we ensure that $H_1$ is well-defined and linear. I will now verify this statement. Suppose $v_1,v_2 \in$ R$(F)$ and $w_1,w_2 \in \mathscr{F}^n$ so that \[Gw_1=H_1v_1 \neq H_1v_2=Gw_2 \text{ where } v_1=Fw_1 \text{ and } v_2 = Fw_2.\] It follows that $G(w_1-w_2) \neq \mathbf{0}$, so that $(w_1 - w_2) \notin$ N$(G)$. Now since N$(F) \subseteq$ N$(G)$, we must then have $(w_1 - w_2) \notin$ N$(F)$, ensuring that $v_1=F(w_1) \neq F(w_2)=v_2$ so that $H_1$ is functional. To verify the linearity of $H_1$, let $w_1,w_2 \in \mathscr{F}^n$, $v_1=Fw_1$, $v_2=Fw_2$ and $c \in \mathscr{F}$. Then by the linearity of $F$ and $G$
\begin{eqnarray} 
\nonumber H_1(v_1+cv_2) &=& H_1(F(w_1 + cw_2)) \\
\nonumber &=& G(w_1+cw_2) \\
\nonumber &=& Gw_1+cGw_2 \\
\nonumber &=& H_1(v_1)+cH_1(v_2), \end{eqnarray} proving that $H_1$ is linear. 

It remains to define $H \restriction \mathcal{V}=H_2: \mathcal{V} \rightarrow \mathcal{W}$. The only restriction in this case is that $H_2$ must be linear, thereby ensuring the linearity of $H$. As the transformation of vectors in $\mathcal{V}$ does not affect the equation $G=HF$, $H \restriction \mathcal{V}=H_2$ can be an arbitrary linear transformation. 

Suppose now we have bases $\alpha_1,\alpha_2$ for R$(F), \mathcal{V}$ and bases $\beta_1, \beta_2$ for R$(G), \mathcal{W}$. Let the (matrix) represention of $H$ in terms of the bases $\alpha=\alpha_1 \cup \alpha_2$ and $\beta=\beta_1 \cup \beta_2$ be $\Phi_{\alpha\beta}(H)$, of $H_1$ in terms of $\alpha_1, \beta_1$ be $\Phi_{\alpha_1\beta_1}(H_1)$, and of $H_2$ in terms of $\alpha_2, \beta_2$ be $\Phi_{\alpha_2\beta_2}(H_2)$. Then, by the preceding discussion, \[H\sim \Phi_{\alpha\beta}(H)=\begin{bmatrix}  \Phi_{\alpha_1\beta_1}(H_1) & \mathbf{0} \\ \mathbf{0} &  \Phi_{\alpha_2\beta_2}(H_2) \end{bmatrix}.\] The rank of $H$ is therefore the sum of the rank of $H_1$ and $H_2$. So then since r$(G) \leq$ r$(F)$ we must have r$(H_1) =$ r$(G)$. Since $H_2$ is arbitrary, its rank can vary between zero and the minimum of its row size and column size, that is $\min(\dim(\mathcal{W}),\dim(\mathcal{V}))$. Now 
\[\min(\dim(\mathcal{W}),\dim(\mathcal{V}))=\min(m-\text{r}(G),k-\text{r}(F))=\min(\text{n}(G^T),\text{n}(F^T))\]
 by the rank-nullity theorem \cite[Theorem 2.11, p.87]{cullen}. This establishes that there exists a class of matrices with ranks varying from the lower to upper bound in \eqref{prop1pt2}.

Verifying the universality of \eqref{prop1pt2} consists of proving that the rank of any matrix $H$ satisfying \eqref{prop1pt1} must be within the bounds as stated in \eqref{prop1pt2}. For the lower bound, since for any matrix $H$ we have col$(HF) \subseteq$ col$(H)$ \cite[Theorem 2.5, p.77]{cullen}, it follows directly that r$(G)=$ r$(HF) \leq$ r$(H)$. For the upper bound, $H$ is first decomposed in terms of its actions on the direct sum components of its domain as established before: $\mathscr{F}^k=$ R$(F)\oplus \mathcal{V}$. For a general matrix $H$ satisfying $G=HF$ it cannot be expected that R$(H)=H($R$(F))\oplus H(\mathcal{V})$ as was the case in the existential part of the proof, but the image of a subspace under a linear transformation is again a subspace so that we have R$(H)=H($R$(F))+ H(\mathcal{V})$. Now it follows that 
\begin{eqnarray}
\nonumber \text{r}(H) &=& \dim(\text{R}(H)) \\
\nonumber &=& \dim(H(\text{R}(F))+ H(\mathcal{V}))\\
 \nonumber &=& \dim(H(\text{R}(F)))+ \dim(H(\mathcal{V})) - \dim(H(\text{R}(F)) \cap H(\mathcal{V})) \\
 \nonumber &=& \dim(\text{R}(G))+ \dim(H(\mathcal{V})) - \dim(H(\text{R}(F)) \cap H(\mathcal{V})) \\
  &=& \text{r}(G)+ \dim(H(\mathcal{V})) - \dim(H(\text{R}(F)) \cap H(\mathcal{V})).\label{eqprp1univh}
\end{eqnarray} 
The third step above follows from Grassmann's theorem (\text{\cite[Proposition 5.16]{golan}}). To conclude the proof it remains to show that
\[ \dim(H(\mathcal{V})) - \dim(H(\text{R}(F)) \cap H(\mathcal{V})) \leq \min(\text{n}(G^T),\text{n}(F^T)).\]
To this end, first notice that
\[ \dim(H(\mathcal{V}))- \dim(H(\text{R}(F)) \cap H(\mathcal{V})) \leq \dim(H(\mathcal{V})) \leq \dim(\mathcal{V}) = k - \text{r}(F),\] so that we have an upper bound for $\dim(H(\mathcal{V}))- \dim(H(\text{R}(F)) \cap H(\mathcal{V}))$. 

Now we can extend a basis for $H(\text{R}(F)) \cap H(\mathcal{V})$ to a basis for $H(\mathcal{V})$, thereby producing the direct sum decomposition $H(\mathcal{V})=\mathcal{U} \oplus (H(\text{R}(F)) \cap H(\mathcal{V}))$. Now since $\mathcal{U} \subseteq H(\mathcal{V})$ and $\mathcal{U} \cap (H(\text{R}(F)) \cap H(\mathcal{V}))=\{\mathbf{0}\}$ it follows that $\mathcal{U} \cap H(\text{R}(F)) = \{\mathbf{0}\}$, and therefore we also have the direct sum $\mathcal{U} \oplus H(\text{R}(F)) \subseteq \mathscr{F}^m$. It follows that 
\begin{eqnarray} 
\nonumber \dim(H(\mathcal{V}))- \dim(H(\text{R}(F)) \cap H(\mathcal{V})) = \dim(\mathcal{U}) &\leq& m - \dim(H(\text{R}(F)))\\
\nonumber &=& m - \dim(\text{R}(G))\\
\nonumber &=& m-\text{r}(G),\end{eqnarray} which is a second upper bound for $\dim(H(\mathcal{V}))- \dim(H(\text{R}(F)) \cap H(\mathcal{V}))$.

Therefore 
\begin{eqnarray} 
\nonumber \dim(H(\mathcal{V}))- \dim(H(\text{R}(F)) \cap H(\mathcal{V})) &\leq& \min(m-\text{r}(G),k-\text{r}(F)) \\
\nonumber &=& \min(\text{n}(G^T),\text{n}(F^T))
\end{eqnarray} by the rank-nullity theorem \cite[Theorem 2.11, p.87]{cullen}. 

Combining these results with \eqref{eqprp1univh} we have 
\[ \text{r}(H) \leq \text{r}(G)+ \min(\text{n}(G^T),\text{n}(F^T)),\] proving the universality of \eqref{prop1pt2}. \hfill $\square$
\subsubsection{Analogous result for left division}
 \begin{quote}
{\bf Proposition \ref{prop:matrixDivision}: Left Division} (Proposition 1 \cite{botha1}). Let $G \in M_{m \times n}(\mathscr{F})$ and $F \in M_{m \times k}(\mathscr{F})$ where $\mathscr{F}$ is a field. The following statements are equivalent:
\begin{enumerate}[(a)]
\item There exists a matrix $H \in M_{k \times n}(\mathscr{F})$ such that $G=FH$.
\item $\text{R}(G) \subseteq \text{R}(F).$
\item $\text{N}(F^T) \subseteq \text{N}(G^T).$
\item $\text{r}(F)=\text{r}\left ( \begin{bmatrix} G & F \end{bmatrix} \right ).$
\end{enumerate}
If any of the statements above are true, then there exists a quotient $H$ of each rank such that
\[\text{r}(G)\leq \text{r}(H ) \leq \text{r}(G) + \min(\text{n}(G),\text{n}(F)).\]
\end{quote}
{\bf Proof.} Note that for left division by $F$, we have $G=FH$ which has transpose $G^T=H^TF^T$, to which we can apply Proposition \ref{prop:matrixDivision} for right division directly. Furthermore $\text{row}(G^T) = \text{R}(G)$, $\text{row}(F^T) = \text{R}(F)$, $\text{r}(F) = \text{r}(F^T)$ and 
\[\text{r}\left (\begin{bmatrix*} G^T \\ F^T \end{bmatrix*}\right ) = \text{r}\left (\begin{bmatrix*} G^T \\ F^T \end{bmatrix*}^T\right ) = \text{r}\left (\begin{bmatrix*} G & F \end{bmatrix*}\right ),\] which proves the equivalence of (a), (b), (c) and (d). The conditions on the rank of $H$ follows directly from Proposition \ref{prop:matrixDivision} for right division, and the fact that the rank of a matrix remains the same under transposition. \hfill $\square$ 
\subsection{Calculating the quotient}
Proposition 2 in the text by Botha \cite{botha1} provides an explicit general formulation of the quotient $H$, which allows it to assume any rank within the bounds determined in Proposition \ref{prop:matrixDivision}. Below I provide a statement and proof following Botha's Proposition 2 for division on the right.
\begin{quote}
\begin{prop}[Proposition 2 (a) \cite{botha1}] Let $G \in M_{m \times n}(\mathscr{F})$ be right divisible by $F \in M_{k \times n}(\mathscr{F})$ where $\mathscr{F}$ is a field and the rank of $F$ is $r$. Then $H \in M_{m \times k}(\mathscr{F})$ is a right quotient of $G$ and $F$ if and only if 
\begin{equation} H=\begin{bmatrix} G & X \end{bmatrix}\begin{bmatrix} F & C \end{bmatrix}^R, \label{prop2pt1} \end{equation}
where $X \in M_{m \times (k-r)}(\mathscr{F}), C \in M_{k \times (k-r)}(\mathscr{F})$ and $\begin{bmatrix} F & C \end{bmatrix}$ is full row rank.

Further, a quotient $H$ of each possible rank according to Proposition \ref{prop:matrixDivision} can be constructed by choosing $X$ of rank 
\[0 \leq \text{r}(X) \leq \min(m-\text{r}(G),k-\text{r}(F))\] such that $\text{R}(X) \cap \text{R}(G) = \{\mathbf{0}\}$.
\label{prop:constructQuotient}
\end{prop}
\end{quote}
{\bf Proof.} Suppose $G=HF$. If $F$ is full row rank, then $F$ has a right inverse and it follows that $H=GF^R$ which completes the proof for this case. Notice that $C$ and $X$ vanish since the rank of $F$ is equal to $k$. Also notice that in this case $F^T$ is full column rank, so that n$(F^T)=0$ by the rank-nullity theorem \cite[Theorem 2.11, p.87]{cullen}, and therefore by Proposition \ref{prop:constructQuotient}, in this special case, r$(H)=$ r$(G)$ and there is no further scope for varying the rank of $H$. 

Now suppose $F$ is not full row rank, then $F$ does not have a right inverse, which poses the problem of how to express $H$ explicitly in terms of the known matrices $F$ and $G$. The problem can be solved by augmenting $F$ with a matrix $C$ of full column rank, in which each column is linearly independent from the columns in $F$: we then have a matrix $\begin{bmatrix}F & C\end{bmatrix}$ which is full row rank. I will now construct a matrix $C$ which satisfies these requirements.

Each matrix is equivalent to a unique matrix of the form $\text{Dg}[I_x,\mathbf{0}]$, where $x$ is the rank of the relevant matrix \cite[Theorem 1.31, p.62]{cullen}. Matrix equivalence can be defined in terms of a decomposition $A=PBQ$, where $A$ and $B$ are the related matrices (with respect to canonical matrix equivalence), and $P$ and $Q$ are square, invertible matrices of appropriate order ($P$ and $Q$ are not necessarily unique) \cite[Theorem 1.29, p. 60]{cullen}. So since $F \sim \text{Dg}[I_r, \mathbf{0}]$ ($\sim$ indicates the relation \emph{matrix equivalence}), we have $F =  Y(\text{Dg}[I_r, \mathbf{0}])Z$, with $Y\in M_{k \times k}(\mathscr{F})$ and $Z\in M_{n \times n}(\mathscr{F})$. 

Now, since $Z$ is invertible, it follows that $Y(\text{Dg}[I_r, \mathbf{0}])=FZ^{-1}$ has the same column space as $F$ \cite[Theorem 2.5, p.77]{cullen}, which is spanned by the first $r$ columns of $Y$: \[Y(\text{Dg}[I_r, \mathbf{0}])=[ \text{col}_1(Y), \text{col}_2(Y), \cdots, \text{col}_r(Y), \mathbf{0}_{k \times (n-r)}]. \] Let $C$ consist of the last $k-r$ columns of $Y$. By invertibility of $Y$ the columns of $C$ are linearly independent from the columns of $F$ and $C$ is full column rank. We can conclude that $\text{col}\left (\begin{bmatrix}F & C\end{bmatrix}\right ) = \text{col}(F) \oplus \text{col}(C)$ and therefore $\text{r}\left (\begin{bmatrix}F & C\end{bmatrix}\right ) = r + k - r = k$, that is $\begin{bmatrix}F & C\end{bmatrix}$ is full row rank, as required.

Now let $X=HC$, then $H\begin{bmatrix}F & C\end{bmatrix} = \begin{bmatrix}G & X\end{bmatrix}$, and since $\begin{bmatrix}F & C\end{bmatrix}$ is full row rank it has a right inverse. It follows that 
\[H = \begin{bmatrix}G & X\end{bmatrix}\begin{bmatrix}F & C\end{bmatrix}^R,\] which proves the necessary part of this proposition.

Now suppose $H$ is of the form given in \eqref{prop2pt1}. Now since $\begin{bmatrix}F & C\end{bmatrix}$ is full row rank it follows that 
\[\text{r}\left (\begin{bmatrix}F & C\end{bmatrix}\right ) = \dim(\mathscr{F}^k) = r + (k-r) = \text{r}(F) + \text{r}(C),\]
and therefore 
\[\mathscr{F}^k=\text{R}(F) \oplus \text{R}(C).\]
This fact makes it possible to define a linear transformation $H_3:\mathscr{F}^k \rightarrow \mathscr{F}^m$ which can be specified in terms of its restriction to $\text{R}(F)$ and its restriction to $\text{R}(C)$. I proceed with construction of such a transformation.

Now since $G$ is right divisible by $F$, there exists a linear transformation $H_1:\mathscr{F}^k \rightarrow \mathscr{F}^m$ such that $G=H_1F$. Let $H_3 \restriction \text{R}(F) = H_1 \restriction \text{R}(F)$. 

Furthermore, by the given conditions in the statement of the proposition $C$ is full column rank, since $\text{r}(C) = k - \text{r}(F)=k-r$, and since $X$ has the same amount of columns as $C$, it follows that the dimension of the space spanned by the columns of $X$ is less than or equal to the dimension of the space spanned by the columns of $C$. Therefore there exists a linear transformation $H_2: \mathscr{F}^k \rightarrow \mathscr{F}^m$ such that $H_2C=X$, which maps the column space of $C$ to the column space of $X$ \cite[Proposition 6.9]{golan}. Let $H_3 \restriction \text{R}(C) = H_2 \restriction \text{R}(C)$. 

Now we have
\[H_3\begin{bmatrix}F & C\end{bmatrix} = \begin{bmatrix}H_1F & H_2C\end{bmatrix}=\begin{bmatrix}G & X\end{bmatrix}.\]
Substituting this result into \eqref{prop2pt1} yields
\[H = H_3\begin{bmatrix}F & C\end{bmatrix}\begin{bmatrix}F & C\end{bmatrix}^R=H_3,\] and so $HF=H_3F = H_1F = G$, and therefore $H$ is a right quotient of $G$ and $F$. 

It remains to show that the rank of $H$ can be varied by choosing $X$ as specified in the statement of the proposition. First, to choose a matrix $X$ with the desired properties we can make use of exactly the same technique as used in the construction of $C$, replacing $F$ with $G$, $r$ with $\text{r}(G)$, $Y\in M_{k \times k}(\mathscr{F})$ with an appropriate invertible matrix $Y_1 \in M_{m \times m}(\mathscr{F})$ and $Z\in M_{n \times n}(\mathscr{F})$ with an appropriate matrix $Z_1 \in M_{k \times k}(\mathscr{F})$. In the previous paragraph it was proved that, with $H$ as given in \eqref{prop2pt1}, we must have 
\[H\begin{bmatrix}F & C\end{bmatrix} = \begin{bmatrix}G & X\end{bmatrix}.\] Now since 
\[\text{R}\left (\begin{bmatrix}F & C\end{bmatrix}\right ) = \mathscr{F}^k = \text{R}(F) \oplus \text{R}(C),\] and $X$ is chosen so that
\[\text{R}\left (\begin{bmatrix}G & X\end{bmatrix}\right ) = \mathscr{F}^m = \text{R}(G) \oplus \text{R}(X),\] 
 it follows that we can choose bases $\alpha_1,\alpha_2$ for R$(F)$, $\text{R}(C)$ and bases $\beta_1, \beta_2$ for R$(G)$, $\text{R}(X)$ such that the (matrix) representation of $H$ in terms of the bases $\alpha=\alpha_1 \cup \alpha_2$ and $\beta=\beta_1 \cup \beta_2$ is \[H\sim \Phi_{\alpha\beta}(H)=\begin{bmatrix}  \Phi_{\alpha_1\beta_1}(H \restriction \text{R}(F)) & \mathbf{0} \\ \mathbf{0} &  \Phi_{\alpha_2\beta_2}(H \restriction \text{R}(C)) \end{bmatrix}.\]
 Therefore, since the range of $H \restriction \text{R}(F)$ is spanned by the columns of $G$ and the range of $H \restriction \text{R}(C)$ is spanned by the columns of $X$, we have $\text{r}(H) = \text{r}(G) + \text{r}(X)$. Furthermore, by the preceding discussion it is easy to see that 
 \[0 \leq \text{r}(X) \leq \min(m-\text{r}(G), \text{r}(C)) = \min(m-\text{r}(G),k-\text{r}(F)),\] which completes the proof. \hfill $\square$ 

\subsubsection{Construction of a left quotient, and the special case $\mathscr{F}=\mathbb{R}$, $\mathscr{F}=\mathbb{C}$}
Note that again, it is easy to verify that the results for right division above are directly applicable to left division of $G$ by $F$.
\begin{quote}
{\bf Proposition \ref{prop:constructQuotient}: Left Division} (Proposition 2 (b) \cite{botha1}). Let $G \in M_{m \times n}(\mathscr{F})$ be left divisible by $F \in M_{m \times k}(\mathscr{F})$ where $\mathscr{F}$ is a field and the rank of $F$ is $r$. Then $H \in M_{k \times n}(\mathscr{F})$ is a left quotient of $G$ and $F$ if and only if 
\begin{equation} H=\begin{bmatrix} F \\ C \end{bmatrix}^L \begin{bmatrix} G \\ X \end{bmatrix}, \label{prop2b}\end{equation}
where $X \in M_{(k-r) \times n}(\mathscr{F}), C \in M_{(k-r) \times k}(\mathscr{F})$ and $\begin{bmatrix} F \\ C \end{bmatrix}$ is full column rank.

Further, a quotient $H$ of each possible rank according to Proposition \ref{prop:matrixDivision} can be constructed by choosing $X$ of rank 
\[0 \leq \text{r}(X) \leq \min(n-\text{r}(G),k-\text{r}(F))\] such that $\text{row}(X) \cap \text{row}(G) = \{\mathbf{0}\}$.
\end{quote}
{\bf Proof.} The result is easy to prove by applying \eqref{prop2pt1} to the transposes of each of the matrices 
\[ \begin{bmatrix} F \\ C \end{bmatrix}, \begin{bmatrix} G \\ X \end{bmatrix}, H,G,F,X \text{ and }C.\] Notice that $(Y^T)^R=(Y^L)^T$ for any matrix $Y$ that is full column rank (that is, $Y^T$ has a right inverse). \hfill $\square$

The results of Proposition \ref{prop:constructQuotient} can be extended in the case where $\mathscr{F}=\mathbb{R}$ or $\mathscr{F}=\mathbb{C}$, as shown in the following corollary.
\begin{quote}
\begin{cor}[Corollary 1 \cite{botha1})] If $\mathscr{F}=\mathbb{C}$ or $\mathscr{F}=\mathbb{R}$, then \eqref{prop2pt1} reduces to
\begin{equation} H = (GF^*+XC^*)(FF^* + CC^*)^{-1},\end{equation} 
and \eqref{prop2b} reduces to
\begin{equation} H = (F^*F+C^*C)^{-1}(F^*G + C^*X).\end{equation}
\label{cor:constructQuotientRC}
\end{cor}
\end{quote}
{\bf Proof.} In this special case the right inverse in \eqref{prop2pt1} can be specified explicitly as 
\[\begin{bmatrix} F & C \end{bmatrix}^* \left (\begin{bmatrix} F & C \end{bmatrix}\begin{bmatrix} F & C \end{bmatrix}^* \right )^{-1},\] and then simplifying. This explicit formulation of the inverse is possible due to the fact that vector spaces over these fields are inner product spaces \cite[Proposition 19.2, p.442]{golan}). In a similar way an explicit left inverse can be specified for left division over $\mathbb{R}$ or $\mathbb{C}$, which concludes the proof. \hfill $\square$

\subsubsection{Example: construction of a right quotient \label{example1}}
I now provide an example of the construction of a right quotient. Let $\mathscr{F}=\mathbb{R}$, and \[ G = \begin{bmatrix} 1 & 2 & 3 & 4 \\ 5 & 6 & 7 & 8 \\ 4 & 8 & 12 & 16 \end{bmatrix} \text{ and } F = \begin{bmatrix} 1 & 4 & 7 & 10 \\ 2 & 5 & 8 & 11 \\ 3 & 6 & 9 & 12 \end{bmatrix}.\] For simplicity I have chosen the field $\mathbb{R}$ with its standard operations, but I will not apply the result in Corollary \ref{cor:constructQuotientRC} in order to illustrate the construction of the quotient in its most general form. Now a feasible criterium for determining whether $F$ is a right divisor of $G$ is (d) of \eqref{prop1pt1}, as this could be accomplished through first performing row reduction of $F$ and $G$ separately (say the resulting matrices in row echelon form are $F_1$ and $G_1$ respectively), and then performing further row reduction on the augmented matrix  \[ T=\begin{bmatrix} G_1 \\ F_1 \end{bmatrix}.\] Suppose $T_1$ is a matrix in row echelon form that is row equivalent to $T$, then it follows that $F$ is a right divisor of $G$ if $T_1$ and $F_1$ have the same amount of nonzero rows. 

Furthermore, during row reduction of $G$ and $F$ we can determine the matrices $C$ and $X$ (as referred to in Proposition \ref{prop:constructQuotient}) simultaneously, by augmenting both matrices before row reduction to form the matrices $\begin{bmatrix} G & : & I_m \end{bmatrix}$ and $\begin{bmatrix} F & : & I_k \end{bmatrix}$ (where it was assumed $G$ has $m$ rows and $F$ has $k$ rows). By applying the row reduction algorithm and a method outlined in \cite[p.683, 684]{lee} it is then possible to determine invertible matrices $Y_1$ and $Y$ (the matrices referred to in the discussion of the proof of Proposition \ref{prop:constructQuotient}), from which we can then select columns for the construction of $X$ and $C$. 

Applying the method as discussed above, we have
\[ \begin{bmatrix} G & : & I_m \end{bmatrix} = \begin{bmatrix} 1 & 2 & 3 & 4 & : & 1 & 0 & 0 \\ 5 & 6 & 7 & 8 & : & 0 & 1 & 0 \\ 4 & 8 & 12 & 16 & : & 0 & 0 & 1 \end{bmatrix}\stackrel{R}{\sim} \begin{bmatrix} 1 & 2 & 3 & 4 & : & 1 & 0 & 0 \\ 0 & 4 & 8 & 12 & : & 5 & -1 & 0 \\ 0 & 0 & 0 & 0 & : & 4 & 0 & -1 \end{bmatrix},\] and
\[\begin{bmatrix} F & : & I_k \end{bmatrix} = \begin{bmatrix} 1 & 4 & 7 & 10 & : & 1 & 0 & 0 \\ 2 & 5 & 8 & 11 & : & 0 & 1 & 0 \\ 3 & 6 & 9 & 12 & : & 0 & 0 & 1 \end{bmatrix} \stackrel{R}{\sim} \begin{bmatrix} 1 & 4 & 7 & 10 & : & 1 & 0 & 0 \\ 0 & 3 & 6 & 9 & : & 2 & -1 & 0 \\ 0 & 0 & 0 & 0 & : & 1 & -2 & 1 \end{bmatrix}.\]
Now the matrices\footnote{These matrices are determined by selecting the columns from $\begin{bmatrix} G & : & I_m \end{bmatrix}$ and $\begin{bmatrix} F & : & I_k \end{bmatrix}$ corresponding to the columns that have leading entries in the row reduced forms of these systems. See \cite[p.683, 684]{lee}.} \[Y_1=\begin{bmatrix}1 & 2 & 1 \\ 5 & 6 & 0 \\ 4 & 8 & 0\end{bmatrix} \text{ and } Y=\begin{bmatrix}1 & 4 & 1 \\ 2 & 5 & 0 \\ 3 & 6 & 0\end{bmatrix} \] are such that $G=Y_1(\text{Dg}[I_2, \mathbf{0}])Z_1$ and $F =  Y(\text{Dg}[I_2, \mathbf{0}])Z$ for some invertible matrices $Z_1, Z \in M_{4 \times 4}(\mathbb{R})$. So at this point it is not necessary to perform any further reduction, and we can proceed to construct the matrix $T$ mentioned above:  
\[T=\begin{bmatrix} G_1 \\ F_1 \end{bmatrix} = \begin{bmatrix} 1 & 2 & 3 & 4 \\ 0 & 4 & 8 & 12 \\ 0 & 0 & 0 & 0 \\ 1 & 4 & 7 & 10 \\ 0 & 3 & 6 & 9 \\ 0 & 0 & 0 & 0 \end{bmatrix} \stackrel{R}{\sim} \begin{bmatrix} 1 & 2 & 3 & 4 \\ 0 & 1 & 2 & 3 \\ 0 & 0 & 0 & 0 \\ 0 & 0 & 0 & 0 \\ 0 & 0 & 0 & 0 \\ 0 & 0 & 0 & 0 \end{bmatrix}=T_1, \] so that $F_1$ and $T_1$ are both in row echelon form and have two nonzero rows, and therefore $F$ is a right divisor of $G$. Furthermore since \[ \text{n}(G^T)= (\text{number of rows of }G)- \text{r}(G)=3-2=1,\] and \[ \text{n}(F^T)= (\text{number of rows of }F)- \text{r}(F)=3-2=1,\]
it follows that the rank of any right quotient of $F$ and $G$ can vary from r$(G)=2$ to \[\text{r}(G)+\min(\text{n}(G^T),\text{n}(F^T))=2+\min(1,1)=3.\]

Now I will proceed to construct two right quotients $H_1,H_2$ with rank two and three respectively. Now since the number of rows of $F$ is three, and the rank of $F$ is two, the matrix $C$ consists of the last $3-2=1$ columns of $Y$, that is \[C=\begin{bmatrix} 1\\0\\0 \end{bmatrix}.\] Now to construct a matrix $H_1$ of rank two, which is the same rank as $G$, let \[X = \begin{bmatrix} 0\\0\\0 \end{bmatrix},\] then it remains to calculate \[H_1=\begin{bmatrix} 1 & 2 & 3 & 4 & 0\\ 5 & 6 & 7 & 8 & 0 \\ 4 & 8 & 12 & 16 & 0 \end{bmatrix}\begin{bmatrix} 1 & 4 & 7 & 10 & 1\\ 2 & 5 & 8 & 11&0 \\ 3 & 6 & 9 & 12 & 0 \end{bmatrix}^R.\] 
To calculate the right inverse above, since $\mathscr{F}=\mathbb{R}$, for a matrix $A$ that is full row rank, $AA^T$ is a square matrix of full rank\footnote{This is easy to see, since if $A$ is full row rank, the linear transformations $A^T$ and $A$ are injective and surjective respectively.} and hence invertible, so that $A^T(AA^T)^{-1}$ is a right inverse of $A$. It follows that 
\begin{eqnarray} \nonumber H_1&=&\begin{bmatrix} 1 & 2 & 3 & 4 & 0\\ 5 & 6 & 7 & 8 & 0 \\ 4 & 8 & 12 & 16 & 0 \end{bmatrix}\left [ {\renewcommand{\arraystretch}{1.2}\begin{array}{rrr} 0 & -1 & \frac{9}{10} \\
 0 & -\frac{1}{2} & \frac{7}{15} \\
 0 & 0 & \frac{1}{30} \\
 0 & \frac{1}{2} & -\frac{2}{5} \\
 1 & -2 & 1 \end{array}} \right ] \\
 \nonumber &=& \frac{1}{3}\begin{bmatrix} 0 & 0 & 1 \\
 0 & -12 & 13 \\
 0 & 0 & 4 \end{bmatrix},
 \end{eqnarray} which is clearly rank two.
 
 To construct a matrix $H_2$ of rank three, we need to select one column from $Y_1$ which is linearly independent from the columns in $G$. From the results discussed before, the only column that satisfies this requirement is the last column of $Y_1$, therefore let \[X=\begin{bmatrix} 1 \\ 0 \\ 0 \end{bmatrix}.\] Now we have 
 \begin{eqnarray} \nonumber H_2&=& \begin{bmatrix} 1 & 2 & 3 & 4 & 1\\ 5 & 6 & 7 & 8 & 0 \\ 4 & 8 & 12 & 16 & 0 \end{bmatrix}\begin{bmatrix} 1 & 4 & 7 & 10 & 1\\ 2 & 5 & 8 & 11&0 \\ 3 & 6 & 9 & 12 & 0 \end{bmatrix}^R \\
\nonumber &=& \begin{bmatrix} 1 & 2 & 3 & 4 & 1\\ 5 & 6 & 7 & 8 & 0 \\ 4 & 8 & 12 & 16 & 0 \end{bmatrix}\left [ {\renewcommand{\arraystretch}{1.2}\begin{array}{rrr} 0 & -1 & \frac{9}{10} \\
 0 & -\frac{1}{2} & \frac{7}{15} \\
 0 & 0 & \frac{1}{30} \\
 0 & \frac{1}{2} & -\frac{2}{5} \\
 1 & -2 & 1 \end{array}} \right ] \\
 \nonumber &=& \frac{1}{3}\begin{bmatrix} 3 & -6 & 4 \\
 0 & -12 & 13 \\
 0 & 0 & 4 \end{bmatrix},
\end{eqnarray}
which is of rank three. It is easy to verify that $G=H_1F=H_2F$.

\section{Main Results: Products of Square-zero Matrices}
Here I present the main results, which I divide into three parts. First, independent from the theory of matrix division presented in preliminary results, I discuss results by Novak \cite{novak} and Botha \cite{botha2} which indicate when a matrix can be factored into a product consisting of two or three square-zero matrices. 

Then, building on the theory developed in preliminary results, I present a general result by Botha \cite{botha2} indicating necessary and sufficient conditions for a matrix $G$, and right divisor $F$, to have a square-zero right quotient. Similar to the preliminary results, the main result is followed by a theorem providing an explicit construction of a square-zero quotient, whenever such a quotient might exist. 

The results developed in part two are then used to prove the results on factorization into two or three square-zero factors as presented in part one (results by Novak \cite{novak} and Botha \cite{botha2}). Besides providing additional insight into the results of part one, this approach provides a coherent theory, moving from the general to the specific. An additional result, that effectively shows that any matrix that can be written as a product of square-zero matrices, can also be written as a product of arbitrarily many square-zero matrices subject to it being more than three factors, follows from this discussion.

\subsection{Square-zero factorization: previous results \label{sec:sqzeroprev}}
When can a matrix be written as a product consisting only of square-zero factors? I start with two results by Novak \cite[p.11,12]{novak} (section 2 of the relevant text), which is concerned with $M_n(\mathscr{F})$, which is the set of $n \times n$ matrices over an \emph{algebraically closed} field $\mathscr{F}$. The generalization by Botha extends the result to matrices over an arbitrary field \cite{botha2}: Theorems 4 and 6 in the text.  

\subsubsection{Factorization into a product of two square-zero matrices}
I follow the order in \cite{botha2} and start with a theorem indicating conditions equivalent to factorization into a product of two square-zero matrices (Theorem 2.2 \cite{novak}, Theorem 4 \cite{botha2}). Note that in the text by Novak the meaning of the symbol $\ominus$ is not defined explicitly. From the context it would seem to indicate an intersection with a complementary space of some kind. For the purposes of this dissertation, in order to avoid difficulties that arise in the proof when considering `complementary' in the light of inner product spaces or symmetric bilinear forms, I provide the following definition for the proof to follow:
\begin{quote}
\begin{defn}
Let $\mathscr{F}$ be an arbitrary field, $\mathcal{V}$ a vector space over $\mathscr{F}$, and $\mathcal{W}$ be an arbitrary subspace of $\mathcal{V}$. Let $\mathcal{U}_1$ be any subspace of $\mathcal{W}$. Define
\[\mathcal{W} \ominus \mathcal{U}_1 = \mathcal{U}_2,\] where $\mathcal{U}_2$ is any subspace of $\mathcal{W}$ satisfying $\mathcal{W} = \mathcal{U}_1 \oplus \mathcal{U}_2$.
\label{defn:ominus}
\end{defn}
\end{quote}
In words, $\mathcal{W} \ominus \mathcal{U}_1$ indicates any complementary subspace of $\mathcal{U}_1$ in $\mathcal{W}$. The proof of Novak's result below (Theorem \ref{thm:novak2FactorSz}) is adapted in accordance with the definition above.
\begin{quote}
\begin{thm}[Specialized version - Theorem 2.2 \cite{novak}] Let $G \in M_n(\mathscr{F})$ be an operator on $\mathscr{F}^n$ where $\mathscr{F}$ is an algebraically closed field. The following statements are equivalent: 
\begin{enumerate}[(a)]
\item $G=HF$, where $H$ and $F$ are square-zero;
\item $\dim (\text{N}(G) \ominus (\text{N}(G) \cap \text{R}(G))) \geq \text{r}(G)$;
\item $G$ is similar to a matrix $S$ such that $\dim(\text{N}(S) \cap \text{N}(S^T)) \geq \text{r}(S)$. [\emph{sic}]
\end{enumerate}
\label{thm:novak2FactorSz}
\end{thm}
\end{quote}
\begin{quote}
\begin{thm}[Generalized version - Theorem 4 \cite{botha2}] Let $G \in M_n(\mathscr{F})$ where $\mathscr{F}$ is a field. The following statements are equivalent: 
\begin{enumerate}[(a)]
\item $G=HF$, where $H$ and $F$ are square-zero; 
\item $\text{r}(G) \leq \text{n}(G) - \dim(\text{N}(G) \cap \text{R}(G))$;
\item $G$ is similar to a matrix $\text{Dg}[\mathbf{0}_{\text{r}(G)}, A]$ for some square matrix $A$ over~$\mathscr{F}$.
\end{enumerate}
\label{thm:botha2FactorSZ1}
\end{thm}
\end{quote}
{\bf Proof.} During the course of the proof I will compare the proofs of Novak and Botha. Note that there is a slight error in statement (c) of the specialized version, which will be pointed out when discussing parts of the proof involving (c) below.

{\bf (a) implies (b).}  The arguments in both versions are essentially the same, relying on construction of complementary subspaces. Theorem \ref{thm:botha2FactorSZ1} assumes existence of a subspace $\mathcal{W}$ such that 
\begin{equation} \text{N}(G)=\mathcal{W} \oplus (\text{N}(G) \cap \text{R}(G)).\label{th1.1}\end{equation}
The existence of $\mathcal{W}$ is justified by the fact that we can take any basis, say $\alpha$, for $\text{N}(G) \cap \text{R}(G)$ and then extend it by adding the (linearly independent set of vectors) $\beta$ so that $\alpha \cap \beta = \emptyset$, each vector in $\beta$ is linearly independent from each vector in~$\alpha$, and $\alpha \cup \beta$ is a basis for $\text{N}(G)$ (which basically refers to the fact that any set of linearly independent vectors can be extended to be a basis \cite[Proposition 5.8, p.68]{golan}). Then we have the desired subspace by letting span($\beta)=\mathcal{W}$.

Now by the direct sum \eqref{th1.1} it follows that 
 \begin{equation} \dim(\mathcal{W})=\text{n}(G)- \dim(\text{N}(G) \cap \text{R}(G)).\label{th1.2}\end{equation} 
 
 In terms of the specialized version, as per Definition \ref{defn:ominus} we can replace $\mathcal{W}$ in \eqref{th1.1} and \eqref{th1.2} with \[\text{N}(G) \ominus (\text{N}(G) \cap \text{R}(G)),\] so that it is clear that (b) in the  specialized version states the same result as (b) in the generalized version.
 
 Having established the above, it remains to use the properties of matrix division and square-zero matrices to prove the desired implication. Since $F$ is square-zero we must have
 \[GF = HF^2 = \mathbf{0},\] and it follows that $\text{R}(F) \subseteq \text{N}(G)$. Now considering $H$ as a linear operator on $\mathscr{F}^n$, which is well-defined functionally, and substituting \eqref{th1.1} we have
 \begin{equation}\text{R}(G)=H(\text{R}(F)) \subseteq H(\text{N}(G))=H(\mathcal{W} \oplus (\text{N}(G) \cap \text{R}(G)))=H(\mathcal{W}).\label{th1.3}\end{equation} To clarify the last step on the right above: Let $w \in \mathcal{W}$ and $u \in \text{N}(G) \cap \text{R}(G)$ so that $w+u \in \mathcal{W} \oplus (\text{N}(G) \cap \text{R}(G))$. Now we have $H(w+u)=H(w)+H(u)$ and since $u \in \text{R}(G)$ it follows that $u=G(v)$ for some $v \in \mathscr{F}^n$. So we must have $H(u)=HG(v)=H^2F(v)=\mathbf{0}$ since $H$ is square-zero, and therefore $H(w+u)=H(w) \in H(\mathcal{W})$. It is easy to see that the reverse inclusion holds ($H(\mathcal{W}) \subseteq H(\mathcal{W} \oplus (\text{N}(G) \cap \text{R}(G)))$). 
 
 By \eqref{th1.3} it follows that
\begin{equation}\text{r}(G) \leq \dim(H(\mathcal{W})) \leq \dim(\mathcal{W}) = \text{n}(G) - \dim(\text{N}(G) \cap \text{R}(G)),\label{th1.4}\end{equation} where the second inequality follows from functionality of $H$ and dimension of $\mathcal{W}$ by \eqref{th1.2}. Again, the specialized version is directly obtained by replacing $\mathcal{W}$ with $\text{N}(G) \ominus (\text{N}(G) \cap \text{R}(G))$ in \eqref{th1.3} and \eqref{th1.4}.

{\bf (b) implies (c).} First notice that the statements in (c) are not equivalent between the two versions of the theorem. Consider the matrix \[S=\begin{bmatrix} 1 & i \\ i & -1 \end{bmatrix} \in M_2(\mathbb{C}); \] it is easy to check that $\dim(\text{N}(S) \cap \text{N}(S^T))=1=\text{r}(S)$, but $S$ is not diagonalizable and therefore not similar to a matrix $\text{Dg}[\mathbf{0}_{\text{r}(S)}, A]$ for some matrix $A$. In fact, when discussing the next implication (\emph{(c) implies (a)}) I will prove that \emph{(c) implies (b)} is false for the specialized version, so that (c) is in fact not equivalent to the other statements in the specialized version as originally claimed. When assuming that $S$ is in Jordan canonical form (as is implicitly done) the proof is valid: note that this is an additional condition that should be explicitly stated in (c) of the specialized version of the theorem.

So, proceeding with the proof, assume that statement (b) holds. Both proofs make use of the invariance of $\text{n}(G)$, $\text{r}(G)$ and $\dim(\text{N}(G) \cap \text{R}(G))$ with respect to similarity. I will verify this detail: since rank and nullity are invariant with respect to similarity\footnote{Similarity is a special case of matrix equivalence, and therefore two matrices that are similar are both equivalent to the same matrix Dg$[I_r,\mathbf{0}]$ \cite[Theorem 1.31, p.62]{cullen}}, we only need to verify that $\dim(\text{N}(G) \cap \text{R}(G))$ is also invariant. 

Suppose $S,T \in M_n(\mathscr{F})$ so that $S=P^{-1}TP$. Let $P(\text{N}(S) \cap \text{R}(S)) = \mathcal{W}$, and consider the  linear transformation \[P \restriction (\text{N}(S) \cap \text{R}(S))=P_1 : \text{N}(S) \cap \text{R}(S) \rightarrow \mathcal{W} \subseteq \mathscr{F}^n \text{ defined by } P_1x=Px.\] 
Since $P$ is invertible, $P_1$ is a bijection onto $\mathcal{W}$, and therefore it remains to verify that $\mathcal{W}=\text{N}(T) \cap \text{R}(T)$. Let $y \in \mathcal{W}$ be nonzero, then there exists a nonzero vector $x \in \text{N}(S) \cap \text{R}(S)$ so that $y=Px$ and $Sx=\mathbf{0}$. But then $\mathbf{0}=P^{-1}TPx=P^{-1}Ty$, which is true if and only if $Ty=\mathbf{0}$ so that $y \in \text{N}(T)$. Furthermore, since $x \in  \text{R}(S)$, there exists nonzero $u \in \mathscr{F}^n$ so that $Su=x$. But then $P^{-1}TPu=x$ and therefore $TPu=Px=y$. So $y \in \text{R}(T)$, proving that every vector in $\mathcal{W}$ is also in $\text{N}(T) \cap \text{R}(T)$. 

For the converse, suppose we have nonzero $z \in \text{N}(T) \cap \text{R}(T)$. Then $Tz=\mathbf{0}$ and therefore $PSP^{-1}z=\mathbf{0}$. Again, this is true if and only if $SP^{-1}z=\mathbf{0}$, so that $P^{-1}z \in \text{N}(S)$. Also, since $z \in  \text{R}(T)$ there exists nonzero $w \in \mathscr{F}^n$ so that $Tw=z$. Therefore $PSP^{-1}w=z$ so that $SP^{-1}w=P^{-1}z$, proving that $P^{-1}z \in \text{R}(S)$. Now since $P^{-1}z \in \text{N}(S) \cap \text{R}(S)$, we have $PP^{-1}z=z$ confirming that $z$ is in the range of $P$. That is $z \in \mathcal{W}$, completing the proof. 

Both proofs then make use of similarity to a diagonal block matrix to complete the proof of \emph{(b) implies (c)}. The specialized version makes use of the Jordan canonical form which is admissible within the context of an algebraically closed field. Note that specifying the Jordan canonical form in the statement of (c) would render all statements of the specialized version equivalent. 

First, in terms of the generalized version, we can make use of Fitting's lemma \cite[Theorem 5.10, p.195]{cullen}: $G$ is similar to a matrix $\text{Dg}[N,B]$, where $N$ is nilpotent and $B$ is nonsingular. Any matrix over an arbitrary field is similar to a matrix in rational canonical form. Since $B$ is nonsingular its characteristic polynomial has no factors of the form $x^e$, whereas $N$ is nilpotent and therefore there exists an integer $0\leq k \leq t$ (where $t$ is the order of $N$) so that $N^k=\mathbf{0}$. By the Cayley-Hamilton theorem \cite[Theorem 5.1, p.176]{cullen} the characteristic polynomial of $N$ is therefore $x^t$. So the elementary divisors\footnote{\cite[p.241]{cullen}} of $N$ are all of the form $x^e$, so that the rational canonical form of $N$ consists of diagonal blocks $H(x^e)=J_e(0)$ where $H(x^e)$ indicates the hypercompanion matrix of $x^e$ \cite[Definition 7.3, p.244]{cullen}, and therefore we can refer to ``...the number of Jordan blocks in $N$ equal to $[0]$..." \cite[p.82]{botha2} within the context of an arbitrary field. 

To conclude this part of the proof, we have to prove that (b) implies there are at least $\text{r}(G)$ simple Jordan blocks $J_1(0)$ in the rational canonical form of $G$. Now Botha states that $\text{n}(G)-\dim(\text{N}(G) \cap \text{R}(G))$ is exactly the number of simple Jordan blocks in $N$ of order $1 \times 1$. To verify this statement, first observe that since $B$ is nonsingular, it follows that $\text{n}(G) = \text{n}(N)$ and $\dim(\text{N}(G) \cap \text{R}(G))=\dim(\text{N}(N) \cap \text{R}(N))$. 

Now it was shown above that $N$ is similar to a matrix consisting entirely of simple Jordan blocks $J_k(0)$. Consider any Jordan block $J_k(0)$ with $k>1$: it has rank $k-1$ and $\text{N}(J_k(0))$ is spanned by $(0,0,\ldots,1)^T \in \mathscr{F}^k$, which is exactly the vector in the second last column of $J_k(0)$.\footnote{Note that I am assuming Jordan blocks with ones below the diagonal as is the custom in Cullen \cite{cullen}.} Therefore, 
\[\text{N}(J_k(0)) \cap \text{R}(J_k(0)) = \text{N}(J_k(0))\] for any simple Jordan block $J_k(0)$, where $k>1$. Furthermore it is easy to see that 
\[\text{N}(J_1(0)) \cap \text{R}(J_1(0)) = \{\mathbf{0}\}.\] 
This proves that $\dim(\text{N}(G) \cap \text{R}(G))=\dim(\text{N}(N) \cap \text{R}(N))$ is equal to the number of blocks $J_k(0)$ with $k > 1$ in $N$. Now it follows directly that 
\[\text{n}(G)-\dim(\text{N}(G) \cap \text{R}(G))=\text{n}(N)-\dim(\text{N}(N) \cap \text{R}(N))\] is the number of Jordan blocks $J_1(0)$ in $N$, which verifies the validity of Botha's statement \cite[p.82]{botha2}. 

To conclude the proof of this part: assuming (b) is true, there are at least $\text{r}(G)$ blocks $J_1(0)$ in $N$ so that we have 
\[G \approx \text{Dg}[N, B] \approx \text{Dg}[\mathbf{0}_{\text{r}(G)}, N_1, B],\] where $N_1$ consists of blocks $J_k(0)$ with $k \geq 1$. By now specifying $A=\text{Dg}[N_1, B]$, the desired result follows. 

Now in terms of the specialized version: it is easy to see that for any Jordan block $J_k(0)$ we have $\text{N}(J_k(0)) = \text{N}(J_k(0)^T)$ if and only if $k = 1$, and $\text{N}(J_k(0)) \cap \text{N}(J_k(0)^T) = \{\mathbf{0}\}$ in all other cases. The latter is also true for all Jordan blocks $\text{N}(J_k(\lambda))$ with $\lambda$ nonzero, since such a block is full rank. By specifying $S$ as the Jordan canonical form of $G$ it is then easy to see (by a similar argument as used above for the generalized version) that 
\[\text{n}(G) - \dim(\text{N}(G) \cap \text{R}(G)) = \dim ( \text{N}(S) \cap \text{N}(S^T)),\] and therefore (c) follows from (b). So the specialized version makes use of the Jordan canonical form, but the rest of the proof is essentially the same.

{\bf (c) implies (a).} Here I will first address the shortcoming of statement (c) in the specialized version, by proving that \emph{(c) implies (b)} is false. Consider the matrix
 \[S=\begin{bmatrix} 1 & i \\ i & -1 \end{bmatrix} \in M_2(\mathbb{C}).\] We have $\text{r}(S)=1, \text{n}(S)=1$ and $(-i,1)^T \in \text{N}(S) \cap \text{R}(S)$ so that $\dim(\text{N}(S)~\cap~\text{R}(S))=1$. Since $S$ is symmetric $\dim(\text{N}(S) \cap \text{N}(S^T))=1=\text{r}(S)$, satisfying the conditions of the specialized version of (c), but we have $\text{n}(S)-\dim(\text{N}(S) \cap \text{R}(S))=0<~\text{r}(S)$, and by the invariance of $\text{n}(S)-\dim(\text{N}(S) \cap \text{R}(S))$ and $\text{r}(S)$ with respect to similarity, the matrix $S$ cannot be similar to any matrix that will satisfy the conditions in~(b), thereby providing a counterexample. Now Novak's proof of \emph{(c) implies (a)} in the specialized version implicitly assumes that $S$ is in Jordan canonical form, which renders the proof correctly. The statement of (c) in the specialized version therefore needs the additional specification that $S$ be in Jordan canonical form, for all statements to be equivalent as claimed. 
 
I will proceed to complete the proof of the last implication. Both proofs state that the diagonal block which contains nonzero entries ($A$ in statement (c) of the general version) is similar to a matrix \[ \begin{bmatrix} \mathbf{0} & X \\ \mathbf{0} & C\end{bmatrix}, \] with $C$ some $\text{r}(G) \times \text{r}(G)$ matrix. I will verify this detail: let $G \approx \text{Dg}[\mathbf{0}_{\text{r}(G)}, A]$. We also have that $G \approx \text{Dg}[N, B]$ by Fitting's lemma \cite[Theorem 5.10, p.195]{cullen}. Now $B$ is full rank and therefore must be of order $r(G) \times r(G)$ or smaller, and therefore $N \approx \text{Dg}[\mathbf{0}_{\text{r}(G)}, \mathbf{0}_{k_1}, N_1]$ for some nonnegative integer $k_1$ and \[N_1=\text{Dg}[J_{e_1}(0),J_{e_2}(0),\ldots,J_{e_k}(0)]\] for some nonnegative integer $k$, with $e_i \geq 2$ for each $i \in \{1,2,\ldots,k\}$.\footnote{Stated explicitly, $N$ is similar to a matrix with Jordan blocks which only has an eigenvalue of zero, in which the blocks are arranged so that all blocks of order 1 appear first on the diagonal. The matrix which consists of all blocks of order greater than 1 is $N_1$.} It follows that the proof will be complete if we can show $N_1 \approx \begin{bmatrix} \mathbf{0} & X_1 \end{bmatrix}$ where $X_1$ is full column rank, for supposing $X_1$ has $i$ columns: let $C = \text{Dg}[Y, B]$ where $Y$ is the bottom $i$ rows of $X_1$, and let $X$ consist of the rows above $C$ with $i$ columns on the left being the remaining top rows of $X_1$ and all other columns zero.

Let $P$ be a permutation matrix corresponding to elementary column operations that move all zero columns in $N_1$ to the left. Then $N_1P=\begin{bmatrix} \mathbf{0} & X_2 \end{bmatrix}$ with $X_2$ full column rank. Now $P^{-1}$ is also a permutation matrix which performs row permutations of $N_1P$ when multiplying on the left, so that zero columns remain at the left in the product, that is $P^{-1}N_1P=\begin{bmatrix} \mathbf{0} & X_1 \end{bmatrix}$ with $X_1$ full column rank, as desired. 

Now the proof is concluded in both versions using the following factorization (I follow the notation and symbols as used in the general version by Botha \cite[p.82]{botha2}):
\[G \approx \begin{bmatrix} \mathbf{0}_{\text{r}(G)} & \mathbf{0} & \mathbf{0}_{\text{r}(G)} \\ 
\mathbf{0} & \mathbf{0} & X \\
\mathbf{0}_{\text{r}(G)} & \mathbf{0} & C \end{bmatrix} =  \begin{bmatrix} \mathbf{0}_{\text{r}(G)} & \mathbf{0} & \mathbf{0}_{\text{r}(G)} \\ 
X & \mathbf{0} & \mathbf{0} \\
C & \mathbf{0} & \mathbf{0}_{\text{r}(G)} \end{bmatrix}   \begin{bmatrix} \mathbf{0}_{\text{r}(G)} & \mathbf{0} & I_{\text{r}(G)} \\ 
\mathbf{0} & \mathbf{0} & \mathbf{0} \\
\mathbf{0}_{\text{r}(G)} & \mathbf{0} & \mathbf{0}_{\text{r}(G)} \end{bmatrix}.\] 
By routine calculation it can be verified that the matrices on the right are square-zero. Denote this product as $H_1F_1$, so that $G=Q^{-1}H_1F_1Q$ for some change-of-basis matrix $Q$. Then we have $G=Q^{-1}H_1F_1Q=Q^{-1}H_1QQ^{-1}F_1Q=HF$. Now since the minimum polynomial of a matrix is invariant with respect to similarity it follows that $H$ and $F$ are square-zero, as desired. \hfill $\square$

\subsubsection{Factorization into a product of three square-zero matrices}
Now I discuss a result providing necessary and sufficient conditions for a matrix to be a product of three square-zero matrices.

\begin{quote}
\begin{thm}[Specialized version - Proposition 2.1 \cite{novak}] A matrix $G \in M_n(\mathscr{F})$, where $\mathscr{F}$ is an algebraically closed field, is a product of three square-zero matrices if and only if $\text{r}(G)\leq n/2$.
\label{novak3FactorsSz}
\end{thm}
\end{quote}

\begin{quote}
\begin{thm}[Generalized version - Theorem 6 \cite{botha2}] Let $G \in M_n(\mathscr{F})$ where $\mathscr{F}$ is an arbitrary field. The matrix $G$ is a product of three square-zero matrices if and only if $\text{r}(G) \leq n/2$ (equivalently $\text{r}(G) \leq \text{n}(G)$).
\label{botha3FactorsSz}
\end{thm}
\end{quote}
{\bf Proof.} Suppose $G$ is a product of three square-zero matrices. For any square-zero matrix $H \in M_n(\mathscr{F})$ we must have $\text{r}(H) \leq n/2$, for suppose this is not so: then $\text{r}(H)>\text{n}(H)$ (by the rank-nullity theorem \cite[Theorem 2.11, p.87]{cullen}), and therefore there exists nonzero $v \in \text{R}(H)$ so that $v \notin \text{N}(H)$. But then $H^2w=H(H(w))=H(v)\neq \mathbf{0}$ for some nonzero $w \in \mathscr{F}^n$, which is a contradiction, establishing that we must have $\text{r}(H) \leq n/2$. Now, if $G=H_1H_2H_3$ with $H_1,H_2,H_3$ square-zero, then $\text{r}(G) \leq \min(\text{r}(H_1),\text{r}(H_2),\text{r}(H_3)) \leq n/2$ (follows from \cite[Theorem 2.5, p.77]{cullen}) as proposed.

Conversely, if $\mathscr{F}$ is considered to be an algebraically closed field, then $G$ is similar to a matrix in Jordan canonical form, which is unique up to the ordering of the simple Jordan blocks on the diagonal. Now since $\text{r}(G) \leq n/2$ it follows that $G \approx \text{Dg}[G_1,G_2,\ldots, G_k, \mathbf{0}]$ for some nonnegative integer $k$, where $G_i=\text{Dg}[J_{n_i}(\lambda_i),\mathbf{0}_{n_i}]$ if $\lambda_i \neq 0$, or $G_i=\text{Dg}[J_{n_i}(0),\mathbf{0}_{n_i-2}]$. Now notice that each $G_i$ can be written as a product of three square-zero matrices:
\begin{eqnarray}
\nonumber \lambda_i \neq 0&:& \qquad \text{Dg}[J_{n_i}(\lambda_i),\mathbf{0}_{n_i}] = \begin{bmatrix} \mathbf{0} & J_{n_i}(\lambda_i) \\ \mathbf{0}_{n_i} & \mathbf{0} \end{bmatrix} \begin{bmatrix} -I_{n_i} & -I \\ I & I_{n_i} \end{bmatrix} \begin{bmatrix} \mathbf{0} & \mathbf{0}_{n_i} \\ I_{n_i} & \mathbf{0} \end{bmatrix} \\
\nonumber \lambda_i=0 \text{ and } n_i=2&:& \qquad J_{2}(0) = \begin{bmatrix} 0 & 0 \\ 1 & 0 \end{bmatrix} \begin{bmatrix} 0 & 1 \\ 0 & 0 \end{bmatrix} \begin{bmatrix} 0 & 0 \\ 1 & 0 \end{bmatrix} \\
\nonumber \lambda_i=0 \text{ and } n_i>2&:& \qquad \text{Dg}[J_{n_i}(0),\mathbf{0}_{n_i-2}] = \begin{bmatrix} E_{2,1} & E-E_{2,1} \\ \mathbf{0} & E_{1,n_i-1} \end{bmatrix} \begin{bmatrix} -I & I \\ -I & I \end{bmatrix} \begin{bmatrix} \mathbf{0} & \mathbf{0} \\ I & \mathbf{0} \end{bmatrix}. 
\end{eqnarray}
Here the result has been modified for the case where the Jordan form is expressed with ones on the sub-diagonal for consistency with the rest of this text (compare with \cite[proposition 2.1]{novak}). The matrix $E_{2,1}$ has a one in entry $(2,1)$ and zeros elsewhere, $E_{1,n_i-1}$ has a one in entry $(1,n_i-1)$ and zeros elsewhere, and $E$ has ones on the sub-diagonal and zeros elsewhere, i.e. $E = J_{n_i-1}(0)$. So this type of partitioning of $\text{Dg}[J_{n_i}(0),\mathbf{0}_{n_i-2}]$ results in blocks of the same order, namely blocks of order $n_i-1$.

When generalizing to an arbitrary field, $G$ does not necessarily have a Jordan decomposition, but we can make use of Fitting's lemma \cite[Theorem 5.10, p.195]{cullen}: $G \approx \text{Dg}[N, B]$ where $N$ is nilpotent and $B$ is nonsingular. Now $N \approx \text{Dg}[N_1, \mathbf{0}_{n_2}]$ where $N_1$ consists of simple Jordan blocks $J_k(0)$ where $k>1$ and $\mathbf{0}_{n_2}$ is a zero matrix of order $n_2 \times n_2$ where $n_2$ is some nonnegative number.

Now by the condition $\text{r}(G) \leq n/2$ and the rank-nullity theorem we must have $\text{r}(G) \leq \text{n}(G)$. Furthermore it is easy to see that $\text{r}(G) = \text{r}(B) + \text{r}(N_1)$ and $\text{n}(G) = \text{n}(N_1) + \text{n}(\mathbf{0}_{n_2})=n_1 + n_2$ where $n_1$ is the number of simple Jordan blocks in $N_1$. Combining these results we have $n_1+n_2 \geq \text{r}(B) + \text{r}(N_1)$, and therefore \[n_2 \geq \text{r}(B) + (\text{r}(N_1)-n_1).\] 
Notice that the second term on the right is exactly the amount of zeros required to create paired diagonal blocks of the form $\text{Dg}[J_k(0),\mathbf{0}_{k-2}]$ for each simple Jordan block $J_k(0)$ with $k>1$ in $N_1$. From this result it follows that $G \approx \text{Dg}[N_3, \mathbf{0}_{\text{r}(B)}, B]$ where $N_3 \approx \text{Dg}[J_1,J_2,\ldots,J_t,\mathbf{0}_m]$ for some nonnegative integers $t,m$, where \[J_i=\text{Dg}[H(x^{e_i}),\mathbf{0}_{e_i-2}]=\text{Dg}[J_{e_i}(0),\mathbf{0}_{e_i-2}]\] for each $i \in \{1,2,\ldots,t\}$. This is the same form as used for the blocks $G_i$ where $\lambda_i=0$ as in the proof of the specialized version, and therefore the same factorization into a product of three nilpotent factors is applicable here.

It remains therefore to provide a factorization of the block $\text{Dg}[\mathbf{0}_{\text{r}(B)}, B]$:
\[\text{Dg}[\mathbf{0}_{\text{r}(B)}, B] = \begin{bmatrix} \mathbf{0}_{\text{r}(B)} & \mathbf{0}_{\text{r}(B)} \\ \mathbf{0}_{\text{r}(B)} & B\end{bmatrix} = \begin{bmatrix} \mathbf{0}_{\text{r}(B)} & \mathbf{0}_{\text{r}(B)} \\ I_{\text{r}(B)} & \mathbf{0}_{\text{r}(B)}\end{bmatrix} \begin{bmatrix} I_{\text{r}(B)} & -I_{\text{r}(B)} \\ I_{\text{r}(B)} & -I_{\text{r}(B)}\end{bmatrix} \begin{bmatrix} \mathbf{0}_{\text{r}(B)} & B \\ \mathbf{0}_{\text{r}(B)} & \mathbf{0}_{\text{r}(B)}\end{bmatrix}.\]
It is easy to see that each of the matrices on the right hand side is square-zero. 

To conclude, it was shown above that $G$ is similar to a product of three square-zero matrices. Therefore, by invariance of the minimum polynomial of a matrix under similarity, $G$ can be written as a product of three square-zero matrices (using an argument similar to the one given in the last paragraph of the proof of Theorem \ref{thm:botha2FactorSZ1}). \hfill $\square$

\subsection{Matrix division with a square-zero quotient \label{seq:sqzerodiv}}
I now present a general result pertaining to matrix division with a square-zero quotient as proved by Botha \cite{botha2}. Since the quotient of interest is square, we must have that the orders of the matrix $G$ and its divisor $F$ are the same. Following the theory of matrix division developed in the preliminary results, I will present the results in a format which will allow the reader to follow easily from the previous results. Specifically, I split the theorem into parts which are equivalent, and agree with the general format of Proposition \ref{prop:matrixDivision}. Parts one and two of this theorem provide several tests (similar to Proposition \ref{prop:matrixDivision}) whereby it is possible to assess whether two matrices of the same order have a square-zero quotient. 
\begin{quote}
\begin{thm}[Part 1. Theorem 1 \cite{botha2}] Let $G,F \in M_{m \times n}(\mathscr{F})$, where $\mathscr{F}$ is a field. The following statements are equivalent:
\begin{enumerate}[(a)]
\item There exists a matrix $H \in M_{m}(\mathscr{F})$ such that $\begin{bmatrix} G & \mathbf{0}_{m \times n} \end{bmatrix}=H\begin{bmatrix} F & G \end{bmatrix}$.
\item $\text{row} \left (\begin{bmatrix} G & \mathbf{0}_{m \times n} \end{bmatrix}\right ) \subseteq \text{row}\left (\begin{bmatrix} F & G \end{bmatrix} \right ).$
\item $\text{N} \left (\begin{bmatrix} F & G \end{bmatrix} \right ) \subseteq \text{N}\left (\begin{bmatrix} G & \mathbf{0}_{m \times n} \end{bmatrix} \right ).$
\item $\text{r}\left (\begin{bmatrix} F & G \end{bmatrix} \right )=\text{r}\left ( \begin{bmatrix} G & \mathbf{0}_{m \times n} \\ F & G  \end{bmatrix} \right ).$
\end{enumerate}
\label{thm:bothaSzDivisorPt1}
\end{thm}
\end{quote}  
{\bf Proof.} Equivalence of all statements follow directly from \eqref{prop1pt1} of Proposition \ref{prop:matrixDivision}. \hfill $\square$

Now, by Proposition \ref{prop:matrixDivision}, it also follows directly that a square matrix $H \in M_{m}(\mathscr{F})$ is a quotient of $G$ and $F$ if and only if $\text{row}(G) \subseteq \text{row}(F)$ if and only if $\text{N}(F) \subseteq \text{N}(G)$ if and only if \[\text{r}(F)=\text{r}\left ( \begin{bmatrix} G  \\ F   \end{bmatrix} \right ).\]
If we have the added condition that $H$ is square-zero, then more can be said, in particular all the statements in Theorem \ref{thm:bothaSzDivisorPt1} above are then also equivalent to the statements above. Furthermore additional conditions can be placed on $\text{R}(G) \cap \text{R}(F)$ and $\text{r}(\begin{bmatrix}G & F\end{bmatrix})$. These results are summarized in part two:
\begin{quote}
\begin{thm}[Part 2. Theorem 1 \cite{botha2}] Let $G,F \in M_{m \times n}(\mathscr{F})$, where $\mathscr{F}$ is a field. The following statements are equivalent:
\begin{enumerate}[(a)]
\item There exists a square-zero matrix $H \in M_{m}(\mathscr{F})$ such that $G=HF$.
\item $\begin{bmatrix} G & \mathbf{0}_{m \times n} \end{bmatrix}$ is divisible by $\begin{bmatrix} F & G \end{bmatrix}$ on the right.
\item $\text{N}(F) \subseteq \text{N}(G)$ and $\text{R}(G) \cap \text{R}(F) \subseteq \text{R}(F \restriction \text{N}(G))$.
\item $\text{r}(F)=\text{r}\left ( \begin{bmatrix} G  \\ F   \end{bmatrix} \right )$ and 
\[\text{r}\begin{bmatrix}G & F \end{bmatrix} = \text{r}(G) + \text{r}\left (\begin{bmatrix} G & FB\end{bmatrix}\right ),\] where $B \in M_{n \times \text{n}(G)}(\mathscr{F})$ is such that $\text{R}(B)=\text{N}(G)$.
\end{enumerate}
\label{thm:bothaSzDivisorPt2}
\end{thm}
\end{quote}  
{\bf Proof.} {\bf (a) implies (b).} If (a) is true, then \[H\begin{bmatrix} F & G \end{bmatrix} = \begin{bmatrix} HF & HG \end{bmatrix}= \begin{bmatrix} HF & H^2F \end{bmatrix}=\begin{bmatrix} G & \mathbf{0}_{m \times n} \end{bmatrix}.\] Notice that it is equivalence of (a) and (b) which asserts that parts one and two of the proof are equivalent. 

{\bf (b) implies (c).} The first part, $\text{N}(F) \subseteq \text{N}(G)$, follows from the fact that there is some $K \in M_m(\mathscr{F})$ so that $\begin{bmatrix} G & \mathbf{0}_{m \times n} \end{bmatrix}=K\begin{bmatrix} F & G \end{bmatrix}$ (implying that $G=KF$), and Proposition \ref{prop:matrixDivision} (c). Now suppose $u \in \text{R}(G) \cap \text{R}(F)$. Then there exist $v,w \in \mathscr{F}^n$ so that $u=Gv=Fw$. Now imposing the condition as stated in (b), there exists some $K \in M_m(\mathscr{F})$ so that $\mathbf{0}=KGv=KFw=Gw$. So $w \in \text{N}(G)$, and therefore $Fw=u \in \text{R}(F \restriction \text{N}(G))$.

{\bf (c) implies (a).} Notice that, given that statement (c) is true, we have $\text{N}(F) \subseteq \text{N}(G)$, which guarantees that $F$ is a right divisor of $G$ by Proposition \ref{prop:matrixDivision}(c). The rest of the proof consists of showing that given the additional condition $\text{R}(G) \cap \text{R}(F) \subseteq \text{R}(F \restriction \text{N}(G))$, there exists a square-zero quotient of $G$ and $F$.

So given the conditions in (c), if we have $Gv+Fw=\mathbf{0}$, it means that $-Gv=Fw=u \in \text{R}(G) \cap \text{R}(F)$, so that (by hypothesis) $u \in \text{R}(F \restriction \text{N}(G))$. So it follows that \begin{equation}w \in \text{N}(G). \label{th3pt2ca} \end{equation} 

Now within this context, consider the direct sum $\mathscr{F}^m=(\text{R}(G)+\text{R}(F)) \oplus V$, which can be constructed by extending a basis for $(\text{R}(G)+\text{R}(F))$ to a basis for $\mathscr{F}^m$. Define \[H_1:(\text{R}(G)+\text{R}(F)) \rightarrow (\text{R}(G)+\text{R}(F)) \text{ by } H_1(Gv+Fw)=Gw.\] Let $u_1=Gv_1+Fw_1$ and $u_2=Gv_2+Fw_2$, which are both in $(\text{R}(G)+\text{R}(F))$. To prove that $H_1$ is functionally well-defined suppose $u_1=u_2$. Then by the linearity of $G$ and $F$ we must have $G(v_1-v_2)=F(w_2-w_1)$. But then $w_2-w_1 \in \text{N}(G)$ by~\eqref{th3pt2ca}, and therefore \[H_1(Gv_1+Fw_1)-H_1(Gv_2+Fw_2) = Gw_1-Gw_2 = \mathbf{0} \Leftrightarrow H_1(Gv_1+Fw_1)=H_1(Gv_2+Fw_2).\]
Furthermore, to prove that $H_1$ is linear, for any $k_1,k_2 \in \mathscr{F}$ it also follows easily by the linearity of $G$ and $F$: 
\begin{eqnarray} 
\nonumber H_1(k_1u_1+k_2u_2) &=& H_1(G(k_1v_1+k_2v_2)+F(k_1w_1+k_2w_2))\\
\nonumber &=&G(k_1w_1 + k_2w_2)\\
\nonumber &=&k_1Gw_1 + k_2Gw_2\\
\nonumber &=&k_1H_1(u_1)+k_2H_1(u_2),\end{eqnarray} and therefore $H_1$ is linear. Finally, for any $u = Gv+Fw \in (\text{R}(G)+\text{R}(F))$ we have \[H_1^2(u)=H_1(Gw)=\mathbf{0},\] so that $H_1$ is square-zero.

Now define $H: \mathscr{F}^m \rightarrow \mathscr{F}^m$ by $H \restriction (\text{R}(G) + \text{R}(F)) = H_1$ and $H \restriction V = H_2: V \rightarrow V$ where $H_2$ is an arbitrary square-zero operator on $V$. At least one such an operator $H_2$ exists, namely the zero transformation, sending all vectors in $V$ to $\mathbf{0}$. This proves existence of a square-zero matrix $H \in M_m(\mathscr{F})$ such that $G=HF$.

Having proved the implications above, notice that at this point we have proved equivalence of all statements in part one to statements (a), (b) and (c) of part two. It remains only to prove that statement (d) of part two is equivalent to all other statements in part two. 

{\bf (c) if and only if (d).} I provide one remark before proceeding with a proof which closely follows the reference text \cite[Theorem 1]{botha2}. The purpose of the matrix~$B$, introduced in this part of the theorem, is to restrict the domain of $F$ to $\text{N}(G)$. The order of $B$ in the theorem is fixed to $n \times \text{n}(G)$, which would force $B$ to be full column rank. This requirement does not seem to be necessary, as for example, it is possible to define $B$ as a $n \times n$ projection, projecting vectors in $\mathscr{F}^n$ onto $N(G)$. It is therefore possible to modify statement (d) so that the only restrictions placed on~$B$ is that $B$ must have $n$ rows and $\text{R}(B)=\text{N}(G)$. 

The first part \[\text{N}(F) \subseteq \text{N}(G) \text{ if and only if } \text{r}(F)=\text{r}\left ( \begin{bmatrix} G  \\ F   \end{bmatrix} \right )\] follows directly from Proposition \ref{prop:matrixDivision}(c), (d). Now notice that \[\text{R}\left ( \begin{bmatrix} G  & F   \end{bmatrix} \right )=\text{R}(G) + \text{R}(F),\] since the column space of $\begin{bmatrix} G  & F   \end{bmatrix}$ consists of the space spanned by the columns of $G$ and $F$, and the columns of $G$ are not necessarily linearly independent from the columns in $F$. So we have by \cite[Proposition 5.16]{golan}
\begin{equation}\text{r}\left ( \begin{bmatrix} G  & F   \end{bmatrix} \right )=\text{r}(G) + \text{r}(F) - \dim(\text{R}(G)\cap \text{R}(F)). \label{th2cd1} \end{equation}
Furthermore, since by Proposition \ref{prop:matrixDivision}, if either (c) or (d) holds, then $\text{N}(F) \subseteq \text{N}(G)$, and therefore 
\begin{equation} \text{r}(F \restriction \text{N}(G))=\text{n}(G) - \text{n}(F). \label{th2cd3}\end{equation} 
Also, by a similar argument as used in formulating \eqref{th2cd1}, and recalling that the purpose of $B$ is to restrict the domain of $F$ to $\text{N}(G)$ we have
\[\text{r}\left ( \begin{bmatrix} G  & FB   \end{bmatrix} \right )=\text{r}(G) + \text{r}(F \restriction \text{N}(G)) - \dim(\text{R}(G)\cap \text{R}(F \restriction \text{N}(G))). \]
Now adding $\text{r}(G)$ to both sides and combining with \eqref{th2cd3} and the rank-nullity theorem \cite[Theorem 2.11]{cullen}, we have
\begin{eqnarray}
\nonumber \text{r}(G)+\text{r}\left ( \begin{bmatrix} G  & FB   \end{bmatrix} \right )&=&\text{r}(G)+\text{r}(G) + \text{r}(F \restriction \text{N}(G)) - \dim(\text{R}(G)\cap \text{R}(F \restriction \text{N}(G)))\\
\nonumber &=& \text{r}(G)+\text{r}(G) + \text{n}(G)-\text{n}(F) - \dim(\text{R}(G)\cap \text{R}(F \restriction \text{N}(G)))\\
\nonumber &=& \text{r}(G)+n-\text{n}(F) - \dim(\text{R}(G)\cap \text{R}(F \restriction \text{N}(G)))\\
&=& \text{r}(G)+\text{r}(F) - \dim(\text{R}(G)\cap \text{R}(F \restriction \text{N}(G))). \label{th2cd2}
\end{eqnarray}
Now suppose (c) holds. Then $\text{R}(G) \cap \text{R}(F) \subseteq \text{R}(F \restriction \text{N}(G))$, which implies that $\text{R}(G) \cap \text{R}(F) = \text{R}(G) \cap \text{R}(F \restriction \text{N}(G))$: if $x \in \text{R}(G) \cap \text{R}(F)$ then $x \in \text{R}(G)$ and $x \in \text{R}(F \restriction \text{N}(G))$ since (c) holds, so $x \in \text{R}(G) \cap \text{R}(F \restriction \text{N}(G))$. The reverse inclusion follows easily from the fact that $\text{R}(F \restriction \text{N}(G)) \subseteq \text{R}(F)$.  Combining this result with \eqref{th2cd1} and \eqref{th2cd2} yields
\[\text{r}\left ( \begin{bmatrix} G  & F   \end{bmatrix} \right ) = \text{r}(G) + \text{r}\left ( \begin{bmatrix} G  & FB   \end{bmatrix} \right ),\] proving that (d) holds.

For the converse we can apply the same argument in reverse, that is, we have by \eqref{th2cd1} and \eqref{th2cd2} that $\dim(\text{R}(G) \cap \text{R}(F)) = \dim(\text{R}(G) \cap \text{R}(F \restriction \text{N}(G)))$ and since $\text{R}(F \restriction \text{N}(G)) \subseteq \text{R}(F)$ we must therefore have \[\text{R}(G) \cap \text{R}(F) = \text{R}(G) \cap \text{R}(F \restriction \text{N}(G))\subseteq \text{R}(F \restriction \text{N}(G)),\] proving that (c) holds. \hfill $\square$

This concludes the proof of part two, so that at this point it has been proved that all statements in part one and two are equivalent. The next part characterizes the square-zero quotient $H$ in terms of its rank. Before proceeding to part three and its proof, we require the following lemma:

\begin{quote}
\begin{lem}[Lemma 1 \cite{botha2}] 
Let $G,F \in M_{m \times n}(\mathscr{F})$, where $\mathscr{F}$ is a field, and let
$\mathscr{F}^m = \text{R}\left (\begin{bmatrix}G & F\end{bmatrix} \right) \oplus \text{R}(C)$ and $\text{R}(G) + \text{R}(F \restriction \text{N}(G)) = \text{R}(G) \oplus \text{R}(B)$,
where both $B \in M_{m \times *}(\mathscr{F})$ and $C \in M_{m \times *}(\mathscr{F})$ are of full column rank. If $H$ is a square-zero right quotient of $G$ and $F$, then
\begin{enumerate}[(a)]
\item $\text{R}\left (\begin{bmatrix}G & F\end{bmatrix} \right )$ is $H$-invariant, and $H \restriction \text{R}\left (\begin{bmatrix}G & F\end{bmatrix} \right )$ is uniquely defined by \[ \left (H \restriction \text{R} \left (\begin{bmatrix}G & F\end{bmatrix}\right ) \right ) (Gv + Fw) = Gw \text{, for all } v,w \in \mathscr{F}^n; \]
\item $\text{N}\left (H \restriction \text{R}\left (\begin{bmatrix}G & F\end{bmatrix}\right ) \right )=\text{R}(G) + \text{R}(F \restriction \text{N}(G))$, \\
$\text{R}\left (H \restriction \text{R}\left (\begin{bmatrix}G & F\end{bmatrix}\right )\right )=\text{R}(G)$, and \[\text{n}\left (H \restriction \text{R} \left (\begin{bmatrix}G & F\end{bmatrix}\right )\right )= \text{r}\left (\begin{bmatrix}G & F\end{bmatrix}\right ) - \text{r}(G);\]
\item $\frac{m}{2}-\left (\text{r}\left (\begin{bmatrix}G & F\end{bmatrix}\right ) - \text{r}(G)\right )=\frac{1}{2}(\text{r}(C)-\text{r}(B)).$
\end{enumerate}
\label{lem:bothaSzDivisorPt3}
\end{lem}
\end{quote} 
 
{\bf Proof.} {\bf (a).} Let $u \in \text{R}\left (\begin{bmatrix}G & F\end{bmatrix} \right )$. Then there exist $v,w \in \mathscr{F}^n$ so that \[u = \begin{bmatrix}G & F\end{bmatrix}\begin{bmatrix}v\\w\end{bmatrix} = Gv + Fw.\] Now if $G=HF$, and $H$ is square-zero, then we must have 
\begin{eqnarray}
\nonumber Hu &=&  H\begin{bmatrix}G & F\end{bmatrix}\begin{bmatrix}v\\w\end{bmatrix}\\
\nonumber &=& \begin{bmatrix}HG & HF\end{bmatrix}\begin{bmatrix}v\\w\end{bmatrix} \\
\nonumber &=& \begin{bmatrix}H^2F & HF\end{bmatrix}\begin{bmatrix}v\\w\end{bmatrix} \\
\nonumber &=& \begin{bmatrix}\mathbf{0}_{m \times n} & G\end{bmatrix}\begin{bmatrix}v\\w\end{bmatrix}\\
\ &=& Gw. \label{lem1a} \end{eqnarray} 
Now $Gw \in \text{R}(G) \subseteq \text{R}\left (\begin{bmatrix}G & F\end{bmatrix}\right )$, and therefore $\text{R}\left (\begin{bmatrix}G & F\end{bmatrix}\right )$ is $H$-invariant. Furthermore, by the equality expressed in \eqref{lem1a} it directly follows that 
\[ \left (H \restriction \text{R} \left (\begin{bmatrix}G & F\end{bmatrix}\right ) \right) (Gv + Fw) = Gw \text{, for all } v,w \in \mathscr{F}^n, \] as required.

{\bf (b).} By \eqref{lem1a} it is easy to see that $u = Gv+Fw \in \text{N}\left (H \restriction \text{R}\left (\begin{bmatrix}G & F\end{bmatrix}\right )\right )$ if and only if $w \in \text{N}(G)$, and therefore $\text{N}\left (H \restriction \text{R}\left (\begin{bmatrix}G & F\end{bmatrix}\right )\right )=\text{R}(G) + \text{R}(F \restriction \text{N}(G))$, proving the first part of (b).   

Now since $H \restriction \text{R}\left (\begin{bmatrix}G & F\end{bmatrix}\right )$ is uniquely defined as specified in (a), we have for any $u \in \text{R}\left (\begin{bmatrix}G & F\end{bmatrix}\right )$ that $Hu \in \text{R}(G)$, that is $\text{R}\left (H \restriction \text{R}\left (\begin{bmatrix}G & F\end{bmatrix}\right )\right ) \subseteq \text{R}(G)$ (this can be seen directly from \eqref{lem1a}). Now let $Gx$ be any vector in $\text{R}(G)$, then $H \restriction \text{R}\left (\begin{bmatrix}G & F\end{bmatrix}\right )(Fx)=Gx$ so that $Gx \in \text{R}\left (H \restriction \text{R}\left (\begin{bmatrix}G & F\end{bmatrix}\right )\right )$ and it follows that $\text{R}(G) \subseteq \text{R}\left (H \restriction \text{R}\left (\begin{bmatrix}G & F\end{bmatrix}\right )\right )$. Combining the preceding arguments, it follows that
\begin{equation}\text{R}\left (H \restriction \text{R}\left (\begin{bmatrix}G & F\end{bmatrix}\right )\right ) = \text{R}(G). \label{lem1ab}\end{equation}

Now since $\text{R}\left (\begin{bmatrix}G & F\end{bmatrix}\right )$ is $H$-invariant, by the rank-nullity theorem \cite[Theorem 2.11, p.87]{cullen} we must have 
\[\text{r}\left (\begin{bmatrix}G & F \end{bmatrix}\right )=\text{r}\left (H \restriction \text{R}\left (\begin{bmatrix}G & F\end{bmatrix}\right )\right )+\text{n}\left (H \restriction \text{R}\left (\begin{bmatrix}G & F\end{bmatrix}\right )\right ),\] and therefore by \eqref{lem1ab} 
\[\text{r}\left (\begin{bmatrix}G & F \end{bmatrix}\right )=\text{r}(G)+\text{n}\left (H \restriction \text{R}\left (\begin{bmatrix}G & F\end{bmatrix}\right )\right ),\] which produces the desired result upon rearrangement:
\[\text{n}\left (H \restriction \text{R}\left (\begin{bmatrix}G & F\end{bmatrix}\right )\right ) = \text{r}\left (\begin{bmatrix}G & F \end{bmatrix}\right ) - \text{r}(G).\]
{\bf (c).} Suppose $B$ and $C$ are defined as in the statement of this lemma, and that the conditions in (b) hold. Then
\begin{eqnarray} \nonumber \text{r}(G)+\text{r}(B)&=&\dim(\text{R}(G) + \text{R}(F \restriction \text{N}(G)))\\
\nonumber &=&\text{n}\left (H \restriction \text{R}\left (\begin{bmatrix}G & F\end{bmatrix}\right )\right )\\
\nonumber &=&\text{r}\left (\begin{bmatrix}G & F \end{bmatrix} \right )-\text{r}(G), \label{lem1bc}\end{eqnarray}
and \begin{equation} \nonumber m=\text{r}\left (\begin{bmatrix}G & F \end{bmatrix}\right )+\text{r}(C). \label{lem1bc2}\end{equation}

Now by the equations above
\begin{eqnarray} 
\nonumber && \frac{m}{2}-\left (\text{r} \left (\begin{bmatrix}G & F \end{bmatrix} \right)-\text{r}(G) \right) \\
\nonumber &=& \frac{1}{2} \left (\text{r}\left(\begin{bmatrix}G & F \end{bmatrix}\right )+\text{r}(C)\right ) - \left (\text{r}(G)+\text{r}(B) \right) \\
\nonumber &=& \frac{1}{2}\left (\text{r}\left(\begin{bmatrix}G & F \end{bmatrix}\right )-\text{r}(G)\right) - \frac{1}{2}\left (\text{r}(G)+\text{r}(B)\right ) + \frac{1}{2}\left (\text{r}(C)-\text{r}(B) \right ) \\
\nonumber &=& \frac{1}{2}(\text{r}(G)+\text{r}(B)) - \frac{1}{2}(\text{r}(G)+\text{r}(B)) + \frac{1}{2}(\text{r}(C)-\text{r}(B)) \\
\nonumber &=& \frac{1}{2}(\text{r}(C)-\text{r}(B)),
\end{eqnarray}
which completes the proof of lemma \ref{lem:bothaSzDivisorPt3}. \hfill $\square$
\begin{quote}
\begin{thm}[Part 3. Theorem 1 \cite{botha2}] If any of the conditions in Theorem \ref{thm:bothaSzDivisorPt1} or \ref{thm:bothaSzDivisorPt2} hold, then there exists a square-zero quotient $H$ of each rank such that $\text{r}(G) \leq \text{r}(H)$ and 
\begin{equation} \text{r}(H) \leq \left \{ {\renewcommand{\arraystretch}{1.2} \begin{array}{ll} \frac{m}{2} & \text{if r}\left ({\renewcommand{\arraystretch}{1}\begin{bmatrix} G & F \end{bmatrix}}\right) - \text{r}(G) \leq \frac{m}{2} \\  
m- \left (\text{r} \left ({\renewcommand{\arraystretch}{1}\begin{bmatrix} G & F \end{bmatrix}}\right ) - \text{r}(G)\right ) & \text{if r}\left ({\renewcommand{\arraystretch}{1}\begin{bmatrix} G & F \end{bmatrix}} \right ) - \text{r}(G) > \frac{m}{2} \end{array}} \right. \label{th3pt3} \end{equation}
\label{thm:bothaSzDivisorPt3}
\end{thm} 
\end{quote}  
This part of the theorem is a modification, or extension, of \eqref{prop1pt2} due to the fact that $H$ is square-zero. As in the proof of \eqref{prop1pt2} Botha makes use of a constructive approach to prove existence of a matrix $H$ of each allowable rank, and then verifies universality by showing there is no matrix $H$ outside the bounds, but in this case since the upper bound depends on an extra condition two cases are considered separately. The proof very much depends on Lemma \ref{lem:bothaSzDivisorPt3}. I will elucidate the proof by explaining in detail how the construction of $H$ is composed of different mappings on $H$-invariant subspaces as given in the reference text. 

{\bf Proof.} If any of the conditions in Theorem \ref{thm:bothaSzDivisorPt1} or \ref{thm:bothaSzDivisorPt2} hold, it follows that there exists a matrix $H$, which is a square-zero right quotient of $G$ and $F$. The lower bound of \eqref{th3pt3} then follows directly from \eqref{prop1pt2}. Now, we will make use of Lemma~\ref{lem:bothaSzDivisorPt3}, by decomposing $\mathscr{F}^m$ into $H$-invariant subspaces which makes it possible to specify the rank of $H$ by its actions on each. So let \[\mathscr{F}^m=\text{R}\left (\begin{bmatrix} G & F \end{bmatrix} \right ) \oplus \text{R}(C)\] where $C$ is a matrix as specified in the statement of Lemma \ref{lem:bothaSzDivisorPt3}. Since $H$ is uniquely defined on $\text{R} \left (\begin{bmatrix} G & F \end{bmatrix} \right)$ by Lemma \ref{lem:bothaSzDivisorPt3}(a), and furthermore by \eqref{lem1ab} we have \[\text{r}\left (H \restriction \text{R}\left (\begin{bmatrix}G & F\end{bmatrix}\right ) \right ) = \text{r}(G),\] there is no scope for varying the rank of $H$ on the subspace $\text{R}\left (\begin{bmatrix} G & F \end{bmatrix} \right )$. 

It remains to vary the rank of $H$, by specifying its action on vectors in $\text{R}(C)$. Since $H$ is square-zero we have the additional constraint that for any nonzero vector $v \in \text{R}(C)$ such that $H(v) \neq \mathbf{0}$, we must have $H^2(v) = \mathbf{0}$. So we must have $\text{R}(H \restriction \text{R}(C)) \subseteq \text{N}(H)$. At this point we consider two possible cases.

\begin{enumerate}[(i)]
\item \[\text{r} \left (\begin{bmatrix} G & F \end{bmatrix} \right ) - \text{r}(G) \leq m/2.\] In this case, by Lemma \ref{lem:bothaSzDivisorPt3}(c) we must have $\text{r}(C) \geq \text{r}(B)$. The core argument of the existential part of the proof is then to show that 
\[\mathscr{F}^m = \mathcal{W} \oplus \mathcal{V}_2 =  \text{R} \left (\begin{bmatrix}G & F\end{bmatrix} \right ) \oplus \mathcal{V}_1 \oplus \mathcal{V}_2\] where $\mathcal{W}$ and $\mathcal{V}_2$ are $H$-invariant subspaces. The particular choice of these subspaces then allows $H$ to be specified in such a way that $\text{r}(H \restriction \mathcal{W})$ is at most $\dim(\mathcal{W})/2$ and $\text{r}(H \restriction \mathcal{V}_2)$ is at most $\dim(\mathcal{V}_2)/2$. Combining these two results we then have the desired upper bound in \eqref{th3pt3}.

Proceeding with the proof, since $\text{r}(C) \geq \text{r}(B)$ we can decompose $\text{R}(C)$ into $\mathcal{V}_1 \oplus \mathcal{V}_2$ where $\dim(\mathcal{V}_1)=\text{r}(B)$. Now specify $H \restriction \mathcal{V}_1$ as an arbitrary linear transformation $H_1: \mathcal{V}_1 \rightarrow \text{R}(B)$ and let \[\mathcal{W} = \text{R}\left (\begin{bmatrix}G & F\end{bmatrix} \right) \oplus \mathcal{V}_1.\]
Now by Lemma \ref{lem:bothaSzDivisorPt3}(b) and specifically \eqref{lem1ab} we have \[\text{R} \left (H \restriction \text{R}\left (\begin{bmatrix}G & F\end{bmatrix}\right ) \right )=\text{R}(G),\] and by the definition of $H \restriction \mathcal{V}_1$ above it follows that $\text{R}(H \restriction \mathcal{V}_1) \subseteq \text{R}(B)$ and therefore
\[\text{R}(H \restriction \mathcal{W})=\text{R}\left (H \restriction \text{R}\left (\begin{bmatrix}G & F\end{bmatrix}\right ) \right ) \oplus \text{R}(H \restriction \mathcal{V}_1) \subseteq \text{R}(G) \oplus \text{R}(B),\] 
and therefore 
\[\text{r}(H \restriction \mathcal{W})=\text{r}\left (H \restriction \text{R}\left (\begin{bmatrix}G & F\end{bmatrix}\right )\right ) +  \text{r}(H \restriction \mathcal{V}_1) \leq \text{r}(G) + \text{r}(B).\]
If we now choose $H_1$ to be a bijection, we have $\text{r}(H \restriction \mathcal{W}) = \text{r}(G)+\text{r}(B) = \text{n}(H \restriction \mathcal{W})$, since 
\[\text{N}(H \restriction \mathcal{W}) = \text{N}\left (H \restriction \text{R}\left (\begin{bmatrix}G & F\end{bmatrix} \right ) \right ) \oplus \text{N}(H \restriction \mathcal{V}_1) = \text{N}\left (H \restriction \text{R}\left (\begin{bmatrix}G & F\end{bmatrix}\right ) \right ) \oplus \{\mathbf{0}\},\]
and by Lemma \ref{lem:bothaSzDivisorPt3}(b)
\[\text{N}\left (H \restriction\text{R}\left (\begin{bmatrix}G & F\end{bmatrix}\right )\right )=\text{R}(G) \oplus \text{R}(B).\]
Therefore by the rank-nullity theorem \cite[Theorem 2.11]{cullen}, $\text{r}(H \restriction \mathcal{W}) = \dim(\mathcal{W})/2$ is the greatest possible rank of $H \restriction \mathcal{W}$. Notice that by construction $\dim(\mathcal{W})$ is even: 
\begin{eqnarray}
\nonumber \dim(\mathcal{W})&=&\text{r}\left (H \restriction \text{R} \left (\begin{bmatrix}G & F\end{bmatrix}\right )\right ) + \text{n}\left (H \restriction \text{R}\left (\begin{bmatrix}G & F\end{bmatrix}\right )\right )+\text{r}(B)\\
\nonumber &=& \text{r}(G) + (\text{r}(G)+\text{r}(B))+\text{r}(B)\\
 &=& 2(\text{r}(G)+\text{r}(B)). \label{eq:dimWsqzQ}\end{eqnarray} 
Furthermore notice that since $\text{R}(H_1) \subseteq \text{R}(B) \subseteq \text{N}(H)$, it follows that $H \restriction \mathcal{W}$ is square zero.

It remains to show that it is possible to define $H \restriction \mathcal{V}_2$ so that the maximum value that $\dim(\text{R}(H \restriction \mathcal{V}_2))$ can achieve is $\dim(\mathcal{V}_2)/2$ with $\mathcal{V}_2$ being an $H$-invariant subspace (and of course $H \restriction \mathcal{V}_2$ must be square-zero). Note that it is also possible that $\mathcal{V}_2=\{\mathbf{0}\}$ ($\text{r}(C)=\text{r}(B)$) in which case this construction is not necessary. Suppose therefore $\mathcal{V}_2 \neq \{\mathbf{0}\}$, and define $\alpha=\{\alpha_1, \alpha_2, \ldots, \alpha_t\}$ as a basis for $\mathcal{V}_2$, and define \[H \restriction \mathcal{V}_2 = H_2 : \mathcal{V}_2 \rightarrow \mathcal{V}_2,\] by choosing a nonnegative integer $c \leq \lfloor \frac{t}{2} \rfloor$ and let 
\[ H_2(\alpha_i) =  \left \{ {\renewcommand{\arraystretch}{1.2} \begin{array}{ll} \alpha_{c+i} & \text{if } 1 \leq i \leq c \\
\mathbf{0} & \text{if } c < i \leq t  
\end{array}} \right.. \]
It is easy to see that $H_2$ is then square-zero and $t/2=\dim(\mathcal{V})/2$ is the supremum of the set of all possible ranks that $H_2$ may assume.

Now, with $H$ defined as 
\begin{enumerate}
\item $H \restriction \text{R}\left (\begin{bmatrix}G & F\end{bmatrix}\right )$ is fixed as per lemma \ref{lem:bothaSzDivisorPt3}(a),
\item $H \restriction \mathcal{V}_1 = H_1$,
\item $H \restriction \mathcal{V}_2 = H_2$,
\end{enumerate} we could choose $H_1$ to be the zero transformation, and $c=0$ which results in $H_2$ being the zero transformation, and then \[\text{r}(H)=\text{r}\left (H \restriction \text{R}\left (\begin{bmatrix}G & F\end{bmatrix}\right )\right )=\text{r}(G),\] achieving the lower bound in \eqref{th3pt3}. Now we can vary the rank of $H$ from this lower bound up to the upper bound in \eqref{th3pt3}, which is achieved by choosing $H_1$ as a bijection between $\mathcal{V}_1$ and $\text{R}(B)$ (as discussed before), and let $c=\lfloor \frac{t}{2} \rfloor$. Then $\text{r}(H \restriction \mathcal{W})=\dim(\mathcal{W})/2$ and $\text{r}(H_2)=\lfloor \frac{\dim(\mathcal{V}_2)}{2} \rfloor$ so that $\text{r}(H) = \lfloor \frac{m}{2} \rfloor$. We can summarize the preceding arguments regarding the upper bound:
\begin{eqnarray}
\nonumber \text{r}(H) &=& \text{r}\left (H \restriction \text{R}\left (\begin{bmatrix}G & F\end{bmatrix}\right )\right ) + \text{r}(H_1) + \text{r}(H_2)\\
\nonumber &=& \text{r}(G) + \text{r}(H \restriction \mathcal{V}_1) + c \\
\nonumber &\leq& \text{r}(G) + \text{r}(B) +  \left \lfloor \frac{\dim(\mathcal{V}_2)}{2} \right \rfloor\\
\nonumber &=& \frac{\dim(\mathcal{W})}{2} +  \left \lfloor \frac{\dim(\mathcal{V}_2)}{2} \right \rfloor \\
\nonumber &=& \left \lfloor \frac{\dim(\mathcal{W}) +  \dim(\mathcal{V}_2)}{2} \right \rfloor \\
\nonumber &=& \left \lfloor \frac{m}{2} \right \rfloor. \end{eqnarray}
This proves that, in the case of  $\text{r}\left (\begin{bmatrix} G & F \end{bmatrix}\right ) - \text{r}(G) \leq m/2$ there exists a class of square-zero right quotients of $G$ and $F$, of every rank within the bounds in~\eqref{th3pt3}. 

The universality of the lower bound as specified in \eqref{th3pt3} again follows directly from \eqref{prop1pt2}. Universality of the upper bound in \eqref{th3pt3} is a direct consequence of $H$ being square-zero: if $H^2=\mathbf{0}$ then $\text{R}(H) \subseteq \text{N}(H)$, and therefore, by the rank-nullity theorem \cite[Theorem 2.11]{cullen}
\[m=\text{r}(H)+\text{n}(H)\geq \text{r}(H)+\text{r}(H) =2\text{r}(H).\]

\item \[\text{r}\left (\begin{bmatrix} G & F \end{bmatrix} \right ) - \text{r}(G) > m/2.\] In this case, by Lemma \ref{lem:bothaSzDivisorPt3}(c) we must have $\text{r}(C) < \text{r}(B)$. Now by further conditions as stated in lemma \ref{lem:bothaSzDivisorPt3}(b), we must have 
\[\text{n}\left (H \restriction \text{R}\left (\begin{bmatrix} G & F \end{bmatrix}\right )\right ) = \text{r}(G) + \text{r}(B) >  \text{r}(G) + \text{r}(C)=\text{r}\left (H \restriction \text{R}\left (\begin{bmatrix} G & F \end{bmatrix}\right )\right ) + \text{r}(C),\] and so regardless of the choice of $H \restriction \text{R}(C)$ we have $\text{n}(H) > \text{r}(H)$. Now there is a natural choice for $H \restriction \text{R}(C)$, which ensures that $H$ is square-zero, and satisfies the upper bound in \eqref{th3pt3}. Explicitly, let $H \restriction \text{R}(C) = H_3: \text{R}(C) \rightarrow \text{R}(B)$ be an arbitrary linear transformation. As for the previous case when considering the subspace $\mathcal{W}$ above, we now have 
\[H(\mathscr{F}^m)=\text{R}\left (H \restriction \text{R}\left (\begin{bmatrix} G & F \end{bmatrix}\right )\right ) + \text{R}(H \restriction \text{R}(C))\subseteq \text{R}(G) \oplus \text{R}(B),\]
and by \eqref{lem1ab} we have $\text{R}\left (H \restriction \text{R}\left (\begin{bmatrix} G & F \end{bmatrix}\right )\right ) = \text{R}(G)$, and by definition $\text{R}(H\restriction \text{R}(C)) \subseteq \text{R}(B)$. So again, it follows that the sum $\text{R}\left (H \restriction \text{R}\left (\begin{bmatrix} G & F \end{bmatrix}\right )\right ) + \text{R}(H\restriction \text{R}(C))$ is in fact a direct sum and therefore 
\[\text{r}(H)=\text{r}\left (H \restriction \text{R}\left (\begin{bmatrix} G & F \end{bmatrix}\right )\right )+\text{r}(H \restriction \text{R}(C)).\] 
So, with $H$ defined as
\begin{enumerate}
\item $H \restriction \text{R}\left (\begin{bmatrix}G & F\end{bmatrix}\right )$ is fixed as per Lemma \ref{lem:bothaSzDivisorPt3}(a),
\item $H \restriction \text{R}(C) = H_3$,
\end{enumerate}
we can specify $H_3$ as the zero transformation in order to have $\text{r}(H)=\text{r}(G)$, achieving the lower bound in \eqref{th3pt3}, and it is possible to vary the rank of $H$ by increasing the rank of $H_3$ up to a maximum of $\text{r}(H_3)=\text{r}(C)=m - \text{r}\left (\begin{bmatrix} G & F \end{bmatrix}\right )$ (by the rank-nullity theorem). To summarize
\begin{eqnarray}
\nonumber \text{r}(H) &=& \text{r}\left (H \restriction \text{R}\left (\begin{bmatrix} G & F \end{bmatrix}\right )\right ) + \text{r}(H_3)\\
\nonumber &\leq&\text{r}(G)+\text{r}(C)\\
\nonumber &=& m-\left (\text{r}\left (\begin{bmatrix} G & F \end{bmatrix}\right )-\text{r}(G)\right ),\end{eqnarray} proving there is a class of square-zero right quotients of $G$ and $F$ of each rank within the bounds as specified in \eqref{th3pt3}.

Now to verify the universality of the upper bound in this second case, we have 
\begin{eqnarray}
\nonumber \text{r}(H) &\leq& \text{r}\left (H \restriction \text{R}\left (\begin{bmatrix} G & F \end{bmatrix}\right )\right ) + \text{r}(H \restriction \text{R}(C)) \\
\nonumber &\leq& \text{r}\left (H \restriction \text{R}\left (\begin{bmatrix} G & F \end{bmatrix}\right )\right ) + \text{r}(C) \\
\nonumber &=&\text{r}(G)+m-\text{r}\left (\begin{bmatrix} G & F \end{bmatrix}\right )\\
\nonumber &=&m-\left (\text{r}\left (\begin{bmatrix} G & F \end{bmatrix}\right )-\text{r}(G)\right ),\end{eqnarray} 
completing the proof. \hfill $\square$
\end{enumerate}
\subsubsection{Analogous result for left division}
As before in Proposition \ref{prop:matrixDivision}, an analogous result holds for left division. 
\begin{quote}
\begin{cor}[Part 1. Corollary 1 \cite{botha2}] Let $G,F \in M_{m \times n}(\mathscr{F})$, where $\mathscr{F}$ is a field. The following statements are equivalent:
\begin{enumerate}[(a)]
\item There exists a matrix $H \in M_{n}(\mathscr{F})$ such that $\begin{bmatrix} G \\ \mathbf{0}_{m \times n} \end{bmatrix}=\begin{bmatrix} F \\ G \end{bmatrix}H$.
\item $\text{R} \left (\begin{bmatrix} G \\ \mathbf{0}_{m \times n} \end{bmatrix} \right ) \subseteq \text{R}\left (\begin{bmatrix} F \\ G \end{bmatrix} \right ).$
\item $\text{N} \left (\begin{bmatrix} F \\ G \end{bmatrix}^T \right ) \subseteq \text{N} \left (\begin{bmatrix} G \\ \mathbf{0}_{m \times n} \end{bmatrix}^T \right ).$
\item $\text{r} \left (\begin{bmatrix} F \\ G \end{bmatrix} \right )=\text{r}\left ( \begin{bmatrix} G  & F \\ \mathbf{0}_{m \times n}& G  \end{bmatrix} \right ).$
\end{enumerate}
\label{cor:bothaSzDivisorPt1}
\end{cor}
\end{quote} 

\begin{quote}
\begin{cor}[Part 2. Corollary 1 \cite{botha2}] Let $G,F \in M_{m \times n}(\mathscr{F})$, where $\mathscr{F}$ is a field. The following statements are equivalent:
\begin{enumerate}[(a)]
\item There exists a square-zero matrix $H \in M_{n}(\mathscr{F})$ such that $G=FH$.
\item $\begin{bmatrix} G \\ \mathbf{0}_{m \times n} \end{bmatrix}$ is divisible by $\begin{bmatrix} F \\ G \end{bmatrix}$ on the left.
\item $\text{R}(G) \subseteq \text{R}(F)$ and $\text{R}(G^T) \cap \text{R}(F^T) \subseteq \text{R}(F^T \restriction \text{N}(G^T))$.
\item $\text{r}(F)=\text{r}\left ( \begin{bmatrix} G  & F   \end{bmatrix} \right )$ and $\text{r} \left (\begin{bmatrix}G \\ F \end{bmatrix} \right ) = \text{r}(G) + \text{r} \left (\begin{bmatrix} G \\ BF\end{bmatrix} \right )$, where $B \in M_{\text{n}(G^T) \times m}(\mathscr{F})$ is such that $\text{R}(B^T)=\text{N}(G^T)$.
\end{enumerate}
\label{cor:bothaSzDivisorPt2}
\end{cor}
\end{quote} 

\begin{quote}
\begin{cor}[Part 3. Corollary 1 \cite{botha2}]  If any of the conditions in Corollary \ref{cor:bothaSzDivisorPt1} or \ref{cor:bothaSzDivisorPt2} hold, then there exists a square-zero quotient~$H$ of each rank such that
\[\text{r}(G) \leq \text{r}(H) \leq \left \{ {\renewcommand{\arraystretch}{1.2} \begin{array}{ll} \frac{n}{2} & \text{if r}\left ({\renewcommand{\arraystretch}{1.1}\begin{bmatrix} G \\ F \end{bmatrix}} \right ) - \text{r}(G) \leq \frac{n}{2} \\ 
n-\left (\text{r}\left ({\renewcommand{\arraystretch}{1.1}\begin{bmatrix} G \\ F \end{bmatrix}} \right ) - \text{r}(G) \right ) & \text{if r} \left ({\renewcommand{\arraystretch}{1.1}\begin{bmatrix} G \\ F \end{bmatrix}} \right ) - \text{r}(G) > \frac{n}{2} \end{array}} \right. .\] 
\label{cor:bothaSzDivisorPt3}
\end{cor}
\end{quote}  
{\bf Proof.} The result follows easily from the following facts: $G=FH$ if and only if $G^T=H^TF^T$, $H$ is square-zero if and only if $H^T$ is square-zero (a consequence of the fact that a matrix is similar to its transpose), rank is invariant under transposition (which is also a consequence of the fact that a matrix is similar to its transpose), and that $\text{N}(F^T) \subseteq \text{N}(G^T)$ if and only if $\text{R}(G) \subseteq \text{R}(F)$. Now applying Theorems \ref{thm:bothaSzDivisorPt1}, \ref{thm:bothaSzDivisorPt2} and \ref{thm:bothaSzDivisorPt3} directly to the matrices $G^T$, $F^T$ and $H^T$ the results above follow. \hfill $\square$
 
\subsection{Calculating a square-zero quotient \label{seq:calcsqz}}
In the same way that Proposition \ref{prop:constructQuotient} provides an explicit formulation of a quotient for matrix division in general, Theorem 2 in the reference text \cite[p.78]{botha2} provides an explicit formulation of a square-zero quotient of every rank within the bounds as specified in \eqref{th3pt3}, whenever such a quotient might exist.
\begin{quote}
\begin{thm}[Theorem 2 \cite{botha2}] Let $G,F \in M_{m \times n}(\mathscr{F})$, where $\mathscr{F}$ is a field, and suppose $G$ is right divisible by $F$ with a square-zero quotient. Let $B$ and $C$ be as in Lemma \ref{lem:bothaSzDivisorPt3}.
\begin{enumerate}[(i)]
\item Suppose $\text{r}\left (\begin{bmatrix}G & F\end{bmatrix}\right )-\text{r}(G)\leq m/2$.

There exists a square-zero right quotient $H \in M_{m \times m}(\mathscr{F})$ of each allowable rank within the bounds of \eqref{th3pt3} that can be expressed in the form
\begin{equation} H = \begin{bmatrix}\mathbf{0}_{m \times n} & G & BX & C_2S\end{bmatrix}\begin{bmatrix} G & F & C_1 & C_2 \end{bmatrix}^R, \label{th4pt1} \end{equation}
where $C = \begin{bmatrix}C_1 & C_2\end{bmatrix}$ is such that $C_1$ is of the same order as $B$, and where $X$ and $S$ are square matrices of appropriate orders and $S$ is square-zero. Conversely, any matrix of this form is a square-zero right quotient of $G$ and $F$.
\item Suppose $\text{r}\left (\begin{bmatrix}G & F\end{bmatrix}\right )-\text{r}(G)> m/2$.

There exists a square-zero right quotient $H \in M_{m \times m}(\mathscr{F})$ of each allowable rank within the bounds of \eqref{th3pt3} that can be expressed in the form
\begin{equation} H=\begin{bmatrix} \mathbf{0}_{m \times n} & G &BY \end{bmatrix}\begin{bmatrix} G & F & C\end{bmatrix}^R, \label{th4pt2} \end{equation} where $Y$ is such that $BY$ is of the same order as $C$. Conversely, any matrix of this form is a square-zero right quotient of $G$ and $F$.
\end{enumerate}
\label{thm:constructSzQuotient}
\end{thm}
\end{quote}
{\bf Proof.} (i) $\text{r}\left (\begin{bmatrix}G & F\end{bmatrix}\right )-\text{r}(G)\leq m/2$.

Suppose it is given $G$ is right divisible by $F$ with a square-zero quotient. The matrix~$C$ can be determined in the same way as the matrix $C$ in Proposition \ref{prop:constructQuotient}, that is, $C$ consists of the last $m-\text{r}\left (\begin{bmatrix}G & F \end{bmatrix}\right )$ columns of an invertible matrix $Y_1$ such that \[\begin{bmatrix}G & F \end{bmatrix} = Y_1 \begin{bmatrix} I_{\text{r}([G \; F])} & \mathbf{0}\\ \mathbf{0} & \mathbf{0} \end{bmatrix} Z,\] for some invertible matrix $Z$. Now the matrix $B$ can be determined by constructing a basis for $\text{R}(F \restriction \text{N}(G))$ and then removing all vectors which are also in the subspace spanned by the columns of $G$; the selected vectors then become the columns of $B$. Let $C_1$ consist of the first $\text{r}(B)$ columns of $C$, that is if $B$ has $k$ columns, we simply select the first $k$ columns from $C$ to be the columns of $C_1$ (since $B$ is full column rank). Now $C_2$ consists of the remainder of the columns of $C$. 

Now by Theorem \ref{thm:bothaSzDivisorPt3} we can construct a square-zero matrix $H$ of each allowable rank, and it can be explicitly shown that $H$ is of the desired form by following the same construction used in the proof of (i) in Theorem \ref{thm:bothaSzDivisorPt3}, which I present below. 
\begin{enumerate}
\item By the proof of Lemma \ref{lem:bothaSzDivisorPt3} we must have that $H \restriction \text{R}\left (\begin{bmatrix} G & F \end{bmatrix}\right )$ is uniquely determined by \[H \restriction \text{R}\left (\begin{bmatrix} G & F \end{bmatrix}\right )(Gv+Fw)=Gw.\] Let $H \restriction \text{R}\left (\begin{bmatrix} G & F \end{bmatrix}\right ) = H_3$. Then by the given condition we must have \[H_3\begin{bmatrix} G & F \end{bmatrix} = \begin{bmatrix} \mathbf{0}_{m \times n} & G \end{bmatrix}.\]
\item Let $\mathcal{V}_1$ be the subspace spanned by the columns of $C_1$. Let $H_1 = H \restriction \mathcal{V}_1$ be defined in such a way that $H_1C_1 = BX$, where $X$ is an arbitrary linear transformation of appropriate order. A natural choice for $X$ is a matrix of the form Dg$[I_{k},\mathbf{0}]$ for a chosen integer $0 \leq k \leq \text{r}(B)$. Since $B$ and $C_1$ have the same order and are both full column rank, and $\text{r}(BX) \leq \text{r}(B)$ we can always find a transformation $H_1$.
\item Let $\mathcal{V}_2$ be the subspace spanned by the columns of $C_2$. Now define $H_2 = H \restriction \mathcal{V}_2$ in such a way that $H_2C_2 = C_2S$ where $S$ is an arbitrary square-zero transformation of appropriate order. The matrix $S$ can be constructed following the construction of $H_2$ as outlined in Theorem \ref{thm:bothaSzDivisorPt3}: first, $S$ is a square matrix of order $\text{r}(C)-\text{r}(B)$. Now choose a nonnegative integer $c \leq \lfloor \frac{\text{r}(C)-\text{r}(B)}{2} \rfloor$, and then let 
\begin{equation} \text{ent}_{ij}(S) =  \left \{ {\renewcommand{\arraystretch}{1.2} \begin{array}{ll} 1 & \text{if } j \leq c \text{ and } i=c+j \\
0 & \text{if } j > c \text{ or } i \neq c + j  
\end{array}} \right..\label{th3pt3h2}\end{equation}
Since $\text{R}(C_2S) \subseteq \text{R}(C_2)$ we can always find a linear transformation $H_2$.
\end{enumerate}
Now $H_3^2\begin{bmatrix} G & F\end{bmatrix}=\mathbf{0}$ so that $H \restriction \begin{bmatrix} G & F\end{bmatrix}$ is square-zero, $H_1$ maps vectors from $\text{R}(C_1)$ to $\text{R}(B)$ which is in the null space of $H \restriction \begin{bmatrix} G & F\end{bmatrix}$ as proved in Lemma \ref{lem:bothaSzDivisorPt3}, and $H_2^2C_2 = H_2C_2S = C_2S^2 = \mathbf{0}$, so that $H_2$ is square-zero. It follows that $H$ is square-zero and we also have
\begin{eqnarray} 
\nonumber H\begin{bmatrix} G & F & C_1 & C_2 \end{bmatrix} &=& \begin{bmatrix} H_3\begin{bmatrix} G & F \end{bmatrix} & H_1C_1 & H_2C_2 \end{bmatrix}\\
\nonumber &=& \begin{bmatrix} \mathbf{0}_{m \times n} & G & BX & C_2S\end{bmatrix}, \end{eqnarray}
and since $\begin{bmatrix} G & F & C_1 & C_2 \end{bmatrix}$ is full row rank it has a right inverse, and the desired form of $H$ is achieved. To show that the rank of $H$ can vary between the limits as specified in \eqref{th3pt3}, let $\alpha_1, \alpha_2, \alpha_3$ be bases for $\text{R}\left (\begin{bmatrix} G & F \end{bmatrix}\right )$, $\text{R}(C_1)$ and $\text{R}(C_2)$ respectively, and let $\beta_1, \beta_2$ be bases for $\mathcal{V}_3$ and $\text{R}(B)$ where $\mathcal{V}_3$ is a subspace such that $\text{R}(G) \subseteq \mathcal{V}_3$ and \[\text{R}\left (\begin{bmatrix}G & F \end{bmatrix}\right ) \oplus \text{R}(C_1) = \mathcal{V}_3 \oplus \text{R}(B).\] Now the matrix representation of $H$ with respect to the bases $\alpha = \alpha_1 \cup \alpha_2 \cup \alpha_3$ and $\beta = \beta_1 \cup \beta_2 \cup \alpha_3$ is
\[\Phi_{\alpha \beta}(H) = \text{Dg} [\Phi_{\alpha_1 \beta_1}(H_3),\Phi_{\alpha_2 \beta_2}(H_1), \Phi_{\alpha_3 \alpha_3}(H_2)], \] and therefore $\text{r}(H) = \text{r}(H_3) + \text{r}(H_1) + \text{r}(H_2)$. By Lemma \ref{lem:bothaSzDivisorPt3} the rank of $H_3$ is fixed to $\text{r}(G)$. The rank of $H_1$ can freely vary from zero to $\text{r}(B)$, by modifying the matrix $X$. The rank of $H_2$ can freely vary from zero to $\lfloor \text{r}(C_2)/2 \rfloor$ by modifying the square-zero matrix $S$. It was previously shown \eqref{eq:dimWsqzQ} that the dimension of $\text{R}\left (\begin{bmatrix}G & F \end{bmatrix} \right ) \oplus \text{R}(C_1)$ is $2(\text{r}(G)+\text{r}(B))$, of which $\text{r}(H_3)+\text{r}(H_1)$ is at most half, which completes the necessary part of the proof. 

For the converse, we need to prove that given
\[H = \begin{bmatrix} \mathbf{0}_{m \times n} & G & BX & C_2S\end{bmatrix}\begin{bmatrix} G & F & C_1 & C_2 \end{bmatrix}^R,\] where $X$ is an arbitrary matrix of appropriate order, and $S$ is an arbitrary square-zero matrix, we must have that $H$ is a square-zero right quotient of $G$ and $F$. Note that in this case we assume $X$ and $S$ are fixed. It is given that $G=KF$ for some matrix $K$ which is square-zero. Now if we define $H'$ as $H' \restriction \text{R}\left (\begin{bmatrix} G & F \end{bmatrix}\right ) = K$, $H' \restriction \text{R}(C_1) = K_1$ where $K_1C_1 = BX$ and $H' \restriction \text{R}(C_2) = K_2$ where $K_2C_2 = C_2S$ following the same construction as outlined in the first part of the proof we have $H'$ is square-zero and 
\[H'\begin{bmatrix} G & F & C_1 & C_2 \end{bmatrix} = \begin{bmatrix} \mathbf{0}_{m \times n} & G & BX & C_2S\end{bmatrix}.\] Now since $\begin{bmatrix} G & F & C_1 & C_2 \end{bmatrix}$ is full row rank, it follows that
\[H' = \begin{bmatrix} \mathbf{0}_{m \times n} & G & BX & C_2S\end{bmatrix}\begin{bmatrix} G & F & C_1 & C_2 \end{bmatrix}^R = H,\] and therefore $H$ is a square-zero right quotient of $G$ and $F$ as desired.

(ii) $\text{r}\left (\begin{bmatrix}G & F\end{bmatrix}\right )-\text{r}(G) > m/2$.

Since the full column rank matrix $C$ is not partitioned into submatrices, this case is simpler. Now following the proof of Theorem \ref{thm:bothaSzDivisorPt3} we can define $H \restriction \text{R}(C)$ to map all vectors in the column space of $C$ into \[\text{R}(B) \subseteq \text{N}\left (H \restriction \left (\begin{bmatrix} G & F \end{bmatrix}\right )\right ).\] Since $C$ is not partitioned there is no need to define the matrix $S$. The matrix $Y$ in this case performs a similar function to $X$ in (i), and is the sole means for adjusting the rank of $H$. The rest of the proof for this case is similar to the proof of (i). \hfill $\square$

I now proceed with an example demonstrating the construction of square-zero matrices of every rank in the case of (i) above.
\subsubsection{Example: construction of a square-zero right quotient. \label{example2}}
Let \[ G= \begin{bmatrix} 120 & 240 & 120 \\
 0 & 0 & 0 \\
 80 & 160 & 80 \\
 0 & 0 & 0 \\
 -160 & -320 & -160 \\
 0 & 0 & 0 \\
 -440 & -880 & -440 \end{bmatrix}, \quad
F= \begin{bmatrix}
11 & 40 & 29 \\
 14 & 40 & 26 \\
 17 & 40 & 23 \\
 20 & 40 & 20 \\
 23 & 40 & 17 \\
 26 & 40 & 14 \\
 29 & 40 & 11 \end{bmatrix} \in M_{7 \times 3}(\mathbb{Q}).\]
Let us first assess whether $G$ is right divisible by $F$ with a square-zero quotient $H$. By row reduction on $\begin{bmatrix} F & G \end{bmatrix}$ we obtain
\[ \begin{bmatrix} F & G \end{bmatrix} \stackrel{R}{\sim} \begin{bmatrix} 1 & 0 & -1 & 0 & 0 & 0 \\
 0 & 1 & 1 & 0 & 0 & 0 \\
 0 & 0 & 0 & 1 & 2 & 1 \\
 0 & 0 & 0 & 0 & 0 & 0 \\
 0 & 0 & 0 & 0 & 0 & 0 \\
 0 & 0 & 0 & 0 & 0 & 0 \\
 0 & 0 & 0 & 0 & 0 & 0 \end{bmatrix}, \] and row reduction on $\begin{bmatrix}G&\mathbf{0} \\ F & G \end{bmatrix}$ yields
 \[ \begin{bmatrix}G&\mathbf{0} \\ F & G \end{bmatrix} \stackrel{R}{\sim} \begin{bmatrix} 1 & 0 & -1 & 0 & 0 & 0 \\
 0 & 1 & 1 & 0 & 0 & 0 \\
 0 & 0 & 0 & 1 & 2 & 1 \\
 0 & 0 & 0 & 0 & 0 & 0 \\
 0 & 0 & 0 & 0 & 0 & 0 \\
 0 & 0 & 0 & 0 & 0 & 0 \\
 0 & 0 & 0 & 0 & 0 & 0 \\
 0 & 0 & 0 & 0 & 0 & 0 \\
 0 & 0 & 0 & 0 & 0 & 0 \\
 0 & 0 & 0 & 0 & 0 & 0 \\
 0 & 0 & 0 & 0 & 0 & 0 \\
 0 & 0 & 0 & 0 & 0 & 0 \\
 0 & 0 & 0 & 0 & 0 & 0 \\
 0 & 0 & 0 & 0 & 0 & 0 \end{bmatrix}, \] so that we have 
 \[\text{r}\left (\begin{bmatrix}G&\mathbf{0} \\ F & G \end{bmatrix}\right ) = 3 = \text{r}\left (\begin{bmatrix} F & G \end{bmatrix}\right ),\] and therefore by Theorem \ref{thm:bothaSzDivisorPt1}, there exists a square-zero matrix $H$ so that $G=HF$. Furthermore, since $\text{r}\left (\begin{bmatrix} F & G \end{bmatrix}\right )-\text{r}(G)=3-1=2 < 7/2 = m/2$ the formula in (i) is applicable.
 
I now proceed with construction of square-zero quotients of each allowable rank as per \eqref{th3pt3}. Now as before in the example given in \ref{example1} we can do row reduction on the augmented matrix $\begin{bmatrix} G & F & : & I_{7 \times 7}\end{bmatrix}$ to determine the full column rank matrix $C$ that is required for the construction. We have
 
$ \begin{bmatrix} G & F & : & I_{7 \times 7}\end{bmatrix} \stackrel{R}{\sim} $
\[ \qquad \left [ {\renewcommand{\arraystretch}{1.2}\begin{array}{rrrrrrcrrrrrrr}
 600 & 1200 & 600 & 0 & 0 & 0 &:& 0 & 0 & 0 & 0 & -1 & 2 & -1 \\
 0 & 0 & 0 & 45 & 0 & -45 &:& 0 & 0 & 0 & 0 & -11 & 7 & 4 \\
 0 & 0 & 0 & 0 & 1800 & 1800 &:& 0 & 0 & 0 & 0 & 286 & -137 & -104 \\
 0 & 0 & 0 & 0 & 0 & 0 &:& 15 & 0 & 0 & 0 & -52 & 14 & 23 \\
 0 & 0 & 0 & 0 & 0 & 0 &:& 0 & 15 & 0 & 0 & -44 & 13 & 16 \\
 0 & 0 & 0 & 0 & 0 & 0 &:& 0 & 0 & 15 & 0 & -31 & 2 & 14 \\
 0 & 0 & 0 & 0 & 0 & 0 &:& 0 & 0 & 0 & 15 & -22 & -1 & 8
\end{array}} \right ],
\]
so that \[Y_1 = \begin{bmatrix} 120 & 11 & 40 & 1 & 0 & 0 & 0 \\
0 & 14 & 40 & 0 & 1 & 0 & 0 \\
80 & 17 & 40 & 0 & 0 & 1 & 0\\
0 & 20 & 40 & 0 & 0 & 0 & 1 \\
-160 & 23 & 40 & 0 & 0 & 0 & 0 \\
0 & 26 & 40 & 0 & 0 & 0 & 0 \\
-440 & 29 & 40 & 0 & 0 & 0 & 0 \end{bmatrix} \] is such that \[\begin{bmatrix}G & F \end{bmatrix} = Y_1 \begin{bmatrix} I_{\text{r}([G \; F])} & \mathbf{0}\\ \mathbf{0} & \mathbf{0} \end{bmatrix} Z,\] for some invertible matrix $Z$. Now the matrix $C$ consists of the last \[m-\text{r}\begin{bmatrix} G & F \end{bmatrix} = 7 - 3=4\] columns of $Y_1$, that is
\[C = \begin{bmatrix} 1 & 0 & 0 & 0 \\
0 & 1 & 0 & 0 \\
0& 0 & 1 & 0 \\
0 & 0 & 0 & 1 \\
0 & 0 & 0 & 0 \\
0 & 0 & 0 & 0 \\
0 & 0 & 0 & 0 \end{bmatrix}. \]
Next, we need to determine $C_1$ and $C_2$. It is easy to see that $\text{r}(G) = 1$, and therefore $\text{r}(B)=\text{r}\begin{bmatrix} G & F \end{bmatrix}-2\text{r}(G)=3-2=1$. So $C_1$ consists of the first column of $C$ and $C_2$ the remaining three columns of $C$.

At this point the matrix $B$ should be determined. Now \[\text{R}(G)=\text{span}\{(3,0,2,0,-4,0,-11)^T\}\] and $\text{N}(G)=\text{span}\{(-1,1,-1)^T,(0,1,-2)^T\}$, so that \[\text{R}(F \restriction \text{N}(G)) = \text{span}\{(-3,-2,-1,0,1,2,3)^T\}.\] It follows that $\text{R}(G) \cap \text{R}(F \restriction \text{N}(G)) = \{\mathbf{0}\}$, and since by Lemma \ref{lem:bothaSzDivisorPt3}, \[\text{N}\left (H \restriction \text{R}\left (\begin{bmatrix}G & F\end{bmatrix}\right )\right )=\text{R}(G) + \text{R}(F \restriction \text{N}(G)),\] it follows that we may choose $\text{R}(B) = \text{R}(F \restriction \text{N}(G))$. Explicitly, 
\[ B= \begin{bmatrix} -3&-2&-1&0&1&2&3 \end{bmatrix}^T. \] 

Now we can proceed to construct square-zero quotients of each allowable rank. First, I will determine $\begin{bmatrix} G & F & C_1 & C_2 \end{bmatrix}^R$, which will be used repeatedly in constructing each quotient. As before in \ref{example1} we have \[\begin{bmatrix} G & F & C_1 & C_2 \end{bmatrix}^R = \begin{bmatrix} G & F & C_1 & C_2 \end{bmatrix}^T\left (\begin{bmatrix} G & F & C_1 & C_2 \end{bmatrix}\begin{bmatrix} G & F & C_1 & C_2 \end{bmatrix}^T\right )^{-1},\] so that  
\[\begin{bmatrix} G & F & C_1 & C_2 \end{bmatrix}^R = 
\frac{1}{3600} \left [ {\renewcommand{\arraystretch}{1.2}\begin{array}{rrrrrrr}
 0 & 0 & 0 & 0 & -1 & 2 & -1 \\
 0 & 0 & 0 & 0 & -2 & 4 & -2 \\
 0 & 0 & 0 & 0 & -1 & 2 & -1 \\
 0 & 0 & 0 & 0 & -396 & 282 & 144 \\
 0 & 0 & 0 & 0 & 88 & 4 & -32 \\
 0 & 0 & 0 & 0 & 484 & -278 & -176 \\
 3600 & 0 & 0 & 0 & -12480 & 3360 & 5520 \\
 0 & 3600 & 0 & 0 & -10560 & 3120 & 3840 \\
 0 & 0 & 3600 & 0 & -7440 & 480 & 3360 \\
 0 & 0 & 0 & 3600 & -5280 & -240 & 1920
\end{array}} \right ].
\]
The matrices $X$ and $S$ can now be determined based on the desired rank of the quotient. Since $B$ has a single column, $X$ is a scalar, and since $C_2$ has three columns $S$ is $3 \times 3$. 
\begin{enumerate}[(i)] 
\item Now to determine a square-zero quotient of minimum rank, let $X=0$ and $S=\mathbf{0}_{3 \times 3}$. Then 
\begin{eqnarray} 
\nonumber H_1&=&\begin{bmatrix} \mathbf{0}_{7 \times 3} & G & \mathbf{0} & \mathbf{0}_{3 \times 3}\end{bmatrix}\begin{bmatrix} G & F & C_1 & C_2 \end{bmatrix}^R\\
\nonumber &=& \frac{1}{15} \left [ {\renewcommand{\arraystretch}{1.2}\begin{array}{rrrrrrr}
0 & 0 & 0 & 0 & 132 & 6 & -48 \\
 0 & 0 & 0 & 0 & 0 & 0 & 0 \\
 0 & 0 & 0 & 0 & 88 & 4 & -32 \\
 0 & 0 & 0 & 0 & 0 & 0 & 0 \\
 0 & 0 & 0 & 0 & -176 & -8 & 64 \\
 0 & 0 & 0 & 0 & 0 & 0 & 0 \\
 0 & 0 & 0 & 0 & -484 & -22 & 176 \end{array}} \right ],
\end{eqnarray} is a square-zero right quotient of $G$ and $F$ with $\text{r}(H_1)=1=\text{r}(G)$.
\item To construct a square-zero quotient of rank two, let $X=1$ and $S= \mathbf{0}_{3 \times 3}$. Then 
\begin{eqnarray} 
\nonumber H_2&=&\begin{bmatrix} \mathbf{0}_{7 \times 3} & G & B & \mathbf{0}_{3 \times 3}\end{bmatrix}\begin{bmatrix} G & F & C_1 & C_2 \end{bmatrix}^R \\
\nonumber &=&\frac{1}{15} \left [ {\renewcommand{\arraystretch}{1.2}\begin{array}{rrrrrrr}-45 & 0 & 0 & 0 & 288 & -36 & -117 \\
 -30 & 0 & 0 & 0 & 104 & -28 & -46 \\
 -15 & 0 & 0 & 0 & 140 & -10 & -55 \\
 0 & 0 & 0 & 0 & 0 & 0 & 0 \\
 15 & 0 & 0 & 0 & -228 & 6 & 87 \\
 30 & 0 & 0 & 0 & -104 & 28 & 46 \\
 45 & 0 & 0 & 0 & -640 & 20 & 245 \end{array}} \right ], \end{eqnarray}
is a square-zero right quotient of $G$ and $F$ with $\text{r}(H_2)=2$.
\item To construct a square-zero quotient of rank three, let $X=1$ and let $c=1$ ($c$ as defined in \eqref{th3pt3h2}). Then \[C_2S=\begin{bmatrix}  0 & 0 & 0 \\
 1 & 0 & 0 \\
 0 & 1 & 0 \\
 0 & 0 & 1 \\
 0 & 0 & 0 \\
 0 & 0 & 0 \\
 0 & 0 & 0 \end{bmatrix}\begin{bmatrix} 0 & 0 & 0 \\
1 & 0 & 0 \\
0 & 0 & 0 \end{bmatrix}=\begin{bmatrix}  0 & 0 & 0 \\
 0 & 0 & 0 \\
 1 & 0 & 0 \\
 0 & 0 & 0 \\
 0 & 0 & 0 \\
 0 & 0 & 0 \\
 0 & 0 & 0 \end{bmatrix}, \] and
\begin{eqnarray} 
\nonumber H_3 &=&\begin{bmatrix} \mathbf{0}_{7 \times 3} & G & B & C_2S \end{bmatrix}\begin{bmatrix} G & F & C_1 & C_2 \end{bmatrix}^R \\
\nonumber &=&\frac{1}{15} \left [ {\renewcommand{\arraystretch}{1.2}\begin{array}{rrrrrrr}-45 & 0 & 0 & 0 & 288 & -36 & -117 \\
 -30 & 0 & 0 & 0 & 104 & -28 & -46 \\
 -15 & 15 & 0 & 0 & 96 & 3 & -39 \\
 0 & 0 & 0 & 0 & 0 & 0 & 0 \\
 15 & 0 & 0 & 0 & -228 & 6 & 87 \\
 30 & 0 & 0 & 0 & -104 & 28 & 46 \\
 45 & 0 & 0 & 0 & -640 & 20 & 245 \end{array}} \right ], \end{eqnarray}
is a square-zero right quotient of $G$ and $F$ with $\text{r}(H_3)=3$. Now notice that by \eqref{th3pt3} the maximum rank of a square-zero right quotient of $G$ and $F$ must be $\lfloor m/2 \rfloor = \lfloor 7/2 \rfloor = 3$, and therefore $H_3$ is a square-zero quotient of maximum rank. 
\end{enumerate}
\subsubsection{Construction of a square-zero left quotient, \\and the special case $\mathscr{F}=\mathbb{R}$, $\mathscr{F}=\mathbb{C}$}
As before, a similar result to Theorem \ref{thm:constructSzQuotient} can be formulated for left division by a square-zero quotient, by replacing each matrix in the statement of the theorem with its transpose, and noting that for any matrix $A$ we have $(A^R)^T = (A^T)^L$. 

A special result can also be formulated for the case where $\mathscr{F}=\mathbb{R}$ or $\mathscr{F}=\mathbb{C}$.
\begin{quote}
\begin{cor}[Corollary 2 \cite{botha2}] If $\mathscr{F}=\mathbb{C}$ or $\mathscr{F}=\mathbb{R}$, then $H$ in \eqref{th4pt1} can be expressed as
\begin{equation} H = (GF^*+BXC_1^*+C_2SC_2^*)(GG^*+FF^* + C_1C_1^*+C_2C_2^*)^{-1}, \label{th4pt1_RC} \end{equation} 
and $H$ in \eqref{th4pt2} can be expressed as
\begin{equation} H = (GF^*+BYC^*)(GG^* + FF^* + CC^*)^{-1}. \label{th4pt2_RC} \end{equation}
\label{cor:SzQuotientRC}
\end{cor}
\end{quote}
{\bf Proof.}  In this case we can specify the right inverse explicitly (which is possible due to the fact that vector spaces over these fields are inner product spaces \cite[Proposition 19.2, p.442]{golan}):
\begin{eqnarray}
\nonumber \begin{bmatrix} G & F & C_1 & C_2 \end{bmatrix}^R &=& \begin{bmatrix} G & F & C_1 & C_2 \end{bmatrix}^* \left (\begin{bmatrix} G & F & C_1 & C_2 \end{bmatrix}\begin{bmatrix} G & F & C_1 & C_2 \end{bmatrix}^*\right )^{-1} \\
\nonumber &=& \begin{bmatrix} G^* \\ F^* \\ C_1^* \\ C_2^* \end{bmatrix}\begin{bmatrix} GG^* & FF^* & C_1C_1^* & C_2C_2^* \end{bmatrix}^{-1},
\end{eqnarray}
and then simplifying \eqref{th4pt1} will produce \eqref{th4pt1_RC}.
In a similar way the right inverse in \eqref{th4pt2} can be specified explicitly when $\mathscr{F}=\mathbb{R}$ or $\mathscr{F}=\mathbb{C}$, and if we then simplify the result is \eqref{th4pt2_RC}. \hfill $\square$

\subsection{Extending the results on matrix factorization}
The results of section \ref{seq:sqzerodiv} and \ref{seq:calcsqz} provide a suitable general theoretical framework for further research into products of square-zero matrices, with the potential to unify previously isolated results. Botha demonstrated this by revisiting the results first proved by Novak \cite{botha2}, \cite{novak}. Following Botha I now revisit the results of section \ref{sec:sqzeroprev}, with the aid of Theorems \ref{thm:bothaSzDivisorPt1}, \ref{thm:bothaSzDivisorPt2} and \ref{thm:bothaSzDivisorPt3}. I will expand upon the arguments used in the proofs, in an attempt to exhaustively unpack the theory underpinning the results.  
\subsubsection{Products of two square-zero matrices}
Let us first consider Theorem \ref{thm:botha2FactorSZ1}: necessary and sufficient conditions for factorization into a product of two square-zero matrices.

\begin{quote}
\begin{thm}[Theorem 3 \cite{botha2}] Let $G \in M_{m}(\mathscr{F})$, where $\mathscr{F}$ denotes an arbitrary field. Then $G=HF$, where $H,F \in M_{m}(\mathscr{F})$ are both square-zero if and only if 
\begin{equation} \text{r}(G) \leq \text{n}(G) - \dim(\text{R}(G) \cap \text{N}(G)). \label{ch_sz_fact_eq1th5} \end{equation}
If this condition holds, then the ranks of $H$ and $F$ can be chosen arbitrarily provided 
\[\text{r}(G) \leq \text{r}(H), \text{r}(F) \leq \frac{m}{2}.\] 
\label{thm:botha2FactorSZ2}
\end{thm}
\end{quote}
{\bf Proof.} The proof (following \cite{botha2} but expanding and elucidating condensed arguments) appeals to (c) of Theorem \ref{thm:bothaSzDivisorPt2} in both the necessary and sufficient parts. First, let $\mathscr{F}^m = \mathcal{V} \oplus \text{N}(G)$, where $\dim(\mathcal{V}) = \text{r}(G)$ (which is possible as a consequence of the rank-nullity theorem  \cite[Theorem 2.11]{cullen}).

Suppose $G=HF$, where $H,F \in M_{m \times m}(\mathscr{F})$ are both square-zero. The proof-strategy is now to show that $\text{R}(F \restriction \mathcal{V}) + (\text{R}(G) \cap \text{N}(G))$ is a direct sum, and that \[\dim(\text{R}(F \restriction \mathcal{V}) \oplus (\text{R}(G) \cap \text{N}(G))) = \text{r}(G) + \dim(\text{R}(G) \cap \text{N}(G)) \leq \text{n}(G),\] which is equivalent to \eqref{ch_sz_fact_eq1th5}. 

Since $F$ is a right divisor of $G$ with a square-zero quotient, by Theorem \ref{thm:bothaSzDivisorPt2}(c) it follows that 
\begin{equation} \text{N}(F) \subseteq \text{N}(G) \text{ and } \text{R}(G) \cap \text{R}(F) \subseteq \text{R}(F \restriction \text{N}(G)), \label{eq4th5} \end{equation} and furthermore from Proposition \ref{prop:matrixDivision}(b) $\text{row}(G) \subseteq \text{row}(F)$ so that we must have $\text{r}(G) \leq \text{r}(F)$. Now we also have that $F$ is square-zero, implying that $\text{R}(F) \subseteq \text{N}(F)$. Combining these results, we must have that 
\begin{equation} \text{R}(F \restriction \mathcal{V}) \subseteq \text{N}(G) \label{eq:FpylVinNG} \end{equation} and $F \restriction \mathcal{V}$ is injective, since every nonzero vector in $\mathcal{V}$ is not in $N(G)$ and hence not in $\text{N}(F)$ (\cite[Theorem 4.1, p.127]{cullen}). So we have
\begin{equation} \text{r}(F \restriction \mathcal{V}) = \text{r}(G). \label{eq2th5} \end{equation}

Furthermore $\text{R}(F \restriction \mathcal{V}) \cap \text{R}(G) = \{0\}$ since if $F(v) \in \text{R}(F \restriction \mathcal{V}) \cap \text{R}(G)$ then $F(v) \in \text{R}(F \restriction \text{N}(G))$ by \eqref{eq4th5}. So there is some $u \in \text{N}(G)$ so that $F(v) = F(u)$. But then $v-u \in \text{N}(G)$, and since $N(G)$ is a subspace and hence closed under vector addition we must then have $v \in \text{N}(G)$. But now, since $\mathcal{V} \cap \text{N}(G) = \{ \mathbf{0} \}$, it follows that $v=\mathbf{0}$, and therefore $F(v)=\mathbf{0}$, proving that $\text{R}(F \restriction \mathcal{V}) \cap \text{R}(G) = \{\mathbf{0}\}$.

Now if $\text{R}(F \restriction \mathcal{V}) \cap \text{R}(G) = \{ \mathbf{0} \}$ then $\text{R}(F \restriction \mathcal{V}) \cap (\text{N}(G) \cap \text{R}(G)) = \{ \mathbf{0} \}$, and therefore the sum $\text{R}(F \restriction \mathcal{V}) + (\text{N}(G) \cap \text{R}(G))$ is a direct sum. Now since this subspace is a direct sum each vector $w \in \text{R}(F \restriction \mathcal{V}) \oplus (\text{N}(G) \cap \text{R}(G))$ can be uniquely expressed as $w=u+v$ where $v \in \text{R}(F \restriction \mathcal{V}) \subseteq \text{N}(G)$ (by \eqref{eq:FpylVinNG}) and $u \in (\text{N}(G) \cap \text{R}(G)) \subseteq \text{N}(G)$, and since $\text{N}(G)$ is a subspace $w= u+v \in \text{N}(G)$. It follows that $\text{R}(F \restriction \mathcal{V}) \oplus (\text{N}(G) \cap \text{R}(G)) \subseteq \text{N}(G)$. Combining this result with \eqref{eq2th5} yields
\begin{eqnarray} 
\nonumber \text{n}(G) &\geq& \dim(\text{R}(F \restriction \mathcal{V}) \oplus (\text{N}(G) \cap \text{R}(G)))\\
\nonumber &=& \text{r}(F \restriction \mathcal{V}) + \dim(\text{N}(G) \cap \text{R}(G)) \\
\nonumber &=& \text{r}(G) + \dim(\text{N}(G) \cap \text{R}(G)),\end{eqnarray} completing the first part of the proof.

The converse is proved by construction: given that \eqref{ch_sz_fact_eq1th5} holds, it is possible to construct a square-zero right divisor $F$, of $G$, with a square-zero right quotient $H$. 

I proceed with the construction: given that \eqref{ch_sz_fact_eq1th5} holds, we can decompose the null space of $G$: $\text{N}(G)=\mathcal{W} \oplus \mathcal{U}$ with $\dim(\mathcal{W})=\text{r}(G)$ and $(\text{N}(G) \cap \text{R}(G)) \subseteq \mathcal{U}$. Now define two linear transformations:
\begin{enumerate}[(i)]
\item Let $L_1 : \mathcal{V} \rightarrow \mathcal{W}$ be an isomorphism ($\mathcal{V}$ as defined before: $\mathscr{F}^m = \mathcal{V} \oplus \text{N}(G)$).
\item Define $L_2 : \text{N}(G) \rightarrow \mathcal{U}$ as 
\[L_2(\mathcal{W})=\{ \mathbf{0} \} \text{ and } L_2(u_i) = \left \{ {\renewcommand{\arraystretch}{1.2} \begin{array}{ll} u_{i+c} & \text{if } 1 \leq i \leq c  \\
\mathbf{0} & \text{if } c < i \leq k  \end{array}} \right.,\] where $\{u_1,u_2,\ldots,u_k\}$ is a basis for $\mathcal{U}$, and $0 \leq c \leq \lfloor k/2 \rfloor$.
\end{enumerate}  
Now define $F$ as 
\begin{enumerate}[(i)]
\item $F \restriction \mathcal{V} = L_1,$
\item $F \restriction \text{N}(G) = L_2.$
\end{enumerate}
It is easy to verify that $F^2u=L_2 \circ L_2 (u) = \mathbf{0}$ for any $u \in \text{N}(G)$ and $F^2v=L_2 \circ L_1(v)=\mathbf{0}$ for any $v \in \mathcal{V}$, and therefore $F$ is square-zero. 

Now $\text{R}(F) \subseteq \text{N}(G)$ and $(\text{N}(G) \cap \text{R}(G)) \subseteq \mathcal{U}$ so that 
\begin{eqnarray} 
\nonumber \text{R}(G) \cap \text{R}(F) &=& \text{N}(G) \cap \text{R}(G) \cap \text{R}(F)\\
\nonumber  &\subseteq& \mathcal{U} \cap \text{R}(F) \\
\nonumber &\subseteq& \text{R}(F \restriction \text{N}(G)).\end{eqnarray}
The last step above can be verified: let $u \in \mathcal{U} \cap \text{R}(F)$. Then $u \notin \mathcal{W}$ and since $u$ is also in $\text{R}(F)$, we must have $u \in F \restriction \text{N}(G)$.
  
It follows by Theorem  \ref{thm:bothaSzDivisorPt2}(c) that $F$ is a right divisor of $G$ with a square-zero quotient $H$, proving the converse.

It remains to verify the ranks of $H$ and $F$, given that the preceding conditions hold. We have $R(F \restriction \mathcal{V}) = \mathcal{W}$ with $\dim(\mathcal{W})=\text{r}(G)$ and $R(F \restriction \text{N}(G) = \{u_{1+c}, \ldots, u_{c+c}\} \subseteq \mathcal{U}$ with $\dim(\text{span} \{u_{1+c}, \ldots, u_{c+c}\})=c$ and $\mathcal{W} \cap \mathcal{U} = \{ \mathbf{0} \}$, so that $\text{r}(F) = \text{r}(G)+c$. By varying $c$, the rank of $F$ can assume any value between $\text{r}(G)$ and $m/2$. 

Also since 
\[\text{r} \left ( \begin{bmatrix} G & F \end{bmatrix} \right ) - \text{r}(G) \leq \text{r}(G) + \text{r}(F) - \text{r}(G) = \text{r}(F) \leq \frac{m}{2}\] it follows from Theorem \ref{thm:bothaSzDivisorPt3} that the rank of $H$ can assume any value between $\text{r}(G)$ and $m/2$.

Combining the results above it follows that the ranks of $F$ and $H$ can be chosen arbitrarily provided
\[\text{r}(G) \leq \text{r}(H), \text{r}(F) \leq \frac{m}{2}.\] \hfill $\square$

As a consequence of Theorem \ref{thm:botha2FactorSZ2}, observe that not all square-zero matrices can be written as a product of two square-zero matrices. For example, let $G \in M_{m \times m}(\mathscr{F})$ be such that $\text{R}(G)=\text{N}(G)$. Then $G$ is square-zero, and $\dim(\text{R}(G) \cap \text{N}(G))=\text{n}(G)$ so that \[\text{n}(G)-\dim(\text{R}(G) \cap \text{N}(G))=0<\text{r}(G),\] and therefore, by Theorem \ref{thm:botha2FactorSZ2}, $G$ is not a product of two square-zero matrices. However, it will be shown in section \ref{arbmany} that any square-zero matrix $G$ can be expressed as a product of three or more square-zero factors.
\subsubsection{Products of three square-zero matrices}
Now we can extend Theorem \ref{botha3FactorsSz} to necessary and sufficient conditions for a matrix to be factored into a product of (arbitrarily many) square-zero matrices. If a matrix can be written as a product of square-zero matrices, no more than three such factors are required.
\begin{quote}
\begin{thm}[Theorem 5 \cite{botha2}] Let $G \in M_{m \times m}(\mathscr{F})$, where $\mathscr{F}$ denotes an arbitrary field. Then $G$ is a product of square-zero matrices if and only if $\text{r}(G) \leq \frac{m}{2}$, which is equivalent to stating $\text{r}(G) \leq \text{n}(G)$. The minimum number of square-zero matrices in such a factorization is at most three, which is a sharp upper bound.

Moreover, each matrix $G$ that satisfies the above condition can be expressed as a product of three square-zero matrices with arbitrary ranks $r_i$ subject to \[\text{r}(G) \leq r_i \leq \frac{m}{2} \text { for } 1 \leq i \leq 3.\]
\label{thm:botha3FactorsSz2}
\end{thm}
\end{quote}
{\bf Proof.} For necessity, the proof employed in Theorem \ref{botha3FactorsSz} (first paragraph of proof) remains valid (for an arbitrary amount of square-zero factors Theorem 2.5 in the text by Cullen \cite[p.77]{cullen} can be applied iteratively).

For the converse, I again follow \cite{botha2} but expand upon and elucidate condensed arguments. The proof consists of constructing a right divisor $F$ with a square-zero quotient, so that $F$ itself is a product of two square-zero matrices. Theorem~\ref{thm:bothaSzDivisorPt2}(c) is used to prove that $F$ is a right divisor with a square-zero quotient, and Theorem~\ref{thm:botha2FactorSZ2} to show that $F$ is itself a product of two square-zero matrices. So $F$ needs to have the following properties:
\begin{eqnarray}
\nonumber \text{N}(F) &\subseteq& \text{N}(G) \quad \text{ and }\quad \text{R}(G) \cap \text{R}(F) \subseteq \text{R}(F \restriction \text{N}(G)), \quad \text{ and}\\
\text{r}(F) &\leq& \text{n}(F) - \dim(\text{R}(F) \cap \text{N}(F)). \label{eq3th6}
\end{eqnarray}

A suitable matrix satisfying these properties is any representation of a projection (idempotent matrix) along the null space  of $G$ onto a subspace that only has the the zero vector in common with $\text{R}(G)$, since supposing $F$ is such a matrix, we have $\text{N}(F)=\text{N}(G)$ and 
\begin{equation} \text{R}(G) \cap \text{R}(F) = \{\mathbf{0}\} = \text{R}(F \restriction \text{N}(F)) = \text{R}(F \restriction \text{N}(G)). \label{eq2th6} \end{equation}
 Furthermore, an idempotent matrix has the property\footnote{Since $F^2 = F$, if we let $\mathcal{W}=\text{R}(F)$, then $F(\mathcal{W})=\mathcal{W}$, so $\mathcal{W} \cap \text{N}(F)=\{\mathbf{0}\}$.} $\text{R}(F) \cap \text{N}(F) = \{\mathbf{0}\}$, and since $\text{N}(F) = \text{N}(G)$ and $\text{r}(G) \leq \text{n}(G)$ we must have 
 \[\text{r}(F) \leq \text{n}(F) = \text{n}(F) - \dim(\text{R}(F) \cap \text{N}(F)).\]
 
I will now discuss the construction of the idempotent matrix $F$. The construction requires that we decompose $\text{R}(G) + \text{N}(G)$ into a direct sum of subspaces $\mathcal{R} \oplus (\text{R}(G) \cap \text{N}(G)) \oplus \mathcal{N}$, where $\mathcal{R} \subseteq \text{R}(G)$ and $\mathcal{N} \subseteq \text{N}(G)$. So $\mathcal{R}$ and $\mathcal{N}$ are subspaces of $\text{R}(G)$ and $\text{N}(G)$ respectively which only has the zero vector in common with $\text{N}(G)$ and $\text{R}(G)$ respectively. 

Is it always possible to construct such a decomposition of $\text{R}(G) + \text{N}(G)$? First suppose that $\dim(\mathcal{R}) > 0$ and $\dim(\text{R}(G) \cap \text{N}(G)) >0$. Then we can extend a basis for $\text{R}(G) \cap \text{N}(G)$ to a basis for $\mathcal{R} \oplus (\text{R}(G) \cap \text{N}(G))$, and since $\text{r}(G) \leq \text{n}(G)$ we must have $\dim(\mathcal{N})>0$ and therefore we can again extend this basis to be a basis for $\mathcal{R} \oplus (\text{R}(G) \cap \text{N}(G)) \oplus \mathcal{N}$ (so that we have produced the desired decomposition).  Now by a similar argument, if $\dim(\mathcal{R})=0$, but $\dim(\text{R}(G) \cap \text{N}(G)) >0$ the direct sum simply reduces to $(\text{R}(G) \cap \text{N}(G)) \oplus \mathcal{N}$, and if $\dim(\text{R}(G) \cap \text{N}(G)) = 0$ we already have $\text{R}(G) + \text{N}(G) = \text{R}(G) \oplus \text{N}(G)=\mathcal{R} \oplus \mathcal{N}$. Finally, if $\mathcal{N}=\{\mathbf{0}\}$, since $\text{r}(G) \leq \text{n}(G)$ the direct sum reduces to the single term $\text{R}(G) \cap \text{N}(G)$.

Having proved that we can always find the desired direct sum decomposition, we can now extend the basis defining the decomposition to be a basis for $\mathscr{F}^m$ so that (in general) we have 
\begin{equation} \mathscr{F}^m = \mathcal{W} \oplus (\text{R}(G) + \text{N}(G)) = \mathcal{W} \oplus \mathcal{R} \oplus (\text{R}(G) \cap \text{N}(G)) \oplus \mathcal{N}. \label{eq1th6} \end{equation}

The matrix $F$ is now constructed to be the matrix with properties: 
\begin{itemize}
\item $\text{N}(F) = \text{N}(G)$ and 
\item $R(F) = \mathcal{W} \oplus \tilde{\mathcal{R}}$ where $\tilde{\mathcal{R}} \subseteq \mathcal{R} \oplus \mathcal{N}$ and $\tilde{\mathcal{R}} \cap \text{R}(G) = \{ \mathbf{0} \}$, $\tilde{\mathcal{R}} \cap \text{N}(G) = \{\mathbf{0}\}$.
\end{itemize} 
Explicitly, we can construct $\tilde{\mathcal{R}}$ in the following way: let $\{r_1, r_2, \ldots r_k\}$ be a basis for $\mathcal{R}$, and $\{n_1, n_2, \ldots n_k\}$ be a linearly independent subset of vectors in $\mathcal{N}$. Now let $\tilde{\mathcal{R}} = \text{span}\{r_1+n_1,r_2+n_2,\ldots,r_k+n_k\}$. Now, to prove that $\tilde{\mathcal{R}} \cap \text{R}(G) =\{\mathbf{0}\}$, suppose $\{r_{k+1}, r_{k+2}, \ldots r_{\text{r}(G)}\}$ completes the basis $\{r_1, r_2, \ldots r_k\}$ to be a basis for $\text{R}(G)$, and that 
\[ \sum_{i=1}^{\text{r}(G)} c_ir_i +\sum_{j=1}^{k} c_{\text{r}(G)+j}(r_j+n_j)) = \mathbf{0}.\]
It follows that
\[\sum_{i=1}^k ((c_i + c_{\text{r}(G)+i})r_i + c_{\text{r}(G)+i}n_i) +  \sum_{j=k+1}^{\text{r}(G)} c_{j}r_j = \mathbf{0}.\]
By the direct sum decomposition \eqref{eq1th6} we have $\{r_1,r_2,\ldots,r_{\text{r}(G)},n_1,n_2,\ldots,n_k\}$ is a linearly independent set, and therefore 
\[(c_i+c_{\text{r}(G)+i}) = c_{\text{r}(G)+i} = c_j = 0 \text{ for each } i \in \{1,2,\ldots,k\}, j \in \{k+1,k+2, \ldots,\text{r}(G)\}.\] 
But then $c_1=c_2=\cdots=c_k=0$, and therefore the spanning vectors of $\tilde{\mathcal{R}}$ and basis vectors of $\text{R}(G)$ are linearly independent, so that $\tilde{\mathcal{R}} \cap \text{R}(G) = \{ \mathbf{0} \}$. Notice that the preceding argument also shows that the set $\{r_1+n_1,r_2+n_2,\ldots,r_k+n_k\}$ is in fact a basis for $\tilde{\mathcal{R}}$ and therefore $\dim(\tilde{\mathcal{R}})=\dim(\mathcal{R})$. Now by a similar argument $\tilde{\mathcal{R}} \cap \text{N}(G) = \{ \mathbf{0}\}$. Therefore we must also have $\tilde{\mathcal{R}} \cap (\text{R}(G) \cap \text{N}(G))=\{\mathbf{0}\}$.

From the preceding discussion it follows that 
\[\text{R}(F) \cap \text{R}(G) = (\mathcal{W} \oplus \tilde{\mathcal{R}}) \cap (\mathcal{R} \oplus (\text{R}(G) \cap \text{N}(G))) = \{\mathbf{0}\}\] as required by \eqref{eq2th6}. It remains to verify that $F$ is indeed a valid idempotent operator, by proving that $\text{R}(F) \oplus \text{N}(F)= \mathscr{F}^m$. We have
\begin{eqnarray}
\nonumber \mathscr{F}^m &=& \mathcal{W} \oplus \mathcal{R} \oplus (\text{R}(G) \cap \text{N}(G)) \oplus \mathcal{N} \\
\nonumber &=& \mathcal{W} \oplus \mathcal{R} \oplus \text{N}(G) \\
\nonumber &=& \mathcal{W} \oplus \tilde{\mathcal{R}} \oplus \text{N}(G) \\
\nonumber &=& \text{R}(F) \oplus \text{N}(F),
\end{eqnarray} 
since $\mathcal{R} \oplus \text{N}(G) = \tilde{\mathcal{R}} \oplus \text{N}(G)$. Therefore $F$ is indeed a valid projection operator on $\mathscr{F}^m$ with the required properties as set out in \eqref{eq3th6}. So $F$ is a right divisor of $G$ with a square-zero quotient, and by Theorem \ref{thm:botha2FactorSZ2} $F$ itself is a product of two square-zero matrices, that is, $G$ is a product of three square-zero matrices. 

Now to prove that three factors is a sharp upper bound. Let $G$ be a matrix such that $\text{r}(G) = \text{n}(G)$ and $\dim(\text{R}(G) \cap \text{N}(G)) > 0$, then $\text{r}(G)> \text{n}(G) - \dim(\text{R}(G) \cap \text{N}(G))$. By Theorem \ref{thm:botha2FactorSZ2} the matrix $G$ is then not a product of two square-zero matrices. So $G$ must have at least three square-zero factors in any factorization consisting of exclusively square-zero matrices. 

Finally, by the preceding discussion each $G$ such that $\text{r}(G) \leq m/2$ can be expressed as a product of three square-zero matrices -- I will now verify the ranks of these matrices. For the right square-zero quotient $H_1$ of $G$ and the idempotent matrix $F$, since 
\[\text{r} \left ( \begin{bmatrix} G & F \end{bmatrix} \right ) - \text{r}(G) \leq \text{r}(G) + \text{r}(F) - \text{r}(G) = \text{r}(F) \leq \frac{m}{2},\] it follows directly from Theorem \ref{thm:bothaSzDivisorPt3} that the rank of $H_1$ can be arbitrary provided $\text{r}(G) \leq H_1 \leq m/2$. Also by Theorem \ref{thm:botha2FactorSZ2}, $F$ is a product of two square-zero matrices (let these matrices be $H_2,H_3$) with arbitrary ranks, provided $\text{r}(F) \leq H_2, H_3 \leq m/2$ and since $\text{r}(F)=\text{r}(G)$ the desired result follows. This concludes the proof of Theorem \ref{thm:botha3FactorsSz2}. \hfill $\square$

It is worth observing that in the following extreme cases the construction of $F$ is still valid. If $\dim(\mathcal{R})=0$ the matrix $G$ is already square-zero. The construction of $F$ is still valid, however, with $\text{R}(F) = \mathcal{W}$, and then the factorization results in three square-zero factors. At the other extreme we have $\text{R}(G) \cap \text{N}(G) = \{\mathbf{0}\}$, and since $\text{r}(G) \leq \text{n}(G)$ the conditions of Theorem \ref{thm:botha2FactorSZ2} are then satisfied, so that $G$ is a product of two square-zero matrices. Again, however, the construction of $F$ is valid, with $\text{R}(F)=\tilde{\mathcal{R}}$, and if we make use of $F$ as a divisor of $G$, a factorization into three square-zero matrices results. 

\subsubsection{Products of arbitrarily many square-zero matrices \label{arbmany}}
The following result extends the result obtained in Theorem \ref{thm:botha3FactorsSz2}, by showing that any matrix that can be expressed as a product of square-zero matrices, can be expressed as a product of arbitrarily many square-zero factors subject to it being at least three factors.

\begin{quote}
\begin{thm}[Theorem 7 \cite{botha2}] Let $\mathscr{F}$ be an arbitrary field and $G \in M_{m \times m}(\mathscr{F})$ be such that $\text{r}(G) \leq m/2$. For any positive integer $k \geq 3$, the matrix $G$ can be expressed as $G=S_1S_2 \ldots S_k$, where $S_i$ is a square-zero matrix with arbitrary rank, provided 
\[\text{r}(G) \leq \text{r}(S_i) \leq \frac{m}{2},\]
for each $1 \leq i \leq k$.
\label{thm:bothaArbManyFactors}
\end{thm}
\end{quote}
{\bf Proof.} Suppose $G \in M_{m \times m}(\mathscr{F})$ is such that $\text{r}(G) \leq m/2$. By Theorem \ref{thm:botha3FactorsSz2}, $G$ can be expressed as a product of three square-zero matrices with arbitrary ranks $r_i$ subject to $\text{r}(G) \leq r_i \leq m/2$ for $1 \leq i \leq 3$. 

Now suppose $k>3$, and the result holds for $k-1$ factors, that is, $G$ can be expressed as a product $S_1S_2 \ldots S_{k-1}$, where each $S_i$ is square-zero and $\text{r}(S_i)$ is arbitrary subject to $\text{r}(G) \leq \text{r}(S_i) \leq m/2$ for $1 \leq i \leq k-1$.

It follows that $\text{r}(S_{k-2}S_{k-1}) \leq m/2$ \cite[Proposition 6.11]{golan}. By Theorem \ref{thm:botha3FactorsSz2} it follows that $S_{k-2}S_{k-1}$ can be expressed as a product of three square-zero matrices $\bar{S}_{k-2}\bar{S}_{k-1}\bar{S}_{k}$ with arbitrary ranks, subject to 
\[ \text{r}(S_{k-2}S_{k-1}) \leq \text{r}(\bar{S}_i) \leq m/2 \text{ for } k-2 \leq i \leq k,\] and since, by the induction hypothesis we can choose $\text{r}(S_{k-2}S_{k-1})$ to be at the lower bound $\text{r}(G)$, it follows that 
\[ \text{r}(G) \leq \text{r}(\bar{S}_i) \leq m/2 \text{ for } k-2 \leq i \leq k.\]

By induction it follows that $G$ can be expressed as a product of $k \geq 3$ square-zero matrices with arbitrary ranks subject to the lower bound $\text{r}(G)$ and upper bound $m/2$, which completes the proof. \hfill $\square$
\begin{quote}
\begin{cor} Suppose $G \in M_{m}(\mathscr{F})$ satisfies one of the following conditions:
\begin{enumerate}
\item $G$ is square-zero,
\item $G$ can be written as a product of two square-zero matrices,
\item $G$ can be written as a product of three square-zero matrices.
\end{enumerate} Then, for any positive integer $k \geq 3$, the matrix $G$ can be expressed as $G=S_1S_2 \ldots S_k$, where $S_i$ is a square-zero matrix with arbitrary rank, provided 
\[\text{r}(G) \leq \text{r}(S_i) \leq \frac{m}{2},\]
for each $1 \leq i \leq k$.
\label{cor:krisjanArbManyFactorSz}
\end{cor}
\end{quote}
{\bf Proof.} If $G$ is square-zero then $\text{R}(G) \subseteq \text{N}(G)$, so that $\text{r}(G) \leq \text{n}(G)$. If $G$ can be written as a product of two square-zero matrices, then by Theorem \ref{thm:botha2FactorSZ2}, 
\[\text{r}(G) \leq \text{n}(G)-\text{dim}(\text{R}(G) \cap \text{N}(G)) \leq \text{n}(G).\] Finally if $G$ can be written as a product of three square-zero matrices then it follows directly from Theorem \ref{thm:botha3FactorsSz2} that $\text{r}(G) \leq \text{n}(G)$. In each of these cases then, $G$ satisfies the conditions of Theorem \ref{thm:bothaArbManyFactors}, and the desired result follows. \hfill $\square$
\section{Conclusion on Products of Square-Zero Matrices}
The central premise of this chapter has been to provide a coherent and comprehensive presentation of currently available results on square-zero factorization, with detailed discussion of the underlying structures and principles involved. To summarize, previous results (Theorem \ref{thm:novak2FactorSz}, \ref{thm:botha2FactorSZ1} and \ref{novak3FactorsSz}, \ref{botha3FactorsSz}) can be unified into a coherent theory of square-zero factorization through the general result presented in Theorems \ref{thm:bothaSzDivisorPt1}, \ref{thm:bothaSzDivisorPt2} and \ref{thm:bothaSzDivisorPt3}: necessary and sufficient conditions for a matrix $G$ to be divisible by $F$ with a square-zero quotient, with $G$ and $F$ of the same order but not necessarily square. Results regarding the rank of such a quotient were also discussed and presented. Theorems \ref{thm:botha2FactorSZ2} and \ref{thm:botha3FactorsSz2} make use of this general result to prove Theorems \ref{thm:botha2FactorSZ1} and \ref{botha3FactorsSz}: necessary and sufficient conditions for a square matrix to be factored into two, three and arbitrarily many factors, also extending previous results by adding conditions on the rank of the factors. Finally, Theorem \ref{thm:bothaArbManyFactors} shows that any matrix that can be expressed as a product of square-zero factors can be expressed as a product of arbitrarily many square-zero factors, provided such a factorization contains at least three factors.

\chapter{Sums of Square-Zero Matrices}
In this chapter extensive use is made of the rational canonical form theory. Given a matrix $A \in M_n(\mathscr{F})$ where $\mathscr{F}$ is an arbitrary field, the \emph{characteristic matrix} of $A$ is $xI_n - A$. The Smith Canonical Form $S(x)$ of such a characteristic matrix is always a full rank diagonal matrix $\text{Dg}[f_1(x),f_2(x),\ldots,f_n(x)]$, where $f_i(x)$ is a monic polynomial for each $i \in \{1,2,\ldots,n\}$ and $f_i(x)$ divides $f_{i+1}(x)$ for each $i \in \{1,2,\ldots,n-1\}$ (it might be that some of these polynomials are constant, i.e. $f_j(x)=1$ for some positive integer $j<n$). These polynomials are called the \emph{invariant factors} of $xI-A$ (i.e. of the characteristic matrix of $A$). In this text I will refer to these factors as the \emph{invariant polynomials} of the matrix $A$. 

The product of the invariant polynomials is the characteristic polynomial of $A$, and the invariant polynomial $f_n(x)$ is the minimum polynomial of $A$ (monic polynomial of least degree whose value at $A$ is the zero matrix).\footnote{In some texts the minimum polynomial is referred to as the minimal polynomial.} The invariant polynomials contain the elementary divisors of $A$, which in turn determine the rational canonical form of $A$. There is also a useful canonical form directly derived from the invariant polynomials: if the nonconstant invariant polynomials of $A$ are $f_1(x),f_2(x),\ldots,f_t(x)$ (where $t \leq n$) then $A \approx \text{Dg}[C(f_1(x)),C(f_2(x)),\ldots,C(f_t(x))]$. The interested reader may consult a source such as the text by Cullen \cite{cullen}, chapters 6 and 7, for a comprehensive treatment of the subject.

In this text an \emph{even (-power)} polynomial refers to a polynomial of the form $p(x)=p'(x^2)$, that is, every coefficient of an odd-power of $x$ in the expanded form of $p(x)$ is zero. Similarly \emph{odd (-power)} polynomials are defined as polynomials in which every coefficient of an even-power of $x$ is zero.
\section{Sums of Two Square-Zero Matrices} \label{sotsm}
A matrix over an arbitrary field is a sum of two square-zero matrices if and only if all of its invariant polynomials are even or odd polynomials. The result was first proved by Botha \cite{botha3}, and an elegant shortened proof was later presented by De Seguins Pazzis \cite{pazzis2}. The proof presented here is a hybrid between these two approaches: for necessity I follow De Seguins Pazzis, but with a slight shift in focus towards pure matrix theory, and for sufficiency I follow parts of the proof by Botha. 
\subsection{Necessity} 
\begin{quote}
\begin{lem}[Lemma A.6 \cite{pazzis2}]
Let $\mathscr{F}$ be an arbitrary field. If $A,B \in M_n(\mathscr{F})$ are square-zero matrices, then $A$ and $B$ commute with $(A+B)^2$.
\label{lem:sum_commute}
\end{lem}
\end{quote}
{\bf Proof.} Notice that 
\begin{eqnarray}
\nonumber A(A+B)^2 &=& A(A^2 + AB + BA + B^2) \\
\nonumber &=& A(AB + BA) \\
\nonumber &=& A^2B + ABA \\
\nonumber &=& ABA,\end{eqnarray} and 
\begin{eqnarray}
\nonumber (A+B)^2A &=& (A^2 + AB + BA + B^2)A \\
\nonumber &=& (AB + BA)A \\
\nonumber &=& ABA + BA^2 \\
\nonumber &=& ABA.\end{eqnarray}
By a similar argument $B$ commutes with $(A+B)^2$. \hfill $\square$

\begin{quote}
\begin{lem}[Corollary A.7 \cite{pazzis2}]
Let $\mathscr{F}$ be an arbitrary field. If $A \in M_n(\mathscr{F})$ is the sum of two square-zero matrices, and $\text{Dg}[N, B]$ is a matrix similar to $A$ as per Fitting's Lemma \cite[Theorem 5.10]{cullen}, that is, $N$ is nilpotent and $B$ nonsingular, then $B$ is the sum of two square-zero matrices.
\label{lem:subsum_inv}
\end{lem}
\end{quote}
{\bf Proof.} Let $A=S + T$ where $S$ and $T$ are square-zero. Since the square-zero property is invariant under similarity we can assume without loss of generality that  $A = \text{Dg}[N,B]$. Partition $S$ and $T$ compatibly so that 
\[A =  \begin{bmatrix} N & \mathbf{0} \\ \mathbf{0} & B \end{bmatrix} = \begin{bmatrix} S_{11} & S_{12} \\ S_{21} & S_{22} \end{bmatrix} + 
 \begin{bmatrix} T_{11} & T_{12} \\ T_{21} & T_{22} \end{bmatrix}\]

Now since $N$ is nilpotent with order of nilpotence no greater than $n$ 
\[A^{2n} = \begin{bmatrix}
\mathbf{0} & \mathbf{0} \\
\mathbf{0} & B^{2n}
\end{bmatrix},\] and therefore 
\begin{eqnarray}
\nonumber SA^{2n} &=& \begin{bmatrix} S_{11} & S_{12} \\ S_{21} & S_{22} \end{bmatrix}
\begin{bmatrix} \mathbf{0} & \mathbf{0} \\
\mathbf{0} & B^{2n} \end{bmatrix} \\
 &=& \begin{bmatrix} \mathbf{0} & S_{12}B^{2n} \\ \mathbf{0} & S_{22}B^{2n} \end{bmatrix}, \label{commult1}
\end{eqnarray} 
and
\begin{eqnarray}
\nonumber A^{2n}S &=&
\begin{bmatrix} \mathbf{0} & \mathbf{0} \\
\mathbf{0} & B^{2n} \end{bmatrix} \begin{bmatrix} S_{11} & S_{12} \\ S_{21} & S_{22} \end{bmatrix} \\
 &=& \begin{bmatrix} \mathbf{0} & \mathbf{0} \\ B^{2n}S_{21} & B^{2n}S_{22} \end{bmatrix}. \label{commult2}
\end{eqnarray}
Now by Lemma \ref{lem:sum_commute} $S$ and $T$ commute with $A^2$ and hence also with $A^{2n}$, and therefore \eqref{commult1} and \eqref{commult2} are equal, and it follows that $S_{12}B^{2n} = \mathbf{0} = B^{2n}S_{21}$. If we substitute this result into \eqref{commult1} and multiply by $S$ on the left 
\begin{eqnarray}
\nonumber S^2A^{2n} &=& \begin{bmatrix} S_{11} & S_{12} \\ S_{21} & S_{22} \end{bmatrix}
\begin{bmatrix} \mathbf{0} & \mathbf{0} \\ \mathbf{0} & S_{22}B^{2n} \end{bmatrix} \\
\nonumber &=& \begin{bmatrix} \mathbf{0} & S_{12}S_{22}B^{2n} \\ \mathbf{0} & S_{22}^2B^{2n} \end{bmatrix}.
\end{eqnarray}
Now $S$ is square-zero and therefore the matrix above is the zero matrix, and furthermore since $B$ is invertible, $B^{2n}$ is also invertible and it follows that $S_{22}$ is square-zero. By a similar argument it can be shown that $T_{22}$ is square-zero, which completes the proof. \hfill $\square$ 

\begin{quote}
\begin{lem}[Lemma A.2 \cite{pazzis2}]
Let $\mathscr{F}$ be an arbitrary field. Let $p \in \mathscr{F}[x]$ be a nonconstant monic polynomial of degree $n$. Then
\[C(p(x^2)) \approx \begin{bmatrix} \mathbf{0} & C(p(x)) \\ I_n & \mathbf{0} \end{bmatrix}.\]
\label{lem:anti_diag_block_1}
\end{lem}
\end{quote}
{\bf Proof.} Let \[Q = [e_1,e_{n+1},e_2,e_{n+2}, \ldots, e_{n},e_{2n}].\] Then 
\[Q^{-1} \begin{bmatrix} \mathbf{0} & C(p(x)) \\ I_n & \mathbf{0} \end{bmatrix} Q = C(p(x^2)).\] \hfill $\square$

\begin{quote}
\begin{lem}[Corollary A.3 \cite{pazzis2}]
Let $\mathscr{F}$ be an arbitrary field. For any $A \in M_n(\mathscr{F})$ all the invariant polynomials of 
\[F(A) = \begin{bmatrix} \mathbf{0} & A \\ I_n & \mathbf{0} \end{bmatrix}\] are even polynomials.
\label{lem:anti_diag_block_2}
\end{lem}
\end{quote}
{\bf Proof.} Consider $\text{Dg}[F(B_1),F(B_2)]$ where $B_1 \in M_n(\mathscr{F})$ and $B_2 \in M_m(\mathscr{F})$. Let 
\[Q = \begin{bmatrix} I_n & \mathbf{0} & \mathbf{0} & \mathbf{0} \\
				\mathbf{0} & \mathbf{0} & I_n & \mathbf{0} \\
				\mathbf{0} & I_m & \mathbf{0} & \mathbf{0} \\
				\mathbf{0} & \mathbf{0} & \mathbf{0} & I_m
\end{bmatrix},\] Then $F \left ( \text{Dg}[B_1,B_2] \right ) = Q^{-1} \text{Dg}[F(B_1),F(B_2)] Q$. By inductive application of this result to $\text{Dg}[B_1,B_2, \ldots, B_k]$ it follows that 
\begin{equation} F \left ( \text{Dg}[B_1,B_2, \ldots, B_k] \right ) \text{ is similar to }\text{Dg}[F(B_1),F(B_2), \ldots,F(B_k)].\footnote{The similarity transform matrix $Q$ can be applied to $\text{Dg}[B_1,\text{Dg}[B_2, \ldots, B_k]]$, and then another similarity transform employing the matrix $\text{Dg}[I_{j},Q_1]$ where $j$ is the order of $F(B_1)$ can again be applied to the resulting matrix, etcetera, until the desired form is achieved.} \label{eq:FpartsToFwhole} \end{equation}

Now by \cite[Theorem 7.2]{cullen} any $A \in M_n(\mathscr{F})$ is similar to a matrix \[C = \text{Dg}[C(f_1(x)), C(f_2(x)), \ldots, C(f_t(x))]\] where $f_1(x), f_2(x), \ldots, f_t(x)$ are the nonconstant invariant polynomials of $A$. Suppose $P$ is a change-of-basis matrix such that $P^{-1}AP=C$, then 
\begin{eqnarray}
\nonumber F(C) &=& \begin{bmatrix}\mathbf{0} & C \\ I_n & \mathbf{0}\end{bmatrix} \\
\nonumber &=&  \begin{bmatrix}P^{-1} & \mathbf{0} \\ \mathbf{0} & P^{-1}\end{bmatrix} \begin{bmatrix}\mathbf{0} & A \\ I_n & \mathbf{0}\end{bmatrix}\begin{bmatrix}P & \mathbf{0} \\ \mathbf{0} & P\end{bmatrix} \\
\nonumber &=& \left (\text{Dg}[P^{-1},P^{-1}] \right )F(A) \left ( \text{Dg}[P,P] \right ), \end{eqnarray}
so that $F(A)$ and $F(C)$ are similar. Combining this result with \eqref{eq:FpartsToFwhole} and Lemma~\ref{lem:anti_diag_block_1} it follows that 
\begin{eqnarray} 
\nonumber F(A) &\approx& \text{Dg}[F(C(f_1(x))), F(C(f_2(x))), \ldots, F(C(f_t(x)))] \\
\nonumber &\approx& \text{Dg}[C(f_1(x^2)), C(f_2(x^2)), \ldots, C(f_t(x^2))]. 
\end{eqnarray}
For each $i \in \{1,2,\ldots,t\}$ the matrix $C(f_i(x^2))$ is nonderogatory, and it follows that $xI-F(A)$ is equivalent to 
\[\text{Dg}[1,\ldots,1,f_1(x^2),f_2(x^2),\ldots,f_t(x^2)].\]
Now since $f_i(x)$ is monic and divides $f_{i+1}(x)$ for each $i \in \{1,2,\ldots,t-1\}$, it follows that $f_i(x^2)$ is also monic and divides $f_{i+1}(x^2)$, and therefore by uniqueness of the Smith canonical form, it follows that $f_1(x^2),f_2(x^2),\ldots,f_t(x^2)$ are the nonconstant invariant polynomials of $F(A)$, hence all the invariant polynomials of $F(A)$ are even polynomials. \hfill $\square$

\begin{quote}
\begin{lem}[Lemma A.8 \cite{pazzis2}]
Let $\mathscr{F}$ be an arbitrary field. If $A \in M_n(\mathscr{F})$ is invertible and also the sum of two square-zero matrices, then all the invariant polynomials of $A$ are even polynomials.
\label{lem:inv_sum_evenpoly}
\end{lem}
\end{quote}
{\bf Proof.} Suppose $A = S + T$ where $S$ and $T$ are square-zero, and hence $\text{r}(S) \leq \text{n}(S)$ and $\text{r}(T) \leq \text{n}(T)$. By the rank-nullity theorem it follows that $\text{n}(S), \text{n}(T) \geq n/2$. Suppose $v \in \text{N}(S)$ is nonzero: since $A$ is full rank we must have $Av = Sv+Tv = Tv$ is nonzero, that is, $v$ is not in the null space of $T$. By a similar argument any nonzero vector $u$ in the null space of $T$ cannot be in the null space of $S$, and therefore
\[\text{N}(S) \cap \text{N}(T) = \{\mathbf{0}\}.\] It follows that we have a direct sum $\text{N}(S) \oplus \text{N}(T) \subseteq \mathscr{F}^n$, and since $\text{n}(S), \text{n}(T) \geq n/2$ we must have  $\text{N}(S) \oplus \text{N}(T) = \mathscr{F}^n$ and 
\begin{equation} \text{n}(S) = \text{r}(S) = n/2 = \text{n}(T) = \text{r}(T). \label{sum2sim} \end{equation}
Let $\nu = \{n_1, n_2, \ldots, n_{n/2}\}$ be a basis for the null space of $T$: then \[\mu = \{Sn_1, Sn_2, \ldots, Sn_{n/2}\}\] is a linearly independent set which is also a basis for the null space of $S$, and therefore $\nu \cup \mu$ is a basis for $\mathscr{F}^n$ and the matrix representation of $S$ and $T$ with respect to this basis is 
\[\begin{bmatrix} \mathbf{0} & \mathbf{0} \\ I_{n/2} & \mathbf{0} \end{bmatrix} \text{ and } \begin{bmatrix} \mathbf{0} & T_1 \\ \mathbf{0} & \mathbf{0} \end{bmatrix}\] respectively, where $T_1$ is some matrix of full rank. So
\[A \approx \begin{bmatrix} \mathbf{0} & T_1 \\ I_{n/2} & \mathbf{0} \end{bmatrix},\] and therefore by Lemma \ref{lem:anti_diag_block_2} all the invariant polynomials of $A$ are even polynomials. \hfill $\square$
\begin{quote}
\begin{cor}
Let $\mathscr{F}$ be an arbitrary field. If $A \in M_n(\mathscr{F})$ is invertible and also the sum of two square-zero matrices, then $n$ is even.
\label{cor:inv_sum_evenpoly}
\end{cor}
\end{quote}
{\bf Proof.} The result follows directly from \eqref{sum2sim}. \hfill $\square$

\begin{quote}
\begin{thm}[Theorem 1 and 2 \cite{botha3}, Theorem 2.6 \cite{pazzis2}]
Let $\mathscr{F}$ be an arbitrary field. If $A \in M_n(\mathscr{F})$ is the sum of two square-zero matrices, then all the invariant polynomials of $A$ are even or odd polynomials.
\label{thm:2_sqsum_nec}
\end{thm}
\end{quote}
{\bf Proof.} By Fitting's Lemma  \cite[Theorem 5.10]{cullen} and invariance of the square-zero property with respect to similarity we may assume $A=\text{Dg}[N,B]$ where $N$ is nilpotent and $B$ is invertible. The invariant polynomials of $N$ are of the form $x^k$ where $k$ is some non-negative integer, and by Lemma \ref{lem:subsum_inv} and Lemma \ref{lem:inv_sum_evenpoly} all the invariant polynomials of $B$ are even polynomials (observe that $x^0=1$ is also even). 

Let $S_1(x)$ be the Smith canonical form equivalent to $(xI-N)$ and $S_2(x)$ the Smith canonical form equivalent to $(xI-B)$. Now consider $\text{Dg}[S_1(x),S_2(x)]$, which is equivalent to $\text{Dg}[(xI-N),(xI-B)]=(xI - A)$: we want to find the Smith canonical form of this matrix as it contains all the invariant polynomials of $A$, and since equivalence over $\mathscr{F}[x]$ is transitive it suffices to find the Smith canonical form equivalent to $\text{Dg}[S_1(x),S_2(x)]$. 

Now since $B$ is invertible, its invariant polynomials are relatively prime with $x^k$ where $k$ is a positive integer. Therefore the general form of the greatest common divisor among subdeterminants of a fixed order $t \times t$ for $\text{Dg}[S_1(x),S_2(x)]$ is $x^kg(x)$, where $x^k$ is either $x^0 = 1$ or some product of the invariant polynomials of $N$, and $g(x)$ is either $g(x)=1$ or some product of invariant polynomials of $B$. 

Now suppose $x^{k_1}g_1(x)$ is the greatest common divisor among subdeterminants of order $(t-1) \times (t-1)$ of the matrix $\text{Dg}[S_1(x),S_2(x)]$. Now it follows that
\begin{equation} \frac{x^k}{x^{k_1}} \frac{g(x)}{g_1(x)}=x^{k-k_1}\frac{g(x)}{g_1(x)} \label{eq:2_sqsum_nec} \end{equation} is the polynomial in entry $(t,t)$ on the diagonal of the Smith canonical form equivalent to $\text{Dg}[S_1(x),S_2(x)]$, and hence an invariant polynomial of $A$.\footnote{The theory related to Smith canonical form, and how equivalence of characteristic matrices over $\mathscr{F}[x]$ relates to similarity, can be referenced in sections 6.5 and 6.6 of the text by Cullen \cite{cullen}.} 

One of the following cases must hold: 
\begin{itemize}
\item $x^kg(x)=1$ in which case we must have $x^{k_1}g_1(x)=1$ and it follows that \eqref{eq:2_sqsum_nec} reduces to 1, which is an even polynomial.
\item $x^kg(x) = x^k$ for some positive integer $k$, and \eqref{eq:2_sqsum_nec} reduces to $x^{k-k_1}$, which is either an odd or an even polynomial. 
\item $x^kg(x) = g(x) \neq 1$ is some product of invariant polynomials of $B$, and \eqref{eq:2_sqsum_nec} reduces to $\frac{g(x)}{g_1(x)}$. Now since $S_2(x)$ is in Smith canonical form $\frac{g(x)}{g_1(x)}$ is an invariant polynomial of $B$ and hence an even polynomial.
\item $x^kg(x)$ is such that $k$ is a positive integer and $g(x) \neq 1$ is some product of invariant polynomials of $B$. As above, it then follows that $\frac{g(x)}{g_1(x)}$ is an invariant polynomial of $B$ and hence an even polynomial. It follows that if $k-k_1$ is even \eqref{eq:2_sqsum_nec} is an even polynomial, otherwise it is odd.   
\end{itemize}
So all the invariant polynomials of $A$ are even or odd, which concludes the proof. \hfill $\square$ 

\subsection{Sufficiency}
\begin{quote}
\begin{lem}[Lemma 2 \cite{botha3}]
Let $\mathscr{F}$ be an arbitrary field. Let $C=C(p(x))$ be the companion matrix of an even or odd monic polynomial $p(x) \in \mathscr{F}[x]$. Then $C$ is the sum of two square-zero matrices.
\label{lem:companion_sufficient_sum}
\end{lem}
\end{quote}
{\bf Proof.} Without loss of generality we can assume 
\[p(x) = x^n - c_{n-1}x^{n-1} - \cdots - c_1x - c_0, \]
and since $p(x)$ is even or odd, $c_{n-k} = 0$ for odd $k$. It follows that 
\[C = \left [ {\renewcommand{\arraystretch}{1.2}\begin{array}{c|c} 
\mathbf{0} &  c_0 \\
\hline
I_{n-1} & \begin{matrix} c_1 \\ \vdots \\ 0 \\ c_{n-2} \\ 0 \end{matrix}
 \end{array}} \right ],\] with $c_1=0$ if $n$ is even and $c_0=0$ if $n$ is odd. Express $C$ as $[u_1,u_2, \ldots, u_n]$ where $u_i$ indicates column $i$ of $C$ for each $i \in \{1,2,\ldots,n\}$. 

Suppose $p(x)$ is even, let 
\[ S = [u_1,\mathbf{0}, u_3, \mathbf{0}, \ldots, u_{n-1}, \mathbf{0}], \] and
\[ T = [\mathbf{0}, u_2, \mathbf{0}, u_4, \ldots, \mathbf{0}, u_n]. \]
Then $C=S+T$ and $S^2 = T^2 = \mathbf{0}$, so that $C$ is the sum of two square-zero matrices.

Suppose $p(x)$ is odd, let 
\[ S = [\mathbf{0}, u_2, \mathbf{0}, u_4, \ldots, u_{n-1}, \mathbf{0}], \] and
\[ T = [u_1,\mathbf{0}, u_3, \mathbf{0}, \ldots, \mathbf{0}, u_{n}]. \]
Then $C=S+T$ and $S^2 = T^2 = \mathbf{0}$, so that $C$ is the sum of two square-zero matrices. \hfill $\square$

\begin{quote}
\begin{thm}[Theorem 1 and 2 \cite{botha3}, Theorem 2.6 \cite{pazzis2}]
Let $\mathscr{F}$ be an arbitrary field, and suppose $A \in M_n(\mathscr{F})$ is such that all of its invariant polynomials are even or odd polynomials. Then $A$ is the sum of two square-zero matrices.
\label{thm:sufficient_sqzero_sum}
\end{thm}
\end{quote}
{\bf Proof.} By Theorem 7.2 in the text by Cullen \cite{cullen} $A$ is similar to \[\text{Dg}[C(p_1(x)),C(p_2(x)), \ldots, C(p_t(x))],\] where $p_1(x),p_2(x), \ldots, p_t(x)$ are the nonconstant invariant polynomials of $A$. Since the invariant polynomials of $A$ are all odd or even, it follows from Lemma \ref{lem:companion_sufficient_sum} that each $C(p_i(x))$ can be written as a sum $S_i + T_i$ where $S_i,T_i$ are square-zero matrices, $i \in \{1,2, \ldots,t\}$. So we have
\[A \approx \text{Dg}[S_1, S_2, \ldots, S_t] + \text{Dg}[T_1, T_2, \ldots, T_t]. \] 
Observe that 
\[ \text{Dg}[S_1, S_2, \ldots, S_t]^2 = \text{Dg}[S_1^2, S_2^2, \ldots, S_t^2] = \mathbf{0}, \] and similarly it can be shown that $\text{Dg}[T_1, T_2, \ldots, T_t]$ is square-zero. Since the square-zero property is invariant with respect to similarity, $A$ is therefore the sum of two square-zero matrices. \hfill $\square$

Notice in particular that a nonzero nilpotent matrix is the sum of two square-zero matrices. 

\subsection{Corollaries}
Combining Theorems \ref{thm:2_sqsum_nec} and \ref{thm:sufficient_sqzero_sum} yields the main result of this section, but there are several results proceeding from the main finding, which are interesting in their own right.

\begin{quote}
\begin{lem}[Lemma 4 \cite{botha3}]
Let $\mathscr{F}$ be an arbitrary field, and suppose 
\[p(x) = x^n - c_{n-1}x^{n-1} - \cdots - c_1x - c_0\] is a polynomial in $\mathscr{F}[x]$. Then the only nonconstant invariant polynomial of $-C(p(x))$ is
\[x^n - (-1)c_{n-1}x^{n-1} - \cdots - (-1)^i c_{n-i} x^{n-i} - \cdots - (-1)^nc_0.\]
\label{lem:poly_minusC}
\end{lem}
\end{quote}
{\bf Proof.} A companion matrix is nonderogatory and therefore has only one nonconstant invariant polynomial, namely its minimum polynomial. Observe that 
\[-C(p(x)) = \left [ {\renewcommand{\arraystretch}{1.2}\begin{array}{c|c} 
\mathbf{0} &  -c_0 \\
\hline
-I_{n-1} & \begin{matrix} -c_1 \\ \vdots \\ -c_{n-2} \\ -c_{n-1} \end{matrix}
 \end{array}} \right ].\] Let $P=[e_1,-e_2,e_3,-e_4,\ldots,(-1)^{n-2}e_{n-1},(-1)^{n-1}e_n]$, then
\[P^{-1}(-C(p(x)))P = \left [ {\renewcommand{\arraystretch}{1.2}\begin{array}{c|c} 
\mathbf{0} &  (-1)^{n}c_0 \\
\hline
I_{n-1} & \begin{matrix} (-1)^{n-1}c_1 \\ \vdots \\ (-1)^{i}c_{n-i} \\ \vdots \\ c_{n-2} \\ -c_{n-1} \end{matrix}
 \end{array}} \right ],\] which is a companion matrix, and therefore the only nonconstant invariant polynomial of $-C(p(x))$ is \[ x^n - (-1)c_{n-1}x^{n-1} - \cdots - (-1)^i c_{n-i} x^{n-i} - \cdots - (-1)^nc_0. \] \hfill $\square$  

\begin{quote}
\begin{cor}
Let $\mathscr{F}$ be a field which is not of characteristic 2, and suppose $p(x) \in \mathscr{F}[x]$ is monic. Then $C(p(x))$ is similar to $-C(p(x))$ if and only if $p(x)$ is an even- or odd-power polynomial.
\label{cor:poly_minusC}
\end{cor}
\end{quote}
{\bf Proof.} If 
\[p(x) = x^n - c_{n-1}x^{n-1} - \cdots - c_1x - c_0\] is even or odd, then $c_{n-i}=0$ for all odd $i \in \{1,2,\ldots, n\}$. Therefore, by Lemma~\ref{lem:poly_minusC}, $C(p(x))$ and $-C(p(x))$ have the same minimum polynomial, which is then also the only nonconstant invariant polynomial for both these matrices, and therefore they are similar.

Conversely, by Lemma \ref{lem:poly_minusC}, if $C(p(x))$ and $-C(p(x))$ are similar we must have
\[x^n - c_{n-1}x^{n-1} - \cdots - c_1x - c_0 = x^n - (-1)c_{n-1}x^{n-1} - \cdots - (-1)^i c_{n-i} x^{n-i} - \cdots - (-1)^nc_0,\] which is true if and only if $c_{n-i} = (-1)^ic_{n-i}$ for each $i \in \{1,2, \ldots, n\}$. Now this requirement is always fulfilled if $i$ is even, but if $i$ is odd we must have $c_{n-i} = -c_{n-i}$, which in a field of characteristic not equal to two implies that $c_{n-i}=0$. Therefore $p(x)$ is an even- or odd-power polynomial. \hfill $\square$

For the proof of the next result we require the following definition.
\begin{quote}
\begin{defn}[\cite{botha3}]
Let $f(x)=x^n + \sum_{i=1}^n c_{n-i}x^{n-i}$ be an element of $\mathscr{F}[x]$, where $\mathscr{F}$ is an arbitrary field. Define $^*:\mathscr{F}[x] \rightarrow \mathscr{F}[x]$ as
\[f^*(x) = (-1)^nf(-x) = x^n + \sum_{i=1}^n (-1)^ic_{n-i}x^{n-i}.\]
\label{defn:starpoly_function}
\end{defn}
\end{quote}
It is easy to check the following properties of the defined $^*$-mapping: for any $f,g \in \mathscr{F}[x]$
\begin{itemize}
\item $(f^*)^* = f$
\item $(fg)^* = f^*g^*$
\item If the characteristic of $\mathscr{F}$ is not two then $f$ is even- or odd-power if and only if $f^* = f$
\end{itemize}
We can employ the properties above to prove that $f \in \mathscr{F}[x]$ is irreducible if and only if $f^*$ is irreducible. Explicitly, suppose $f(x)=q(x)r(x)$ (that is, $f$ is reducible), then $q^*(x)r^*(x) = (q(x)r(x))^* = f^*(x)$ so that $f^*$ is also reducible. Now suppose $f^*$ is reducible, say $f^*(x)=q(x)r(x)$, then $q^*(x)r^*(x) = (q(x)r(x))^* = (f^*(x))^* = f(x)$ so that $f$ is also reducible.

Also notice that by Lemma \ref{lem:poly_minusC} the only invariant polynomial of $-C(f(x))$ is $f^*(x)$.
\begin{quote}
\begin{cor}[Theorem 2 \cite{botha3}, Theorem 2.11 \cite{ww}]
Let $\mathscr{F}$ be a field which is not of characteristic 2. The matrix $A \in M_n(\mathscr{F})$ is a sum of two square-zero matrices if and only if $A$ is similar to $\text{Dg}[N,X,-X,C]$ where $N$ is nilpotent, $X$ is invertible and $C$ is similar to a block diagonal matrix with diagonal blocks of the form $C(p(x)^e)$ where $p(x)$ is an irreducible even-power polynomial with nonzero constant term.
\label{cor:X_and_minusX}
\end{cor}
\end{quote}
{\bf Proof.} By Fitting's Lemma \cite[Theorem 5.10]{cullen} any matrix $A$ is similar to $\text{Dg}[N,B]$ where $N$ is nilpotent and $B$ is nonsingular. 

Suppose first that $A$ is the sum of two square-zero matrices. By Lemma~\ref{lem:subsum_inv} $B$ is then also a sum of two square-zero matrices, and since $B$ is invertible, by Lemma~\ref{lem:inv_sum_evenpoly} all the invariant polynomials of $B$ are even. So for each invariant polynomial $p(x)$ of $B$, we have $p(x) = p^*(x)$ where $^*$ is the mapping as defined in Definition \ref{defn:starpoly_function}. Now 
\[p(x) = d_1(x)d_2(x)\ldots d_k(x),\] where $d_i(x) = (p_i(x))^{e_i}$ for each $i \in \{1,2,\ldots,k\}$ where the set $\{p_1,p_2,\ldots,p_k\}$ consists of distinct, monic, irreducible polynomials. Now by the properties of the $^*$-mapping as listed before the statement of this corollary 
\[p^*(x) = d_1^*(x)d_2^*(x)\ldots d_k^*(x) = (p_1^*(x))^{e_1}(p_2^*(x))^{e_2}\ldots (p_k^*(x))^{e_k},\] and $p_1^*,p_2^*,\ldots,p_k^*$ are distinct, monic and irreducible since $p_1,p_2,\ldots,p_k$ are distinct, monic and irreducible. Now since $p(x)=p^*(x)$ we must then have, either $d_i(x)=d_i^*(x)$, or there exists a pair $(d_i(x), d_j(x)) = (d_i(x), d_i^*(x))$; if this is the case $d_i(x)$ and $d_i^*(x)$ must be relatively prime. This is so because if $d_i(x) = (p_i(x))^{e_i} \neq (p_i^*(x))^{e_i} = d_i^*(x)$, we must have $p_i(x) \neq p_i^*(x)$, and since $p_i(x)$ and $p_i^*(x)$ are also irreducible, they are relatively prime.

Now by Theorem 7.3 in the text by Cullen \cite{cullen}
\[C(p(x)) = \text{Dg}[C(d_1(x)),C(d_2(x)), \ldots, C(d_k(x))].\] Let $C'$ consist of all companion matrices associated with $d_i(x)$ where $d_i(x)=d_i^*(x)$. Now since $B$ is invertible, each block in $C'$ must be invertible, and hence $d_i(x)$ must have a non-zero constant term. Furthermore, let $X'$ consist of companion matrices associated with $d_i(x)$ where $d_i(x)$ is such that it is the first component of a pair $(d_i(x),d_i^*(x))$. By Lemma \ref{lem:poly_minusC}: $C(d_i^*(x)) \approx -C(d_i(x))$. Note that $X'$ is invertible by invertibility of $B$. Combining these results we have 
\begin{equation} C(p(x)) \approx \text{Dg}[X',-X',C'].\label{eq:comp_x-xc} \end{equation} 

Let $f_1(x),f_2(x), \ldots, f_t(x)$ be the non-constant invariant polynomials of $B$. Now since $B$ is similar to \[\text{Dg}[C(f_1(x)),C(f_2(x)), \ldots, C(f_t(x))],\] we can apply \eqref{eq:comp_x-xc} to each matrix $C(f_i(x))$, and followed by a simple permutation to rearrange the blocks of $B$, we find that $B$ is similar to a matrix of the form $\text{Dg}[X,-X,C]$ where $X$ and $C$ are as defined in the statement of this corollary.

It follows that $A$ is similar to a matrix $\text{Dg}[N,X,-X,C]$, which is of the desired form.

Conversely, suppose $A$ is similar to $\text{Dg}[N,X,-X,C]$. Now $N$ consists of invariant polynomials which are constant or of the form $x^t$, so that all the invariant polynomials of $N$ are either even or odd. By Theorem \ref{thm:sufficient_sqzero_sum} $N$ is a sum of two square zero matrices. Furthermore, by Lemma \ref{lem:companion_sufficient_sum} and invariance of the square-zero property with respect to similarity $C$ is a sum of two square-zero matrices. Finally observe that 
\[\text{Dg}[X,-X] = \frac{1}{2}\begin{bmatrix} X & -X \\ X & -X \end{bmatrix} +
\frac{1}{2}\begin{bmatrix} X & X \\ -X & -X \end{bmatrix},\] and it is easy to verify that each of the matrices on the right-hand side are square-zero. Combining all of these results, we have $A$ is a sum of two square-zero matrices as desired. \hfill $\square$

\begin{quote}
\begin{cor}[Theorem 2 \cite{botha3}, Corollary 2.7 \cite{pazzis2}, Theorem 2.11 \cite{ww}]
Let $\mathscr{F}$ be a field which is not of characteristic 2. The matrix $A \in M_n(\mathscr{F})$ is a sum of two square-zero matrices if and only if $A$ is similar to $-A$.
\label{cor:A_sim_minusA}
\end{cor}
\end{quote}
{\bf Proof.} By Theorem 7.2 in the text by Cullen \cite{cullen} $A$ is similar to \[\text{Dg}[C(p_1(x)),C(p_2(x)), \ldots, C(p_t(x))],\] where $p_1(x),p_2(x), \ldots, p_t(x)$ are the nonconstant invariant polynomials of $A$. 

Suppose $A$ is a sum of two square-zero matrices. By Theorem \ref{thm:2_sqsum_nec} all the invariant polynomials of $A$ are even or odd polynomials. Now $-A$ is similar to 
\[\text{Dg}[-C(p_1(x)),-C(p_2(x)), \ldots, -C(p_t(x))],\] and since $p_1(x),p_2(x), \ldots, p_t(x)$ are all even or odd, by Corollary \ref{cor:poly_minusC} we must have $C(p_i(x))$ is similar to $-C(p_i(x))$ for each $i \in \{1,2, \ldots, t\}$. It follows that
\begin{equation} \text{Dg}[-C(p_1(x)),-C(p_2(x)), \ldots, -C(p_t(x))] \approx \text{Dg}[C(p_1(x)),C(p_2(x)), \ldots, C(p_t(x))], \label{eq:minusc_sim_c} \end{equation} and therefore, since similarity is an equivalence relation, $A$ is similar to $-A$.

Now suppose $A$ is similar to $-A$. By Theorem 7.5 in the text by Cullen \cite{cullen} $A$ is similar to \[B = \text{Dg}[C(d_1(x)), C(d_2(x)), \ldots, C(d_s(x))]\] where $d_1(x), d_2(x), \ldots, d_s(x)$ are the elementary divisors of $A$, so that we can write $d_i(x) = (p_i(x))^{e_i}$ where $p_i(x)$ is monic and irreducible and $e_i$ is a positive integer. Now since $A \approx -A$ we must have 
\begin{eqnarray} 
\nonumber B &\approx& -\text{Dg}[C(d_1(x)), C(d_2(x)), \ldots, C(d_s(x))] \\
\nonumber &=& \text{Dg}[-C(d_1(x)), -C(d_2(x)), \ldots, -C(d_s(x))] \\
\nonumber &=& \text{Dg}[C(d_1^*(x)), C(d_2^*(x)), \ldots, C(d_s^*(x))],
\end{eqnarray} 
where the last step follows by Lemma \ref{lem:poly_minusC}, that is the only invariant polynomial of $-C(f(x))$ is $f^*(x)$. Now by the properties of the $^*$-mapping it follows that $d_i^*(x) = (p_i^*(x))^{e_i}$, and since $p_i(x)$ is monic and irreducible, so is $p_i^*(x)$. It follows that $d_i^*(x)$ is an elementary divisor of $-A$ for each $i \in \{1,2,\ldots,s\}$.

By Theorem 7.4 in the text by Cullen \cite{cullen}, two matrices are similar if and only if they have the same elementary divisors, and therefore we must have either $d_i(x) = d_i^*(x)$, or among the elementary divisors of $A$ there exists a pair $(d_i(x), d_i^*(x))$ where $d_i(x)$ and $d_i^*(x)$ must be relatively prime (the same argument as applied in the proof of Corollary \ref{cor:X_and_minusX} is valid). To summarize these results: we can apply a permutation similarity transformation to $\text{Dg}[C(d_1(x)), C(d_2(x)), \ldots, C(d_s(x))]$, and reindex the diagonal blocks in such a way as to show that $A$ is similar to
\[\text{Dg}[C(d_{1}(x)),C(d_{1}^*(x)),\ldots, C(d_{m}(x)),C(d_{m}^*(x)), C(d_{m+1}(x)), \ldots, C(d_{s-m}(x))],\] where $d_i(x) = d_i^*(x)$ for $m+1 \leq i \leq s-m$. Now by the remarks that follow Definition \ref{defn:starpoly_function} each of the blocks $C(d_{m+1}(x)), C(d_{m+2}(x)), \ldots, C(d_{s-m}(x))$ are even- or odd-power and are therefore sums of two square-zero matrices by Theorem \ref{thm:sufficient_sqzero_sum}. 

It remains to prove that each block $\text{Dg}[C(d_i(x)), C(d_i^*(x))]$ $(1 \leq i \leq m)$ is a sum of two square-zero matrices. Now in this case $d_i(x)$ must have a nonzero constant term, for if this is not so, since $d_i(x)$ is a power of some irreducible polynomial, we must have $d_i(x) = x^{e_i}$, in which case $d_i^*(x) = x^{e_i}$ which is a contradiction. It follows that $C(d_i(x)) = X$ where $X$ is an invertible matrix, and by Lemma \ref{lem:poly_minusC} it follows that $C(d_i^*(x)) = -X$. By Corollary \ref{cor:X_and_minusX} it follows that $\text{Dg}[C(d_i(x)), C(d_i^*(x))]$ is a sum of two square-zero matrices, which concludes the proof. \hfill $\square$

\begin{quote}
\begin{cor}[Theorem 2 \cite{botha3}, Theorem 2.11 \cite{ww}]
Let $\mathscr{F}$ be a field which is not of characteristic 2. The matrix $A \in M_n(\mathscr{F})$ is a sum of two square-zero matrices if and only if there exists an involution~$V$ such that $AV=-VA$.
\label{cor:AV_eq_minusVA}
\end{cor}
\end{quote}
{\bf Proof.} If $AV=-VA$ for some involution $V$ then $A$ is similar to $-A$ and therefore, by Corollary \ref{cor:A_sim_minusA} $A$ is the sum of two square-zero matrices.

Suppose $A$ is the sum of two square-zero matrices. As before, by Theorem 7.2 in the text by Cullen \cite{cullen} $A$ is similar to \[C = \text{Dg}[C(p_1(x)),C(p_2(x)), \ldots, C(p_t(x))],\] where $p_1(x),p_2(x), \ldots, p_t(x)$ are the nonconstant invariant polynomials of $A$, and since $A$ is a sum of two square-zero matrices the similarity indicated in \eqref{eq:minusc_sim_c} holds (recall that the matrix on the left is similar to $-A$). 

Now consider the matrix $P$ employed in the proof of Lemma \ref{lem:poly_minusC}: this matrix is an involution with the property $P^{-1}(-C(p(x)))P$ is a companion matrix for any monic polynomial $p(x)$. Now for each $i \in \{1,2,\ldots,t\}$ let $P_i$ be an involution of this form and of appropriate order such that $P_i^{-1}(-C(p_i(x)))P_i$ is defined. Note that this product is then a companion matrix as shown in Lemma \ref{lem:poly_minusC}, but now since \eqref{eq:minusc_sim_c} holds we also have that $-C(p_i(x))$ is similar to the companion matrix $C(p_i(x))$. By uniqueness of this form it follows that \[C(p_i(x))P_i^{-1} = -P_i^{-1}C(p_i(x)).\] Let $Q = \text{Dg}[P_1^{-1}, P_2^{-1}, \ldots, P_t^{-1}]$, then it follows that $Q$ is an involution and $CQ = -QC$. Furthermore if $R$ is the change-of-basis matrix such that $R^{-1}AR = C$,  then we must have $ARQR^{-1} = -RQR^{-1}A$. It is easy to check that the involution property is invariant with respect to similarity, and therefore $RQR^{-1}$ is an involution, completing the proof. \hfill $\square$

\begin{quote}
\begin{cor}[Corollary 2.8 \cite{pazzis2}]
Let $\mathscr{F}$ be a field which is of characteristic 2. The matrix $A \in M_n(\mathscr{F})$ is a sum of two square-zero matrices if and only if, in some algebraic closure of $\mathscr{F}$, the simple Jordan blocks of $A$ associated with nonzero eigenvalues all have even order.
\label{cor:2sum_char2_evenBlocks}
\end{cor}
\end{quote}
{\bf Proof.} Suppose $A$ is a sum of two square-zero matrices. By Fitting's Lemma \cite[Theorem 5.10]{cullen} $A$ is similar to $\text{Dg}[N,B]$ where $N$ is nilpotent and $B$ is invertible. By Lemma \ref{lem:subsum_inv} $B$ is a sum of two square-zero matrices, and therefore by Lemma \ref{lem:inv_sum_evenpoly} all the invariant polynomials of $B$ are even polynomials. 

Let $f(x)$ be an invariant polynomial of $B$. Then $f(x) = g(x^2)$ for some monic polynomial $g(x)$. Now over an algebraic closure of $\mathscr{F}$ (denote this field as $\mathscr{F}'$) we can write
\[g(x) = (x+\lambda_1^2)(x+\lambda_2^2)\cdots(x+\lambda_{m/2}^2),\] where $m$ is the degree of $f(x)$, and $\lambda_1, \lambda_2, \cdots, \lambda_{m/2}$ are non-zero elements in $\mathscr{F}'$ (which are not necessarily distinct). It follows that 
\[f(x) = g(x^2) = (x^2+\lambda_1^2)(x^2+\lambda_2^2)\cdots(x^2+\lambda_{m/2}^2),\] and since $\mathscr{F}'$ is a field of characteristic two \[(x-\lambda)^2 = x^2 - 2 \lambda x +\lambda^2 = x^2 + \lambda^2\] for any $\lambda \in \mathscr{F}'$. So
\[f(x) = (x-\lambda_1)^2(x-\lambda_2)^2\cdots(x-\lambda_{m/2})^2,\] and therefore, over $\mathscr{F}'$, all the simple Jordan blocks of $B$ are even order, which concludes the first part of the proof.

Consider the converse: suppose in some algebraic closure of $\mathscr{F}$ (let this field be $\mathscr{F}'$), the simple Jordan blocks of $A$ associated with nonzero eigenvalues all have even order. Again, by Fitting's Lemma \cite[Theorem 5.10]{cullen} $A$ is similar to $\text{Dg}[N,B]$ where $N$ is nilpotent and $B$ is invertible. Since $N$ is nilpotent it is a sum of two square-zero matrices, and the matrix $B$ will contain the Jordan blocks of $A$ associated with nonzero eigenvalues (over $\mathscr{F}'$).

It follows that every elementary divisor of $B$ is of the form $(x-\lambda)^{2e}$ over $\mathscr{F}'$. The binomial expansion of an elementary divisor of the form $(x-\lambda)^{2e}$, with $\lambda \neq 0$, is
\[\sum_{k=0}^{2e} \binom{2e}{k} \lambda^kx^{2e-k}.\]
Now by Lucas's Theorem \cite{fine} if $k$ is odd it follows that $\binom{2e}{k}$ is even, and therefore each term of odd degree has an even binomial coefficient. Now since $\mathscr{F}'$ is of characteristic 2, each term of odd degree is therefore annihilated in the expansion of $(x-\lambda)^{2e}$, and therefore each simple Jordan block with nonzero eigenvalue is associated with an even-power elementary divisor. It follows that each invariant polynomial of $B$ is a product of even-power polynomials, and hence is also even.

Now by Theorem 7.20 in the text by Roman \cite{roman} the invariant polynomials of~$B$ must be polynomials over $\mathscr{F}$, and therefore all the invariant polynomials of $B$ over  $\mathscr{F}$ must be even. By Theorem \ref{thm:sufficient_sqzero_sum} $B$ is a sum of two square-zero matrices.

By invariance of the square-zero property with respect to similarity, $A$ is a sum of two square-zero matrices. \hfill $\square$

The results obtained in this section provide a complete characterization of sums of two square-zero matrices. I will now proceed to sums of four square-zero matrices, as this result is also completely resolved, before I address sums of three square-zero matrices which still requires a complete characterization. Notice that although not explicitly stated, a sum of two square-zero matrices must have trace zero by Theorem~\ref{thm:nilpotentSums}. It is also easy to see that if the invariant polynomials of a matrix are odd or even then such a matrix has trace zero: the companion matrix of an odd or even polynomial has only zero on the diagonal, and similar matrices have the same trace. So we see that a zero trace is also a necessary condition for sums of two square-zero matrices, but more can be said regarding this requirement in the next section.

\section{Sums of Four Square-Zero Matrices}
This result was first obtained by Wang and Wu over the complex field, and later shown by Pazzis to be true over an arbitrary field.
\begin{quote}
\begin{thm}[Theorem 1.1 \cite{pazzis2}, Theorem 3.6 \cite{ww}]
Let $\mathscr{F}$ be an arbitrary field. The matrix $A \in M_n(\mathscr{F})$ is a sum of four square-zero matrices if and only if the trace of $A$ is zero. 
\label{thm:trace_zero_sum_of_4}
\end{thm}
\end{quote}
{\bf Proof.} By Theorem \ref{thm:nilpotentSums}, if $A$ is a sum of square-zero matrices, the trace of $A$ is zero.

Now suppose the trace of $A$ is zero. First suppose that $A$ is not scalar. Then by Theorem \ref{thm:nilpotentSums} $A$ is a sum of two nilpotent matrices, and since each of the invariant polynomials of a nilpotent matrix is either odd- or even-power it follows by Theorem \ref{thm:sufficient_sqzero_sum} that each of these nilpotent matrices is a sum of two square-zero matrices, concluding the proof for this case.

Now suppose $A$ is scalar, that is $A=\lambda I_n$. In this case the characteristic of $\mathscr{F}$ (let this characteristic be denoted by $c$) divides $n$. I will now prove that $I_c$ is a sum of four square-zero matrices, from which the desired result will follow. 

First I will show that if the characteristic of $\mathscr{F}$ is $c$ it is always possible to choose a monic polynomial $p(x)$ of degree $c$ such that $p(x)$ and $p(x-1)$ are both even- or odd-power. If $c=2$ let $p(x) = x^2$, then 
\[p(x-1) = x^2 - 2x +1 = x^2 + 1.\] It follows that both $p(x)$ and $p(x-1)$ are even-power polynomials, so that $C(p(x))$ and $C(p(x-1))$ are both sums of two square-zero matrices by Theorem \ref{thm:sufficient_sqzero_sum}. If $c \neq 2$ let $p(x) = x^c - x$. Now notice that in the binomial expansion of $(x-1)^c$ each term $\binom{c}{k} (-1)^kx^{c-k}$ with $0<k<c$ is equal to zero. Explicitly, since $c$ is a prime it follows by Lucas's Theorem \cite{fine} that the coefficient $\binom{c}{k}$ is divisible by $c$, and therefore must be zero. It follows that  
\[p(x-1) = (x-1)^c - (x-1) = x^c - 1 - x + 1 = x^c - x = p(x),\] so that both $p(x)$ and $p(x-1)$ are odd-power polynomials, and therefore $C(p(x)) = C(p(x-1))$ is a sum of two square-zero matrices (again by Theorem \ref{thm:sufficient_sqzero_sum}).  

Now for any monic polynomial $p(x)$ we have 
\[ \det (xI - (I+C(p(x)))) = \det ((x-1)I-C(p(x))) = p(x-1), \] and it follows that the characteristic polynomials of the matrices $I+C(p(x))$ and $C(p(x-1))$ are the same. Furthermore $I+C(p(x))$ is nonderogatory, for consider the $(c-1) \times (c-1)$ submatrix obtained by deleting the first row and last column of $xI-(I + C(p(x)))$: this submatrix has determinant $(-1)^{c-1}$. Therefore the greatest common divisor among the $(c-1) \times (c-1)$ subdeterminants of $xI-(I + C(p(x)))$ is 1, and it follows that the Smith canonical form equivalent to $xI-(I + C(p(x)))$ is
\[ \text{Dg}[1,1,\ldots , 1, p(x-1)] .\]  
Now $xI-C(p(x-1))$ has the same Smith canonical form as above, and therefore by Theorem 6.16 in the text by Cullen \cite{cullen} $C(p(x-1))$ and $I+C(p(x))$ are similar. 

Now applying the preceding results where $C(p(x))$ and $C(p(x-1))$ are both sums of square-zero matrices (by appropriate choice of $p(x)$ this is always possible over a field of characteristic $c$ as shown above), it follows by invariance of the square-zero property with respect to similarity that $I_c$ is a sum of four square-zero matrices. Explicitly, since $I_c+C(p(x)) \approx C(p(x-1))$, it follows that $I_c = -C(p(x)) + P^{-1} C(p(x-1)) P$ (where $P$ is some change of basis matrix), and with appropriate choice of $p(x)$ (which is always possible since the characteristic of $\mathscr{F}$ is $c$) both $-C(p(x))$ and $P^{-1} C(p(x-1)) P$ are sums of square-zero matrices.

Finally, since $c$ divides $n$ we have
\[ A = \lambda I_n = \lambda \text{Dg}[I_c, I_c, \ldots, I_c],\] where each $I_c$ is a sum of four square-zero matrices, and it follows that $A$ is a sum of four square-zero matrices as desired. \hfill $\square$ 

The necessary condition in the theorem above indicates that a matrix which has trace not equal to zero cannot be a sum consisting of only square zero matrices. Adding sufficiency, this theorem is the most general requirement for a matrix to be written as a sum of square-zero matrices. 

\section{Sums of Three Square-Zero Matrices}
In the previous section it was shown that any matrix that can be written as a sum of square-zero matrices, can be written as a sum of four square-zero matrices. And in the preceding section it was shown that sums of two square-zero matrices have invariant polynomials that are even- or odd-power. It remains to characterize sums of three square-zero matrices, which at the time of writing has not yet been fully resolved. 

Historically, Wang and Wu \cite{ww} gave a necessary condition for a matrix to be a sum of three square-zero matrices, as well as some special cases for sufficiency. They conjectured that the necessary condition might also be sufficient. Takahashi \cite{takahashi} then proved that there is a class of matrices that satisfy the necessary condition as set forth by Wang and Wu, but that is not a sum of three square-zero matrices, and then proved a special case for matrices with a minimum polynomial of degree two. It is worth noting that the results presented in section \ref{sotsm} were not fully available to Wang and Wu, and Takahashi (Wang and Wu did prove some of these results), and the new results could significantly aid in characterizing sums of three square-zero matrices. In fact de Seguins Pazzis \cite{pazzis2} employed these results to prove that for the special case of matrices over a field of characteristic two, any trace zero matrix can be written as a sum of three square-zero matrices. De Seguins Pazzis also characterised sums of three square-zero matrices where the matrix in question may be augmented with rows and columns of zero vectors.

We can appreciate the complexity of the problem by considering an example. Let $S_1,S_2,S_3$ be square-zero. By the preceding sections
\[ S_1 + S_2 + S_3 = C + S_3 \] where $C$ is a matrix with only even- or odd-power invariant polynomials. By invariance of the square-zero property with respect to similarity we can assume \[C=\text{Dg}[C(f_1(x)),C(f_2(x)), \ldots, C(f_t(x))],\] where $f_i(x)$ is a nonconstant even- or odd-power invariant polynomial of $C$ for each $i \in \{1,2,\ldots,t\}$. Let 
\[C=C(x^5 - x^3 - x),\] and 
\[S_3 = \begin{bmatrix}
1 & 1 & 0 & 0 & 0 \\
-1 & -1 & 0 & 0 & 0 \\
0 & 0 & 0 & 0 & 0 \\
0 & 0 & 0 & 0 & 0 \\
0 & 0 & 0 & 0 & 0
\end{bmatrix},\] so that $C+S_3$ is the sum of three square-zero matrices. Note that $C+S_3$ is nonderogatory\footnote{The determinant of the matrix obtained by deleting the third row and first column of $xI-(C+S_3)$ is 1.} with characteristic polynomial $x^5-2x^3+1$ and we can therefore be sure that $C+S_3$ is not a sum of two square-zero matrices. Now the purpose of this example is to show that $C$ is a nonderogatory matrix with odd-power minimum polynomial, but adding the relatively simple square-zero matrix $S_3$ alters the similarity structure in a way that may be difficult to characterize. It is difficult to say more about the problem based on this single example, but it is clear that the fact that $S_3$ could be any square-zero matrix increases the complexity of the problem significantly.  

Can we answer the following question: What specific characterizable properties does the sum of an arbitrary square-zero matrix and a  matrix with only even- or odd-power invariant polynomials have? The rational canonical form of such a matrix, and therefore its invariant polynomials and elementary divisors should reveal its unique structure, but at the time of writing the \emph{specific} properties remained elusive.

Let us start with a necessary and sufficient condition which is an easy corollary to the preceding work, and will be useful in much of what is to follow.
\begin{quote}
\begin{cor}
Let $\mathscr{F}$ be an arbitrary field. The matrix $A \in M_n(\mathscr{F})$ is a sum of three square-zero matrices if and only if it can be written as a sum $C+S$, where $C$ is a matrix with only odd- or even-power invariant polynomials, and $S$ is an arbitrary square-zero matrix of compatible order. 
\label{cor:3sqz_is_sqzPlusC}
\end{cor}
\end{quote}
{\bf Proof.} The result follows directly from the fact that $C$ is a sum of two square-zero matrices if and only if all the invariant polynomials of $C$ are odd- or even-power polynomials (Theorem \ref{thm:2_sqsum_nec} and Theorem \ref{thm:sufficient_sqzero_sum}). \hfill $\square$

The following result provides a mechanism whereby we can construct a nonderogatory matrix with a characteristic polynomial of choice, by adding certain elements to a matrix in the canonical form as described in Theorem 7.2 in the text by Cullen \cite{cullen}.
\begin{quote}
\begin{prop}[Lemma 11 \cite{pazzis3}]
Let $\mathscr{F}$ be an arbitrary field, and let 
\[p_1(x) = x^k - \sum_{i=0}^{k-1} a_ix^i,\] and 
\[p_2(x) = x^{n-k} - \sum_{i=0}^{n-k-1} b_ix^i\] where each coefficient $a_i, b_i$ is in $\mathscr{F}$. 

Choose any polynomial 
\[p(x) = x^{n} - \sum_{i=0}^{n-1} c_ix^i,\] subject to the constraint $c_{n-1} = a_{k-1} + b_{n-k-1}$. Let $N \in M_{(n-k) \times k}(\mathscr{F})$ be the matrix with a one in entry $(1,k)$ and zeros elsewhere. Then there exists a matrix $D=(d_{ij}) \in M_{k \times (n-k)}(\mathscr{F})$ such that $C(p(x))$ is similar to
\[ M(D) = \begin{bmatrix} 
C(p_1(x)) & D \\
N & C(p_2(x))
\end{bmatrix}. \]
\label{prop:choice_of_polynomial}
\end{prop}
\end{quote}
{\bf Proof.} First, notice that $M(D)$ is nonderogatory, since the determinant of the $(n-1) \times (n-1)$ submatrix obtained from $xI_n - M(D)$ by deleting the first row and last column, is $(-1)^{n-1}$. Therefore, if we can prove that the characteristic polynomial of $M(D)$ is $p(x)$, the proof will be complete.

Let us determine the characteristic polynomial of $M(D)$ in terms of the polynomials $p_1(x), p_2(x), p(x)$. Apply the following elementary row operation algorithm:
\begin{enumerate}
\item Let $i=k-1$.
\item Replace $\text{row}_i(M(D)-xI)$ by $\text{row}_{i}(M(D)-xI)+x \cdot \text{row}_{i+1}(M(D)-xI)$.
\item If $i=1$ stop, else let $i=i-1$, and repeat the algorithm from step 2 above.
\end{enumerate}
It follows that
\[M(D) - xI_n \stackrel{R}{\sim} 
\left [ {\renewcommand{\arraystretch}{1.2}\begin{array}{c|c|c} 
\mathbf{0} &  -p_1(x) & \begin{matrix} d_1(x) & d_2(x) & \cdots & d_{n-k}(x) \end{matrix} \\ 
\hline
I_{k-1} & \begin{matrix}  
(-p_1(x)-a_0)/x \\  
(-p_1(x)-a_1x-a_0)/x^2 \\ 
\vdots \\ 
-x^2 + xa_{k-1} + a_{k-2} \\ 
-x+a_{k-1} \end{matrix}
 & D_1(x) \\ \hline
\mathbf{0} & e_1 & C(p_2(x))-xI_{n-k} 
 \end{array}} \right ],\] where the row vector $\begin{bmatrix} d_1(x) & d_2(x) & \ldots & d_{n-k}(x)\end{bmatrix}$ consist of polynomial components, each with maximum degree $k-1$, and $D_1(x)$ is some matrix resulting from the sequence of row operations as described above. By virtue of the type of row operation performed, it follows that the determinant of the matrix on the right-hand side equals $\det(M(D)-xI_n)$, and therefore if we develop this determinant inductively along the first column
\begin{eqnarray}
\nonumber \det(M(D)-xI_n) &=& (-1)^{k-1} \det \left [ {\renewcommand{\arraystretch}{1.2}\begin{array}{c|c} 
 -p_1(x) & \begin{matrix} d_1(x) & d_2(x) & \cdots & d_{n-k}(x) \end{matrix} \\ 
\hline
 e_1 & C(p_2(x))-xI_{n-k} 
 \end{array}} \right ] \\
 \nonumber &=& (-1)^{n} \det \left [ {\renewcommand{\arraystretch}{1.2}\begin{array}{c|c} 
 p_1(x) & \begin{matrix} -d_1(x) & -d_2(x) & \cdots & -d_{n-k}(x) \end{matrix} \\ 
\hline
 -e_1 & xI_{n-k}-C(p_2(x)) 
 \end{array}} \right ].
 \end{eqnarray}  
It follows that the characteristic polynomial of $M(D)$ is 
\[p_1(x)p_2(x) - \sum_{i=1}^{n-k} d_i(x) \left (x^{n-k-i}-\sum_{j=i}^{n-k-1}b_jx^{j-i}\right ).\]
Explicitly, the determinant of the characteristic matrix is developed inductively along the first column. The first step in this inductive development is 
\[p_1(x) \det(xI_{n-k}-C(p_2(x))) + 1\cdot \det \left ( \left [{\renewcommand{\arraystretch}{1.2}\begin{array}{c} 
 \begin{matrix} -d_1(x) & -d_2(x) & \cdots & -d_{n-k}(x) \end{matrix} \\ 
\hline
 C 
 \end{array}} \right ] \right ),\] where the matrix $C$ is obtained from $xI_{n-k}-C(p_2(x))$ by deleting the first row. The determinant on the right is then calculated by continuing the process of inductive development along the first column.

Now to require that the characteristic polynomial of $M(D)$ is $p(x)$, is to require
\[p(x) - p_1(x)p_2(x) = -\sum_{i=1}^{n-k} d_i(x) \left (x^{n-k-i}-\sum_{j=i}^{n-k-1}b_jx^{j-i}\right ).\] Therefore it remains to verify that polynomials $d_1(x), d_2(x), \ldots, d_{n-k}(x)$ which can satisfy this requirement always exist. Now since both $p(x)$ and $p_1(x),p_2(x)$ are monic polynomials, and furthermore $c_{n-1} = a_{k-1} + b_{n-k-1}$ we must have that $p(x)-p_1(x)p_2(x)$ is a polynomial of degree at most $n-2$. Let 
\[ p_{(2,i)}(x) = \left ( x^{n-k-i}-\sum_{j=i}^{n-k-1}b_jx^{j-i} \right ),\] for each $i \in \{1,2, \ldots, n-k\}$. Notice that $p_{(2,1)},p_{(2,2)},\ldots,p_{(2,(n-k))}$ is a basis for the vector space of all polynomials of maximum degree $n-k-1$. Furthermore the polynomials $d_1(x), d_2(x), \ldots, d_{n-k}(x)$ can be chosen to be any polynomials in $\mathscr{F}_{k-1}[x]$, the vector space of polynomials of maximum degree $k-1$. Since these polynomials can be chosen freely, it follows that any vector in $\mathscr{F}_{n-2}[x]$ can be constructed from $\sum_{i=1}^{n-k} d_i(x) p_{(2,i)}$ (and hence also from $-\sum_{i=1}^{n-k} d_i(x) p_{(2,i)}$). 

Explicitly, suppose $p(x) - p_1(x)p_2(x) = -\sum_{i=0}^{n-2} e_ix^i$, start by setting $d_1(x) = e_{n-2}x^{k-1}$. Then 
\[d_1(x)p_{(2,1)}(x) = e_{n-2}x^{n-2} - e_{n-2} b_{n-k-1} x^{n-3} - e_{n-2}x^{k-1}\sum_{j=1}^{n-k-2}b_jx^{j-1}.\] Next set $d_2(x) = (e_{n-3}+ e_{n-2}b_{n-k-1})x^{k-1}$,  then
\[d_2(x)p_{(2,2)}(x) = (e_{n-3}+e_{n-2}b_{n-k-1})x^{n-3} - (e_{n-3} b_{n-k-1}+e_{n-2}b_{n-k-1}^2) x^{n-4} \]
\[\qquad \qquad \qquad \qquad \qquad \qquad - (e_{n-3}+e_{n-2}b_{n-k-1})x^{k-1}\sum_{j=2}^{n-k-2}b_jx^{j-2}.\]
Now notice that 
\[d_1(x)p_{(2,1)}(x) + d_2(x)p_{(2,2)}(x) = e_{n-2}x^{n-2} + e_{n-3}x^{n-3} + \sum_{j=k-1}^{n-4} f_jx^j, \]
where $f_{k-1},\ldots,f_{n-4} \in \mathscr{F}$ are some coefficients resulting from the addition. Now inductively continuing this argument for $i=3, 4, \ldots, n-k-1$ we have 
\[ \sum_{i=1}^{n-k-1} d_i(x)p_{(2,i)}(x) = \left ( \sum_{i=k}^{n-2} e_ix^i \right ) + fx^{k-1},\] where $f \in \mathscr{F}$ is some coefficient resulting from the sum. Now $p_{(2,n-k)}=1$ and therefore if we let 
\[d_{n-k}(x) = (e_{k-1}-f)x^{k-1} + \sum_{i=0}^{k-2} e_i x^i \] we have 
\[-\sum_{i=1}^{n-k} d_i(x) p_{(2,i)} = -\sum_{i=0}^{n-2} e_ix^i = p(x) - p_1(x)p_2(x)\] as desired. \hfill $\square$

In what is to follow we will require the following proposition, which is a corollary to Roth's theorem \cite{roth}.
\begin{quote}
\begin{prop}
Let $\mathscr{F}$ be an arbitrary field, and let $A \in M_n(\mathscr{F}), B \in M_m(\mathscr{F})$. If the characteristic polynomials of $A$ and $B$ are relatively prime, 
\[ \begin{bmatrix} A & D \\ \mathbf{0} & B \end{bmatrix} \text{ is similar to } \begin{bmatrix} A & \mathbf{0} \\ \mathbf{0} & B \end{bmatrix} \] for any matrix $D \in M_{n \times m}(\mathscr{F})$.
\label{prop:coprime_chars_sim}
\end{prop}
\end{quote}
{\bf Proof.} Define $S:M_{n \times m}(\mathscr{F}) \rightarrow M_{n \times m}(\mathscr{F})$ as $S(X)=AX-XB$. It is easy to verify that $S$ is a linear transformation. Therefore if we can prove that $S(X)=\mathbf{0}$ implies $X=\mathbf{0}$, we will have proved that $S$ is an automorphism of $M_{n \times m}$. Consequently a unique solution exists for the equation $S(X)=D$ where $D$ is any matrix in $M_{n \times m}(\mathscr{F})$, and the result then follows directly by Roth's theorem, since then
\[\begin{bmatrix} I & X \\ \mathbf{0} & I \end{bmatrix}^{-1} 
\begin{bmatrix} A & \mathbf{0} \\ \mathbf{0} & B \end{bmatrix} 
\begin{bmatrix} I & X \\ \mathbf{0} & I \end{bmatrix} = 
\begin{bmatrix} I & -X \\ \mathbf{0} & I \end{bmatrix} 
\begin{bmatrix} A & \mathbf{0} \\ \mathbf{0} & B \end{bmatrix} 
\begin{bmatrix} I & X \\ \mathbf{0} & I \end{bmatrix} = \begin{bmatrix} A & D \\ \mathbf{0} & B \end{bmatrix}.\]

To this end, let $S(X)=\mathbf{0}$, then $AX=XB$. Let the characteristic polynomial of~$A$ be $p(x)$ and the characteristic polynomial of $B$ be $q(x)$. If these polynomials are relatively prime there exist polynomials $f(x),g(x)$ such that 
\[f(x)p(x) + g(x)q(x) = 1\] Now since $A$ satisfies its own characteristic equation $p(x)=0$ (by the Cayley-Hamilton theorem) we have 
\[f(A)p(A)+g(A)q(A)=g(A)q(A)=I.\] Now it follows that 
\[X=IX=g(A)q(A)X=g(A)Xq(B)=\mathbf{0},\] where the second last step follows from repeated application of $AX=XB$, and the last step is due to the fact that $B$ satisfies its own characteristic equation. \hfill $\square$
 
Now I present the definition of a \emph{well-partitioned} matrix, and some related results, all of which are due to de Seguins Pazzis \cite{pazzis2}. The definition and results are central to proving the main result of sums of three square-zero matrices over a field of characteristic two.

\begin{quote}
\begin{defn}[Definition 2.1 \cite{pazzis2}]
A square matrix $A$ is \emph{well-partitioned} if there are positive integers $s,t$ and monic polynomials \[p_1(x), p_2(x), \ldots, p_s(x), q_1(x), q_2(x), \ldots, q_t(x)\] such that:
\begin{enumerate}
\item $A=\text{Dg}[C(p_1(x)), \ldots, C(p_s(x)),C(q_1(x)), \ldots, C(q_t(x))]$,
\item the degree of $p_i(x)$ is greater than or equal to 2 for all $i \in \{2,3,\ldots,s\}$,
\item the degree of $q_j(x)$ is greater than or equal to 2 for all $j \in \{1,2,\ldots,t-1\}$,
\item for each $i \in \{1, 2, 3,\ldots,s\}$ the polynomial $p_i(x)$ is relatively prime to each polynomial in $\{q_1(x),q_2(x),\ldots,q_t(x)\}$. \hfill $\square$
\end{enumerate} 
\label{defn:well_partitioned}
\end{defn}
\end{quote}

\begin{quote}
\begin{prop}[Lemma 2.5 \cite{pazzis2}]
Let $\mathscr{F}$ be an arbitrary field, and let $A \in M_n(\mathscr{F})$ be a well-partitioned matrix. Choose any monic polynomial $p(x) = x^n - \sum_{i=0}^{n-1} a_{i}x^i$, subject to the constraint $a_{n-1} = \text{tr}(A)$. Then there exists a square-zero matrix $S$ such that $A$ is similar to $C(p(x))+S$. 
\label{prop:well_partitioned_is_similar_to_polyOfChoice}
\end{prop}
\end{quote}
{\bf Proof.} Let 
\[A= \text{Dg}[C(p_1(x)), C(p_2(x)), \ldots, C(p_s(x)),C(q_1(x)), C(q_2(x)), \ldots, C(q_t(x))],\] where $p_1(x), p_2(x), \ldots, p_s(x),$ $q_1(x), q_2(x), \ldots, q_t(x)$ are as described in Definition~\ref{defn:well_partitioned}, and let $m_1, m_2, \ldots, m_{s+t}$ be the respective degrees of each of these polynomials. Let 
$C_1 = C(p_1(x)), C_2 = C(p_2(x)), \ldots, C_{s+t} = C(q_t(x))$. Now we can write
\begin{eqnarray}
\nonumber A - S'  &=& \begin{bmatrix} C_1 & \mathbf{0} & \mathbf{0} & \cdots & \cdots & \mathbf{0} \\
\mathbf{0} & C_2 & \mathbf{0} & \cdots & \cdots & \mathbf{0} \\
\mathbf{0} & \mathbf{0} & C_3 & \ddots &  & \vdots \\
\vdots & & \ddots & \ddots & & \vdots \\
& & & \ddots & \ddots & \mathbf{0} \\
\mathbf{0} & \cdots & & & \mathbf{0} & C_{s+t}
\end{bmatrix} - 
\begin{bmatrix} \mathbf{0}  & \mathbf{0} & \mathbf{0} & \cdots & \cdots & \mathbf{0} \\
-N_1 & \mathbf{0}  & \mathbf{0} & \cdots & \cdots & \mathbf{0} \\
\mathbf{0} & -N_2 & \mathbf{0}  & \ddots &  & \vdots \\
\vdots & & \ddots & \ddots & & \vdots \\
& & & \ddots & \ddots & \mathbf{0} \\
\mathbf{0} & \cdots & & & -N_{s+t-1} & \mathbf{0} 
\end{bmatrix} \\
\nonumber &=& \begin{bmatrix} C_1 & \mathbf{0} & \mathbf{0} & \cdots & \cdots & \mathbf{0} \\
N_1 & C_2 & \mathbf{0} & \cdots & \cdots & \mathbf{0} \\
\mathbf{0} & N_2 & C_3 & \ddots &  & \vdots \\
\vdots & & \ddots & \ddots & & \vdots \\
& & & \ddots & \ddots & \mathbf{0} \\
\mathbf{0} & \cdots & & & N_{s+t-1} & C_{s+t}
\end{bmatrix}, \end{eqnarray}
 where $N_i$ is an $m_{i+1} \times m_{i}$ matrix with a one in entry $(1,m_i)$ and zeros elsewhere, for each $i \in \{1,2,\ldots,s+t-1\}$, and $S'$ is square-zero by virtue of $A$ being well-partitioned. 

If we designate $A_1 = \text{Dg}[C_1, C_2, \ldots, C_s]$ and $A_2= \text{Dg}[C_{s+1},C_{s+2},\ldots,C_{s+t}]$, and then define 
\[ A_1-S_1 = \begin{bmatrix} C_1 & \mathbf{0} & \mathbf{0} & \cdots & \cdots & \mathbf{0} \\
N_1 & C_2 & \mathbf{0} & \cdots & \cdots & \mathbf{0} \\
\mathbf{0} & N_2 & C_3 & \ddots &  & \vdots \\
\vdots & & \ddots & \ddots & & \vdots \\
& & & \ddots & \ddots & \mathbf{0} \\
\mathbf{0} & \cdots & & & N_{s-1} & C_{s}
\end{bmatrix},\] and
\[ A_2-S_2 = \begin{bmatrix} C_{s+1} & \mathbf{0} & \mathbf{0} & \cdots & \cdots & \mathbf{0} \\
N_{s+1} & C_{s+2} & \mathbf{0} & \cdots & \cdots & \mathbf{0} \\
\mathbf{0} & N_{s+2} & C_{s+3} & \ddots &  & \vdots \\
\vdots & & \ddots & \ddots & & \vdots \\
& & & \ddots & \ddots & \mathbf{0} \\
\mathbf{0} & \cdots & & & N_{s+t-1} & C_{s+t}
\end{bmatrix},\] then it is easy to see that $A_1-S_1,A_2-S_2$ are nonderogatory\footnote{The determinant of the submatrix obtained from $xI-(A_i-S_i)$ by deleting the top row and right-most column is 1.} and that 
\[A-S' = \begin{bmatrix} A_1-S_1 & \mathbf{0} \\
N'_{s} & A_2-S_2 \end{bmatrix}, \] where $N'_s$ is order $(\sum_{i=s+1}^{s+t} m_i) \times (\sum_{j=1}^{s} m_j)$ with a 1 in entry $(1,(\sum_{j=1}^{s} m_j))$ and zeros elsewhere.  

By Proposition \ref{prop:choice_of_polynomial} there exists a matrix $D$, such that 
\[ \begin{bmatrix} A_1-S_1 & D \\
N'_{s} & A_2-S_2 \end{bmatrix} \approx C(p(x)),\] since $\text{tr}(A_1-S_1)+\text{tr}(A_2-S_2) = \text{tr}(A) = c_{n-1}$. By virtue of $A$ being well-partitioned, we have that the characteristic polynomials of $A_1$ and $A_2$ are relatively prime, and therefore by Proposition \ref{prop:coprime_chars_sim} 
\[  A = \begin{bmatrix} A_1 & \mathbf{0} \\ \mathbf{0} & A_2 \end{bmatrix} \approx \begin{bmatrix} A_1 & D \\ \mathbf{0} & A_2 \end{bmatrix} = \begin{bmatrix} A_1-S_1 & D \\
N'_{s} & A_2-S_2 \end{bmatrix} + S' \approx C(p(x)) + S,\] where $S$ is square-zero by invariance of the square-zero property with respect to similarity. \hfill $\square$

\begin{quote}
\begin{prop}[Lemma 4.2 \cite{pazzis2}]
Let $\mathscr{F}$ be an arbitrary field, and let $A \in M_n(\mathscr{F})$. If the rational canonical form of $A$ contains at most one block $H(x-\lambda)$ associated with each of its characteristic values $\lambda$, and its minimum polynomial is not a power of some irreducible polynomial, then $A$ is similar to a well-partitioned matrix. 
\label{prop:one_block_per_eigen_similar_to_wp}
\end{prop}
\end{quote}
{\bf Proof.} If the rational canonical form of $A$ has at most one block of order $1 \times 1$ associated with each characteristic value, then for each elementary divisor of $A$ which is a power of an irreducible polynomial of degree one (that is, of the form $(x-\lambda)^k$), there is at most one invariant polynomial where $k=1$. Requiring that the minimum polynomial of $A$ is not a power of some irreducible polynomial, is equivalent to requiring that the minimum polynomial of $A$ must contain at least two elementary divisors. 

Let us fix one of these elementary divisors $(p(x))^{e_t}$ where $p(x)$ is an irreducible monic polynomial. By Theorem 7.5  in the text by Cullen \cite{cullen} $A$ is similar to 
\[\text{Dg}[C((p(x))^{e_1}),C((p(x))^{e_2}),\ldots,C((p(x))^{e_t}),C(d_2(x)),C(d_3(x)),\ldots,C(d_s(x))],\]
that is, $A$ is similar to a diagonal block matrix, where the diagonal blocks are companion matrices of its elementary divisors. Set \[A_1 = \text{Dg}[C((p(x))^{e_1}),C((p(x))^{e_2}),\ldots,C((p(x))^{e_t})]\] where we may assume (by a permutation of the diagonal blocks of $A_1$ if necessary) that $1 \leq e_1 \leq e_2 \leq \cdots \leq e_t$ and if $e_1 = 1$ we must have $e_1 < e_2$ (by the first requirement). Also, set \[A_2=[C(d_2(x)),C(d_3(x)),\ldots,C(d_s(x))],\] where we can assume (by a permutation of diagonal blocks of $A_2$ if needed) that all blocks of order $1 \times 1$ can be placed next to each other on the diagonal, if we relabel the blocks of this matrix so that $d_1(x),d_2(x),\ldots,d_{s_1}(x)$ are elementary divisors of degree 2 or higher in $A_2$, and $x-\lambda_1,x-\lambda_2,\ldots,x-\lambda_{s_2}$ are elementary divisors of degree 1 in $A_2$, then $A_2$ is similar to
\[\text{Dg}[C(d_1(x)),C(d_2(x)),\ldots,C(d_{s_1}(x)),C(x-\lambda_1),C(x-\lambda_2),\ldots,C(x-\lambda_{s_2})].\] By the conditions placed on $A$, we must have $x-\lambda_1,x-\lambda_2,\ldots,x-\lambda_{s_2}$ are distinct polynomials, and therefore $\text{Dg}[C(x-\lambda_1),C(x-\lambda_2),\ldots,C(x-\lambda_{s_2})]$ is a nonderogatory matrix which is similar to $C((x-\lambda_1)(x-\lambda_2)\ldots(x-\lambda_{s_2}))$. So we can write
\[A_2 \approx A'_2 = \text{Dg}[C(d_1(x)),C(d_2(x)),\ldots,C(d_{s_1}(x)),C((x-\lambda_1)(x-\lambda_2)\ldots(x-\lambda_{s_2}))].\] Now $A$ is similar to $\text{Dg}[A_1,A'_2]$, which is a well-partitioned matrix. \hfill $\square$    

\subsection{Fields of characteristic two}
\begin{quote}
\begin{lem}[Lemma 4.1 \cite{pazzis2}]
Let $\mathscr{F}$ be a field of characteristic two. For any $\lambda \in \mathscr{F}$, the matrix $\lambda I_2$ is the sum of three square-zero matrices. 
\label{lem:lambdaI2_is_sum_of_three}
\end{lem}
\end{quote}
{\bf Proof.} Notice that
\[\lambda I_2 = J_2(\lambda) - J_2(0) \approx C((x-\lambda)^2) - J_2(0). \]
Now $J_2(0)$ is square-zero, and within the context of a field of characteristic two, $(x-\lambda)^2 = x^2 + \lambda^2$, which is an even-power polynomial. By Theorem \ref{thm:sufficient_sqzero_sum} we therefore have $C((x-\lambda)^2)$ is a sum of two square-zero matrices. By invariance of the square-zero property with respect to similarity the result follows. \hfill $\square$ 

\begin{quote}
\begin{lem}
Let $\mathscr{F}$ be a field of characteristic two. Let $f(x) = a_mx^m + \cdots + a_0 \in \mathscr{F}[x]$ be an arbitrary polynomial. Then $f(x)^n$ is an even-power polynomial for any even integer $n$.
\label{lem:char_2_n_even_poly}
\end{lem}
\end{quote}
{\bf Proof.} We may assume without loss of generality that $m \geq 1$ and $n \geq 2$. Since the product of two even-power polynomials is again an even-power polynomial, the result will follow from $f(x)^n = (f(x)^2)^{n/2}$, by showing that $f(x)^2$ is an even-power polynomial. Let us denote $f(x)^2 = b_{2m}x^{2m} + \cdots + b_0$. Now calculating the coefficient~$b_k$ by expanding $(a_mx^m + \cdots + a_0)^2$, we find
\[b_k = a_ka_0 + a_{k-1}a_1 + \cdots + a_{1}a_{k-1} + a_0a_k,\] and if $k$ is odd we must have an even number of terms on the right-hand side. It follows that 
\[b_k = 2(a_ka_0 + a_{k-1}a_1 + \cdots + a_{(k+1)/2}a_{(k-1)/2}) = 0.\]
It follows that $f(x)^2$ is an even-power polynomial. \hfill $\square$

\begin{quote}
\begin{lem}[Lemma 4.3 \cite{pazzis2}]
Let $\mathscr{F}$ be a field of characteristic two. Let the minimum polynomial of $A \in M_n(\mathscr{F}$) be a power of some irreducible polynomial. If the trace of $A$ is zero then $A$ is a sum of three square-zero matrices.
\label{lem:minpoly_power_irred_is_sum_of_three}
\end{lem}
\end{quote}
{\bf Proof.} Let the minimum polynomial of $A$ be $p(x)^t$ where $p(x) = x^k - \sum_{i=0}^{k-1} c_ix^i$ is an irreducible polynomial. By Theorem 7.2 in the text by Cullen \cite{cullen} $A$ is similar to \[\text{Dg}[C(p(x)^{t_1}),C(p(x)^{t_2}),\ldots,C(p(x)^{t_m}),C(p(x)^{t})]\] where $1 \leq t_1 \leq t_2 \leq \cdots \leq t_m \leq t$.

If $c_{k-1}=0$, then for any positive integer $s$, the coefficient of $x^{k\cdot s-1}$ in the expansion of $p(x)^s$ is $s \cdot c_{k-1} = 0$. The reader may now refer to Corollary \ref{cor:tr0_nond_3sqz} where it is proved that a diagonal block matrix over an arbitrary field, where each diagonal block is a trace-zero companion matrix, is a sum of three square-zero matrices. 

It remains to prove that the result holds if $c_{k-1} \neq 0$. In the multinomial expansion of $p(x)^s$, we must have that the coefficient of $x^{k\cdot s-1}$ is $s \cdot c_{k-1}$ for any positive integer $s$. Now since $\mathscr{F}$ is characteristic two, $s \cdot c_{k-1}=0$ if $s$ is even, and furthermore since the trace of $A$ is zero it then follows that $A$ must have an even amount of odd-degree invariant polynomials. Applying Theorem 7.2 in the text by Cullen \cite{cullen}, followed by a simple permutation of diagonal blocks, we find that $A$ is similar to 
\[\text{Dg}[C(p(x)^{d_1}), C(p(x)^{d_2}), \ldots, C(p(x)^{d_{s_1}}),C(p(x)^{e_1}),C(p(x)^{e_2}), \ldots, C(p(x)^{e_{s_2}})],\] where $d_i$ is even for each $i \in \{1,2,\ldots,s_1\}$ and $e_j$ is odd for each $j \in \{1,2,\ldots,s_2\}$, and by the preceding arguments we must have $s_2$ is even. Now by Lemma \ref{lem:char_2_n_even_poly} it follows that $p(x)^{d_i}$ is even-power, and therefore each block $C(p(x)^{d_i})$ is a sum of two square-zero matrices. Now for $j \in \{1,3, \ldots, s_2-1\}$ we have 
\[ \text{Dg}[C(p(x)^{e_j}),C(p(x)^{e_{j+1}})] = \begin{bmatrix} C(p(x)^{e_j}) & \mathbf{0} \\
N & C(p(x)^{e_{j+1}}) \end{bmatrix} + \begin{bmatrix} \mathbf{0} & \mathbf{0} \\
-N & \mathbf{0}  \end{bmatrix} = C + S,\] where $N$ is a matrix with a one in entry $(1,ke_j)$ and zeros elsewhere. It is easy to verify that $S$ is square-zero. 

It is also easy to verify that $C$ is nonderogatory, and its characteristic polynomial is the product of the characteristic polynomials of its diagonal blocks (this follows from \cite[Proposition 11.12]{golan}), that is, its characteristic polynomial is $p(x)^{e_j+e_{j+1}}$. Now since $e_j$ and $e_{j+1}$ are both odd, the sum $e_j+e_{j+1}$ is even. By Lemma \ref{lem:char_2_n_even_poly} it follows that $p(x)^{e_j+e_{j+1}}$ is an even-power polynomial.



So by Theorem \ref{thm:sufficient_sqzero_sum} it follows that $C$ is a sum of two square-zero matrices. It follows that $\text{Dg}[C(p(x)^{e_j}),C(p(x)^{e_{j+1}})]$ is a sum of three square-zero matrices for each $j \in \{1,3, \ldots, s_2-1\}$. 

To summarize, and conclude on the results above: $A$ is similar to a diagonal block matrix, where each diagonal block is a sum of three square-zero matrices, and therefore $A$ is similar to a sum of three square-zero matrices. By invariance of the square-zero property with respect to similarity $A$ is a sum of three square-zero matrices. \hfill $\square$ 

The main result of this section can now be proved.
\begin{quote}
\begin{thm}[Theorem 1.3 \cite{pazzis2}]
Let $\mathscr{F}$ be a field of characteristic two. Then $A \in M_n(\mathscr{F}$) is a sum of three square-zero matrices if and only if the trace of $A$ is zero.
\label{thm:fieldchar2_sum_of_three_main_result}
\end{thm}
\end{quote}
{\bf Proof.} The trace of any square-zero matrix is zero, since it is similar to a matrix in Jordan canonical form which has only zeros on the diagonal, and trace is invariant with respect to similarity. Therefore the sum of three square-zero matrices has trace zero.

For the converse, suppose $A$ is a matrix with trace zero. Observe that the rational canonical form of $A$ is similar to 
\[\text{Dg}[A_1, \lambda_1I_2, \lambda_2 I_2, \ldots, \lambda_t I_2],\] where $A_1$ is in rational canonical form and contains at most one block $H(x-\lambda)$ per distinct characteristic value $\lambda$ of $A$. Note that the characteristic values $\lambda_i$ are not necessarily distinct, that is, it is possible that $\lambda_i = \lambda_j$ for $i \neq j$. To describe explicitly, if in the rational canonical form of $A$ we have a repeated block $H(x-\lambda)$ for any characteristic value $\lambda$, such blocks are paired to form the block $\lambda I_2$ and moved to the last position on the diagonal. In this way it is ensured that $A_1$ has at most one block $H(x-\lambda)$ per distinct $\lambda$.

Now since $\mathscr{F}$ is of characteristic two and the trace of $A$ is zero it follows that the trace of $A_1$ is also zero. Furthermore, by Lemma \ref{lem:lambdaI2_is_sum_of_three}, for each $i \in \{1,2, \ldots,t\}$ the block $\lambda_i I_2$ is a sum of three square-zero matrices. The proof will therefore be complete if we can prove that $A_1$ is a sum of three square-zero matrices.
 
Now if the minimum polynomial of $A_1$ is the power of some irreducible polynomial the result follows from Lemma \ref{lem:minpoly_power_irred_is_sum_of_three}.  Suppose this is not the case. Then by Proposition \ref{prop:one_block_per_eigen_similar_to_wp} $A_1$ is similar to a well-partitioned matrix. Suppose $A_1$ is order $k \times k$: choose $p(x)=x^k$. Then by Proposition \ref{prop:well_partitioned_is_similar_to_polyOfChoice} $A_1$ is similar to $C(x^k)+S$ where $S$ is a square-zero matrix. By Theorem \ref{thm:sufficient_sqzero_sum} $C(x^k)$ is a sum of two square-zero matrices, and therefore by invariance of the square-zero property with respect to similarity it follows that $A_1$ is a sum of three square-zero matrices. \hfill $\square$

\subsection{Fields which are not of characteristic two}
\subsubsection{Necessary Conditions}
The following necessary condition for matrices over the complex field is due to Wang and Wu, and remains valid over a field which is not of characteristic two.
\begin{quote}
\begin{prop}[Theorem 3.1 \cite{ww}]
Let $\mathscr{F}$ be a field which is not of characteristic two. If $A \in M_n(\mathscr{F})$ is the sum of three square-zero matrices, then $\text{n}(A - \lambda I) \leq 3n/4$ for any nonzero $\lambda \in \mathscr{F}$. 
\label{prop:3sum_nec_ww}
\end{prop}
\end{quote}
{\bf Proof.} By Corollary \ref{cor:3sqz_is_sqzPlusC} $A=C+S$ where $C$ is a matrix with only odd or even invariant polynomials, and $S$ is square-zero. Let $m = \text{n}(A-\lambda I)$. By Grassmann's Theorem \cite[Proposition 5.16]{golan} (rearranging terms so that the dimension of the sum of the subspaces and the dimension of the intersection of subspaces are exchanged):
\[  \dim(\text{N}(A-\lambda I) \cap \text{N}(S)) = \text{n}(A-\lambda I) + \text{n}(S) - \dim(\text{N}(A-\lambda I) + \text{N}(S)), \] and since $\text{n}(S) \geq n/2$ it follows that
\[ \dim(\text{N}(A-\lambda I) \cap \text{N}(S)) \geq m + \frac{n}{2} - n = m - \frac{n}{2}.\]

Let $u \in \text{N}(A-\lambda I) \cap \text{N}(S)$. \\Then we have $(A-\lambda I)u= Au - \lambda u=\mathbf{0}$ and $Su=\mathbf{0}$ so that
\[ \lambda u = Au = (C+S)u = Cu,\] and therefore $(C - \lambda I)u = \mathbf{0}$. We therefore have 
\[ \text{N}(A-\lambda I) \cap \text{N}(S) \subseteq \text{N}(C - \lambda I).\] 
From these arguments it follows that $\text{n}(C - \lambda I) \geq m - \frac{n}{2}$. Now by Corollary \ref{cor:A_sim_minusA} $C$ is similar to $-C$, so that if $u \in \text{N}(C - \lambda I)$ and $P$ is an invertible matrix such that $P^{-1}(-C)P = C$ then 
\[ (C - \lambda I)u = (P^{-1}(-C)P - \lambda P^{-1}P)u = -P^{-1}(C + \lambda I)Pu = \mathbf{0}, \] that is $Pu \in \text{N}(C + \lambda I)$. Since $P$ is injective, it follows that $\text{n}(C - \lambda I) \leq \text{n}(C + \lambda I)$. By a similar argument it can be shown that for $v \in \text{N}(C + \lambda I)$ we have $P^{-1}v \in \text{N}(C - \lambda I)$ and therefore \[\text{n}(C + \lambda I) = \text{n}(C - \lambda I) \geq m - \frac{n}{2}.\]
Now again by Grassmann's Theorem \cite[Proposition 5.16]{golan}
\[  \dim(\text{N}(A-\lambda I) \cap \text{N}(C + \lambda I)) = \text{n}(A-\lambda I) + \text{n}(C + \lambda I) - \dim(\text{N}(A-\lambda I) + \text{N}(C + \lambda I)), \] and by the preceding arguments it follows that 
\[ \dim(\text{N}(A-\lambda I) \cap \text{N}(C + \lambda I)) \geq m + \left ( m - \frac{n}{2} \right ) - n = 2m - \frac{3}{2}n.\]

Let $v \in \text{N}(A-\lambda I) \cap \text{N}(C + \lambda I)$, then
\[ (A-\lambda I)v = (C + S - \lambda I)v = \mathbf{0}, \] and furthermore $Cv = -\lambda v$ so that substituting in the equation above we have
\[ (S - 2 \lambda I)v = \mathbf{0}. \] So $v \in \text{N}(S - 2 \lambda I)$, and it follows that $\text{N}(A-\lambda I) \cap \text{N}(C + \lambda I) \subseteq \text{N}(S - 2 \lambda I)$. Therefore $\text{n}(S - 2 \lambda I) \geq 2m - 3n/2$.

Now $S$ is square-zero and therefore its characteristic polynomial splits over an arbitrary field, and its only eigenvalue is zero, which implies that $S - 2 \lambda I$ is invertible for any $\lambda \neq 0$.\footnote{Note that it is at this point where the result becomes specific to fields which are not of characteristic two.} Therefore $\text{n}(S - 2 \lambda I) = 0$, and it follows that 
\[ m \leq \frac{3}{4}n, \] which is the desired result. \hfill $\square$   

\subsubsection{Sufficient Conditions}
A necessary condition for a sum of three square-zero matrices is that its trace must be zero. That this condition is not sufficient for an arbitrary field follows easily from Proposition \ref{prop:3sum_nec_ww}: as an example it may be verified that $A = \text{Dg}[-4,1,1,1,1] \in M_5(\mathbb{C})$ has trace zero, but 
\[\text{n}(A-I) = 4 > \frac{15}{4}=\frac{3\cdot 5}{4}.\]

Wang and Wu \cite{ww} proposed that Proposition \ref{prop:3sum_nec_ww} might also be a sufficient condition (considered for matrices over the complex field), but Takahashi \cite{takahashi} proved that this is not the case. A part of this section will be devoted to a presentation of these results, but first I will present some sufficient results which are true over an arbitrary field.

\begin{quote}
\begin{prop}
Let $\mathscr{F}$ be an arbitrary field. If the trace of the matrix $A \in M_n(\mathscr{F})$ is zero and $A$ is nonderogatory then $A$ is a sum of three square-zero matrices. 
\label{prop:tr0_nond_3sqz}
\end{prop}
\end{quote}
{\bf Proof.} The matrix $A$ is similar to a companion matrix $C(x^n - c_{n-2}x^{n-2} - c_{n-3}x^{n-3} - \cdots - c_0)$, where $x^n - c_{n-2}x^{n-2} - c_{n-3}x^{n-3} - \cdots - c_0$ is the characteristic polynomial of $A$ (notice that the term $c_{n-1}x^{n-1}$ is absent by virtue of $A$ having trace zero). Now $C(x^n - c_{n-2}x^{n-2} - c_{n-3}x^{n-3} - \cdots - c_0)$ is equal to
\[C(x^n - \sum_{k=1}^{\lfloor n/2 \rfloor} c_{n-2k}x^{n-2k}) + \begin{bmatrix} 
\vdots & \vdots &  & \vdots & \vdots\\
0 & 0 & \cdots & 0 & c_{n-5}\\
0 & 0 & \cdots & 0 & 0\\
0 & 0 & \cdots & 0 & c_{n-3}\\
0 & 0 & \cdots & 0 & 0\\
0 & 0 & \cdots & 0 & 0 
\end{bmatrix},\] which is the sum of a matrix with one nonconstant even- or odd-power invariant polynomial and a square-zero matrix. By Corollary \ref{cor:3sqz_is_sqzPlusC} and the invariance of the square-zero property with respect to similarity $A$ is a sum of three square-zero matrices. \hfill $\square$

\begin{quote}
\begin{cor}
Let $\mathscr{F}$ be an arbitrary field. Let $A \in M_n(\mathscr{F})$ be similar to a diagonal matrix with distinct (non-repeating) values on the diagonal. If the trace of $A$ is zero then $A$ is a sum of three square-zero matrices. 
\label{cor:tr0_diag_distinct_3sqz}
\end{cor}
\end{quote}
{\bf Proof.} The characteristic polynomial of $A$ is 
\[ \prod_{i=1}^n (x-\lambda_i)\] where $\lambda_i \neq \lambda_j$ if $i \neq j$, and since each distinct linear factor in the characteristic polynomial of a matrix must also be in its minimum polynomial \cite[Corollary~3]{cullen}, the matrix $A$ must be nonderogatory. Therefore the result follows directly from Proposition 
\ref{prop:tr0_nond_3sqz}. \hfill $\square$

\begin{quote}
\begin{cor}
Let $\mathscr{F}$ be an arbitrary field. Suppose the matrix $A \in M_n(\mathscr{F})$ has $t$ nonconstant invariant polynomials. Denote each polynomial as $f_i(x) = x^{k_i} - \sum_{j=1}^{k_i} c_{k_i-j}x^{k_i-j}$ for $i \in \{1,2,\ldots,t\}$. If $c_{k_i-1} = 0$ for each $i \in \{1,2,\ldots,t\}$ then $A$ is a sum of three square-zero matrices.
\label{cor:tr0_nond_3sqz}
\end{cor}
\end{quote}
{\bf Proof.} By Theorem 7.2 in the text by Cullen \cite{cullen} $A$ is similar to 
\[C = \text{Dg}[C(f_1(x)),C(f_2(x)), \ldots, C(f_t(x))],\] and by Proposition \ref{prop:tr0_nond_3sqz} each of the diagonal blocks is a sum of three square-zero matrices. \hfill $\square$
 
Suppose we have an invariant polynomial $f(x)=x^n - \sum_{j=1}^{n} c_{n-j}x^{n-j}$. Let us designate such an invariant polynomial as a trace zero invariant polynomial whenever $c_{n-1}=0$. This would agree with the fact that $C(f(x))$ has trace zero. The preceding Corollary shows that a matrix which has only trace zero invariant polynomials is a sum of three square-zero matrices. 

The following result is a corollary to Propositions \ref{prop:well_partitioned_is_similar_to_polyOfChoice} and \ref{prop:one_block_per_eigen_similar_to_wp}; this result was also employed as part of the proof of Theorem \ref{thm:fieldchar2_sum_of_three_main_result}, which is the main result for sums of three square-zero matrices over a field of characteristic two.

\begin{quote}
\begin{cor}
Let $\mathscr{F}$ be an arbitrary field, and let $A \in M_n(\mathscr{F})$. If the trace of $A$ is zero, and the rational canonical form of $A$ contains at most one block $H(x-\lambda)$ associated with each of its characteristic values $\lambda$, and its minimum polynomial is not a power of some irreducible polynomial, then $A$ is a sum of three square-zero matrices.
\label{cor:atmost_onedegreeone_blockpereigen_3sqz}
\end{cor}
\end{quote}
{\bf Proof.} By Proposition \ref{prop:one_block_per_eigen_similar_to_wp} $A$ is similar to a well-partitioned matrix, and by Proposition \ref{prop:well_partitioned_is_similar_to_polyOfChoice} $A$ is then similar to $C(p(x)) + S$ where $S$ is square-zero and $p(x)$ is an arbitrary monic polynomial of degree $n$, subject to $\text{tr}(C(p(x)))=\text{tr}(A)=0$. Since the trace of $p(x)$ is zero, we can choose $p(x) = x^n$, so that $C(p(x))$ is a sum of two square-zero matrices. Therefore $A$ is a sum of three square-zero matrices. \hfill $\square$

Let us consider the following results due to Takahashi: necessary and sufficient conditions for a matrix with minimum polynomial $(x-a)(x-b)$ where $a \neq b$ to be a sum of three square-zero matrices. The results are only valid over a field of characteristic zero. As a consequence of these results we will be able to prove that Proposition \ref{prop:3sum_nec_ww} is not a sufficient condition for a matrix over a field not of characteristic two to be a sum of three square-zero matrices. The following four lemmas are required to prove the main result.

\begin{quote}
\begin{lem}[Lemma 3 \cite{takahashi}]
Let $\mathscr{F}$ be an arbitrary field. Let $A \in M_n(\mathscr{F})$ have minimum polynomial $(x-a)(x-b)$ and let $S \in M_n(\mathscr{F})$ be square-zero. If $\lambda$ is a characteristic value of $A+S$, and $\lambda \neq a$, $\lambda \neq b$, then $a + b - \lambda$ is also a characteristic value of $A+S$. 
\label{lem:eigens_of_A_plus_S}
\end{lem}
\end{quote}
{\bf Proof.}\footnote{The proof as given here is due to Botha.} The result is proved in three parts. First it is proved that the result holds for $a = 1$ and $b = 0$, that is, for $A$ idempotent.

Without loss of generality we may assume that in this case \[A = \text{Dg}[I_r, \mathbf{0}_{n-r}] \text{ and } S = \begin{bmatrix} S_1 & S_2 \\ S_3 & S_4 \end{bmatrix} \] where $S_1$ is order $r \times r$ (so $S$ and $A$ are compatibly partitioned). Since $S$ is square-zero, it follows that 
\begin{eqnarray}
\nonumber S_1^2 + S_2S_3 &=& \mathbf{0} \\
\nonumber S_1S_2 + S_2S_4 &=& \mathbf{0} \\
\nonumber S_3S_1 + S_4S_3 &=& \mathbf{0} \\
\nonumber S_3S_2 + S_4^2 &=& \mathbf{0} \\
\label{eqset:sqzeroS}
\end{eqnarray}
Now suppose that $(A + S) v = \lambda v$ where $\lambda \neq 1, \lambda \neq 0$, $v \neq \mathbf{0}$ and $v = (v_1, v_2)^T$ is a compatibly partitioned column vector. Then
\[\begin{bmatrix} I_r + S_1 & S_2 \\ S_3 & S_4 \end{bmatrix} \begin{bmatrix} v_1 \\ v_2 \end{bmatrix} = \lambda \begin{bmatrix} v_1 \\ v_2 \end{bmatrix}\] and it follows that
\begin{equation}
S_1v_1 + S_2v_2 = (\lambda - 1)v_1 \label{eqset:AplusSeigenTop}
\end{equation}
\begin{equation}
S_3v_1 + S_4v_2 = \lambda v_2.  \label{eqset:AplusSeigenBottom}
\end{equation}
Now multiplying \eqref{eqset:AplusSeigenTop} by $S_1$ on the left, and \eqref{eqset:AplusSeigenBottom} by $S_2$ on the left, adding the resulting equations together, and then applying \eqref{eqset:sqzeroS} we have
\begin{equation}
(\lambda - 1)S_1v_1 + \lambda S_2v_2 = \mathbf{0}. \label{eq:S1S2v1v2}
\end{equation}
Multiplying \eqref{eqset:AplusSeigenTop} by $S_3$ on the left, and \eqref{eqset:AplusSeigenBottom} by $S_4$ on the left, adding the resulting equations together, and then applying \eqref{eqset:sqzeroS} we have
\begin{equation}
(\lambda - 1)S_3v_1 + \lambda S_4v_2 = \mathbf{0}. \label{eq:S3S4v1v2}
\end{equation}
Now by \eqref{eqset:AplusSeigenTop} $S_1v_1 = (\lambda - 1)v_1 - S_2v_2$ and substituting into \eqref{eq:S1S2v1v2} and rearranging terms as necessary we have 
\[ S_2v_2 = -(\lambda - 1)^2 v_1.\] Substituting this result back into \eqref{eqset:AplusSeigenTop} and simplifying as necessary we have 
\[ S_1v_1 = \lambda(\lambda - 1) v_1.\] Similarly it follows from \eqref{eqset:AplusSeigenBottom} and \eqref{eq:S3S4v1v2} that
\[S_3v_1 = \lambda^2 v_2 \text{ and } S_4v_2 = -\lambda(\lambda-1)v_2.\]
Now the proof will be complete if we can show there exists a nonzero vector $w = (w_1, w_2)^T$ such that
\[(A+S)w = \begin{bmatrix} I_r + S_1 & S_2 \\ S_3 & S_4 \end{bmatrix} \begin{bmatrix} w_1 \\ w_2 \end{bmatrix} = \begin{bmatrix} w_1 + S_1w_1 + S_2w_2 \\ S_3w_1 + S_4w_2 \end{bmatrix}  = (1-\lambda) \begin{bmatrix} w_1 \\ w_2 \end{bmatrix},\]
that is, we must have 
\[S_1w_1 + S_2w_2 = -\lambda w_1 \text{ and } S_3w_1 + S_4w_2 = (1-\lambda) w_2.\] But now, by the preceding results, it can be verified that setting \[w_1 = \frac{1}{\lambda^2}v_1 \text{ and } w_2 = \frac{1}{(\lambda-1)^2}v_2,\] results in a vector $w$ with the required properties. The result therefore holds for $a = 1, b = 0$.

Now consider the case where $a$ and $b$ are arbitrary subject to the constraint $a \neq b$. Since the square-zero property is invariant with respect to similarity, we can assume without loss of generality that $A = \text{Dg}[a I_r, b I_{n-r}]$. In this case we must also have $b - a \neq 0$, and therefore 
\[\frac{1}{b - a}(A - a I) = \frac{1}{b - a}\text{Dg}[\mathbf{0}_r, (b-a) I_{n-r}] = \text{Dg}[\mathbf{0}_r, I_{n-r}],\] which is an idempotent matrix. 

Let $\lambda \neq a, \lambda \neq b$ be a characteristic value of $A+S$ where $S$ is an arbitrary square-zero matrix. It follows that $(\lambda - a)/(b - a)$ is a characteristic value of $\frac{1}{b - a}(A - a I + S).$ It is easy to check that the square-zero property is preserved under scalar multiplication, and therefore $\frac{1}{b - a}(A - a I + S)$ is the sum of an idempotent matrix and a square-zero matrix, and therefore we can apply the result from the first part. Explicitly, since $(\lambda - a)/(b - a)$ is not equal to zero or one, we must have $1-(\lambda - a)/(b - a)$ is a characteristic value of $\frac{1}{b - a}(A - a I + S)$, and therefore
\[(b - a)\left (1-\frac{\lambda - a}{b - a} \right ) + a = a + b - \lambda\] is a characteristic value of $A+S$.

It remains to prove that the result holds for $a = b$. By invariance of the square-zero property with respect to similarity we may assume that $A = \text{Dg}[a I_r, A_1]$ where $A_1$ is a matrix consisting exclusively of blocks $J_2(a)$. It is easy to verify that $A - a I$ is then a square-zero matrix, and hence $A - a I + S$ is a sum of two square-zero matrices. 

Now suppose $(A+S)v = \lambda v$ for some nonzero column vector $v$, and $\lambda \neq a$. Then $(A - a I + S)v = (\lambda - a)v$. If $\mathscr{F}$ is not characteristic two it follows from Corollary~\ref{cor:A_sim_minusA} that $A - a I + S$ is similar to $-(A - a I + S)$, and therefore we must have $-(\lambda - a)$ is also a characteristic value of $(A - a I + S)$. Suppose $\mathscr{F}$ is characteristic two: then $-(\lambda - a) = \lambda - a$, and again it follows that $-(\lambda - a)$ is a characteristic value of $(A - a I + S)$.

But then $-(\lambda - a) + a = 2 a - \lambda$ is a characteristic value of $A+S$. \hfill $\square$

\begin{quote}
\begin{lem}[Lemma 3.2 \cite{ww}]
Let $\mathscr{F}$ be an arbitrary field. Let $C(p(x)) \in M_n(\mathscr{F})$ be the companion matrix of the polynomial $p(x)=x^n - c_{n-1}x^{n-1} - \cdots - c_0$, where $\{c_0, c_1, \ldots, c_{n-1}\} \subseteq \mathscr{F}$. Then 
\[C(p(x)) = S + C(x^n-c_{n-1}x^{n-1} - b_{n-2}x^{n-2} - \cdots - b_0)\] where $S$ is square-zero, and $b_0, b_1, \ldots, b_{n-2}$ are arbitrary elements in $\mathscr{F}$. As a consequence, it is possible to choose $n$ arbitrary values $\lambda_1, \lambda_2, \ldots, \lambda_n \in \mathscr{F}$ subject to $\lambda_1 + \lambda_2 + \ldots + \lambda_n = c_{n-1}$ so that 
\[C(p(x)) = S + C((x-\lambda_1)(x-\lambda_2)\cdots(x-\lambda_n)).\]
\label{lem:construct_3_mechanism}
\end{lem}
\end{quote}
{\bf Proof.} We have 
\begin{eqnarray}
\nonumber C(p(x)) = \left [ {\renewcommand{\arraystretch}{1.2}\begin{array}{c|c} 
\mathbf{0} &  c_0 \\
\hline
I_{n-1} & \begin{matrix} c_1 \\ \vdots \\ c_{n-2} \\ c_{n-1} \end{matrix}
 \end{array}} \right ] &=& \left [ {\renewcommand{\arraystretch}{1.2}\begin{array}{c|c} 
\mathbf{0} &  c_0-b_0 \\
\hline
\mathbf{0} & \begin{matrix} c_1-b_1 \\ \vdots \\ c_{n-2}-b_{n-2} \\ 0 \end{matrix}
 \end{array}} \right ] + \left [ {\renewcommand{\arraystretch}{1.2}\begin{array}{c|c} 
\mathbf{0} &  b_0 \\
\hline
I_{n-1} & \begin{matrix} b_1 \\ \vdots \\ b_{n-2} \\ c_{n-1} \end{matrix}
 \end{array}} \right ]\\
\nonumber &=& S + C(x^n-c_{n-1}x^{n-1} - b_{n-2}x^{n-2} - \cdots - b_0),
\end{eqnarray} and it is easy to verify that $S$ is square-zero. Now the last part follows from the fact that we can choose $b_0, b_1, \ldots, b_{n-2}$ to be the coefficients of $x^0, x , x^2, \ldots, x^{n-2}$ in the expansion of $(x-\lambda_1)(x-\lambda_2)\cdots(x-\lambda_n)$, which completes the proof. \hfill $\square$ 

\begin{quote}
\begin{lem}[Lemma 1 \cite{takahashi}]
Let $\mathscr{F}$ be a field with at least four elements, and let $a,b,c,d \in \mathscr{F}$ be such that $a \neq b$, $c \neq d$ and $a+b = c+d$. Then there exists a square-zero matrix $S$ such that
\[\text{Dg}[a,b,a,b] \approx S + \text{Dg} \left [ J_2(c), J_2(d)\right ]\]
\label{lem:lem1takahashi}
\end{lem}
\end{quote}
{\bf Proof.} It is easy to see that the result is true in the trivial case where $a=c$ or $a=d$: we can apply a similarity transformation with a suitable permutation matrix to 
\[\text{Dg} \left [ J_2(c), J_2(d) \right ] - \text{Dg}[J_2(0),J_2(0)]\] to see that the result holds.

Assume therefore that $a \neq c,d$ (then we must have $b \neq c,d$). Let 
\[P = \begin{bmatrix} 0 & \frac{b - d}{(a - c) (c - b)} & 0 &  
  \frac{1}{a - c}\\
  \frac{b - d}{(a - c) (c - b)} & 0 &  \frac{1}{a - c} & 0 \\ 
  0 & 1 & 0 & 1\\
  1 &0 & 1 & 0 \end{bmatrix}, \] then 
  \[\text{Dg}[a,b,a,b] \approx \text{Dg}[b,b,a,a] = P^{-1}\begin{bmatrix} cI_2 & I_2\\
  (a - c) (c-b)I_2 & dI_2 \end{bmatrix}P. \] Furthermore, let 
  \[Q = \text{Dg} [1, J_2(1),1],\] then 
  \[ Q^{-1} \begin{bmatrix} cI_2 & I_2\\
  (a - c) (c-b)I_2 & dI_2 \end{bmatrix} Q =  
\begin{bmatrix}
 c & 1 & 1 & 0 \\
 0 & c & 0 & 1 \\
 (a-c) (c-b) & d-c & d & -1 \\
 0 & (a-c) (c-b) & 0 & d \\
\end{bmatrix}.\] Summarizing the above: $\text{Dg}[a,b,a,b]$ is similar to 
\[ \begin{bmatrix}
 0 & 0 & 0 & 0 \\
 0 & 0 & 0 & 0 \\
 (a-c) (c-b) & d-c & 0 & 0 \\
 0 & (a-c) (c-b) & 0 & 0 \\
\end{bmatrix} + \begin{bmatrix}
 c & 1 & 1 & 0 \\
 0 & c & 0 & 1 \\
 0 & 0 & d & -1 \\
 0 & 0 & 0 & d \\
\end{bmatrix}.\] By direct computation it is easy to verify that the matrix on the left is square-zero, and since $c \neq d$ it follows by Proposition \ref{prop:coprime_chars_sim} that the matrix on the right is similar to 
\[\begin{bmatrix}
 c & 1 & 0 & 0 \\
 0 & c & 0 & 0 \\
 0 & 0 & d & -1 \\
 0 & 0 & 0 & d \\
\end{bmatrix} \approx \text{Dg}[J_2(c),J_2(d)].\] Since the square-zero property is invariant with respect to similarity the desired result follows. \hfill $\square$

The next lemma was proved for an arbitrary matrix over the complex field by Wang and Wu, but only for the case where all the invariant polynomials are of degree higher than one. For the special case of matrices with minimum polynomial $(x-a)(x-b)$ the result was extended by Takahashi. This result will allow us to prove the sufficient part of the main result.
\begin{quote}
\begin{lem}[Proposition 3.3 \cite{ww}, Lemma 2 \cite{takahashi}]
Let $\mathscr{F}$ be a field of characteristic zero. Suppose $A \in M_n(\mathscr{F})$ has minimum polynomial $(x-a)(x-b)$ where $a \neq b$ and that $A$ has at most two invariant polynomials of degree one. If the trace of $A$ is zero, then $A$ is a sum of three square-zero matrices.   
\label{lem:invariant_factors_order_2}
\end{lem}
\end{quote}
{\bf Proof.}  By Theorem 7.2 in the text by Cullen \cite{cullen} $A$ is similar to 
\[B = \text{Dg}[C(f_1(x)),C(f_2(x)), \ldots, C(f_t(x))],\] where $f_1(x), f_2(x), \ldots, f_t(x)$ are the nonconstant invariant polynomials of $A$. Since the minimum polynomial of $A$ is $(x-a)(x-b)$ we can assume without loss of generality that $f_i(x) = (x-a)(x-b)$ or $f_i(x) = (x-a)$, and therefore the trace of $f_i(x)$ is either $a + b$ or $a$ for each $i \in \{1,2, \ldots,t\}$.

Now suppose first that all the invariant polynomials of $A$ have degree greater than one. In this case each invariant polynomial of $A$ is $(x-a)(x-b)$. Then since the trace of $A$ is zero we must have $t(a + b) = 0$, and therefore since the characteristic of $\mathscr{F}$ is zero and $t \neq 0$ we must have $a + b = 0$, that is the trace of each invariant polynomial of $A$ is zero. By Corollary \ref{cor:tr0_nond_3sqz} it follows that $A$ is a sum of three square-zero matrices.

Now suppose $A$ has one invariant polynomial of degree one. Then we must have $f_1(x) = x-a$ is the invariant polynomial of degree one. Since the trace of $A$ is zero, $a \neq b$, and the characteristic of $\mathscr{F}$ is zero, we must have $a \neq 0$. Now we can apply Lemma \ref{lem:construct_3_mechanism} to each matrix $C(f_i(x))$ for $i \geq 2$, that is $C(f_i(x))$ is the sum of a square-zero matrix and an arbitrary companion matrix $C_i$ subject only to the condition that $\text{tr}(C(f_i(x))) = \text{tr}(C_i)$. Let
\[C_i = C \left ( \left (x- \left (-a-(i-1)(a+b) \right ) \right ) \left (x- \left (i(a+b)+a \right ) \right ) \right ).\]
Notice that the trace of $C_i$ is $a+b$ as required. I will now verify that the factors $(x-  (-a-(i-1)(a+b) )$ and $(x- (i(a+b)+a )$ are distinct: first $(t-1)b + ta = \text{tr}(A) = 0$, and since $t \geq 2$ it follows that
\begin{equation} b = -\frac{t}{(t-1)}a. \label{prop:invariant_factors_order_2_eq1}\end{equation} Now for any $i \in \{2,3,\ldots,t\}$, suppose that  $-a-(i-1)(a+b) = i(a+b)+a,$ then substituting $b$ from \eqref{prop:invariant_factors_order_2_eq1} and simplifying we have
\[(2(i-t)+1)a=0.\] Now since $2(i-t)+1$ is odd, and $\mathscr{F}$ is characteristic zero it follows that $a=0$ which is a contradiction. Therefore the factors $(x-  (-a-(i-1)(a+b) )$ and $(x- (i(a+b)+a )$ must be distinct.

It follows that for $i \geq 2$ the block $C_i$ has two characteristic values, and by construction one of the characteristic values of $C_i$ is also a characteristic value of $C_{i-1}$ (we set $C_1 = C(f_1(x)) = C(x-a)$), and the block $C_{t}$ has only one nonzero characteristic value, since
\[ t(a+b)+a = \text{tr}(A) = 0. \] 
Applying a similarity transformation with a permutation matrix to the rational canonical form of $\text{Dg}[C_1,C_2,\ldots,C_t]$, we find that it is similar to $\text{Dg}[N,X,-X]$, and therefore by Corollary \ref{cor:X_and_minusX} it follows that it is a sum of two square-zero matrices. By invariance of the square-zero property with respect to similarity, it follows that $A$ is a sum of three square-zero matrices.

It remains to prove the result for a matrix $A$ which has two invariant polynomials of degree one. First suppose that $A$ has an even number (denote this number as $m$) of invariant polynomials of degree two: in this case $A$ is similar to 
$\text{Dg}[A_1,A_1]$ where \[A_1 = \text{Dg}[\stackrel{m/2 \text{ blocks}}{\overbrace{C((x-a)(x-b)), \ldots, C((x-a)(x-b))}}, a].\] It follows that $2\text{tr}(A_1) = \text{tr}(A) =0$ and since $\mathscr{F}$ is of characteristic zero, $\text{tr}(A_1)=0$. By the preceding part it follows that $A_1$ is a sum of three square-zero matrices, and hence $A$ is a sum of three square-zero matrices, which concludes the proof in this case.

Now suppose that $A$ has an odd number of invariant polynomials of degree two: denote this number as $2k+1$. Notice that, since the trace of $A$ is zero, we must have $(2k+1)(a+b) + 2a = 0$, and therefore $a+b = -2a/(2k+1)$. Let us change the index of the invariant polynomials of $A$ so that the invariant polynomials of degree two are denoted as $f_1(x), f_2(x), \ldots, f_{2k+1}(x)$. Now for $i \in \{1,2, \ldots, k\}$, by Lemma~\ref{lem:lem1takahashi} we can write
\begin{eqnarray} 
\nonumber \text{Dg}[C(f_{2i-1}(x)), \text{Dg}[C(f_{2i}(x))] &\approx& \text{Dg}[a,b,a,b] \\
\nonumber &\approx& S_i + \frac{1}{2k+1}\text{Dg} [J_2(-(2i+1)a),J_2((2i-1)a)],
\end{eqnarray}
where $S_i$ is a square-zero matrix. Furthermore, by Lemma \ref{lem:construct_3_mechanism} we can write
\begin{eqnarray}
\nonumber C(f_{2k+1}(x)) &\approx& \text{Dg}[a,b] \\
\nonumber &\approx& S_{2k+1} + C\left ( \left (x+\frac{a}{2k+1}\right ) \left (x+\frac{a}{2k+1} \right) \right )\\
\nonumber &\approx& S_{2k+1} + \frac{1}{2k+1}J_2(-a),
\end{eqnarray} where $S_{2k+1}$ is square-zero. Now consider the two invariant polynomials of degree one; we can write:
\[\text{Dg}[C(x-a),C(x-a)] = J_2(a) - J_2(0) = \frac{1}{2k+1}J_2((2k+1)a) - J_2(0).\] Taking all of the above together, and applying a simple permutation similarity transformation as needed, $A$ is similar to the sum $S+B_1$ where 
\[S= \text{Dg}[S_{2k+1}, S_1, S_2, \ldots, S_{2k}, -J_2(0)],\] which is square-zero since each of its diagonal blocks are square-zero, and
\[B_1  = \frac{1}{2k+1}\text{Dg}[J_2(-a), J_2(a), J_2(-3a), J_2(3a), \ldots, J_2(-(2k+1)a), J_2((2k+1)a)].\] Now we can apply a simple permutation similarity transformation to see that $B \approx \text{Dg}[X,-X]$ where $X$ is invertible since each of the diagonal blocks of $B$ is invertible, and therefore by Corollary \ref{cor:X_and_minusX} $B$ is a sum of two square-zero matrices. It follows that $A$ is a sum of three square-zero matrices, concluding the proof. \hfill $\square$

\begin{quote}
\begin{prop}[Proposition 1 \cite{takahashi}]
Let $\mathscr{F}$ be a field of characteristic zero. Suppose $A \in M_n(\mathscr{F})$ has trace zero and minimum polynomial $(x-a)(x-b)$ where $a \neq b$ and that $A$ has at least one invariant polynomial of degree one. Let $m$ be the number of invariant polynomials of degree two and $r$ the number of invariant polynomials of degree one. Then $A$ is a sum of three square-zero matrices if and only if $r$ divides $2m$.   
\label{prop:takahashi_prop_1}
\end{prop}
\end{quote}
{\bf Proof.} Without loss of generality we may assume that the invariant polynomials of degree one are $(x-a)$.

Suppose $rs = 2m$ for some integer $s$. If $r$ is odd, then $s$ must be even, and we can write $r(s/2) = m$ where $s/2$ is an integer. It follows that $A$ is similar to a diagonal block matrix with $r$ identical blocks 
\[ B = \text{Dg}[ \stackrel{s/2 \text{ blocks}}{\overbrace{C((x-a)(x-b)), \cdots, C((x-a)(x-b))}}, a]\] on the diagonal. Now $\text{tr}(A) = r \cdot \text{tr}(B) = 0$ and since $r > 0$ and $\mathscr{F}$ is of characteristic zero, $\text{tr}(B) = 0$. By Lemma \ref{lem:invariant_factors_order_2} $B$ is therefore a sum of three square-zero matrices and consequently $A$ is a sum of three square-zero matrices. 

If $r$ is even then we can write $(r/2)s=m$ where $r/2$ is an integer. It follows that $A$ is similar to a diagonal block matrix with $r/2$ identical diagonal blocks
\[ B = \text{Dg}[ \stackrel{s \text{ blocks}}{\overbrace{C((x-a)(x-b)), \cdots, C((x-a)(x-b))}}, a,a]\] on the diagonal, and $\text{tr}(A) = r/2 \cdot \text{tr}(B) = 0$ and since $r > 0$ and $\mathscr{F}$ is of characteristic~zero, $\text{tr}(B) = 0$. By Lemma \ref{lem:invariant_factors_order_2} $B$ is therefore a sum of three square-zero matrices and consequently $A$ is a sum of three square-zero matrices. 

Now suppose that $A$ is a sum of three square-zero matrices. Then there exists a square-zero matrix $S$ such that $A+S$ is a sum of two square-zero matrices. Now $A$ is similar to
\[\text{Dg}[\stackrel{m \text{ blocks}}{\overbrace{C((x-a)(x-b)), \ldots, C((x-a)(x-b))}}, \stackrel{r \text{ blocks}}{\overbrace{C(x-a), \ldots, C(x-a)}}],\] and therefore $n=2m+r$ and $\text{n}(A-aI) = m+r$. By the first equation $m=(n-r)/2$, and substituting this result into the second equation we have $\text{n}(A-aI) = (n+r)/2$. Therefore, since $r \geq 1$ it follows that $\text{n}(A-aI)>n/2$. Also, since $S$ is square-zero $\text{n}(S) \geq n/2$. Therefore there exists a nonzero vector $v \in \mathscr{F}^n$ such that $v \in \text{n}(A-aI) \cap \text{n}(S)$. But then $((A+S) - aI)v = \mathbf{0}$ and it follows that $a$ is a characteristic value of $A+S$. By Corollary \ref{cor:X_and_minusX} $-a$ must also be a characteristic value of $A+S$. 

Now since $a \neq b$, and 
\begin{equation} \text{tr}(A)=(r+m)a + mb=0, \label{eq:trAnul_invariant_factors_order_2} \end{equation} and $r \geq 1$ we must have $a \neq -b$. 
Therefore $ka+kb \neq 0$ for every integer $k \neq 0$, and consequently $-(ka + (k-1)b)\neq b$ for every integer $k \neq 0$. Furthermore, since $\mathscr{F}$ is characteristic zero and $a \neq b$ we cannot have $a=0$, and hence $a \neq -a$.
Now since $-a \neq a,b$ it follows by Lemma \ref{lem:eigens_of_A_plus_S} that $2a+b$ must also be a characteristic value of $A+S$, and hence $-(2a+b)$ is a characteristic value of $A+S$. Now if $-(2a+b)=a$, then $3a=-b$ and by \eqref{eq:trAnul_invariant_factors_order_2} we then have $r=2m$, which is the desired result. 

Suppose therefore that $-(2a+b) \neq a$. Since $-(ka + (k-1)b)\neq b$ for every integer $k \neq 0$ as shown above, we must have $-(2a+b) \neq b$. Therefore by Lemma \ref{lem:eigens_of_A_plus_S} we must then have $3a+2b$ is a characteristic value of $A+S$. Now we can repeat the preceding argument inductively to prove the desired result. Explicitly, suppose $-(ka+(k-1)b) \neq a$ for all $k \geq 2$, then since we also have that $-(ka+(k-1)b) \neq b$ we must have by Lemma \ref{lem:eigens_of_A_plus_S} that $-(ka+(k-1)b)$ is a characteristic value of $A+S$ for all $k$ which is impossible since $A+S$ has finite order $n$. Therefore we must have that $-(ka+(k-1)b) = a$ for some $k \geq 2$. But then by \eqref{eq:trAnul_invariant_factors_order_2} we must have $(k-1)r = 2m$, which proves the desired result. \hfill $\square$

Now by Proposition \ref{prop:takahashi_prop_1} we can prove that Proposition \ref{prop:3sum_nec_ww} is not a sufficient condition for a matrix (over a field which is not of characteristic two) to be a sum of three square-zero matrices. As an example it may be verified that \[A = \text{Dg}[-5,2,-5,2,2,2,2] \in M_7(\mathbb{C})\] has trace zero and that the maximum nullity of $A-\lambda I$ is achieved when $\lambda = 2$ and then we have $\text{n}(A-2I) = 5 < 3 \cdot 7/4$, so that $A$ has the desired properties as set forth in Proposition \ref{prop:3sum_nec_ww}. But now since $A$ has $r=3$ invariant polynomials of degree one, and $m=2$ invariant polynomials of degree two, it follows that $r$ does not divide $2m$, and therefore by Proposition \ref{prop:takahashi_prop_1} $A$ is not a sum of three square-zero matrices.

Finally, let us consider some sufficient conditions for the special case where the characteristic of the underlying field $\mathscr{F}$ is zero. 

\begin{quote}
\begin{prop}[Proposition 3.3 \cite{ww}]
Let $\mathscr{F}$ be a field which has characteristic zero. If all the invariant polynomials of $A \in M_n(\mathscr{F})$ have degree greater than one, and the trace of $A$ is zero, then $A$ is a sum of three square-zero matrices.   
\label{prop:invariant_factors_order_2_general}
\end{prop}
\end{quote}
{\bf Proof.}
By Theorem 7.2 in the text by Cullen \cite{cullen} $A$ is similar to 
\[B = \text{Dg}[C(f_1(x)),C(f_2(x)), \ldots, C(f_t(x))],\] where $f_1(x), f_2(x), \ldots, f_t(x)$ are the nonconstant invariant polynomials of $A$. Let the trace of $C(f_i(x))$ be $c_i$, and the degree of $f_i$ be $n_i > 1$ for each $i \in \{1,2, \ldots, t\}$. 

First, suppose $c_i =0$ for all $i \in \{1,2, \ldots, t\}$: then the result follows directly from Corollary \ref{cor:tr0_nond_3sqz}. Let us assume therefore that $c_i \neq 0$ for some $i \in \{1,2, \ldots, t\}$. Now since $\mathscr{F}$ is characteristic zero we can choose $c \in \mathscr{F}$ so that 
\begin{itemize}
\item $c \neq 0$,
\item $2c \neq \sum_{j=1}^{i-1} c_j + \sum_{j=1}^{i} c_j$ for each $i \in \{1,2, \ldots, t\}$,
\item and $c \neq \sum_{j=1}^{i} c_j$ for each $i \in \{1,2, \ldots, t\}$.
\end{itemize}
 
Now applying Lemma \ref{lem:construct_3_mechanism} to each diagonal block of $B$ we can write 
\[ B = S + \text{Dg}[C_1,C_2, \ldots, C_t],\]
where $S$ is square-zero and
\begin{eqnarray}
\nonumber C_i &=& C \left ( \left (x- \left (c-\sum_{j=1}^{i-1}c_j \right ) \right ) \left (x- \left (\sum_{j=1}^{i}c_j-c \right ) \right ) x^{n_i-2} \right ) \\
 \nonumber &=& C \left ( \left (x- \lambda_{i1} \right ) \left (x- \lambda_{i2} \right ) x^{n_i-2} \right ), \end{eqnarray} for each $i \in \{1,2,\ldots,t\}$ (notice that the trace of $C_i$ is $c_i$ as required by Lemma \ref{lem:construct_3_mechanism}). Furthermore, by the choice of $c$ it follows that $\lambda_{i1} = c-\sum_{j=1}^{i-1}c_j \neq \sum_{j=1}^{i}c_j-c=\lambda_{i2}$. Explicitly, suppose for some $i$ we have $c-\sum_{j=1}^{i-1}c_j = \sum_{j=1}^{i}c_j-c$: then we must have $2c = \sum_{j=1}^{i-1}c_j + \sum_{j=1}^{i}c_j$ which by our choice of $c$ is a contradiction. Furthermore $\lambda_{i1}, \lambda_{i2} \neq 0$ since $c \neq 0$ and $c \neq \sum_{j=1}^{i} c_j$ for each $i \in \{1,2, \ldots, t\}$.

It follows that the rational canonical form of $C_i$ is 
\[C'_i = \text{Dg}\left [ c-\sum_{j=1}^{i-1}c_j,  \sum_{j=1}^{i}c_j-c, N_i \right] = \text{Dg} [ \lambda_{i1},  \lambda_{i2}, N_i ], \] where $N_i$ is some nilpotent matrix and $\lambda_{i1},  \lambda_{i2} \neq 0$.
 
So we have $\text{Dg}[C_1,C_2, \ldots, C_t]$ is similar to $\text{Dg}[C'_1,C'_2, \ldots, C'_t]$, where each $C'_i$ is in the form as described above. Since $\mathscr{F}$ is not of characteristic 2, notice that by a permutation of diagonal blocks $\text{Dg}[C'_1,C'_2, \ldots, C'_t]$ is similar to $\text{Dg}[N,X,-X]$ where $X$ is invertible. Explicitly, for each pair $C'_i, C'_{i+1}$ we have \[\lambda_{i2} = \sum_{j=1}^{i}c_j-c = -\left (c-\sum_{j=1}^{(i+1)-1}c_j \right ) = -\lambda_{(i+1)1}\] and furthermore $\lambda_{11} = c-\sum_{j=1}^{0}c_j = c$ is a characteristic value of $C'_1$, whereas $\lambda_{t2} = \sum_{j=1}^{t}c_j-c=-c$ (since the trace of $A$ is zero) is a characteristic value of $C'_t$. 
 
 It follows by Corollary \ref{cor:X_and_minusX} that $\text{Dg}[C'_1,C'_2, \ldots, C'_t]$ is a sum of two square-zero matrices. By invariance of the square-zero property with respect to similarity, it follows that $A$ is a sum of three square-zero matrices. \hfill $\square$

The final result is a sufficient condition in terms of a similar criterium as employed in Proposition \ref{prop:3sum_nec_ww}. I present this result in two propositions, where the second requires the first. The first proposition requires the following ``choice of characteristic polynomial" lemma.
\begin{quote}
\begin{lem}[Lemma 5 \cite{takahashi}]
Let $\mathscr{F}$ be an algebraically closed field. Suppose $A = \text{Dg}[B,aI_r]$ where $B \in M_n(\mathscr{F})$ is nonderogatory, $a \in \mathscr{F}$, and $r \leq n-2$. Then there exists a square-zero matrix $S$ so that the characteristic polynomial of $A+S$ can be arbitrarily chosen subject to the usual trace condition, that is the coefficient of $x^{n+r-1}$ must be equal to $-\text{tr}(A) = -(\text{tr}(B) + \text{tr}(aI_r))$.
\label{lem:choice_of_poly_takahashi}
\end{lem}
\end{quote}
{\bf Proof.} If $A$ is nonderogatory (i.e. $r=0$, or $(x-a)$ is not a factor of the characteristic polynomial of $B$ and $r=1$) then the result follows directly from Lemma \ref{lem:construct_3_mechanism}. 

Let us therefore assume $(x-a)$ is a factor of $B$ and $r \geq 1$ (if $r \geq 2$ and $(x-a)$ is not a factor of $B$ we can replace $B$ with $\text{Dg}[B,a]$). Let $(x-a) \cdot \prod_{i=1}^{n-1}(x-b_i)$ be the characteristic polynomial of $B$ (where the values $b_i$ are not necessarily distinct or different from $a$). Then $B$ is similar to 
\[B_1 = \begin{bmatrix}
a & 0 & \cdots & & 0\\
1 & b_1 & 0 & \cdots & 0\\
0 & 1 & b_2 & \ddots & \vdots \\
\vdots & \ddots & \ddots & \ddots & 0 \\
0 & \cdots & 0 & 1 & b_{n-1}
\end{bmatrix},\] which is easy to see since the determinant of the submatrix obtained by deleting the first row and last column of $xI_n-B_1$ is 1, and $B_1$ and $B$ have the same characteristic polynomial. Let $P$ be an invertible matrix such that $B_1 = P^{-1} B P$, and let 
\[Q = [e_2, e_3, \ldots, e_{n}, e_1, e_2 + e_{n+1}, e_3 + e_{n+2},\ldots, e_{r+1} + e_{n+r}],\] then a similarity transformation of $A$ via the change-of-basis matrix $\text{Dg}[P,I] Q$ shows that 
\[A \approx \begin{bmatrix}
B_2 & C \\ \mathbf{0} & aI_{r+1}
\end{bmatrix},\] where 
\[B_2 = \begin{bmatrix}
b_1 & 0 & \cdots & & 0\\
1 & b_2 & 0 & \cdots & 0\\
0 & 1 & b_3 & \ddots & \vdots \\
\vdots & \ddots & \ddots & \ddots & 0 \\
0 & \cdots & 0 & 1 & b_{n-1}
\end{bmatrix} \text{ and } 
C = \begin{bmatrix}
1 & b_1-a & 0 & \cdots &  0\\
0 & 1 & b_2-a & \ddots &   \vdots\\
 \vdots &  & \ddots &  \ddots &   0 \\
 &  & &  1& b_{r}-a \\
0 & \cdots &  & 0 & 1 \\
0 & \cdots &  & 0 & 0 \\
\vdots &  &  & \vdots & \vdots
\end{bmatrix}.\] 
Note that this form of $C$ is ensured, since $r+1 \leq n-1$. The proof will be complete if we can prove that there exists a matrix $N$ so that \[\begin{bmatrix}
B_2 & C \\ N & aI_{r+1}
\end{bmatrix}\] has an arbitrary characteristic polynomial, subject to the trace condition as specified in the statement of this lemma. Notice that
\[\begin{bmatrix}
I & \mathbf{0} \\ -X & I
\end{bmatrix}
\begin{bmatrix}
B_2 & C \\ N & aI_{r+1}
\end{bmatrix}
\begin{bmatrix}
I & \mathbf{0} \\ X & I
\end{bmatrix} = \begin{bmatrix}
B_2 + CX & C \\ 
N-XB_2-XCX+aX & aI_{r+1}-XC
\end{bmatrix}.
\] Now if we choose $N = XB_2+XCX-aX$, the matrix on the right above is block upper triangular, and therefore its characteristic polynomial is then the product of the characteristic polynomials of $B_2+CX$ and $aI_{r+1}-XC$. And since 
\[\begin{bmatrix}
B_2 & C \\ N & aI_{r+1}
\end{bmatrix}\] is similar to this matrix, its characteristic polynomial is then also the product of the characteristic polynomials of $B_2+CX$ and $aI_{r+1}-XC$. The proof therefore reduces to proving the existence of a matrix $X$ such that the product of these characteristic polynomials is arbitrary, subject to the trace condition. To this end, let 
\[ X = \begin{bmatrix}
d_1 & s_1 & s_2 & \cdots &  s_r & s_{r+1} & \cdots & s_{n-2}\\
0 & d_2 & 0 &  &   0& 0 & & 0\\
 \vdots &  & \ddots &  \ddots &   \vdots & \vdots & & \vdots\\
 &  & &  & 0 & 0& \cdots &\\
0 & \cdots &  & 0 & d_{r+1} & 0 & \cdots & 0
\end{bmatrix}. \] By this choice $aI_{r+1}-XC$ is an upper triangular matrix with entries $a-d_1, a-d_2, \ldots, a-d_{r+1}$ on the diagonal: we can therefore fix the characteristic polynomial of $aI_{r+1}-XC$ to be an arbitrary monic polynomial of degree $r+1$ by specifying $d_1, d_2, \ldots, d_{r+1}$ accordingly.

Now 
\[B_2 + CX = \begin{bmatrix}
e_1 & s_1+f_1 & s_2 & \cdots&  & s_{n-2}\\
1 & e_2 & f_2 & 0 & \cdots & 0\\
0 & 1 & e_3 & \ddots & \ddots & \vdots \\
\vdots & \ddots & \ddots & \ddots & & 0 \\
& & & & e_{n-2} & f_{n-2} \\
0 & \cdots & & 0 & 1 & e_{n-1}
\end{bmatrix} \] where 
\[e_i = \left \{ {\renewcommand{\arraystretch}{1.2} \begin{array}{ll} b_i + d_i & \text{ if } i \leq r+1 \\
b_i & \text{ if } i > r+1\end{array}} \right.\] and 
\[f_i = \left \{ {\renewcommand{\arraystretch}{1.2} \begin{array}{ll} 
(b_i-a)d_{i+1} & \text{ if } i \leq r \\
0 & \text{ if } i > r
\end{array}} \right. .\] 
Note that $e_i$ and $f_i$ are fixed for all $i \in \{1,2, \cdots,n-1\}$, as the matrices $B$ and $aI_{r}$ are prescribed, and $d_i$ was fixed in assigning the characteristic polynomial of $aI_{r+1}-XC$. It remains to prescribe the characteristic polynomial of $B_2 + CX$ through the choice of $s_1, s_2, \ldots, s_{n-2}$: to this end consider the determinant of 
\[xI-(B_2 + CX) = \begin{bmatrix}
(x-e_1) & -s_1-f_1 & -s_2 & \cdots&  & -s_{n-2}\\
-1 & (x-e_2) & -f_2 & 0 & \cdots & 0\\
0 & -1 & (x-e_3) & \ddots & \ddots & \vdots \\
\vdots & \ddots & \ddots & \ddots & & 0 \\
& & & & (x-e_{n-2}) & -f_{n-2} \\
0 & \cdots & & 0 & -1 & (x-e_{n-1})
\end{bmatrix}. \]
Developing the determinant inductively along the first column yields
\begin{equation} \det(xI-(B_2 + CX)) = (x-e_1)p_2(x) - f_1p_3(x) - \sum_{i=1}^{n-2} s_i p_{i+2}(x) \label{detb2cx} \end{equation} where $p_i(x)$ is some monic polynomial of degree $n-i$. Now $(x-e_1)p_2(x)$ is monic and of degree $n-1$, whereas all the other terms in \eqref{detb2cx} are degree less than or equal to $n-3$. It follows that the trace of $\det(xI-(B_2 + CX))$ is determined by the coefficient of $x^{n-2}$ in the expansion of $(x-e_1)p_2(x)$, and it is easy to see that this coefficient is fixed to $-\sum_{i=1}^{n-1} e_i$. Now since $\{p_3(x), p_4(x), \ldots, p_{n}(x)\}$ is a basis for the vector space of polynomials of maximum degree $n-3$, the coefficients of $1,x,x^2,\ldots,x^{n-3}$ in $\det(xI-(B_2 + CX))$ are completely determined in terms of $s_1, s_2, \ldots, s_{n-2}$.

Suppose now that 
\[ \det(xI-(B_2 + CX)) = x^{n-1} - \left (\sum_{i=1}^{n-1} e_i \right )x^{n-2} -\sum_{i=0}^{n-3} c_ix^i,\] where $c_i$ is a coefficient of choice as determined in the previous step. Then the characteristic polynomial of  \[\begin{bmatrix}
B_2 & C \\ N & aI_{r+1}
\end{bmatrix}\] is $\det(xI-(B_2 + CX)) \cdot \det(aI_{r+1} -XC)$, which is equal to 
\[\left ( x^{n-1} - \left (\sum_{i=1}^{n-1} e_i \right )x^{n-2} -\sum_{i=0}^{n-3} c_ix^i \right ) \left ( \prod_{i=1}^{r+1} (x-(a-d_i)) \right ). \]
The trace of this polynomial is 
\begin{eqnarray} 
\nonumber \sum_{i=1}^{n-1} e_i + \sum_{i=1}^{r+1} (a-d_i) &=& \sum_{i=1}^{n-1} b_i + \sum_{i=1}^{r+1} d_i  +\sum_{i=1}^{r+1}  (a-d_i)\\
\nonumber  &=& \left ( \sum_{i=1}^{n-1} b_i + a \right ) + ra \\
\nonumber &=& \text{tr} (B) + \text{tr} (aI_{r}),
\end{eqnarray} and the coefficients of $1,x,x^2,\ldots,x^{n+r-2}$ can be freely chosen in terms of the entries $c_i$ and $a_i-d_i$, which in turn is determined by the entries $d_1, d_2, \ldots, d_{r+1}$ and $s_1, s_2, \ldots, s_{n-2}$ of $X$ as proved above. The matrix $X$ determines the matrix $N$, and therefore
\[ S = \text{Dg}[P,I] Q \begin{bmatrix} \mathbf{0} & \mathbf{0} \\ N & \mathbf{0} \end{bmatrix} (\text{Dg}[P,I] Q)^{-1} \] is a square-zero matrix such that the characteristic polynomial of $A+S$ can be arbitrarily chosen subject to the coefficient of $x^{n+r-1}$ being $-\text{tr}(A)$. \hfill $\square$

\begin{quote}
\begin{prop}[Proposition 3 \cite{takahashi} (part 1)]
Let $\mathscr{F}$ be an algebraically closed field which has characteristic zero. Let $A \in M_n(\mathscr{F})$ have trace zero. Then $A$ is a sum of three square-zero matrices if $\text{n}(A-\lambda I) \leq n/2$ for all $\lambda \in \mathscr{F}$.   
\label{prop:tak_suf_1}
\end{prop}
\end{quote}
{\bf Proof.} Choose $\lambda$ so that $\text{n}(A-\lambda I)$ is at its maximum, and let $l, m$ and $r$ be the number of invariant polynomials of degree at least three, degree two, and degree one respectively. Notice that if $\lambda = 0$, $A$ is similar to $\text{Dg}[B,0,\ldots,0]$ where all the invariant polynomials of $B$ have degree greater than one, and therefore by Proposition \ref{prop:invariant_factors_order_2_general} $A$ is a sum of three square-zero matrices. 

Let us therefore assume that $\lambda \neq 0$. By the choice of $\lambda$, and Theorem 7.2 in the text by Cullen \cite{cullen} $A$ is similar to 
\[B = \text{Dg}[\lambda I_r, C(f_{1}(x)),C(f_{2}(x)), \ldots, C(f_{l+m}(x))],\] where we may assume $f_{1}(x), f_{2}(x), \ldots, f_{l}(x)$ are the nonconstant invariant polynomials of $A$ of degree at least three, and $f_{l+1}(x), f_{l+2}(x), \ldots, f_{l+m}(x)$ are the nonconstant invariant polynomials of $A$ of degree two (apply a simple permutation similarity transformation as needed). Let the degree of $f_{i}$ be $n_i$. Notice that by the rational canonical form theory $(x-\lambda)$ must be a factor of each invariant polynomial $f_i(x)$ and therefore
\[\text{n}(A-\lambda I) = l + m + r, \] and furthermore
\[n = \sum_{i=1}^l n_i + 2\cdot m + r.  \]

Therefore if we require that $\text{n}(A-\lambda I) \leq n/2$, we require
\[l+m+r \leq \frac{\sum_{i=1}^l n_i + 2\cdot m + r}{2},\] and therefore that
\[r \leq \sum_{i=1}^l (n_i-2).\] Choose $r_1, r_2, \ldots, r_l$ so that $r = \sum_{i=1}^l r_i$ and $r_i \leq n_i -2$. By a simple permutation similarity transformation $B$ is similar to 
\[\text{Dg}\left [\begin{bmatrix} C(f_{1}) & \mathbf{0} \\ \mathbf{0} & \lambda I_{r_1}\end{bmatrix}, \begin{bmatrix} C(f_{2}) & \mathbf{0} \\ \mathbf{0} & \lambda I_{r_2}\end{bmatrix}, \ldots, \begin{bmatrix} C(f_{l}) & \mathbf{0} \\ \mathbf{0} & \lambda I_{r_l}\end{bmatrix},C(f_{l+1}), \ldots, C(f_{l+m}) \right ].\] Set \[A_i = \begin{bmatrix} C(f_{i}) & \mathbf{0} \\ \mathbf{0} & \lambda I_{r_i}\end{bmatrix}\] if $i \leq l$ and set $A_i = C(f_{i})$ if $i > l$. Let $t_{i}$ be the trace of $A_i$. Now as in the proof of Proposition \ref{prop:invariant_factors_order_2_general} choose $c \in \mathscr{F}$ so that \begin{itemize}
\item $c \neq 0$,
\item $2c \neq \sum_{j=1}^{i-1} t_j + \sum_{j=1}^{i} t_j$ for each $i \in \{1,2, \ldots, l+m\}$,
\item and $c \neq \sum_{j=1}^{i} t_j$ for each $i \in \{1,2, \ldots, l+m\}$.
\end{itemize}

Now apply Lemma \ref{lem:choice_of_poly_takahashi} to select square-zero matrices $S_i$ such that the characteristic polynomial of $A_i+S_i$ is
\[\left (x- \left (c-\sum_{j=1}^{i-1}t_j \right )\right )  \left (x- \left (\sum_{j=1}^{i}t_j-c \right ) \right ) x^{n_i+r_i-2} =  \left (x- \lambda_{i1} \right ) \left (x- \lambda_{i2} \right ) x^{n_i+r_i-2} , \] for each $i \in \{1,2,\ldots,l+m\}$. As shown before in the proof of Proposition \ref{prop:invariant_factors_order_2_general}, by the choice of $c$, and over a field of characteristic zero: $\lambda_{i1}, \lambda_{i2} \neq 0$, $\lambda_{i1} \neq \lambda_{i2}$, $\lambda_{i2} = -\lambda_{(i+1)1}$, and furthermore $\lambda_{11} = c$ and by the zero trace condition $\lambda_{(l+m)2} = -c$. Combining all these properties, 
\[ \text{Dg}[A_1, A_2, \ldots, A_{l+m}] + \text{Dg}[S_1, S_2, \ldots, S_{l+m}] \approx A + S \] is similar to a matrix $\text{Dg}[N,X,-X]$ where $N$ is nilpotent and $X$ is invertible. By Corollary \ref{cor:X_and_minusX} it follows that $A+S$ is a sum of two square-zero matrices, and therefore $A$ is a sum of three square-zero matrices. \hfill $\square$

The following lemma is required for the final result; I follow Takahashi's proof.
\begin{quote}
\begin{lem}[Lemma 6 \cite{takahashi}]
Let $\mathscr{F}$ be an algebraically closed field. Suppose $A = \text{Dg}[B,aI_{n-2}]$ where $B \in M_n(\mathscr{F})$ is nonderogatory, and $a \in \mathscr{F}$ is a characteristic value of $B$, and let $s$ be a positive integer such that $s\leq n-1$. Then there exists a square-zero matrix $S$ so that $A+S$ is similar to $\text{Dg}[A_1,A_2+\lambda I_s]$ where $A_2$ is a square-zero matrix, $\lambda \in \mathscr{F}$ is arbitrary, and $A_1$ is a matrix with an arbitrary characteristic polynomial subject to $\text{tr}(A_1) + s \lambda = \text{tr}(A)$, and $\lambda$ is not one of its characteristic values.
\label{lem:choice_of_poly_2_takahashi}
\end{lem}
\end{quote}
{\bf Proof.} We can follow the first part of the proof of Lemma \ref{lem:choice_of_poly_takahashi} to see that 
\[A \approx \begin{bmatrix}
B_2 & C \\ \mathbf{0} & aI_{n-1}
\end{bmatrix}\] where $B_2$ and $C$ have the same form as in the proof of Lemma \ref{lem:choice_of_poly_takahashi}. In this case $C$ is a square matrix of order $n-1$, and hence is invertible. It follows that 
\[A \approx  \begin{bmatrix}
C^{-1} & \mathbf{0} \\ \mathbf{0} & I_{n-1}
\end{bmatrix}  \begin{bmatrix}
B_2 & C \\ \mathbf{0} & aI_{n-1}
\end{bmatrix}  \begin{bmatrix}
C & \mathbf{0} \\ \mathbf{0} & I_{n-1}
\end{bmatrix} =  \begin{bmatrix}
C^{-1}B_2C & I_{n-1} \\ \mathbf{0} & aI_{n-1}
\end{bmatrix},\] and since $B_2$ is nonderogatory, $C^{-1}B_2C$ is also nonderogatory. 

Now suppose we choose the characteristic polynomial of $A_1$ to be 
\[\prod_{i=1}^{2(n-1)-s} (x-\lambda_i),\] where the $\lambda_i$ are arbitrary subject to $s\lambda + \sum_{i=1}^{2(n-1)-s} \lambda_i = \text{tr}(A)$, and $\lambda_i \neq \lambda$ for all $i \in \{1,2,\ldots,2(n-1)-s\}$.

Let $D = C^{-1}B_2C + aI_{n-1}$; it is easy to see that $D$ is also nonderogatory: $C^{-1}B_2C + aI_{n-1} = C^{-1}B_2C + aC^{-1}C = C^{-1}(B_2 + aI_{n-1})C$ and since $B_2 + aI_{n-1}$ is nonderogatory\footnote{Determinant of the submatrix obtained by deleting the first row and last column is 1.} the result follows. Hence there exists a vector $v$ such that the set $\{v, Dv, D^2v, \ldots, D^{n-2}v\}$ is linearly independent. Let
\[d_i = \left \{ {\renewcommand{\arraystretch}{1.2} \begin{array}{ll} 
\lambda_i+\lambda_{n-1+i} & \text{ if } 1 \leq i \leq n-1-s \\
\lambda_i+\lambda & \text{ if } n-1-s < i \leq n-1\end{array}} \right.,\] and let 
\[\pi_i = \left (\prod_{j=1}^{i-1} (D-d_j I_{n-1})\right )v \text{ for } 1 \leq i \leq n-1. \]
Then the matrix $P=[\pi_1, \pi_2, \ldots, \pi_{n-1}]$ is a change-of-basis matrix and 
\[P^{-1} D P = \begin{bmatrix}
d_1 & 0 & \cdots & 0 & c_0 \\
1 & d_2 & \ddots & \vdots & c_1 \\
0 & \ddots & \ddots & 0 & \vdots \\
\vdots & \ddots & &d_{n-2} & c_{n-3} \\
0 & \cdots & 0 & 1 & d_{n-1}
\end{bmatrix}\] where $c_0, c_1, \ldots, c_{n-3}$ are some scalars. The last column is obtained by expressing $D^{n-1}v$ as a linear combination of the columns of $P$ (which is a basis for the column space of $D$) and using the fact that $d_1+...+d_{n-1}=\text{tr}(A)=\text{tr}(D)$. Set
\[G = \begin{bmatrix}
\lambda_1 & 0 & \cdots &  & 0 \\
1 & \lambda_2 & \ddots &  & \vdots \\
0 & \ddots & \ddots &  &  \\
\vdots & \ddots & & \lambda_{n-2}& 0 \\
0 & \cdots & 0 & 1 & \lambda_{n-1}
\end{bmatrix} \] and $H = P^{-1} D P - G$. Partition $H$ so that 
\[H = \begin{bmatrix} 
H_{11} & H_{12} \\ H_{21} & H_{22}
\end{bmatrix}\] where $H_{11} = \text{Dg}[\lambda_n, \lambda_{n+1}, \ldots, \lambda_{2(n-1)-s}]$, $H_{12}$ is the matrix with last column \[(c_0,c_1, \ldots, c_{n-2-s})^T\] and zeros elsewhere, $H_{21} = \mathbf{0}$ and 
\[H_{22} = \begin{bmatrix}
\lambda & 0 & \cdots & 0 & c_{n-1-s} \\
0 & \lambda & \ddots & \vdots & c_{n-s} \\
\vdots & \ddots & \ddots & 0 & \vdots \\
 &  & &\lambda & c_{n-3} \\
0 & \cdots &  & 0 & \lambda
\end{bmatrix}.\] Now since $\lambda \neq \lambda_i$ for all $i \in \{1,2,\ldots,2(n-1)-s\}$ the characteristic polynomials of $H_{11}$ and $H_{22}$ are relatively prime and therefore we can apply Proposition \ref{prop:coprime_chars_sim} to obtain $H \approx \text{Dg}[H_{11},H_{22}]$. Note that $H_{22} = \lambda I_s + A_2$ where $A_2$ is square-zero, and therefore it remains to show that 
\[\begin{bmatrix}
C^{-1}B_2C & I \\ \mathbf{0} & aI_{n-1}
\end{bmatrix}\] is similar to $\text{Dg}[A_1,H_{22}]$ where $A_1$ is of the desired form. To this end, let $B_3=C^{-1}B_2C$ and apply a similarity transformation employing the change-of-basis matrix 
\[\begin{bmatrix} I & \mathbf{0} \\ PGP^{-1}-B_3 & I \end{bmatrix} 
\begin{bmatrix} P & \mathbf{0} \\ \mathbf{0} & P \end{bmatrix} = \begin{bmatrix} P & \mathbf{0} \\ (PGP^{-1}-B_3)P & P \end{bmatrix}.\] Now
\[ \begin{bmatrix} P^{-1} & \mathbf{0} \\ P^{-1}(B_3-PGP^{-1}) & P^{-1} \end{bmatrix} 
\begin{bmatrix} B_3 & I \\ \mathbf{0} & aI \end{bmatrix}  
\begin{bmatrix} P & \mathbf{0} \\ (PGP^{-1}-B_3)P & P \end{bmatrix}\]  is equal to
\[\begin{bmatrix} G & I \\ HG - aP^{-1}B_3P & H \end{bmatrix}=A'.\]

Set 
\[S' = - \begin{bmatrix} \mathbf{0} & \mathbf{0} \\ HG - aP^{-1}B_3P & \mathbf{0} \end{bmatrix}, \] then $S'$ is square-zero and $A'+S'$ is similar to 
\[\begin{bmatrix} G & I \\
\mathbf{0} & \text{Dg}[H_{11},H_{22}]
\end{bmatrix}.\] Now again, since the characteristic polynomials of $\text{Dg}[G, H_{11}]$ and $H_{22}$ are relatively prime, we can apply Proposition \ref{prop:coprime_chars_sim} to obtain 
\[A'+S' \approx \begin{bmatrix} G & J & \mathbf{0}  \\
\mathbf{0} & H_{11} & \mathbf{0} \\
\mathbf{0} & \mathbf{0} & H_{22}
\end{bmatrix},\] where $J = \begin{bmatrix}I_{n-1-s}\\ \mathbf{0}\end{bmatrix}$. Now 
\[A_1 = \begin{bmatrix} G & J  \\
\mathbf{0} & H_{11} \end{bmatrix} \] is a matrix of the desired form, and by the invariance of the square-zero property with respect to similarity $A+S \approx A'+S'$, which completes the proof. \hfill $\square$ 

\begin{quote}
\begin{prop}[Proposition 3 \cite{takahashi} (part 2)]
Let $\mathscr{F}$ be an algebraically closed field which has characteristic zero. Let $A \in M_n(\mathscr{F})$ have trace zero, and let $m$ be the number of invariant polynomials of $A$ which is of degree two. Then $A$ is a sum of three square-zero matrices if $\text{n}(A-\lambda I) \leq (2n-m)/3$ for all $\lambda \in \mathscr{F}$.   
\label{prop:tak_suf_2}
\end{prop}
\end{quote}
{\bf Proof.} As in the proof of Proposition \ref{prop:tak_suf_1}, choose $\lambda$ so that $\text{n}(A-\lambda I)$ is at its maximum, and let $l$ and $r$ be the number of invariant polynomials of degree at least three, and degree one respectively. As also shown before if $\lambda = 0$ the result follows easily from Proposition \ref{prop:invariant_factors_order_2_general}. Let us therefore assume that $\lambda \neq 0$. 

By the specific choice of $\lambda$, and Theorem 7.2 in the text by Cullen \cite{cullen} $A$ is similar to 
\[B = \text{Dg}[\lambda I_r, C(f_{1}(x)),C(f_{2}(x)), \ldots, C(f_{l+m}(x))],\] where we may assume $f_{1}(x), f_{2}(x), \ldots, f_{l}(x)$ are the nonconstant invariant polynomials of $A$ of degree at least three, and $f_{l+1}(x), f_{l+2}(x), \ldots, f_{l+m}(x)$ are the nonconstant invariant polynomials of $A$ of degree two (apply a simple permutation similarity transformation as needed). Let the degree of $f_{i}$ be $n_i$. 

If $r \leq \sum_{i=1}^l (n_i-2)$, then it follows as in the proof of Proposition \ref{prop:tak_suf_1} that $A$ is a sum of three square-zero matrices. Let us therefore assume $r > \sum_{i=1}^l (n_i-2)$, and let $s = r - \sum_{i=1}^l (n_i-2)$. Now since $\text{n}(A-\lambda I) \leq (2n-m)/3$, and as shown before $\text{n}(A-\lambda I) = l+m+r$ and $n = \sum_{i=1}^l n_i + 2\cdot m + r$, we must have
\[l+m+r \leq \frac{2(\sum_{i=1}^l n_i + 2\cdot m + r)-m}{3}\] which is true if and only if 
\[r \leq \sum_{i=1}^l (2n_i - 3).\] It follows that 
\[s \leq \sum_{i=1}^l (2n_i - 3) - \sum_{i=1}^l (n_i-2) = \sum_{i=1}^l (n_i-1),\] and therefore we can take integers $s_1, s_2, \ldots, s_t$ where $1 \leq t \leq l$ so that $s = \sum_{i=1}^t s_i$ and $1 \leq s_i \leq n_i-1$. Now let 
\[c_i = \left \{ {\renewcommand{\arraystretch}{1.2} \begin{array}{ll} 
\text{tr}(\text{Dg}[C(f_i(x)), \lambda I_{n_i-2}]) + \lambda s_i & \text{ if } 1 \leq i \leq t \\
\text{tr}(\text{Dg}[C(f_i(x)), \lambda I_{n_i-2}]) & \text{ if } t < i \leq l+m \end{array}} \right.\] and choose $c \in \mathscr{F}$ so that 
\begin{itemize}
\item $c \neq 0$,
\item $2c \neq \sum_{j=1}^{i-1} c_j + \sum_{j=1}^{i} c_j$, and $c \neq -\lambda + \sum_{j=1}^{i-1} c_j$, and $c \neq \lambda + \sum_{j=1}^{i} c_j$ for each $i \in \{1,2, \ldots, l+m\}$,
\item $c \neq \sum_{j=1}^{i} c_j$ for each $i \in \{1,2, \ldots, l+m\}$, and $c \neq -\lambda$.
\end{itemize} 
It is easy to verify that this choice of $c$ ensures that $-\lambda$, $c-\sum_{j=1}^{i-1} c_j$ and $\sum_{j=1}^{i} c_j - c$ are nonzero and mutually different for each $i \in \{1,2,\ldots,l+m\}$. 

Now for $1\leq i \leq t$ apply Lemma \ref{lem:choice_of_poly_2_takahashi}: there is a square-zero matrix $S_i$ such that $\text{Dg}[C(f_i(x)), \lambda I_{n_i-2}] + S_i$ is similar to $\text{Dg}[C_i, D_i-\lambda I_{s_i}]$ where we choose the characteristic polynomial of $C_i$ to be 
\[\left (x- \left (c-\sum_{j=1}^{i-1}c_j \right )\right )  \left (x- \left (\sum_{j=1}^{i}c_j-c \right ) \right ) x^{2n_i-s_i-4} =  \left (x- \lambda_{i1} \right ) \left (x- \lambda_{i2} \right ) x^{2n_i-s_i-4} , \]
and $D_i$ is square-zero. Furthermore, for $t < i \leq l+m$ we can apply Lemma \ref{lem:choice_of_poly_takahashi} so that there is a square-zero matrix $S_i$ such that $\text{Dg}[C(f_i(x)),\lambda I_{n_i-2}] + S_i$ is similar to $C_i$ where we choose the characteristic polynomial of $C_i$ to be
\[\left (x- \left (c-\sum_{j=1}^{i-1}c_j \right )\right )  \left (x- \left (\sum_{j=1}^{i}c_j-c \right ) \right ) x^{2n_i-4} =  \left (x- \lambda_{i1} \right ) \left (x- \lambda_{i2} \right ) x^{2n_i-4} .\]

At this stage we have that $B$ (and therefore also $A$) is similar to 
\begin{eqnarray} 
\nonumber B' &=& \text{Dg}[\lambda I_s, \text{Dg}[C(f_1(x)), \lambda I_{n_1-2}], \text{Dg}[C(f_2(x)), \lambda I_{n_2-2}], \ldots \\ 
\nonumber && \qquad \ldots, \text{Dg}[C(f_{l+m}(x)), \lambda I_{n_{(l+m)}-2}]]
\end{eqnarray} and
\begin{eqnarray}
\nonumber B'+S' &\approx& \text{Dg}[\lambda I_s, \text{Dg}[C_1, D_1- \lambda I_{s_1}], \ldots, \text{Dg}[C_t, D_t- \lambda I_{s_t}], C_{t+1}, \ldots, C_{l+m}]
\end{eqnarray} where $S' = \text{Dg}[\mathbf{0}_s, S_1, S_2, \ldots, S_{l+m}]$.

Now let $S'' = \text{Dg}[D_1, D_2, \ldots, D_t,\mathbf{0}_{n-s}];$\footnote{Note that $\text{Dg}[D_1, D_2, \ldots, D_t]$ is a matrix of order $s$} then 
\[S' + S'' = \text{Dg}[D_1, D_2, \ldots, D_t, S_1, S_2, \ldots, S_{l+m}]\] is square-zero, and 
\begin{eqnarray}
\nonumber B' + (S' + S'') &\approx& \text{Dg}[D_1+\lambda I_{s_1}, \ldots, D_t + \lambda I_{s_t}, \text{Dg}[C_1, D_1-\lambda I_{s_1}], \ldots \\
\nonumber && \qquad \ldots, \text{Dg}[C_t, D_t-\lambda I_{s_t}], C_{t+1}, \ldots, C_{l+m}] \\
\nonumber &\approx& \text{Dg}[D_1+ \lambda I_{s_1}, D_1- \lambda I_{s_1}, \ldots, D_t +  \lambda I_{s_t}, D_t- \lambda I_{s_t}, C_1, \ldots, C_{l+m}].
\end{eqnarray}
As shown before in the proof of Proposition \ref{prop:invariant_factors_order_2_general}, by the choice of $c$, and over a field of characteristic zero, $\lambda_{i1}, \lambda_{i2} \neq 0$, $\lambda_{i1} \neq \lambda_{i2}$, $\lambda_{i2} = -\lambda_{(i+1)1}$, and furthermore $\lambda_{11} = c$ and by the zero trace condition $\lambda_{(l+m)2} = -c$. It follows that the submatrix $\text{Dg}[C_1, C_2, \ldots, C_{l+m}]$ is similar to a matrix $\text{Dg}[N,X,-X]$ where $N$ is nilpotent and $X$ is invertible, and therefore by Corollary \ref{cor:X_and_minusX} it is a sum of two square-zero matrices. 

It remains to prove that each block $\text{Dg}[D_i+\lambda I_{s_i}, D_i- \lambda I_{s_i}]$ is a sum of two square-zero matrices, which is true if and only if $-(D_i+\lambda I_{s_i})$ is similar to $D_i- \lambda I_{s_i}$. To this end, notice that since $D_i$ is square-zero, all of its invariant polynomials are even-power and therefore by Corollary \ref{cor:A_sim_minusA} $D_i$ is similar to $-D_i$. Suppose $P$ is a change-of-basis matrix such that $P^{-1}(-D_i)P = D_i$; then 
\[P^{-1}(-D_i-\lambda I_{s_i})P = P^{-1}(-D_i)P + P^{-1}(-\lambda I_{s_i})P = D_i - \lambda P^{-1}P = D_i-\lambda I_{s_i}.\] It follows that $B' + (S' + S'')$ is a sum of two square-zero matrices, and by invariance of the square-zero property with respect to similarity $B' + (S' + S'') \approx A + S$ where $S$ is square-zero. And therefore $A$ is a sum of three square-zero matrices. \hfill $\square$

\begin{quote}
\begin{cor}[Corollary 1 \cite{takahashi}]
Let $\mathscr{F}$ be an algebraically closed field which has characteristic zero. If the trace of $A \in M_n(\mathscr{F})$ is zero and $\text{n}(A-\lambda I) \leq n/2 + 1$ for all $\lambda \in \mathscr{F}$, then $A$ is a sum of three square-zero matrices.   
\label{tak_cor_1}
\end{cor}
\end{quote}
{\bf Proof.} Choose $\lambda \in \mathscr{F}$ so that $\text{n}(A-\lambda I)$ is at its maximum. Let $l, m$ and $r$ be the number of invariant polynomials of degree at least three, degree two, and degree one respectively. By the particular choice of $\lambda$ and Theorem 7.2 in the text by Cullen \cite{cullen} $A$ is similar to 
\[B = \text{Dg}[\lambda I_r, C(f_{1}(x)),C(f_{2}(x)), \ldots, C(f_{l+m}(x))],\] where we may assume $f_i(x)$ is a polynomial of degree greater than two for $i \in \{1,2,\ldots,l\}$, otherwise it is of degree two. Let $n_i$ be the degree of $f_i(x)$. Then $n = \sum_{i=1}^l n_i + 2m + r$ and $\text{n}(A-\lambda I) = l + m + r$ since, by the Smith canonical form theory the factor $(x-\lambda)$ must also be a factor of invariant polynomials of degree greater than one. So the condition $\text{n}(A-\lambda I) \leq n/2+1$ is equivalent to 
\[l+m+r \leq \frac{\sum_{i=1}^l n_i + 2m + r}{2} + 1,\] which is true if and only if
\[r \leq \sum_{i=1}^l (n_i - 2) + 2.\] 
Now if $l=0$ this equation reduces to $r \leq 2$, and therefore in this case, by Lemma \ref{lem:invariant_factors_order_2} $A$ is a sum of three square-zero matrices. Suppose $l > 0$; then 
\[r \leq \sum_{i=1}^l (n_i - 2) + 2 \leq \sum_{i=1}^l (n_i - 2) + \sum_{i=1}^l (n_i - 1) = \sum_{i=1}^l (2n_i - 3),\] and therefore it follows as in the proof of Proposition \ref{prop:tak_suf_2} that $A$ is a sum of three square-zero matrices. \hfill $\square$
\section{Conclusion on Sums of Square-Zero Matrices}
In this chapter a full characterization of sums of two and four square-zero matrices was presented from the body of existing research. For two square-zero matrices we have a succinct result, by combining Theorems \ref{thm:2_sqsum_nec} and \ref{thm:sufficient_sqzero_sum}: 
\begin{quote} Let $\mathscr{F}$ be an arbitrary field. The matrix $A \in M_n(\mathscr{F})$ is a sum of two square-zero matrices if and only if all of its invariant polynomials are even-, or odd-power polynomials.\end{quote} 
An even-, or odd-power polynomial is defined to be a polynomial \[p(x) = \sum_{i=0}^n c_{n-i}x^{n-i}\] where $c_{n-i}=0$ if $i$ is odd. The polynomial is even-power if $n$ is even, else it is odd.

Several useful corollaries also proceed from this result. Corollaries \ref{cor:A_sim_minusA}, \ref{cor:AV_eq_minusVA} and \ref{cor:X_and_minusX} provide the following results for matrices over a field which is not of characteristic two:

\begin{quote}
Let $\mathscr{F}$ be a field which is not of characteristic two. Then the following statements are equivalent for a matrix $A \in M_n(\mathscr{F})$:
\begin{enumerate}
\item $A$ is a sum of two square-zero matrices.
\item $A$ is similar to $-A$.
\item There exists an involution $V  \in M_n(\mathscr{F})$ such that $AV=-VA$.
\item $A$ is similar to $\text{Dg}[N,X,-X,C]$, where $N$ is nilpotent, $X$ is invertible, and $C$ is similar to a block diagonal matrix with diagonal blocks of the form $C(p(x)^e)$ where $p(x)$ is an irreducible even-power polynomial with nonzero constant term. 
\end{enumerate}
\end{quote}
Corollary \ref{cor:2sum_char2_evenBlocks} provides us with the following result for matrices over a field of characteristic two:
\begin{quote}
Let $\mathscr{F}$ be a field which is of characteristic 2. The matrix $A \in M_n(\mathscr{F})$ is a sum of two square-zero matrices if and only if, in some algebraic closure of $\mathscr{F}$, the simple Jordan blocks of $A$ associated with nonzero eigenvalues all have even order.
\end{quote}

The most general question we can ask regarding sums of square-zero matrices, is whether a matrix can be written as a sum of square-zero matrices, regardless of the number of matrices in the sum. Again, in this case there is a succinct answer in Theorem \ref{thm:trace_zero_sum_of_4}, and it turns out that four is a sufficient number in the sum whenever it is possible to write a matrix as a sum of square-zero matrices.
\begin{quote} Let $\mathscr{F}$ be an arbitrary field. The matrix $A \in M_n(\mathscr{F})$ is a sum of four square-zero matrices if and only if its trace is zero.\end{quote}

Finally results were presented on sums of three square-zero matrices. For matrices over a field of characteristic two we have a complete characterization in Theorem~\ref{thm:fieldchar2_sum_of_three_main_result}:
\begin{quote} Let $\mathscr{F}$ be a field of characteristic two. Then $A \in M_n(\mathscr{F}$) is a sum of three square-zero matrices if and only if the trace of $A$ is zero. \end{quote}
For matrices over a field which is not of characteristic two, the problem remains open (at the time of writing). The following necessary condition was presented in Theorem \ref{prop:3sum_nec_ww}:
\begin{quote}  Let $\mathscr{F}$ be a field which is not of characteristic two. If $A \in M_n(\mathscr{F})$ is the sum of three square-zero matrices, then $\text{n}(A - \lambda I) \leq 3n/4$ for any nonzero $\lambda \in \mathscr{F}$. \end{quote}
By Proposition \ref{prop:takahashi_prop_1} we know that this condition is not sufficient.  

The final results are only valid over certain fields. First, we have a sufficient condition in Proposition \ref{prop:invariant_factors_order_2_general} which is valid over any field which has characteristic zero.
\begin{quote} 
Let $\mathscr{F}$ be a field which has characteristic zero. If all the invariant polynomials of $A \in M_n(\mathscr{F})$ have degree greater than one, and the trace of $A$ is zero, then $A$ is a sum of three square-zero matrices.   
\end{quote}
The following sufficient condition, which is further restricted to algebraically closed fields, was presented in Proposition \ref{prop:tak_suf_2}:
\begin{quote}
Let $\mathscr{F}$ be an algebraically closed field which has characteristic zero. Let $A \in M_n(\mathscr{F})$ have trace zero, and let $m$ be the number of invariant polynomials of $A$ which are of degree two. Then $A$ is a sum of three square-zero matrices if $\text{n}(A-\lambda I) \leq (2n-m)/3$ for all $\lambda \in \mathscr{F}$.   
\end{quote}
A useful corollary to this result which references the matrix structure only in terms of its order and its maximal eigenspace was presented in Corollary \ref{tak_cor_1}:
\begin{quote}
Let $\mathscr{F}$ be an algebraically closed field which has characteristic zero. If the trace of $A \in M_n(\mathscr{F})$ is zero and $\text{n}(A-\lambda I) \leq n/2 + 1$ for all $\lambda \in \mathscr{F}$, then $A$ is a sum of three square-zero matrices.  
\end{quote}

\chapter{Conclusion}
In this text I have presented a consolidation of results on sums and products of square-zero matrices. Let us summarize the key findings. 

Firstly within the broader context of nilpotent matrices, 
\begin{quote} A matrix $A \in M_n(\mathscr{F})$ ($\mathscr{F}$  an arbitrary field) is a product of (two) nilpotent matrices if and only if it is singular and not a nonzero nilpotent matrix of order $2 \times 2$. \end{quote}

\begin{quote} The matrix $A$ is a sum of nilpotent matrices if and only if its trace is zero, and if this is the case and $A$ is a nonzero scalar matrix, then $A$ can be written as a sum of three square-zero matrices, which is best possible. If $A$ is the zero matrix or non-scalar then $A$ can be written as a sum of two nilpotent matrices. \end{quote}

 An open problem within the broader context is the characterization of (pairs of) matrices in $M_{m \times n}(\mathscr{F})$ which have a nilpotent divisor or quotient (matrix division is discussed in detail in Section \ref{prelim}).

Which matrices are products of square-zero matrices? A full characterization of (pairs of) matrices in $M_{m \times n}(\mathscr{F})$ which have a square-zero divisor or quotient is given in Theorems \ref{thm:bothaSzDivisorPt1}, \ref{thm:bothaSzDivisorPt2} and \ref{thm:bothaSzDivisorPt3}. The last of these provides a method for the construction of a square-zero quotient for each allowable rank, whenever it is possible for two matrices to have such a quotient. Perhaps, if one was to single out one useful characterization from these theorems: 
\begin{quote} For matrices $G,F \in M_{m \times n}(\mathscr{F})$ the equation $G=HF$ where $H$ is a square-zero matrix will hold if and only if $\text{N}(F) \subseteq \text{N}(G)$ and $\text{R}(G) \cap \text{R}(F) \subseteq \text{R}(F \restriction \text{N}(G))$. \end{quote}
Whenever this condition holds a matrix $H$ of each allowable rank within the following range may be constructed:
\begin{quote}
\[\text{r}(G) \leq \text{r}(H) \leq \left \{ {\renewcommand{\arraystretch}{1.2} \begin{array}{ll} \frac{m}{2} & \text{if r}({\renewcommand{\arraystretch}{1}\begin{bmatrix} G & F \end{bmatrix}}) - \text{r}(G) \leq \frac{m}{2} \\ 
m-(\text{r}({\renewcommand{\arraystretch}{1}\begin{bmatrix} G & F \end{bmatrix}}) - \text{r}(G)) & \text{if r}({\renewcommand{\arraystretch}{1}\begin{bmatrix} G & F \end{bmatrix}}) - \text{r}(G) > \frac{m}{2} \end{array}} \right. .\]
\end{quote}
Similar results hold where the roles of $G$ and $F$ are interchanged, i.e. in the case of left division. The theory developed above can be applied in proving results for square matrices. Theorem \ref{thm:botha2FactorSZ2} presents the main result for products of two square-zero matrices:
\begin{quote}
Let $G \in M_{m}(\mathscr{F})$, where $\mathscr{F}$ denotes an arbitrary field. Then $G=HF$, where $H,F \in M_{m}(\mathscr{F})$ are both square-zero if and only if 
\[ \text{r}(G) \leq \text{n}(G) - \dim(\text{R}(G) \cap \text{N}(G)). \label{eq1th5} \]
If this condition holds, then the ranks of $H$ and $F$ can be arbitrary provided 
\[\text{r}(G) \leq \text{r}(H), \text{r}(F) \leq \frac{m}{2}.\] 
\end{quote}
Products of three (or more) square-zero matrices are addressed in Theorem \ref{thm:botha3FactorsSz2}:
\begin{quote}
Let $G \in M_{m}(\mathscr{F})$, where $\mathscr{F}$ denotes an arbitrary field. Then $G$ is a product of square-zero matrices if and only if $\text{r}(G) \leq \frac{m}{2}$, which is equivalent to stating $\text{r}(G) \leq \text{n}(G)$. The minimum number of square-zero matrices in such a factorization is at most three, which is a sharp upper bound.

Moreover, each matrix $G$ that satisfies the above condition can be expressed as a product of three square-zero matrices with arbitrary ranks~$r_i$ subject to \[\text{r}(G) \leq r_i \leq \frac{m}{2} \text { for } 1 \leq i \leq 3.\]
\end{quote}

 Finally, which matrices are sums of square-zero matrices?
 \begin{quote}
 In general, a matrix can be written as a sum of square-zero matrices if and only if its trace is zero, and if this is the case the least number of square-zero matrices required is at most four. (Theorem \ref{thm:trace_zero_sum_of_4}).\end{quote}
 Which matrices are sums of two square-zero matrices? Combining Theorems \ref{thm:2_sqsum_nec} and \ref{thm:sufficient_sqzero_sum} we have
 \begin{quote}
 A matrix $A \in M_n(\mathscr{F})$ is a sum of two square-zero matrices if and only if all its invariant polynomials are either odd or even.
 \end{quote}
 An even, or odd(-power) polynomial is defined to be a polynomial \[p(x)=\sum_{i=0}^n c_{n-i}x^{n-i}\] where $c_{n-i}=0$ if $i$ is odd. The polynomial is even-power if $n$ is even, else it is odd.
 
 Which matrices are sums of three square-zero matrices? In general this problem remains open, but we have the following result for matrices over a field of characteristic two (Theorem \ref{thm:fieldchar2_sum_of_three_main_result}):
 \begin{quote}
 A matrix $A \in M_n(\mathscr{F})$ is a sum of three square-zero matrices if and only if its trace is zero.
 \end{quote}
 For fields not of characteristic two we have the following necessary conditions (Theorem \ref{prop:3sum_nec_ww}):
 \begin{quote}
 If the matrix $A \in M_n(\mathscr{F})$ is a sum of three square-zero matrices then the trace of $A$ is zero and $\text{n}(A - \lambda I) \leq 3n/4$ for any nonzero $\lambda \in \mathscr{F}$.
 \end{quote}
 These conditions are unfortunately not suffcient, as shown in Proposition \ref{prop:takahashi_prop_1}. 
 
 The final results are only valid over certain fields. First, we have a sufficient condition in Proposition \ref{prop:invariant_factors_order_2_general} which is valid over any field which has characteristic zero.
\begin{quote} 
Let $\mathscr{F}$ be a field which has characteristic zero. If all the invariant polynomials of $A \in M_n(\mathscr{F})$ have degree greater than one, and the trace of $A$ is zero, then $A$ is a sum of three square-zero matrices.   
\end{quote}
The following sufficient condition, which is further restricted to algebraically closed fields, was presented in Proposition \ref{prop:tak_suf_2}:
\begin{quote}
Let $\mathscr{F}$ be an algebraically closed field which has characteristic zero. Let $A \in M_n(\mathscr{F})$ have trace zero, and let $m$ be the number of invariant polynomials of $A$ which are of degree two. Then $A$ is a sum of three square-zero matrices if $\text{n}(A-\lambda I) \leq (2n-m)/3$ for all $\lambda \in \mathscr{F}$.   
\end{quote}
A useful corollary to this result which references the matrix structure only in terms of its order and its maximal eigenspace was presented in Corollary \ref{tak_cor_1}:
\begin{quote}
Let $\mathscr{F}$ be an algebraically closed field which has characteristic zero. If the trace of $A \in M_n(\mathscr{F})$ is zero and $\text{n}(A-\lambda I) \leq n/2 + 1$ for all $\lambda \in \mathscr{F}$, then $A$ is a sum of three square-zero matrices.  
\end{quote}
\addcontentsline{toc}{chapter}{Bibliography}
\printbibliography

\end{document}